\pdfoutput=1
\documentclass[a4paper,leqno]{compositio}
\usepackage[francais]{babel}
\usepackage[latin1]{inputenc}
\usepackage{amsthm}
\usepackage{amsmath}
\usepackage{amssymb}
\usepackage{amscd}
\usepackage[mathscr]{euscript}
\usepackage{rotating}
\usepackage[all]{xy}
\usepackage{bm}

\theoremstyle{plain}

\newtheorem{lemme}[subsection]{Lemme}
\newtheorem{proposition}[subsection]{Proposition}
\newtheorem{corollaire}[subsection]{Corollaire}
\newtheorem*{utheoreme}{Théorème}
\newtheorem*{ucorollaire}{Corollaire}

\newtheorem{subtheoreme}[subsubsection]{Théorème}
\newtheorem{sublemme}[subsubsection]{Lemme}
\newtheorem{subproposition}[subsubsection]{Proposition}
\newtheorem{subcorollaire}[subsubsection]{Corollaire}
\newtheorem{subfait}[subsubsection]{Fait}

\theoremstyle{definition}

\newtheorem{notation}[subsection]{Notation}

\newtheorem{subdefinition}[subsubsection]{Définition}

\theoremstyle{remark}

\newtheorem{remarque}[subsection]{Remarque}
\newtheorem{exemple}[subsection]{Exemple}

\newtheorem{subremarques}[subsubsection]{Remarques}
\newtheorem{subremarque}[subsubsection]{Remarque}
\newtheorem{subexemple}[subsubsection]{Exemple}

\newenvironment{preuve}{\vspace{0.3cm}\emph{Démonstration.}}
{\hfill $\square$ \vspace{0.3cm}}
\newenvironment{preuvet}{\vspace{0.3cm}\emph{Démonstration du théorème.}}
{\hfill $\square$ \vspace{0.3cm}}
\newenvironment{preuvep}{\vspace{0.3cm}\emph{Démonstration de la proposition.}}
{\hfill $\square$ \vspace{0.3cm}}
\newenvironment{preuvepn}[1]{\vspace{0.3cm}\emph{Démonstration de la
proposition #1.}} {\hfill $\square$ \vspace{0.3cm}}

\newcommand\C{\mathbb{C}}
\newcommand\Nat{\mathbb{N}}
\newcommand\R{\mathbb{R}}
\newcommand\Z{\mathbb{Z}}
\newcommand\Q{\mathbb {Q}}
\newcommand\Fi{\mathbb{F}}
\newcommand\Af{\mathbb{A}_f}
\newcommand\Of{\mathcal{O}}
\newcommand\Ade{\mathbb{A}}

\newcommand\SD{\mathbb{S}}
\newcommand\T{{\bf T}}
\newcommand\G{{\bf G}}
\renewcommand\H{{\bf H}}
\newcommand\GU{{\bf GU}}
\newcommand\U{{\bf U}}

\newcommand\GSp{{\bf GSp}}
\newcommand\Sp{{\bf Sp}}
\newcommand\GSpin{{\bf GSpin}}
\newcommand\Spin{{\bf Spin}}
\newcommand\GSO{{\bf GSO}}
\newcommand\GO{{\bf GO}}
\newcommand\PSO{{\bf PSO}}
\renewcommand\O{{\bf O}}
\newcommand\GL{{\bf GL}}
\newcommand\Pa{{\bf P}}
\newcommand\QP{{\bf Q}}

\newcommand\B{{\bf B}}
\newcommand\Gr{\mathbb{G}}
\newcommand\N{{\bf N}}
\newcommand\M{{\bf M}}
\newcommand\Le{{\bf L}}

\newcommand\A{{\bf A}}

\newcommand\K{\mathrm {K}}

\newcommand\PGL{{\bf PGL}}

\newcommand\SO{{\bf SO}}

\newcommand\F{{\mathcal F}}
\newcommand\Gf{{\mathcal G}}

\newcommand\Kgoth{{\mathfrak K}}
\newcommand\Sgoth{{\mathfrak S}}

\DeclareMathOperator{\Ad}{Ad}

\DeclareMathOperator{\Aut}{Aut}
\DeclareMathOperator{\Cent}{Cent}
\DeclareMathOperator{\Gal}{Gal}
\DeclareMathOperator{\Ho}{H}
\DeclareMathOperator{\Int}{Int}
\DeclareMathOperator{\Ker}{Ker}

\DeclareMathOperator{\Nor}{Nor}

\DeclareMathOperator{\Tr}{Tr}
\DeclareMathOperator{\vol}{vol}
\newcommand\Bo{\mathcal{B}}

\newcommand\Dcal{\mathcal{D}}

\newcommand\Ell{\mathcal{E}}

\newcommand\Hecke{{\mathcal H}}

\newcommand\Levi{{\mathcal L}}

\newcommand\Norme{{\mathcal N}}

\DeclareMathOperator{\oubli}{oub}

\newcommand\Par{{\mathcal P}}
\newcommand\Qpar{{\mathcal Q}}
\newcommand{\quash}[1]{}

\newcommand\sous{\setminus}

\newcommand\X{{\mathcal X}}

\newcommand\fl{\longrightarrow}
\newcommand\fle{\longmapsto}
\newcommand\iso{\stackrel {\sim} {\fl}}

\newcommand\ungras{1\mkern -5mu\mathrm{l}}
\newcommand\ddotsinv{\begin{turn}{45}\large\ldots\end{turn}}

\title{Cohomologie d'intersection des variétés modulaires de Siegel, suite}
\author{Sophie Morel}
\email{morel@math.harvard.edu}
\address{Department of Mathematics, Harvard University\\
One Oxford Street\\
Cambridge, MA 02138, USA}

\classification{11G18 (primary), 14F20, 20G35 (secondary)}
\keywords{Siegel modular varieties, intersection cohomology, discrete
automorphic representations of symplectic groups}

\thanks{Ce texte a été écrit pendant que j'étais employée par le 
Clay Mathematics Institute en tant que Clay Research Fellow, et accueillie en
tant que membre à l'Institute for Advanced Study à Princeton. Il a été révisé
pendant mon séjour à l'université Harvard en tant que visiteur, puis en tant
que professeur.
De plus, j'ai bénéficié
du soutien financier de la NSF à travers les contrats DMS-0111298 et
DMS-0635607.}

\begin{document}

\begin{abstract}

In this work, we study the intersection cohomology of Siegel modular varieties.
The goal is to express the trace of a Hecke operator
composed with a power of the Frobenius endomorphism (at a good place)
on this cohomology in terms of the geometric side of Arthur's invariant trace
formula for well-chosen test functions.

Our main tools are the results of Kottwitz about the contribution of the
cohomology with compact support and about the stabilization of the trace
formula, Arthur's $L^2$ trace formula and the
fixed point formula of \cite{M2}. We ``stabilize'' this last formula, ie
express it as a sum of stable distributions on the general symplectic groups
and its endoscopic groups, and obtain the formula conjectured by Kottwitz
in \cite{K-SVLR}.

Applications of the results of this article have already been given by
Kottwitz, assuming Arthur's conjectures. Here, we give weaker unconditional
applications in the cases of the groups $\GSp_4$ and $\GSp_6$.

\end{abstract}

\maketitle

\tableofcontents

\section*{Introduction}

Cet article est la suite de l'article \cite{M2}, et il est parallèle au livre
\cite{M3}, dont il suit globalement la structure et utilise certains des
résultats.
Son but est de continuer l'étude de la cohomologie d'intersection
de la compactification de Satake-Baily-Borel des variétés modulaires de Siegel
(l'objet de \cite{M3} était le cas des variétés de Shimura associées aux
groupes unitaires sur $\Q$).

La méthode utilisée est celle développée par Ihara, Langlands et Kottwitz :
comparaison de la formule des points fixes de Grothendieck-Lefschetz et
de la formule des traces d'Arthur-Selberg. Le premier pas,
c'est-à-dire le calcul de la trace sur la cohomologie d'intersection
d'une correspondance de Hecke composée avec une puissance (assez grande)
du morphisme de Frobenius en une bonne place, était l'objet de
l'article \cite{M2} (complété par le chapitre 1 de \cite{M3}). On exprime
ici la formule obtenue en fonction du côté géométrique de la
formule des traces stable sur le groupe général symplectique et ses
groupes endoscopiques. Notons que cette ``stabilisation'' utilise
le lemme fondamental et certains cas du lemme fondamental tordu; ces
résultats sont désormais disponibles
grâce aux travaux de Laumon-Ngo (\cite{LN}),
Ngo (\cite{Ng}) et Waldspurger (\cite{Wa1}, \cite{Wa2},
\cite{Wa3}).

Pour la partie elliptique de la formule des points fixes (c'est-à-dire la
partie qui provient de la cohomologie à supports compacts de la variété
de Shimura), la stabilisation modulo les lemmes fondamentaux est due à Kottwitz
(\cite{K-SVLR}). On applique ici la méthode de Kottwitz aux sous-groupes
de Levi du groupe général symplectique pour obtenir la stabilisation
des termes non elliptiques; pour faire le lien entre les sous-groupes
de Levi des groupes endoscopiques et les groupes endoscopiques des
sous-groupes de Levi, on utilise aussi des méthodes de Kottwitz (cf
\cite{K-NP}). Le seul point qui est vraiment plus difficile que dans
le cas des variétés de Shimura compactes est la stabilisation de la
partie à l'infini (qui prend d'ailleurs la moitié de cet article).

Le théorème principal est le théorème \ref{th:adieu_partie_lineaire_GSp}.
Son corollaire \ref{cor:stab_FT_IC} est la formule conjecturée par
Kottwitz dans \cite{K-SVLR} (10.1) (pour les variétés modulaires de
Siegel) :

\begin{utheoreme} Soit $S^\K$ la variété modulaire de Siegel associée
au groupe $\G=\GSp_{2n}$ et
de niveau $\K\subset\G(\Af)$ net. On note $IH^*(S^\K,V)$ la cohomologie
d'intersection de la compactification de Satake-Baily-Borel de  $S^K$ à
coefficients dans une représentation algébrique $V$ de $\G$. Soit
$p$ un nombre premier tel que $\K=\G(\Z_p)\K^p$, avec
$\K^p\subset\G(\Af^p)$; on note $Frob_p\in\Gal(\overline{\Q}/\Q)$ un
relèvement du Frobenius géométrique en $p$. Alors, pour toute fonction
$f^\infty\in C_c^\infty(\K\sous\G(\Af)/\K)$
telle que $f^\infty=f^{\infty,p}\ungras_{\G(\Z_p)}$
et pour tout entier $j$ assez grand, on a
\[\Tr(Frob_p^j\times f^\infty,IH^*(S^\K,V))=\sum_{(\H,s,\eta_0)}\iota(\G,\H)
ST^H(f_\H^{(j)}),\]
où la somme est sur les classes d'équivalence de triplets endoscopiques
elliptiques $(\H,s,\eta_0)$ de $\G$ tels que $\H_\R$ admette un tore
maximal elliptique, les $f_\H^{(j)}$ sont des fonctions
sur $\H(\Ade)$ déterminées par $f^\infty$, $j$ et $V$ et $ST^H$ est le
côté géométrique de la formule des traces stable sur $\H$.

\end{utheoreme}

Les notations précises sont expliquées dans \ref{stabilisation2}.
Signalons tout de même ce que l'on entend par ``côté
géométrique de la formule des traces stable''. Pour les fonctions
qui apparaissent lorsque l'on calcule la cohomologie des variétés de
Shimura (c'est-à-dire les fonctions cuspidales stables à l'infini),
Arthur a donné dans \cite{A-L2} (formule (3.5) et théorème 6.1) une
expression simple de la formule des traces invariante. Kottwitz a
stabilisé le côté géométrique de cette expression dans \cite{K-NP}; notons
que cette stabilisation n'utilise que le lemme fondamental (et pas le
lemme fondamental pondéré). On utilise ici la formule de Kottwitz.
Malheureusement, l'article \cite{K-NP} n'est pas publié, et il ne considère
que le cas des groupes à groupe dérivé simplement connexe, alors que
les groupes endoscopiques des groupes symplectiques ne vérifient pas
tous cette condition. Cependant, il ne fait aucun doute que les méthodes
de \cite{K-NP} s'adapteraient au cas général, avec plus de complications
techniques. C'est la formule de \cite{K-NP} qui est utilisée dans cet
article.

Si l'on admet les conjectures d'Arthur sur la description du spectre
discret des groupes symplectiques et orthogonaux-symplectiques de
\ref{formule_points_fixes1}, alors on peut déduire d'une formule
comme celle du corollaire \ref{cor:stab_FT_IC} la description complète
de la cohomologie d'intersection. Cela a été fait par Kottwitz dans les
section 8 à 10 de \cite{K-SVLR}. Arthur a annoncé une preuve de ses
conjectures pour les groupes orthogonaux et symplectiques (en
admettant la stabilisation de la formule des traces tordue). Cependant,
on aurait besoin ici des conjectures d'Arthur pour les groupes
\emph{généraux} symplectiques (et certains groupes généraux
orthogonaux-symplectiques). Plutôt que d'utiliser
les résultats annoncés par Arthur pour tenter de donner des applications
aussi générales que possible du corollaire \ref{cor:stab_FT_IC}, on a
choisi de donner des applications inconditionnelles
pour les groupes $\GSp_4$ et
$\GSp_6$, où la situation est particulièrement simple (principalement
parce que le groupe $\PSO_4$ est isomorphe à $\PGL_2\times\PGL_2$, cf
la proposition \ref{prop:ST=T}). On obtient par exemple le résultat suivant
sur la fonction $L$ du complexe d'intersection
(proposition \ref{prop:fonction_L}) :

\begin{ucorollaire} On garde les notations du théorème ci-dessus, et
on suppose que $\G=\GSp_4$ ou $\G=\GSp_6$. Alors
\begin{flushleft}$\displaystyle{
\log L_p(s,IH^*(S^\K,V))=\sum_{\pi_f}c_\G(\pi_f)\dim(\pi_f^\K)\log
L(s-\frac{d}{2},\pi_p,r_{-\mu})
}$\end{flushleft}
\begin{flushright}$\displaystyle{
+\sum_{\pi_{H,f}}c_\H(\pi_{H,f})
\Tr(\pi_{H,f}((\ungras_\K)^\H))
\sum_{\mu_H\in M_H}\varepsilon(\mu_H)\log L(s-\frac{d}{2},\pi_{H,p},r_{-\mu_H})
.
}$\end{flushright}
Dans cette formule, $\H$ est égal à $\GSO_4$ si $\G=\GSp_4$ et à $\G(\Sp_2
\times\SO_4)$ si $\G=\GSp_6$, la première (resp. deuxième) somme est sur
l'ensemble des classes d'isomorphisme de représentations irréductibles
admissibles $\pi_f$ (resp. $\pi_{H,f}$) de $\G(\Af)$ (resp. $\H(\Af)$), les
$c_\G(\pi_f)$ et $c_\H(\pi_{H,f})$ sont des coefficients définis dans
la section \ref{applications}, $(\ungras_\K)^\H$ est un transfert de
$\ungras_\K$,
$M_H$ est un ensemble de cocaractères de
$\H$, les $\varepsilon(\mu_h)$ sont des signes,
$r_{-\mu}$ (resp. $r_{-\mu_H}$) est la représentation algébrique
de $\widehat{\G}$ (resp. $\widehat{\H}$) de plus haut poids $-\mu$ (resp.
$-\mu_H$) et $d=\dim(S^\K)$.

\end{ucorollaire}

Signalons enfin que le cas du groupe $\GU(2,1)$ est traité dans le livre
\cite{LR} (plus complètement, puisque les conjectures d'Arthur sont connues
dans ce cas par les travaux de Rogawski, cf \cite{Ro1}), celui
des groupes unitaires sur $\Q$ est traité dans
le livre \cite{M3} (dont on utilisera souvent les lemmes techniques) et que,
dans le cas du groupe $\GSp_4$, Gérard Laumon
a appliqué les mêmes méthodes à l'étude de la cohomologie à supports
compacts (cf \cite{Lau}, \cite{Lau2}).

Donnons une description rapide des différentes sections.

La section \ref{formule_points_fixes}
est consacrée à des définitions et des rappels sur la
formule des points fixes. La formule des points fixes que l'on utilise
ici (théorème \ref{th:points_fixes_moi}) est celle de \cite{M3} 1.7 où
on a un peu réarrangé les termes pour rendre la stabilisation plus facile.

La section \ref{endoscopie}
contient des rappels sur l'endoscopie et le calcul des données
endoscopiques elliptiques du groupes $\GSp_{2n}$ et de ses sous-groupes de
Levi. En particulier, on rappelle dans \ref{endoscopie2} quelques
définitions non standard de \cite{K-NP}.

Les sections \ref{infini} et \ref{partie_en_p} prouvent les résulats locaux
qui sont utilisés dans \ref{stabilisation2}. La section \ref{partie_en_p}
est consacrée aux calculs à la place $p$ et la section \ref{infini} à ceux à la
place infinie (cette dernière section, avec son appendice
\ref{lemmes_combinatoires}, contient la partie la plus technique et pénible de
l'article, et celle où la différence avec le cas des groupes unitaires est la
plus grande).

La section \ref{stabilisation} énonce le résultat de Kottwitz (\cite{K-NP})
sur le côté géométrique de la formule des traces stable et donne la
stabilisation de la formule des points fixes.

Enfin, la section \ref{applications} donne des applications des résultats
de \ref{stabilisation2} dans le cas des groupes $\GSp_4$ et $\GSp_6$. On
y prouve une formule pour la trace d'une puissance du Frobenius sur
les composantes isotypiques (pour l'action de $\GSp_{2n}(\Af)$) de la
cohomologie d'intersection (théorème \ref{th:composantes_isotypiques}).

Je remercie vivement Robert Kottwitz, qui m'a apporté une aide précieuse en
corrigeant certaines de mes idées fausses sur l'endoscopie et en me permettant
de lire son manuscrit \cite{K-NP}, ainsi que Gérard Laumon. Je remercie aussi 
les autres mathématiciens qui ont répondu à mes questions ou m'ont signalé des
simplifications, en particulier Pierre-Henri Chaudouard, Laurent Fargues,
Günter Harder, Colette Moeglin, Bao Chau Ngo,
Sug Woo Shin et Marie-France Vignéras. Enfin, je remercie le rapporteur anonyme
qui a signalé plusieurs erreurs et inexactitudes dans les versions précédentes
de ce texte, et m'a patiemment aidée à débusquer une erreur tenace dans
l'appendice.

\vspace{1cm}

Dans tout cet article, on utilisera la notation suivante :
Soient $F$ un corps,
$\G$ un groupe algébrique sur $F$ et $\gamma\in\G(F)$. Alors
on note $\G_\gamma=\Cent_\G(\gamma)^0$.

\section{La formule des points fixes}
\label{formule_points_fixes}

Le but de cette section est de rappeler la formule des points fixes du chapitre
1 de \cite{M3} (cf aussi \cite{M2} pour le cas d'une correspondance de Hecke
triviale), sous une forme un peu plus commode pour le processus de
stabilisation de la section \ref{stabilisation}.

\subsection{Définition des groupes et des données de Shimura}
\label{formule_points_fixes1}

Dans ce paragraphe, on définit les groupes symplectiques
et leurs données de Shimura, et on rappelle la description de
leurs sous-groupes paraboliques. On introduit aussi les
groupes orthogonaux, car ceux-ci
apparaissent comme groupes endoscopiques
des groupes symplectiques (cf \ref{endoscopie1}).

Pour $n\in\Nat^*$, on note
\[I=I_n=\left(\begin{array}{ccc}1 & & 0 \\ & \ddots & \\ 0 & & 1\end{array}
\right)\in\GL_n(\Z)\]
\[A_n=\left(\begin{array}{ccc}0 & & 1 \\ & \begin{turn}{45}\large\ldots
\end{turn}& \\ 1 & & 0\end{array}\right)\in\GL_n(\Z)\]
et
\[B_n=\left(\begin{array}{cc}0 & A_n \\
-A_n & 0\end{array}\right)\in\GL_{2n}(\Z).\]

Soient $n\in\Nat^*$ et $J\in\GL_n(\Q)$ une matrice symétrique. On définit des
groupes algébriques $\GO(J)$ et $\GSp_{2n}$ sur $\Q$ en posant, pour
toute $\Q$-algèbre $A$ :
\[\GO(J)(A)=\{g\in\GL_n(A)|{}^tgJg=c(g)J,c(g)\in A^\times\}\]
\[\GSp_{2n}(A)=\{g\in\GL_{2n}(A)|{}^tgB_ng=c(g)B_n,c(g)\in A^\times\}.\]
On a des morphismes de groupes algébriques sur $\Q$:
\[c:\GO(J)\fl\Gr_m\mbox{ et }\det:\GO(J)\fl\Gr_m\]
\[c:\GSp_{2n}\fl\Gr_m\mbox{ et }\det:\GSp_{2n}\fl\Gr_m.\]
On note $\O(J)=\Ker(c)$ et $\SO(J)=\Ker(c)\cap\Ker(\det)$
(resp. $\Sp_{2n}=\Ker(c)=
\Ker(c)\cap\Ker(\det)$). Le groupe $\GO(J)$ n'est pas forcément connexe (il a
deux composantes connexes si $n$ est pair); on note $\GSO(J)=\GO(J)^0$.

Le groupe $\SO(J)$ (resp. $\Sp_{2n}$)
est le groupe dérivé de $\GSO(J)$ (resp.
de $\GSp_{2n}$).
Les groupes $\GSO(J)$ et $\GSp_{2n}$ sont réductifs
connexes, et le groupe
$\Sp_{2n}$ est semi-simple simplement connexe.

On note $\GSO_n=\GSO(A_n)$.
Le groupe $\GSO_n$ est la forme intérieure quasi-déployée
du groupe $\GSO(J)$,
pour toute $J\in\GL_n(\Q)$ symétrique, et il est déployé de rang semi-simple
$\lfloor n/2\rfloor$.
D'autre part, le groupe $\GSp_{2n}$ est déployé sur $\Q$, de rang semi-simple
$n$.

Enfin, on note $\GSp_0=\GSO_0=\Gr_m$, et
$(c:\Gr_m\fl\Gr_m)=id$.

\begin{subremarque}\label{rq:groupes_sur_Z}
Si la matrice $J$ est dans $\GL_n(\Z)$,
on peut étendre $\GO(J)$ en un schéma en groupes $\Gf'$ sur $\Z$ en posant, 
pour toute $\Z$-algèbre $A$,
\[\Gf'(A)=\{g\in\GL_n(A)|{}^tgJg=c(g)J,c(g)\in A^\times\}.\]
Alors $\Gf:=(\Gf')^0$ est un schéma en groupes sur $\Z$ qui étend
$\GSO(J)$ et,
pour tout nombre premier $\ell$, $\Gf_{\Fi_\ell}$ est
un groupe algébrique réductif connexe sur $\Fi_\ell$.

On a évidemment une construction similaire pour les groupes
$\GSp_{2n}$.

\end{subremarque}

\vspace{1cm}

On définit maintenant les données de Shimura.
On note $\SD=R_{\C/\R}\Gr_m$.
On fixe
$n\in\Nat^*$, et on note $\G=\GSp_{2n}$.
On note $\X_{n}^+$ ou $\X^+$ (resp. $\X_n^-$ ou $\X^-$) l'ensemble des
morphismes $h:\SD\fl\G_{\R}$ qui induisent une structure de Hodge pure de type 
$\{(0,1),(1,0)\}$ sur $\Q^{2n}$, et tels que la forme bilinéaire $\R^{2n}\times
\R^{2n}\fl\R$, $(v,w)\fle {}^tvB_nh(i)w$ soit symétrique définie positive
(resp. négative).
Le groupe $\G(\R)$ agit transitivement (par conjugaison) sur $\X=\X^+\cup\X^-$,
et le morphisme
$$h:z=a+ib\fle\left(\begin{array}{cc}aI_n & -bA_n \\
bA_n & aI_n\end{array}\right)$$
est dans $\X^+$, 
donc $\X\simeq\G(\R)/Stab_{\G(\R)}(h)$.
Le triplet $(\G,\X,h)$ est une donnée de Shimura pure au sens de \cite{P1} 2.1.
Remarquons au passage que le cocaractère $\mu:\Gr_{m,\C}\fl\G_\C$ associé à
$h$ comme dans (par exemple) \cite{K-SVLR} \S1 est $z\fle diag(z I_n,I_n)$;
en particulier, il est défini sur $\Q$, donc le corps reflex de la donnée
$(\G,\X,h)$ est $\Q$.

\vspace{1cm}

Rappelons la description des sous-groupes paraboliques des groupes
symplectiques.

Un tore maximal de $\G$ est 
\[\T=\left\{\left(\begin{array}{ccc}\lambda_1 &  & 0 \\  & \ddots &  
\\ 0 &  & \lambda_{2n}\end{array}\right),\lambda_1,\dots,\lambda_n\in
\Gr_m,\lambda_1\lambda_{2n}=\lambda_2\lambda_{2n-1}\dots=\lambda_n\lambda_{n+1}
\right\}.\]
Il est déployé.
L'intersection $\Pa_0$ de $\G$ avec le sous-groupe des matrices triangulaires
supérieures de $\GL_{2n}$ est un sous-groupe de Borel de $\G$ contenant $\T$.
Les sous-groupes paraboliques standard de $\G$
sont indexés par les sous-ensembles de $\{1,\dots,n\}$ de la manière 
suivante.

Soit $S\subset\{1,\dots,n\}$. On écrit $S=\{r_1,r_1+r_2,\dots,r_1+\dots+r_m\}$
avec $r_1,\dots,r_m\in\Nat^*$ et on note $r=r_1+\dots+r_m$. Le sous-groupe
parabolique standard $\Pa_S$ correspondant est l'intersection avec $\G$
du groupe
\[\left(\begin{array}{ccccccc}\GL_{r_1} &  &  &  &  &  & * \\
 & \ddots &  &  &  &  &  \\
 &  & \GL_{r_m} &  &  &  &  \\
 &  &  & \GSp_{2(n-r)} &  &  &  \\
 &  &  &  & \GL_{r_m} &  &  \\
 &  &  &  &  & \ddots &  \\
0 &  &  &  &  &  & \GL_{r_1}\end{array}\right).\]
En particulier, les sous-groupes paraboliques maximaux standard de $\G$ 
sont les
\[\Pa_r:=\Pa_{\{r\}}=\left(\begin{array}{ccc}\GL_r &  & * \\
 & \GSp_{2(n-r)} &  \\
0 &  & \GL_r\end{array}\right)\cap\G\]
pour $r\in\{1,\dots,n\}$, et on a $\Pa_S=\bigcap\limits_{r\in S}\Pa_r$.

On note $\N_S$ (ou $\N_{P_S}$) le radical unipotent de $\Pa_S$, $\M_S$
(ou $\M_{P_S}$) le sous-groupe de Levi
évident (formé des matrices diagonales par blocs) et $\A_{M_S}$ le sous-tore
déployé maximal du centre de $\M_S$.
Si on écrit comme plus
haut $S=\{r_1,\dots,r_1+\dots+r_m\}$ et $r=r_1+\dots+r_m$, on a un
isomorphisme
\[\begin{array}{ccl}\M_S & \iso & \GL_{r_1}\times\dots\times\GL_{r_m}
\times\GSp_{2(n-r)} \\
diag(g_1,\dots,g_m,g,h_m,\dots,h_1) & \fle & (c(g)^{-1}g_1,\dots,c(g)^{-1}g_m,
g)\end{array}.\]
L'image réciproque par cet isomorphisme de $\GL_{n_1}\times\dots\times
\GL_{r_m}$ est appelée \emph{partie linéaire} de $\M_S$; on la note
$\Le_S$ ou $\M_{S,l}$. L'image réciproque de $\GSp_{2(n-r)}$ est appelée
\emph{partie hermitienne} de $\M_S$, et notée $\G_r$ ou $\M_{S,h}$.
On voit en particulier que les sous-groupes paraboliques maximaux
de $\G$ vérifient l'hypothèse de la section 1.1 de \cite{M3}.
Remarquons aussi que, gr\^ace au choix de l'isomorphisme entre $\M_S$ et
$\GL_{r_1}\times\dots\times\GL_{n_r}\times\GSp_{2(n-r)}$, l'image du cocaractère
$\mu:\Gr_{m,\C}\fl\G_\C$ est contenue dans la partie hermitienne $\G_r$; le
cocaractère de $\G_r$ que l'on obtient est l'analogue de $\mu$ pour
$\G_r\simeq\GSp_{2(n-r)}$.

\subsection{Formule des points fixes}
\label{formule_points_fixes2}

Les références pour tous les faits énoncés dans ce paragraphe sont données
dans le chapitre 1 de \cite{M3}. Comme dans \cite{M3}, on fixe une fois
pour toutes des clôtures algébriques $\overline{\Q}$ de $\Q$ et
$\overline{\Q}_p$ de $\Q_p$, pour tout nombre premier $p$, et des inclusions
$\overline{\Q}\subset\overline{\Q}_p$ et $\overline{\Q}\subset\C$.

On fixe $n\in\Nat^*$, et on note $\G=\GSp_{2n}$.
Soit $\K$ un sous-groupe compact ouvert net de $\G(\Af)$ (cf \cite{P1} 0.6).
On note
$S^\K$ la variété de Shimura de niveau $\K$ définie par la donnée
$(\G,\X,h)$ de \ref{formule_points_fixes1}; c'est une variété algébrique
quasi-projective sur $\Q$, dont l'ensemble des points complexes est donné
par la formule
\[S^\K(\C)=\G(\Q)\sous(\X\times\G(\Af)/\K).\]
Soit $j:S^\K\fl\overline{S}^\K$ la compactification de Satake-Baily-Borel
de $S^\K$; la variété $\overline{S}^\K$ est aussi définie sur $\Q$, et elle
est projective et normale, mais elle n'est pas lisse si $n\geq 2$. Soit
$V$ une représentation algébrique de $\G$. Comme $\G$ est déployé,
$V$ est définie sur $\Q$. On fixe un nombre premier $\ell$ et un
isomorphisme $\overline{\Q}_\ell\simeq\C$.
Alors il est bien connu qu'on peut associer à $V$ un faisceau
$\ell$-adique lisse $\F^{\K}V$ sur $S^\K$ de fibre $V\otimes_\Q\Q_\ell$
en tout point géométrique de $S^\K$ (voir par exemple la section 5.1 de
\cite{P2}). Le complexe d'intersection sur $\overline{S}^\K$ à coefficients
dans $V$ est par définition
\[IC^\K V=(j_{!*}(\F^\K V[d]))[-d],\]
où $d$ est la dimension de $S^\K$ (donc $d=n(n+1)/2$).
La cohomologie de $IC^{\K}V$ est une représentation (graduée) de
$\Gal(\overline{\Q}/\Q)\times\Hecke_\K$, où $\Hecke_\K=\Hecke(\G(\Af),\K)$
est l'algèbre de Hecke des fonctions de $C_c^\infty(\G(\Af),\C)$
qui sont bi-invariantes par $\K$.

On fixe $f^\infty\in\Hecke_\K$ et un nombre premier $p\not=\ell$ tel que
$\K=\K^p\G(\Z_p)$
avec $\K^p\subset\G(\Af^p)$ et $f^\infty=f^{\infty,p}\ungras_{\G(\Z_p)}$ avec
$f^{\infty,p}\in C_c^\infty(\G(\Af^p))$.
Soit $Frob_p\in\Gal(\overline{\Q}/\Q)$ un
relèvement du Frobenius géométrique en $p$. On va rappeler la formule de
\cite{M3} 1.7 pour la trace de $Frob_p^j\times f$ ($j>>0$) sur
la cohomologie de $IC^\K V$.

Rappelons (cf \cite{K-NP} 1.5) qu'un sous-groupe de Levi $\M$ de $\G$ est
dit \emph{cuspidal} si $(\M/\A_M)_\R$ admet un tore maximal anisotrope
(sur $\R$), où $\A_M$ est la partie déployée du centre de $\M$.
Il est facile de voir que les sous-groupes de Levi cuspidaux
de $\G$ sont ceux isomorphes à un groupe $\Gr_m^r\times\GL_2^t\times\GSp_{2(
n-r-2t)}$, avec $r,t\in\Nat$ et $r+2t\leq n$. 

Soit $\M$ un sous-groupe de Levi cuspidal
standard de $\G$ et soit $\M_h$ sa partie hermitienne.
Alors $\mu:\Gr_m\fl\G$ se factorise par $\M_h$, donc on
peut utiliser $\mu$ pour définir l'ensemble des triplets de Kottwitz
associé à $\M_h$ et à $j\in\Nat^*$ (les définitions sont celles de \cite{K-SVLR}
\S2 et \S3) : Soit
$j\in\Nat^*$. On note $L$ l'extension non ramifiée de degré $j$ de $\Q_p$ dans
$\overline{\Q}_p$ et $\sigma\in\Gal(L/\Q_p)$ le relèvement du Frobenius
arithmétique. Alors un triplet de Kottwitz est un triplet $(\gamma_0;
\gamma,\delta)$, où $\gamma_0\in\M_h(\Q)$, $\gamma=(\gamma_v)_{\not=p,\infty}
\in\M_h(\Af^p)$ et $\delta\in
\M_h(L)$, qui vérifie un certain ensemble de conditions, notées (C) dans
\cite{M3} 1.6 (on utilisera la même notation ici). On renvoie à
\cite{K-SVLR} \S2
pour l'énoncé précis de ces conditions: notons simplement
ici qu'en particulier, on demande que $\gamma_0$ soit semi-simple et
elliptique dans $\M_h(\R)$, que $\gamma_v$ et $\gamma_0$ soient stablement
conjugués pour toute place $v\not=p,\infty$ de $\Q$ et que $N\delta:=
\delta\sigma(\delta)\dots\sigma^{j-1}(\delta)$ et $\gamma_0$ soient
stablement conjugués. On dit que deux triplets de Kottwitz $(\gamma_0;\gamma,
\delta)$ et $(\gamma_0';\gamma',\delta')$ sont équivalents, et on note
$(\gamma_0;\gamma,\delta)\sim(\gamma'_0;\gamma,\delta)$, si $\gamma_0$ et
$\gamma'_0$ sont stablement conjugués, $\gamma$ et $\gamma'$ sont
conjugués dans $\M_h(\Af^p)$ et $\delta$ et $\delta'$ sont $\sigma$-conjugués
dans $\M_h(L)$ (ie il existe $x\in\M_h(L)$ tel que $\delta'=x\delta\sigma(x)^
{-1}$). Enfin, dans \cite{K-SVLR} \S2, Kottwitz associe à tout
$(\gamma_0;\gamma,\delta)$ un invariant cohomologique $\alpha(\gamma_0;
\gamma,\delta)\in\Kgoth(I_0/\Q)^D$ ($I_0=\M_{h,\gamma_0}$, cf loc. cit. pour
les autres notations), qui ne dépend que de la classe d'équivalence de
$(\gamma_0;\gamma,\delta)$ (en un sens précisé dans loc. cit.). Si
$(\gamma_0;\gamma,\delta)$ est tel que $\alpha(\gamma_0;\gamma,
\delta)=0$, Kottwitz définit dans \cite{K-SVLR} \S3 une forme intérieure
$I$ de $I_0$ telle que $I/\A_{M_h}$ soit anisotrope sur $\R$,
que $I_{\Q_v}\simeq\M_{h,\gamma_v}$ pour $v\not=p,\infty$ et que
$I(\Q_p)\simeq\M_h(L)_{\delta}^\sigma$ (le $\sigma$-centralisateur de $\delta$
dans $\M_h(L)$, ie l'ensemble des $x\in\M_h(L)$ tels que $x\delta\sigma(x)^{-1}
=\delta$); on note alors
\[c(\gamma_0;\gamma,\delta)=\vol(I(\Q)\sous I(\Af))\]
(l'autre facteur dans la définition de \cite{K-SVLR} \S3 disparaît
grâce au lemme \ref{lemme:Tamagawa}).
Comme dans \cite{M3} 1.6, on note $C_{M_h,j}$ l'ensemble des classes
d'équivalence de triplets de Kottwitz $(\gamma_0;\gamma,\delta)$ tels que
$\alpha(\gamma_0;\gamma,\delta)=0$.
Si $\M_h$ est un tore, on a un sous-ensemble
$C'_{M_h,j}$ de $C_{M_h,j}$ défini dans la remarque 1.6.5 de \cite{M3} :
c'est le sous-ensemble des classes d'équivalence de
triplets $(\gamma_0;\gamma,\delta)$ tels que
$c(\gamma_0)>0$. Si $\M_h$ n'est pas un tore, on note
$C'_{M_h,j}=C_{M_h,j}$.
Rappelons aussi que l'on définit une fonction
$\phi_j^{M_h}$ de $\Hecke(\M_h(L),\M_h(\Of_L))$ par
\[\phi_j^{M_h}=\ungras_{\M_h(\Of_L)\mu(\varpi_L^{-1})\M_h(\Of_L)}\in\Hecke(\M_h(L),
\M_h(\Of_L)),\]
où $\varpi_L$ est une uniformisante de $L$.

Soit $\M$ un sous-groupe de Levi standard de $\G$ (pas forcément cuspidal).
On notera dans la suite $\M_l$ et $\M_h$ les parties linéaire et hermitienne
de $\M$. Soit
$\Nor'_\G(\M)=\Nor_\G(\M)\cap\Nor_\G(\M_l)\cap\Nor_\G(\M_h)$ (le sous-groupe
de $\Nor_\G(\M)$ qui préserve la décomposition de $\M$ en sa partie linéaire
et sa partie hermitienne). 
On note $\Par(\M)$ l'ensemble des paires $(\QP,g)$, où $\QP$ est un sous-groupe
parabolique standard de $\G$ et $g$ est un élément de $\G(\Q)$ tel que
$g\M_h g^{-1}=\G_Q$ et $g\M_l g^{-1}$ soit
un sous-groupe de Levi de $\Le_Q$.
On considère la relation d'équivalence suivante sur $\Par(\M)$ :
$(\QP,g)\sim (\QP',g')$ si et seulement si $\QP=\QP'$ et il existe
$h_1\in\M_Q(\Q)$ et $h_2\in\M(\Q)$ tels que $g'=h_1gh_2$.
Pour tout $(\QP,g)\in\Par(\M)$, la classe de $\M_Q(\Q)$-conjugaison du
sous-groupe de Levi $g\M g^{-1}$ de $\M_Q$ ne dépend que de la classe
d'équivalence de $(\QP,g)$. Soit $\Pa$ le sous-groupe parabolique standard
de $\G$ de sous-groupe de Levi $\M$. On écrit $\Pa=\Pa_S$, avec
$S\subset\{1,\dots,n\}$. Pour tout $(t_s)_{s\in S}\in\Z^S$, on a défini
dans \cite{M2} 4.2 les complexes de cohomologie tronquée $\Ho^*(Lie(\N_S),V)_
{<t_s,s\in S}$ et $\Ho^*(Lie(\N_S),V)_{>t_s,s\in S}$; ce sont des
sous-$\M$-représentations (graduées) de $\Ho^*(Lie(\N_S),V)$. Si
$t_s=s(s+1)/2-n(n+1)/2$ pour tout $s\in S$, on note
$\Ho^*(Lie(\N_S),V)_{<0}$ et $\Ho^*(Lie(\N_S),V)_{>0}$ au lieu de
$\Ho^*(Lie(\N_S),V)_{<t_s,s\in S}$ et $\Ho^*(Lie(\N_S),V)_{>t_s,s\in S}$.

Enfin, pour tout sous-groupe parabolique $\Pa$ de $\G$, tous
sous-groupes de Levi $\M\subset\M'$ de $\G$ et toute place $v$ de $\Q$,
si $\gamma\in\Pa(\Q_v)$ et $\gamma_M\in\M(\Q_v)$, on note
\[\delta_{P(\Q_v)}(\gamma)=|\det(\Ad(\gamma),Lie(\N_P))|_v\]
\[D_M^{M'}(\gamma_M)=\det(1-\Ad(\gamma),Lie(\M')/Lie(\M))\]
\[n_M^G=|\Nor_\G(\M)(\Q)/\M(\Q)|.\]

Soit $\M$ un sous-groupe de Levi standard cuspidal de $\G$. 
Pour tout $\gamma_M\in\M(\R)$ semi-simple elliptique, on pose
\begin{flushleft}$\displaystyle{
L_M(\gamma_M)=\sum_{(\QP,g)\in\Par(\M)/\sim}|(\Nor'_\G(\M)/(\Nor'_\G(\M)
\cap g^{-1}\M_Qg))(\Q)|^{-1}m_Q(-1)^{\dim(\A_{gMg^{-1}}/\A_{M_Q})}
}$\end{flushleft}
\begin{flushright}$\displaystyle{
(n_{gMg^{-1}}^{M_Q})^{-1}|D_{gMg^{-1}}^{M_Q}(g\gamma_M
g^{-1})|^{1/2}\delta_{\QP(\R)}^{1/2}(g\gamma_Mg^{-1})\Tr(g\gamma_Mg^{-1},
\Ho^*(Lie(\N_Q),V)_{>0}),}$\end{flushright}
où $m_Q$ est égal à $2$ si la partie hermitienne de $\M_Q$ est un tore, et
à $1$ sinon.
On note
\[\Tr_M(f^\infty,j)=\sum_{\gamma_L}\sum_{(\gamma_0;\gamma,\delta)\in C'_{M_h,j}}
\chi(\M_{l,\gamma_L})c(\gamma_0;\gamma,\delta)
\delta_{P(\Q_p)}^{1/2}(\gamma_L\gamma_0)O_{\gamma_L\gamma}(f^{\infty,p}_M)
O_{\gamma_L}(\ungras_{\M_l(\Z_p)})TO_\delta
(\phi_j^{M_h})L_M(\gamma_L\gamma_0),\]
où $\gamma_L$ parcourt l'ensemble des classes de conjugaison de $\M_l(\Q)$,
$\chi(\M_{l,\gamma_L})$ est comme dans \cite{GKM} 7.10 et
$\Pa$ est le sous-groupe parabolique standard de sous-groupe de Levi
$\M$.
Les définitions des intégrales orbitales $O_\gamma$, des intégrales
orbitales tordues $TO_\delta$ et du terme constant $f^{\infty,p}_M$
sont rappelées
dans \cite{M3} 1.6 et 1.7. On normalise les mesures de Haar comme dans
\cite{M3} 1.6 (c'est la convention de \cite{K-SVLR} \S3) :
On utilise sur $\M_l(\Q_p)$ (resp. $\M_h(L)$) la mesure de Haar telle
que le volume de $\M_l(\Z_p)$ (resp. $\M_h(\Of_L)$) soit égal à
$1$. On prend sur $I(\Af^p)$ (resp. $I(\Q_p)$, resp.
$\M_l(\Af^p)_{\gamma_M}$, resp. $\M_l(\Q_p)_{\gamma_L}$) une mesure de Haar telle
que les volumes des sous-groupes ouverts compacts soient des nombres
rationnels, et on utilise les torseurs intérieurs de \cite{K-SVLR} \S3
pour transporter les deux premières mesures sur $\M_h(\Af^p)_\gamma$ et
$\M_h(L)_\delta^\sigma$.

\begin{subtheoreme}\label{th:points_fixes_moi} Si $j$ est assez grand, alors
\[\Tr(Frob_p^j\times f^\infty,H^*(\overline{S}^{\K}_{\overline{\Q}},IC^{\K}V))
=\sum_M\Tr_M(f^\infty,j),\]
où la somme est sur l'ensemble des sous-groupes de Levi cuspidaux standard
de $\G$.
De plus, si $f^\infty=\ungras_\K$, alors le théorème est vrai pour tout 
$j\in\Nat^*$.

\end{subtheoreme}

On aurait bien entendu un résultat similaire pour les espaces de
cohomologie pondérée de \cite{M2}. Ce théorème se déduit facilement du
théorème 1.7.1 de \cite{M3}.

\section{Endoscopie}
\label{endoscopie}

\subsection{Groupes endoscopiques des groupes symplectiques}
\label{endoscopie1}

Dans cette section, on s'intéresse aux triplets endoscopiques 
elliptiques $(\H,s,\eta_0)$ des groupes $\G$ définis dans
\ref{formule_points_fixes1}.
On utilise la définition des triplets endoscopiques et des isomorphismes de
triplets endoscopiques de \cite{K-STF:CTT} 7.4 et 7.5.

On commence par rappeler la définition des groupes généraux spinoriels (cf
\cite{Bo} 2.2(5)). Soit $n\in\Nat$. Le groupe $\{\pm 1\}$ s'identifie au
sous-groupe central $Z$ de $\Spin_n$ tel que $\Spin_n/Z=\SO_n$ (avec la
convention $\Spin_0=\{\pm 1\}$). On pose
\[\GSpin_n=(\Gr_m\times\Spin_n)/\{\pm 1\},\]
où $\{\pm 1\}$ est plongé diagonalement dans $\Gr_m\times\Spin_n$. Le
caractère $\Gr_m\times\Spin_n\fl\Gr_m$, $(z,g)\fle z^2$, passe au quotient
et donne un caractère $c:\GSpin_n\fl\Gr_m$, dont le noyau s'identifie
à $\Spin_n$. 
De même, si $n_1,\dots,n_r\in\Nat$, on pose
\[\G(\Spin_{n_1}\times\dots\times\Spin_{n_r})=(\Gr_m\times\Spin_{n_1}\times
\dots\Spin_{n_r})/\{\pm 1\},\]
où $\{\pm 1\}$ est plongé diagonalement. Le groupe $\G(\Spin_{n_1}\times
\dots\times\Spin_{n_r})$ s'identifie à $\{(g_1,\dots,g_r)\in\GSpin_{n_1}\times
\dots\times\GSpin_{n_r}|c(g_1)=\dots=c(g_r)\}$.

Soit $n\in\Nat^*$. On note $\G=\GSp_{2n}$. Soient $\T$ le tore diagonal de $\G$
(c'est un tore maximal de $\G$) et $\B$ le sous-groupe de $\G$ formé
des matrices triangulaires supérieures (c'est un sous-groupe de Borel de
$\G$).
On a un isomorphisme
\[\begin{array}{rcl}\T & = & \{(\lambda_1,\dots,\lambda_{2n})|\exists
\lambda\in\Gr_m,\forall i\in\{1,\dots,n\}\lambda_{i}\lambda_{2n+1-i}=
\lambda\}\\
& \simeq & \Gr_m\times\Gr_m^n\end{array}\]
donné par la formule $(\lambda_1,\dots,\lambda_{2n})\fle
(\lambda_1\lambda_{2n},(\lambda_1,\dots,\lambda_n))$.
Pour tout $i\in\{1,\dots,n\}$, on note $e_i$ le caractère de $\T$ défini par
\[e_i(\lambda_1,\dots,\lambda_{2n})=\lambda_i.\]
Alors le groupe des caractères de $\T$ est
\[X^*(\T)=\Z c\oplus\bigoplus_{i=1}^n\Z e_i.\]
Donc le tore dual de $\T$ est
\[\widehat{\T}=\C^\times\times(\C^\times)^n,\]
avec l'action triviale de $\Gal(\overline{\Q}/\Q)$.
Pour tout $i\in\{1,\dots,n\}$, on note $\widehat{e}_i$ le caractère
$(\lambda,(\lambda_1,\dots,\lambda_n))\fle\lambda_i$ de $\widehat{\T}$.
L'ensemble des racines de $\T$ dans $Lie(\G)$ est
\[\Phi=\Phi(\T,\G)=\{\pm (e_i-\frac{c}{2})\pm (e_j-\frac{c}{2}),1\leq i<j\leq
n\}\cup\{\pm(2e_i-c),1\leq i\leq n\}.\]
Le sous-ensemble de racines simples déterminé par $\B$ est
\[\Delta=\{\alpha_i:=e_i-e_{i+1},1\leq i\leq n-1\}\cup\{\alpha_n:=2e_n-c\}.\]
L'ensemble des coracines est
\[\widehat{\Phi}=\{\pm \widehat{e}_i\pm \widehat{e}_j,1\leq i<j\leq n\}\cup
\{\pm \widehat{e}_i,1\leq i\leq n\}.\]
En particulier, $Z(\widehat{\G})=\bigcap\limits_{\widehat{\alpha}\in\widehat{
\Phi}}\Ker(\widehat{\alpha})=\C^\times\times\{1\}$.
Le groupe de Weyl de $\T$ dans $\G$ est $\{\pm 1\}^n\rtimes\Sgoth_n$, agissant
sur $\T\simeq\Gr_m\times\Gr_m^n$ par
\[((\varepsilon_1,\dots,\varepsilon_n)\rtimes\sigma,(\lambda,(\lambda_1,
\dots,\lambda_n)))\fle (\lambda,(\lambda\lambda_{\sigma^{-1}(1)}^{\varepsilon_1},
\dots,\lambda\lambda_{\sigma^{-1}(n)}^{\varepsilon_n}))\]
(et sur $\widehat{\T}$ par une formule similaire).

Enfin, le groupe dual de $\G$ est
\[\widehat{\G}=\GSpin_{2n+1}(\C),\]
avec l'action triviale de $\Gal(\overline{\Q}/\Q)$ (cf \cite{Bo} 2.2(5)).

\begin{subproposition}\label{prop:groupes_endoscopiques_GSp}
Soit $K$ un corps local ou global de caractéristique $0$.
Un triplet endoscopique elliptique $(\H,s,\eta_0)$ de $\G_K$ est uniquement
déterminé par la donnée de $s$. On peut supposer que $s\in\widehat{\T}$, et on
a forcément $s\in Z(\widehat{\G})(\{1\}\times\{\pm 1\}^n)$ dans ce cas.
Un ensemble de représentants des classes d'équivalence
de $s$ est
\[\{s_{n_1}:=(\overbrace{1,\dots,1}^{n_1},\overbrace{-1,\dots,-1}^{n-n_1})|
0\leq n_1\leq n,n_1\not=n-1\},\]
et le groupe endoscopique correspondant à $s_{n_1}$ est
\[\H=\G(\Sp_{2n_1}\times\SO_{2(n-n_1)})_K:=\{(g_1,g_2)\in\GSp_{2n_1,K}\times
\GSO_{2(n-n_1),K}|c(g_1)=c(g_2)\}.\]
Le groupe $\Lambda(\H,s,\eta_0)$ de \cite{K-STF:CTT} 7.5 est égal à
$\{\pm 1\}$ si $n-n_1\geq 2$ et à $\{1\}$ si $n=n_1$.

De plus, si $K=\Q$ ou $\R$, alors $\H_\R$ admet un tore maximal elliptique
si et seulement si $n-n_1$ est pair.

\end{subproposition}

Soit $(\H,s,\eta_0)$ un triplet endoscopique elliptique de $\G_K$. Comme
$\G$ et $\H$ sont déployés sur $K$, on a un prolongement évident de
$\eta_0:\widehat{\H}\fl\widehat{\G}$ en un $L$-morphisme $\eta:{}^L\H\fl{}^L
\G$, qui est $\eta_0\times id_{W_K}$. Si $\H=\G(\Sp_{2n_1}\times\SO_{2n_2})_K$
(avec $n_1+n_2=n$ et $n_2\not=1$), alors $\widehat{\H}=\G(\Spin_{2n_1+1}\times
\Spin_{2n_2})$.

\begin{preuve} Comme $\G$ est déployé sur $K$ et de centre connexe, tous ses
groupes endoscopiques sont déployés sur $K$ (cf la définition 1.8.1 de
\cite{Ng}). Ceci implique qu'un triplet endoscopique $(\H,s,\eta_0)$ est
uniquement déterminé par $s$.

Soit $(\H,s,\eta_0)$ un triplet endoscopique elliptique de $\G_K$. Comme
$\H$ et $\G$ sont déployés sur $K$, la condition d'ellipticité s'écrit
simplement $Z(\widehat{\H})^0\subset Z(\widehat{\G})$.
Quitte à
remplacer $(\H,s,\eta_0)$ par un triplet équivalent, on peut supposer
que $s\in\widehat{\T}$ et
\[s=(1,(\overbrace{s_1,\dots,s_1}^{n_1},\dots,\overbrace{s_r,\dots,s_r}^{n_r})
),\]
avec $s_1,\dots,s_r\in\C^\times$ deux à deux distincts et $n_1,\dots,n_r\in
\Nat^*$ tels que $n=n_1+\dots+n_r$. Supposons qu'il existe $t\in\{1,\dots,r\}$
tel que $s_t\not\in\{\pm 1\}$. Quitte à permuter les $s_i$ (ce qui
remplace $(\H,s,\eta_0)$ par un triplet équivalent), on peut supposer que
$t=r$; on note $n'=n_1+\dots+n_{r-1}$.
Alors le système de coracines $\widehat{\Phi}_H$ de $\H$, qui est
l'ensemble des coracines de $\G$ annulant $s$, vérifie 
\[\widehat{\Phi}_H\subset\{\pm\widehat{e_i}\pm\widehat{e_j},1\leq i<j\leq n'\}
\cup\{\pm\widehat{e_i},1\leq i\leq n'\}\cup\{\pm(\widehat{e_i}-\widehat{e_j}),
n'+1\leq i<j\leq n\},\]
donc
\[Z(\widehat{\H})=\bigcap_{\widehat{\alpha}\in\widehat{\Phi}_H}\Ker(\widehat{
\alpha})\supset\C^\times\times\{1\}^{n'}\times(\C^\times)^{n_r}.\]
Ceci contredit le fait que $Z(\widehat{\H})^0\subset Z(\widehat{\G})$. On en
déduit que $s_1,\dots,s_r\in\{\pm 1\}$ (donc que $r\leq 2$). Un argument
similaire (calcul de $Z(\widehat{\H})$) montre que la multiplicité de
$-1$ dans $s$ ne peut pas être égale à $1$. Ceci prouve l'assertion sur
la description des classes d'équivalence de $s$ possibles.

Soit $n_1\in\{1,\dots,n\}-\{n-1\}$, et soit $\H$ le groupe endoscopique
associé comme dans l'énoncé. Alors le système de coracines de $\H$ est
\[\widehat{\Phi}_H=\{\pm\widehat{e_i}\pm\widehat{e_j},1\leq i<j\leq n_1\}
\cup\{\pm\widehat{e_i},1\leq i\leq n_1\}\cup\{\pm\widehat{e_i}\pm\widehat{
e_j},n_1+1\leq i<j\leq n\},\]
donc le système de racines de $\H$ est isomorphe à
\[\{\pm(e_i-\frac{c}{2})\pm(e_j-\frac{c}{2}),1\leq i<j\leq n_1\}\cup
\{\pm(2e_i-c),1\leq i\leq n_1\}\cup
\{\pm(e_i-\frac{c}{2})\pm(e_j-\frac{c}{2}),n_1+1\leq i<j\leq n\},\]
c'est-à-dire de type $C_{n_1}\times D_{n-n_1}$.
La formule pour $\H$ résulte facilement de ceci.
Pour calculer $\Lambda(\H,s,\eta_0)$, on remarque que, si $n\not=n_1$,
alors il existe, à translation
par $\widehat{\H}$ près, un unique élément $g\in\widehat{\G}-\widehat{\H}$
qui normalise $\widehat{\H}$ : c'est un relèvement de l'élément de
longueur maximale du groupe de Weyl de $B_n$.

Enfin, si $K=\Q$ ou
$\R$, $\H$ admet un
$\R$-tore maximal elliptique si et seulement s'il existe un élément de son
groupe de Weyl qui agit par $-1$ sur les éléments d'un ensemble de racines
simples. Le système de racines $C_{n_1}$ vérifie cette propriété pour tout
$n_1$, mais $D_{n-n_1}$ ne vérifie cette propriété que si $n-n_1$ est pair. 

\end{preuve}

\subsubsection*{Un calcul de nombres de Tamagawa}

\begin{sublemme}\label{lemme:Tamagawa}\begin{itemize}
\item[(i)] Soit $\G$ un groupe réductif connexe déployé sur $\Q$. Alors
$\Ker^1(\Q,\G)=\{1\}$.
\item[(ii)] On a
\[\tau(\G(\Sp_{2n_1}\times\SO_{2n_2}))=\left\{\begin{array}{ll}1 & \mbox{ si }
n_2=0 \\
2 & \mbox{ si }n_2\geq 2\end{array}\right..\]
\item[(iii)] Soient $F$ une extension finie de $\Q$ et $\Le=R_{F/\Q}\GL_{n,F}$,
avec $n\in\Nat^*$. Alors $\tau(\Le)=1$.

\end{itemize}
\end{sublemme}

\begin{preuve}
Rappelons d'abord que, d'après \cite{K-STF:CTT} 4.2.2 et 5.1.1, \cite{K-TN} et
\cite{C},
pour tout groupe algébrique $\G$ réductif connexe sur $\Q$, on a
\[\tau(\G)=|\pi_0(Z(\widehat{\G})^{\Gal(\overline{\Q}/\Q)})|.|\Ker^1(\Q,\G)|
^{-1}.\]

\begin{itemize}
\item[(i)] D'après \cite{K-STF:CTT} 4.2, on a une bijection canonique
$\Ker^1(\Q,\G)\fl\Ker^1(\Q,Z(\widehat{\G}))^D$, donc il suffit de
montrer que $\Ker^1(\Q,Z(\widehat{\G}))=\{1\}$. Comme $Z(\widehat{\G})$
est commutatif et que $\Gal(\overline{\Q}/\Q)$ agit trivialement sur
$Z(\widehat{\G})$, cela résulte du théorème de densité de \v{C}ebotarev.

\item[(ii)] D'après le rappel ci-dessus et (i), il suffit de calculer
$|\pi_0(Z(\widehat{\G(\Sp_{2n_1}\times\SO_{n_2})}))|$. En raisonnant comme
dans la preuve de la proposition \ref{prop:groupes_endoscopiques_GSp}, on
voit que
\[Z(\widehat{\G(\Sp_{2n_1}\times\SO_{n_2})})\simeq\left\{\begin{array}{ll}
\C^\times & \mbox{ si }n_2=0 \\
\C^\times\times\{\pm 1\} & \mbox{ si }n_2\geq 2\end{array}\right..\]
La conclusion de (ii) en résulte.

\item[(iii)] C'est le (ii) du lemme 2.3.3 de \cite{M3}.

\end{itemize}
\end{preuve}

\subsection{Sous-groupes de Levi et groupes endoscopiques}
\label{endoscopie2}

Dans ce paragraphe, on rappelle quelques notions de la section 7 de
\cite{K-NP}. On utilise les notations et les définitions de la section 7
de \cite{K-STF:CTT}.

Soit $\G$ un groupe réductif connexe sur un corps local ou global $F$.
On note $\Ell(\G)$ l'ensemble
des classes d'isomorphisme de triplets endoscopiques elliptiques de $\G$
(au sens de \cite{K-STF:CTT} 7.4) et $\Levi(\G)$ l'ensemble des
classes de $\G(F)$-conjugaison de sous-groupes de Levi de $\G$.
Soit $\M$ un sous-groupe de Levi de
$\G$. On a un plongement canonique $\Gal(\overline{F}/F)$-équivariant
$Z(\widehat{\G})\fl Z(\widehat{\M})$.

\begin{subdefinition}(\cite{K-NP} 7.1) Un \emph{$\G$-triplet endoscopique}
de $\M$
est un triplet endoscopique $(\M',s_M,\eta_{M,0})$ de $\M$ tel que :
\begin{itemize}
\item[(i)] l'image de $s_M$ dans $Z(\widehat{\M'})/Z(\widehat{\G})$ soit
fixe par $\Gal(\overline{F}/F)$;
\item[(ii)] l'image de $s_M$ dans $\Ho^1(F,Z(\widehat{\G}))$ (via le
morphisme $(Z(\widehat{\M'})/Z(\widehat{\G}))^{\Gal(\overline{F}/F)}
\fl\Ho^1(F,Z(\widehat{\G}))$ induit par la suite exacte
$1\fl Z(\widehat{\G})\fl Z(\widehat{\M})\fl Z(\widehat{\M})/Z(\widehat{\G})\fl
1$) est dans $\Ker^1(F,
Z(\widehat{\G}))$.

\end{itemize}
On dit que $(\M',s_M,\eta_{M,0})$ est elliptique
s'il est elliptique en tant que triplet
endoscopique de $\M$. 

Soient $(\M'_1,s_1,\eta_{1,0})$ et $(\M'_2,s_2,\eta_{2,0})$ deux $\G$-triplets
endoscopiques de $\M$. Un \emph{isomorphisme de $\G$-triplets endoscopiques}
de $(\M'_1,s_1,\eta_{1,0})$ sur $(\M'_2,s_2,\eta_{2,0})$ est un isomorphisme
$\alpha:\M'_1\fl\M'_2$ de triplets endoscopiques de $\M$ (au sens de
\cite{K-STF:CTT} 7.5) tel que les images de $s_1$ et $\widehat{\alpha}(s_2)$
dans $Z(\widehat{\M'_1})/Z(\widehat{\G})$ soient égales.

\end{subdefinition}

Dans \cite{K-NP} 3.7 et 7.4, Kottwitz explique comment associer
une classe d'isomorphisme de triplets endoscopiques de $\G$
à un $\G$-triplet endoscopique de $\M$ (cette construction
est rappelée dans \cite{M3} 2.4).
On note $\Ell_\G(\M)$ l'ensemble des classes d'isomorphisme de $\G$-triplets
endoscopiques elliptiques de $\M$ telles que la classe d'isomorphisme de
triplets endoscopiques de $\G$ associée soit elliptique.
On a des applications évidentes
$\Ell_\G(\M)\fl\Ell(\M)$ et $\Ell_\G(\M)\fl\Ell(\G)$.
Pour tout $(\M',s_M,\eta_{M,0})\in\Ell_\G(\M)$, on
note $\Aut(\M',s_M,\eta_{M,0})$ le groupe des $\G$-automorphismes
de $(\M',s_M,\eta_{M,0})$ et $\Lambda_\G(\M',s_M,\eta_{M,0})=
\Aut(\M',s_M,\eta_{M,0})/\M'_{ad}(\Q)$ le groupe des $\G$-automorphismes
extérieurs; si $\M=\G$, on omet l'indice $\G$.

Rappelons que l'on note $n_M^G=|\Nor_\G(\M)(\Q)/\M(\Q)|$.

Le lemme \ref{lemme:Levi_endoscopiques} ci-dessous est un cas particulier du
lemme 7.2 de \cite{K-NP}. Comme \cite{K-NP} n'est pas publié, on le prouve
directement par un calcul pour les groupes considérés dans ce texte.
On suppose désormais que $\G=\GSp_{2n}$, $n\in\Nat$.
Si $(\M',s_M,\eta_{M,0})\in\Ell_\G(\M)$
et $(\H,s,\eta_0)$ est son image dans $\Ell(\G)$, alors on voit
facilement que $\M'$ détermine une classe de $\H(\Q)$-conjugaison de
sous-groupes de Levi de $\H$. (Cette remarque est vraie en général et
prouvée dans \cite{K-NP}, mais
on n'a pas besoin de ce fait; cf la note 3 de \cite{M3} 2.4.)

\begin{sublemme}\label{lemme:Levi_endoscopiques}
Soit $\varphi:\coprod\limits
_{(\H,s,\eta_0)\in\Ell(\G)}\Levi(\H)\fl\C$. Alors
\begin{flushleft}$\displaystyle{\sum_{(\H,s,\eta_0)\in\Ell(\G)}|\Lambda(\H,s,
\eta_0)|^{-1}\sum_{\M_H\in\Levi(\H)}(n_{M_H}^H)^{-1}\varphi(\H,\M_H)}$
\end{flushleft}
\begin{flushright}$\displaystyle{=\sum_{\M\in\Levi(\G)}(n_M^G)^{-1}
\sum_{(\M',s_M,\eta_{M,0})\in\Ell_\G(\M)}|\Lambda_G(\M',s_M,\eta_{M,0})|^{-1}
\varphi(\H,\M_H),}$\end{flushright}
où, dans la deuxième somme, $(\H,s,\eta_0)$ est l'image de $(\M',s_M,\eta_
{M,0})$ dans $\Ell(\G)$ et $\M_H$ est l'élément de $\Levi(\H)$ associé
à $(\M',s_M,\eta_{M,0})$.

\end{sublemme}

(On écrira parfois $\M'$ au lieu de $\M_H$, car $\M'$ et $\M_H$ sont
isomorphes.)

Nous n'utiliserons ce lemme que pour des fonctions $\varphi$ qui s'annulent dès
que leur deuxième argument n'est pas un sous-groupe de Levi cuspidal.
Dans ce cas, le lemme résulte facilement du lemme
ci-dessous, qui se prouve comme la proposition
\ref{prop:groupes_endoscopiques_GSp}.
(Il serait facile aussi mais plus fastidieux de montrer le lemme
pour $\varphi$ quelconque par un calcul explicite.)

\begin{sublemme}\label{lemme:Levi_endoscopiques_GSp} Soit $n\in\Nat^*$; on
note $\G=\GSp_{2n}$. On note $\T$ le tore diagonal de $\G$, et on identifie
$\widehat{\T}$ à $\C^\times\times(\C^\times)^n$ comme dans \ref{endoscopie1}.
Soit $\M$ un sous-groupe de Levi cuspidal de $\G$. Alors $\M$ est isomorphe
à $\Gr_m^r\times\GL_2^t\times\GSp_{2m}$, avec $m,r,t\in\Nat$ tels que $n=m+r+
2t$. Soit $\T_M$ le tore diagonal de $\M$.
On choisit un isomorphisme $\widehat{\T}_M\simeq\widehat{\T}$ tel que
l'ensemble des coracines du facteur $\GSp_{2m}$ de $\M$ soit
$\{\pm\widehat{e_i}\pm\widehat{e_j},r+2t+1\leq i<j\leq n\}\cup
\{\pm\widehat{e_i},r+2t+1\leq i\leq n\}$ et celui du facteur
$\GL_2^t$ soit $\{\pm(\widehat{e}_{r+2i-1}-\widehat{e}_{r+2i}),1\leq i\leq t\}$
(les notations sont celles de \ref{endoscopie1}).

Alors un élément $(\M',s_M,\eta_{M,0})$ de $\Ell_\G(\M)$ est
uniquement déterminé par la donnée de $s_M$ et si on suppose (comme on peut
le faire) que $s_M\in\widehat{\T}_M\simeq\widehat{\T}$, on a nécessairement
$s_M\in Z(\widehat{\G})(\{1\}\times\{\pm 1\}^n)$.
Pour tous $A\subset\{1,\dots,r\}$, $B\subset\{1,\dots,t\}$ et $m_1,m_2\in
\Nat$ tels que $m=m_1+m_2$ et $m_2\not=1$, on note
\[s_{A,B,m_1,m_2}=(s_1,\dots,s_r,s'_1,\dots,s'_{2t},\overbrace{1,\dots,1}^
{m_1},\overbrace{-1,\dots,-1}^{m_2}),\]
avec $s_i=-1$ si $i\in A$ et $s_i=1$ si $i\not\in A$, et $s'_{2i-1}=s'_{2i}=-1$
si $i\in B$ et $s'_{2i-1}=s'_{2i}=1$ si $i\not\in B$. Alors l'ensemble des
$(1,s_{A,B,m_1,m_2})$ est
un ensemble de représentants des classes
d'équivalence de $s_M$ possibles.

Soit $s_M=(1,(s_1,\dots,s_n))\in\{1\}\times\{\pm 1\}^n$ 
dans l'ensemble de représentants ci-dessus. Soit $(\M',s_M,\eta_{M,0})$
l'élément de $\Ell_\G(\M)$ associé à $s_M$, et $(\H,s,\eta_0)$
son image dans $\Ell(\G)$. Soient $n_1=|\{i\in\{1,\dots,n\}|s_i=1\}$,
$n_2=n-n_1$, $m_1=|\{i\in\{r+2t+1,\dots,n\}|s_i=1\}|$, $m_2=m-m_1$, $r_1=
|\{i\in\{1,\dots,r\}|s_i=1\}|$,$r_2=r-r_1$, $t=\frac{1}{2}|\{i\in\{r+1,\dots,
r+2t+1\}|s_i=1\}|$, $t_2=t-t_1$ (grâce aux hypothèses sur $s_M$,
$t_1$ et $t_2$ sont entiers et $n_2$ et $m_2$ différents de $1$).
Alors $\H=\G(\Sp_{2n_1}\times\SO_{2n_2})$, $\M'=\Gr_m^r\times\GL_2^t\times\G(
\Sp_{2m_1}\times\SO_{2m_2})$, et $n_{M'}^H=2^{r+t}(r_1)!(t_1)!(r_2)!(t_2)!$.
Enfin, $|\Lambda_G(\M',s_M,\eta_{M,0})|=1$ si $\M\not=\G$.

\end{sublemme}

On finit cette section en rappelant le résultat de \cite{K-NP} 7.3.
On suppose à nouveau que $\G$ est un groupe réductif connexe sur un corps
local ou global $F$. Soit $\M$ un sous-groupe de Levi de $\G$.

\begin{subdefinition} Soit $\gamma\in\M(F)$ semi-simple.
Un \emph{$\G$-quadruplet endoscopique} de $(\M,\gamma)$ est
un quadruplet $(\M',s_M,\eta_{M,0},\gamma')$, où $(\M',s_M,\eta_{M,0})$ est
un $\G$-triplet endoscopique de $\M$ et $\gamma'\in\M'(F)$ est un élément
semi-simple $(\M,\M')$-régulier tel que $\gamma$ soit une image de
$\gamma'$ (cf \cite{K-STF:EST} 3 pour la définition de ces termes).
Un \emph{isomorphisme de $\G$-quadruplets endoscopiques}
$\alpha:(\M'_1,s_{M,1},\eta_
{M,0,1},\gamma'_1)\fl(\M'_2,s_{M,2},\eta_{M,0,2},\gamma'_2)$ est un
isomorphisme de $\G$-triplets endoscopiques $\alpha:\M'_1\fl\M'_2$ tel que
$\alpha(\gamma'_1)$ et $\gamma'_2$ soient stablement conjugués.

\end{subdefinition}

Soit $I$ un sous-groupe réductif connexe de $\G$ qui contient un tore
maximal de $\G$. On a une inclusion
canonique $\Gal(\overline{F}/F)$-équivariante $Z(\widehat{\G})\subset
Z(\widehat{I})$. On note $\Kgoth_\G(I/F)$ l'ensemble des éléments de
$(Z(\widehat{I})/Z(\widehat{\G}))^{\Gal(\overline{F}/F)}$ dont l'image par
le morphisme $(Z(\widehat{I})/Z(\widehat{\G}))^{\Gal(\overline{F}/F)}\fl
\H^1(F,Z(\widehat{\G}))$ est triviale si $F$ est
local, dans $\Ker^1(F,Z(\widehat{\G}))$ si $F$ est global.

\begin{subremarque}\label{rq:Kgoth_G_M}(\cite{K-NP} (7.2.1))
Si $I$ est inclus dans $\M$, alors
on a une suite exacte
\[1\fl(Z(\widehat{\M})/Z(\widehat{\G}))^{\Gal(\overline{F}/F)}\fl\Kgoth_\G
(I/F)\fl\Kgoth_\M(I/F)\fl 1.\]

\end{subremarque}

Cette remarque est démontrée dans la preuve du lemme 6.3.4 de \cite{M3}.

On fixe $\gamma\in\M(F)$ semi-simple, et on note $I=\M_\gamma$.
Soit $(\M',s_M,\eta_{M,0},\gamma')$ un $\G$-quadruplet endoscopique de
$(\M,\gamma)$. On note $I'=\M'_{\gamma'}$. Comme $\gamma'$ est
$(\M,\M')$-régulier, $I'$ est une forme intérieure de $I$ (cf
\cite{K-STF:EST} 3), donc on a un isomorphisme canonique $Z(\widehat{I})\simeq
Z(\widehat{I'})$. On note $\kappa(\M',s_M,\eta_{M,0},\gamma')$ l'image de
$s_M$ par le morphisme $Z(\widehat{\M'})\subset Z(\widehat{I'})\simeq
Z(\widehat{I})$.

\begin{sublemme}\label{lemme:param_groupes_endoscopiques} L'application
$(\M',s_M,\eta_{M,0},\gamma')\fle\kappa(\M',s_M,\eta_{M,0},\gamma')$ induit
une bijection de l'ensemble des classes d'isomorphisme de $\G$-quadruplets
endoscopiques de $(\M,\gamma)$ sur $\Kgoth_\G(I/F)$. De plus, tout
automorphisme d'un $\G$-quadruplet endoscopique de $(\M,\gamma)$ est
intérieur.

\end{sublemme}

Ce lemme est le lemme 7.1 de \cite{K-NP}. Il s'agit d'une généralisation du
lemme 9.7 de \cite{K-STF:EST}, et il se prouve de la même manière que ce
lemme.

\section{Calculs en la place infinie}
\label{infini}

\subsection{Notations et rappels}
\label{infini1}

Soit $\G$ un groupe algébrique réductif connexe sur $\R$.
Dans cette 
section, on forme les $L$-groupes en utilisant le groupe de Weil $W_\R$.
Rappelons que $W_\R=W_\C\sqcup W_\C\tau$, avec $W_\C=\C^\times$,
$\tau^2=-1$ 
et, pour tout $z\in\C^\times$, $\tau z\tau^{-1}=\overline{z}$, et que
$W_\R$ agit sur $\widehat{\G}$ via son quotient 
$\Gal(\C/\R)\simeq W_\R/W_\C$.
On note $\Pi(\G(\R))$ (resp. $\Pi_{temp}(\G(\R))$) l'ensemble des classes 
d'équivalence de représentations admissibles
irréductibles (resp. et tempérées) de 
$\G(\R)$. Pour tout $\pi\in\Pi(\G(\R))$, on note
$\Theta_\pi$ le caractère de Harish-Chandra de $\pi$ (c'est une fonction
analytique réelle sur $\G_{reg}(\R)$).

On suppose maintenant que $\G(\R)$ a une série discrète. Soient $\A_G$
le sous-tore déployé (sur $\R$) maximal du centre de $\G$ et $\overline{\G}$
la forme intérieure de $\G$ telle que $\overline{\G}/\A_G$ soit anisotrope
sur $\R$.
On pose $q(G)=\dim(X)/2$, où $X$ est l'espace symétrique de $\G(\R)$.
On note $\Pi_{disc}(\G(\R))\subset\Pi(\G(\R))$
l'ensemble des classes d'équivalence de représentations de la série discrète.

L'ensemble $\Pi_{disc}(\G(\R))$ est réunion disjointe de sous-ensembles finis
tous de même cardinal, qui sont appelés $L$-paquets et paramétrés par les
classes d'équivalence de paramètres de Langlands elliptiques $\varphi:W_\R\fl
{}^L\G$ ou, ce qui revient au même, par les représentations irréductibles de
$\overline{\G}(\R)$. On note $\Pi(\varphi)$ (resp. $\Pi(E)$) le $L$-paquet 
associé au paramètre $\varphi$ (resp. à la représentation irréductible $E$
de $\overline{\G}(\R)$), et
$d(\G)$ le cardinal des $L$-paquets de $\Pi_{disc}(\G)$. Rappelons que, si $E$
est une représentation irréductible de $\overline{\G}(\R)$, alors $\Pi(E)$
est l'ensemble des éléments de $\Pi_{disc}(\G(\R))$ qui ont le même
caractère central et le même caractère infinitésimal que $E$.

Soit $\pi\in\Pi_{disc}(\G(\R))$. 
On note $f_\pi$ un pseudo-coefficient de $\pi$ (cf \cite{CD}).

Pour tout paramètre de Langlands elliptique $\varphi:W_\R\fl{}^L\G$, on note
\[S\Theta_\varphi=\sum_{\pi\in\Pi(\varphi)}\Theta_\pi.\]

On va maintenant calculer l'entier $d(\G)$ dans le cas
des groupes symplectiques.

Soit $n\in\Nat^*$. On note $\G=\GSp_{2n}$ et $\T$ le tore diagonal de $\G$.
On a $\A_G=\Gr_mI_{2n}$.
Soit
\[\T_{ell}=\left\{\left(\begin{array}{cc}
\begin{array}{ccc}a_1 & & 0 \\ & \ddots & \\ 0 & & a_n\end{array} 
& \begin{array}{ccc}0 & & b_1 \\ & \ddotsinv & \\ b_n & & 0 \end{array} \\
\begin{array}{ccc}0 & & - b_n \\ & \ddotsinv & \\ -b_1 & & 0 \end{array}
& \begin{array}{ccc}a_n & & 0 \\ & \ddots & \\ 0 & & a_1\end{array}
\end{array}\right),a_1^2+b_1^2=\dots=a_n^2+b_n^2\not=0\right\}.\]
Alors $\T_{ell}$ est un tore maximal de $\G$, et $\T_{ell}/\A_G$ est anisotrope
sur $\R$. Notons
\[u_G=\frac{1}{\sqrt{2}}\left(\begin{array}{cc}I_n & iA_n \\
iA_n &  I_n\end{array}
\right)\in\Sp_{2n}(\C).\]
La conjugaison par $u_G^{-1}$ induit un isomorphisme $\alpha:\T_{ell,\C}\iso
\T_\C$. On utilise $\alpha$ pour identifier $\widehat{\T}_{ell}$ et
$\widehat{\T}$.
On note $\Omega_\G=W(\T_{ell}(\C),\G(\C))$ et $\Omega_{\G(\R)}=W(\T_{ell}(\R),
\G(\R))$
les groupes de Weyl de $\T_{ell}$ sur $\C$ et sur $\R$. On a
$\Omega_\G\simeq W(\T(\C),\G(\C))\simeq\{\pm 1\}^n\rtimes\Sgoth_n$, où
$\Sgoth_n$ agit sur $\{\pm 1\}^n$ par $(\sigma,(\varepsilon_1,\dots,\varepsilon
_n))\fle(\varepsilon_{\sigma^{-1}(1)},\dots,\varepsilon_{\sigma^{-1}(n)})$.
Le sous-groupe $\Omega_{\G(\R)}$ de $\Omega_{\G}$ est le groupe engendré
par $(-1,\dots,-1)\in\{\pm 1\}^n$ et par $\Sgoth_n$,
donc $d(\G)=2^{n-1}$.

\begin{subremarque}\label{rq:T_ell}
Soit $E=\Q[\sqrt{-1}]$. On note $\U(1)$ le groupe
des éléments de norme $1$ de $E$ et $\G(\U(1)^n)=\{(x_1,\dots,x_n)\in
E^\times,\overline{x}_1x_1=\dots=\overline{x}_nx_n\}$.
Alors le tore $\T_{ell}$ est isomorphe à
$\G(\U(1)^n)$ par le morphisme
\[\beta_\G:\left(\begin{array}{cc}
diag(a_1,\dots,a_n) & diag(b_1,\dots,b_n)A_n \\
-A_n diag (b_1,\dots,b_n) & diag(a_n,\dots,a_1)
\end{array}\right)\fle (a_1+ib_1,\dots,a_n+ib_n).\]

\end{subremarque}

On aura aussi besoin de connaître un tore maximal elliptique d'un groupe
orthogonal.
Soit $n\in\Nat^*$. On note $\G=\GSO_{2n}$ et $\T$ le tore diagonal de $\G$.
On a $\A_G=\Gr_m I_{2n}$. Si $n$ est impair, alors $\G_\R$ n'a pas de tore
maximal elliptique. On suppose à partir de maintenant que $n$ est pair, et
on note $q=n/2$.
On note $\T_{ell}$ le sous-groupe de $\G$ formé des matrices
\[g=\left(\begin{array}{cccccc}
\begin{array}{cc}a_1 & b_1 \\ -b_1 & a_1\end{array} & & 0 & 0 & &
\begin{array}{cc}d_1 & c_1 \\ c_1 & -d_1\end{array} \\
& \ddots & & & \ddotsinv & \\
0 & & \begin{array}{cc}a_q & b_q \\ -b_q & a_q\end{array} &
\begin{array}{cc}d_q & c_q \\ c_q & -d_q\end{array} & & 0 \\
0 & & \begin{array}{cc}-d_q & c_q \\ c_q & d_q\end{array} &
\begin{array}{cc} a_q & -b_q \\ b_q & a_q\end{array} & & 0 \\
& \ddotsinv & & & \ddots  & \\
\begin{array}{cc}-d_1 & c_1 \\ c_1 & d_1\end{array} & & 0 & 0 & &
\begin{array}{cc}a_1 & -b_1 \\ b_1 & a_1\end{array} 
\end{array}\right),\]
avec $a_r,b_r,c_r,d_r\in\Gr_a$ tels que $a_rc_r+b_rd_r=0$ pour tout $r\in
\{1,\dots,q\}$ et $a_1^2+b_1^2+c_1^2+d_1^2=\dots=a_q^2+b_q^2+c_q^2+d_q^2\not=
0$ (on a $c(g)=a_1^2+b_1^2+c_1^2+d_1^2$).
Alors $\T_{ell}$ est un tore maximal de $\G$, et $\T_{ell}/\A_G$ est anisotrope
sur $\R$ (donc $\G$ est cuspidal). Notons
\[u_G=\frac{1}{2}
\left(\begin{array}{cccccc}
\begin{array}{cc}1 & i \\ i & -1\end{array} & & 0 & 0 & &
\begin{array}{cc}i & 1 \\ 1 & -i\end{array} \\
& \ddots & & & \ddotsinv & \\
0 & & \begin{array}{cc}1 & i \\ i & -1\end{array} &
\begin{array}{cc} i & 1 \\ 1 & -i\end{array} & & 0 \\
0 & & \begin{array}{cc} i & 1 \\ 1 & -i\end{array} &
\begin{array}{cc} -1 & -i \\ -i & 1\end{array} & & 0 \\
& \ddotsinv & & & \ddots  & \\
\begin{array}{cc} i & 1 \\ 1 & -i\end{array} & & 0 & 0 & &
\begin{array}{cc} -1 & -i \\ -i & 1\end{array} 
\end{array}\right)\in\SO_{2n}(\C).\]
La conjugaison par $u_G^{-1}$ induit un isomorphisme $\alpha:\T_{ell,\C}\iso
\T_\C$.
On utilise $\alpha$ pour identifier $\widehat{\T}_{ell}$ et
$\widehat{\T}$.

\begin{subremarque} Le tore $\T_{ell}$ est isomorphe à $\G(\U(1)^n)$ (ce groupe
est défini dans la remarque \ref{rq:T_ell})
par le morphisme $\beta_G$ qui envoie la matrice $g$ ci-dessus sur
\[(a_1+ib_1+c_1+id_1,a_1+ib_1-c_1-id_1,\dots,a_q+ib_q+c_q+id_q,a_q+ib_q-c_q
-id_q).\]

\end{subremarque}

\begin{subremarque}\label{rq:T_Sp_SO}
Soient $n_1,n_2\in\Nat$ et $\G=\G(\Sp_{2n_1}\times\SO_{2n
_2})$. Alors $\G_\R$ ne peut avoir un tore maximal elliptique (sur $\R$) que
si $n_2$ est pair. On suppose que $n_2$ est pair.
En utilisant les exemples ci-dessus, on construit facilement un
tore maximal elliptique $\T_{ell}$ de $\G$ (en particulier, $\G$ est cuspidal)
et un élément $u_G\in\G_{der}(\C)$
tel que $\Int(u_G^{-1})$ envoie $\T_{ell}$ sur le tore diagonal.
On obtient
aussi un isomorphisme $\beta_G:\T_{ell}\iso\G(\U(1)^{n_1+n_2})$.

\end{subremarque}

\vspace{1cm}

Rappelons une construction d'Arthur et Shelstad.

Soit $\G$ un groupe réductif connexe sur $\R$. Un \emph{caractère virtuel} 
$\Theta$
sur $\G(\R)$ est une combinaison linéaire à coefficients dans $\Z$ de fonctions
$\Theta_\pi$, $\pi\in\Pi(\G(\R))$. On dit que $\Theta$ est \emph{stable} si
$\Theta(\gamma)=\Theta(\gamma')$ pour tous $\gamma,\gamma'\in\G_{reg}(\R)$
stablement conjugués.

Soit $\T$ un tore maximal de $\G$. On note $\A$ le sous-tore déployé maximal
de $\T$ et $\M=\Cent_\G(\A)$ (un sous-groupe de Levi de $\G$). Pour tout 
$\gamma\in\M(\R)$, on note (comme dans \ref{formule_points_fixes2})
\[D_M^G(\gamma)=\det(1-\Ad(\gamma),Lie(\G)/Lie(\M)).\]

\begin{sublemme}\label{lemme:Phi_M}
(\cite{A-L2} 4.1, \cite{GKM} 4.1) Soit $\Theta$ un caractère
virtuel stable sur $\G(\R)$. Alors la fonction
\[\gamma\fle |D_M^G(\gamma)|^{1/2}\Theta(\gamma)\]
sur $\T_{reg}(\R)$ s'étend en une fonction continue sur $\T(\R)$, qu'on
notera $\Phi_M(.,\Theta)$ ou $\Phi_M^G(.,\Theta)$.

\end{sublemme}

Dans la suite, on considérera parfois $\Phi_M(.,\Theta)$ comme une fonction
sur $\M(\R)$ définie de la manière suivante : si $\gamma\in\M(\R)$ est conjugué
à $\gamma'\in\T(\R)$, alors $\Phi_M(\gamma,\Theta)=\Phi_M(\gamma',\Theta)$; si
$\gamma\in\M(\R)$ n'a aucun conjugué dans $\T(\R)$, alors
$\Phi_M(\gamma,\Theta)=0$.

\begin{subremarque}\label{rq:Phi_M}
La fonction $\Phi_M(.,\Theta)$ sur $\M(\R)$ est invariante par 
conjugaison par $\Nor_\G(\M)(\R)$ (car $\Theta$ et $D_M^G$ le sont).

\end{subremarque}

Soit $V$ une représentation algébrique irréductible de $\G_\C$. On peut
voir $V$ comme une représentation irréductible de $\overline{\G}(\R)$, donc
on a un $L$-paquet $\Pi(V)$ associé à $V$. On note $\varphi$ un paramètre de
Langlands de $\Pi(V)$.
On va rappeler la formule pour $\Phi_M^G(.,S\Theta_\varphi)$ qui est donnée
dans \cite{GKM} \S4.

On rappelle que, si $(X,X^*,R,R^\vee)$ est un système de racines tel que $R$
engendre $X$ et que $-1$ soit dans le groupe de Weyl de $R$, alors
on a une fonction $\overline{c}_R:X_{reg}\times X^*_{reg}\fl\Z$, dont les
propriétés sont données (par exemple) dans \cite{GKM} \S3. On a noté ici
$X_{reg}$ (resp. $X^*_{reg}$)
l'ensemble des éléments réguliers de $X$ (resp. $X^*$).
\footnote{La définition d'un élément régulier de $X^*$ est celle de la section
3 de \cite{GKM}.}

On note $\Bo(\T)$ l'ensemble des sous-groupes de Borel de $\G_\C$ contenant
$\T$, et on fixe $\B\in\Bo(\T)$. Pour tout $\B'\in\Bo(\T)$, on pose
$\rho_{\B'}=\frac{1}{2}\sum\limits_{\alpha\in\Phi(\T,\B')}\alpha\in X^*(\T)
\otimes_\Z\Q$, $\Delta_{\B'}=\prod\limits_{\alpha\in\Phi(\T,\B')}(1-\alpha^{-1
})$, et on note $\lambda_{\B'}$ le plus haut poids de $V$ relativement à $\B'$.
On note $\Phi^+=\Phi(\T,\B)$.
Soit $\Omega$ le groupe de Weyl de $\T(\C)$ dans
$\G(\C)$. Pour tout $\omega\in\Omega$, on note $\Phi(\omega)=\Phi^+\cap
(-\omega\Phi^+)$, $\ell(\omega)=|\Phi(\omega)|$ la longueur
de $\omega$, $\varepsilon(\omega)=(-1)^{\ell(\omega)}$.
Le groupe $\Omega$ agit simplement transitivement sur $\Bo(\T)$. Si $\B'\in
\Bo(\T)$ est égal à $\omega\B$, alors $\rho_{\B'}=\omega\rho_\B$ et
$\lambda_{\B'}=\omega\lambda_\B$.

Soit $\T(\R)_1$ le sous-groupe compact maximal de $\T(\R)$. On a
$\T(\R)\simeq\A(\R)^0\times\T(\R)_1$. On note $p:X^*(\T)\otimes_\Z\R\fl
X^*(\A)\otimes_\Z\R$ le morphisme induit par la restriction. Comme $V$ est
irréductible, le tore central $\A_G$ de $\G$ agit sur $V$ par un caractère, que
l'on note $\lambda_0$. Le $\Q$-espace vectoriel $X^*(\A_G)\otimes_\Z\Q$ est de
manière naturelle un facteur direct de $X^*(\T)\otimes_\Z\Q$, donc on
peut aussi voir $\lambda_0$ comme un élément de $X^*(\T)\otimes_\Z\Q$.

On note $R$ l'ensemble des racines réelles dans $\Phi$.
Soit $\gamma\in\T(\R)$.
On note $R_\gamma=\{\alpha\in R|\alpha(\gamma)>0\}$, $R^+_\gamma=\{\alpha\in
R|\alpha(\gamma)>1\}$ et $\varepsilon_R(\gamma)=(-1)^{\Phi^+\cap(-R_\gamma^+)
}$. On écrit
$\gamma=\exp(x)\gamma_1$, avec $x\in X_*(\A)\otimes_\Z\R=Lie(\A)$ et
$\gamma_1\in\T(\R)_1$, et on suppose que $\gamma\in\T_{reg}(\R)$.
On définit des entiers
$n(\gamma,\B')$, $\B'\in\Bo(\T)$, de la manière suivante : Si
$\gamma$ n'est pas dans $Z(\G)(\R)Im(\G_{sc}(\R)\fl\G(\R))$, alors
$n(\gamma,\B')=0$ pour tout $\B'\in\Bo(\T)$ ($Z(\G)$ est le centre de $\G$, et
$\G_{sc}\fl\G_{der}$ est le rêvetement universel du groupe dérivé
$\G_{der}$ de $\G$).
Sinon, alors $-1$ est dans le
groupe de Weyl de $R_\gamma$ (cf \cite{GKM} \S4), donc $R_\gamma$ donne
une fonction
\[\overline{c}=\overline{c}_{R_\gamma}:(X_*(\A/\A_G)\otimes_\Z\R)_{reg}\times
(X^*(\A/\A_G)\otimes_\Z\R)_{reg}\fl\Z,\]
et, pour tout $\B'\in\Bo(\T)$, on pose
\[n(\gamma,\B')=\overline{c}(x,p(\lambda_{\B'}+\rho_{\B'}-\lambda_0)).\]

Enfin, on note $\Pa$ le sous-groupe parabolique de $\G_\R$ de sous-groupe de
Levi $\M$ qui contient $\B$ et $\B_M=\B\cap\M$ (un sous-groupe de Borel
de $\M$).
Comme le tore maximal $\T$ de $\M$ est tel que $\T/\A$ est elliptique
(sur $\R$), $\Pa$ est défini sur $\Q$ (cf \cite{GKM} \S 5).

\begin{subfait}\label{fait:formule_Phi_M_G} On a
\[\Tr(\gamma,V)=\sum_{\B'\in\Bo(\T)}\lambda_{\B'}(\gamma)^{-1}\Delta_{\B'}
(\gamma)^{-1}=\Delta_\B(\gamma)^{-1}\sum_{\omega\in\Omega}\varepsilon(\omega)
(\omega\lambda_\B)(\gamma)\prod_{\alpha\in\Phi(\omega)}\alpha^{-1}(\gamma),\]
\[S\Theta_\varphi(\gamma)=(-1)^{q(\G)}\sum_{\B'\in\Bo(\T)}n(\gamma,\B')
\lambda_{\B'}(\gamma)\Delta_{\B'}(\gamma)^{-1}\]
et $\Phi_M^G(\gamma,S\Theta_\varphi)$ est égal à
\[(-1)^{q(\G)}\varepsilon_R(\gamma)
\delta_{\Pa(\R)}^{1/2}(\gamma)\Delta_{\B_M}(\gamma)^{-1}\sum_{\omega\in
\Omega}\varepsilon(\omega)n(\gamma,\omega\B)(\omega\lambda_\B)(\gamma)
\prod_{\alpha\in\Phi(\omega)}\alpha^{-1}(\gamma).\]

\end{subfait}

\begin{preuve} La première formule pour $\Tr(\gamma,V)$ est celle
de \cite{GKM} \S4 (c'est la formule du caractère de Weyl);
la deuxième formule s'en déduit immédiatement.
La formule pour $S\Theta_\varphi$ est celle de \cite{GKM} \S4. 
Pour en déduire la formule pour $\Phi_M^G(\gamma,S\Theta_\varphi)$, il suffit
de montrer que
\[|D_M^G(\gamma)|^{1/2}=\delta_{\Pa(\R)}^{1/2}(\gamma)\Delta_\B(\gamma)
\Delta_{\B_M}(\gamma)^{-1}(-1)^{|\Phi^+\cap(-R_\gamma^+)|}.\]
Or on a
\[D_M^G(\gamma)=\det(1-\Ad(\gamma),Lie(\G)/Lie(\M))=\prod_{\alpha\in
\Phi^+-\Phi(\T,\B_M)}(-\alpha(\gamma))(1-\alpha^{-1}(\gamma))^2,\]
donc
\[|D_M^G(\gamma)|^{1/2}=\prod_{\alpha\in\Phi^+-\Phi(\T,\B_M)}|\alpha(\gamma)|
^{1/2}|1-\alpha^{-1}(\gamma)|=\delta_{\Pa(\R)}^{1/2}(\gamma)|\Delta_\B(\gamma)
\Delta_{\B_M}(\gamma)^{-1}|.\]
Toutes les racines imaginaires de $\T$ sont dans $\Phi(\T,\M)$, donc
une racine dans $\Phi^+-\Phi(\T,\B_M)$ est soit complexe, soit réelle. Soit
$\alpha\in\Phi^+-\Phi(\T,\B_M)$. Si $\alpha$
est complexe, alors il existe $\alpha'\in
\Phi^+-\Phi(\T,\B_M)$, $\alpha'\not=\alpha$, telle que
$(1-\alpha^{-1}(\gamma))(1-{\alpha'}^{-1}(\gamma))=|1-\alpha^{-1}(\gamma)|^2
>0$. Si $\alpha$ est réelle, alors $1-\alpha^{-1}(\gamma)<0$ si et seulement
si $\alpha\in -R_\gamma^+$. Ceci finit la preuve.

\end{preuve}

\begin{subremarque}\label{rq:serie_discrete2}
Si $\G=\GSp_{2n}$, alors le groupe dérivé de $\G$ est simplement connexe, donc
$Im(\G_{sc}(\R)\fl\G(\R))$ est simplement $\G_{der}(\R)$. On voit
facilement qu'un élément $g$ de $\G(\R)$ est dans $Z(\G)(\R)\G_{der}(\R)$ si
et seulement si $c(g)>0$.

\end{subremarque}

\subsection{Transfert}
\label{infini2}

Rappelons d'abord quelques définitions de \cite{K-SVLR} \S7.

Soit $\G$ un groupe algébrique réductif connexe sur $\Q$. Pour tout tore
maximal $\T$ de $\G$, on note $\Bo_G(\T)$ l'ensemble des sous-groupes
de Borel de $\G_\C$ contenant $\T$. On suppose que
$\G$ a un tore maximal $\T_G$ tel que $(\T_G/\A_G)_\R$ soit anisotrope, et
on note $\overline{\G}$ une forme intérieure de $\G$ sur $\R$ telle que
$\overline{\G}/\A_{G,\R}$ soit anisotrope.
On note $\Omega_G=W(\T_G(\C),\G(\C))$.
Soit $\varphi:W_\R\fl{}^L\G$ un paramètre de Langlands elliptique. 

Soit $(\H,s,\eta_0)$ un triplet endoscopique elliptique de $\G$. On suppose
qu'il existe
un $L$-morphisme $\eta:{}^L\H\fl{}^L\G$ qui prolonge $\eta_0:\widehat{\H}\fl
\widehat{\G}$ (ceci est vrai si $\G_{der}$ est simplement connexe
d'après la proposition 1 de \cite{L2}), 
et on note $\Phi_H(\varphi)$ l'ensemble des classes 
d'équivalence de paramètres de Langlands $\varphi_H:W_\R\fl{}^L\H$ tels que
$\eta\circ\varphi_H$ et $\varphi$ soient équivalents.
On suppose que le tore $\T_G$ provient d'un tore maximal
$\T_H$ de $\H$, et on fixe un isomorphisme admissible $j:\T_H\iso\T_G$ (ie
obtenu comme à la fin de \cite{L3} II.4).
On note
$\Omega_H=W(\T_H(\C),\H(\C))$. On a 
$j_*(\Phi(\T_H,\H))\subset\Phi(\T_G,\G)$, donc $j$ induit une application
$j^*:\Bo_G(\T_G)\fl\Bo_H(\T_H)$ et un morphisme injectif $\Omega_H\fl
\Omega_G$, qu'on utilise pour identifier $\Omega_H$ à un sous-groupe de
$\Omega_G$.

Soient $\B\in\Bo_G(\T_G)$ et $\B_H=j^*(\B)$. On note
\[\begin{array}{rcl}\Omega_* & = & \{\omega\in\Omega_G|j^*(\omega(\B))=\B_H\}\\
& = & \{\omega\in\Omega_G|\omega^{-1}(j_*(\Phi(\T_H,\B_H)))\subset\Phi(\T_G,
\B)\}.\end{array}\]
Alors, pour tout $\omega\in\Omega_G$, il existe un unique couple
$(\omega_H,\omega_*)\in\Omega_H\times\Omega_*$ tel que $\omega=\omega_H
\omega_*$. De plus, on a une bijection $\Phi_H(\varphi)\iso\Omega_*$ : à un
élément $\varphi_H\in\Phi_H(\varphi)$, on associe l'unique $\omega_*(\varphi_H)
\in\Omega_*$ tel que $(\omega_*(\varphi_H)^{-1}\circ j,\B,\B_H)$ soit aligné
avec $\varphi_H$ (au sens de \cite{K-SVLR} \S7 p 184). 

Le sous-groupe de Borel $\B$ définit aussi un $L$-morphisme
$\eta_B:{}^L\T_G\fl{}^L\G$, unique à conjugaison par $\widehat{\G}$ près
(cf \cite{K-SVLR} p 183). Ce morphisme vérifie $\eta_B(W_\C)\subset
{}^L\G_{ad}\subset{}^L\G$ et, pour tout $z\in W_\C=\C^\times$,
$\eta_B(z)=\eta_B(z^\rho\overline{z}^{-\rho})\rtimes z$, où $\rho=\rho_B$ et
où on utilise les conventions de \cite{K-SVLR} p 183 (ou \cite{Bo} 9.1)
pour noter les morphismes $\C^\times\fl\widehat{\T}_G$.

Rappelons la normalisation
des facteurs de transfert utilisée dans \cite{K-SVLR} \S7.

\begin{subdefinition}\label{def:facteur_transfert}
Avec les notations ci-dessus, on pose, pour tout
$\gamma_H\in\T_H(\R)$,
\[\Delta_{j,B}(\gamma_H,\gamma)=(-1)^{q(\G)+q(\H)}\chi_B(\gamma)\prod_
{\alpha\in\Phi(\T_G,\B)-j_*(\Phi(\T_H,\B_H))}(1-\alpha(\gamma^{-1})),\]
où $\gamma=j(\gamma_H)$ et $\chi_B$ est le quasi-caractère de $\T_G(\R)$ 
associé
au $1$-cocyle $a:W_\R\fl\widehat{\T}_G$ tel que $\eta\circ\eta_{B_H}\circ
\widehat{j}$ et $\eta_B.a$ soient conjugués par $\widehat{\G}$.

\end{subdefinition}

\begin{subremarques}\label{rq:facteur_transfert}\begin{itemize}
\item[(1)] Soit $\varphi_H\in\Phi_H(\varphi)$ tel que $\omega_*(\varphi_H)=1$.
Quitte à conjuguer $\varphi$ (resp. $\varphi_H$) par $\widehat{\G}$ 
(resp. $\widehat{\H}$), on peut écrire $\varphi=\eta_B\circ\varphi_B$ (resp.
$\varphi_H=\eta_{B_H}\circ\varphi_{B_H}$), où $\varphi_B$ (resp. 
$\varphi_{B_H}$) est un paramètre de Langlands pour $\T_G$ (resp. $\T_H$).
On note $\chi_{\varphi,B}$ (resp. $\chi_{\varphi_H,B_H}$) le quasi-caractère
de $\T_G(\R)$ (resp. $\T_H(\R)$) associé à $\varphi_B$ (resp. $\varphi_{B_H}$).
Alors $\chi_B=\chi_{\varphi,B}(\chi_{\varphi_H,B_H}\circ j^{-1})^{-1}$.
\item[(2)] Soit $\omega\in\Omega_G$. On écrit $\omega=\omega_H\omega_*$, avec
$\omega_H\in\Omega_H$ et $\omega_*\in\Omega_*$. On note, comme dans
\ref{infini1}, $\Phi(\omega)=\Phi(\T_G,\B)\cap(-\omega\Phi(\T_G,\B))$
et $\Phi_H(\omega_H)=\Phi(\T_H,\B_H)\cap(-\omega_H\Phi(\T_H,\B_H))$. Alors
\[\chi_B\chi_{\omega(B)}^{-1}=\left(\prod_{\alpha\in j_*\Phi_H(\omega_H)}\alpha
\right)\left(\prod_{\alpha\in\Phi(\omega)}\alpha\right)^{-1}.\]
On en déduit en particulier que $\Delta_{j,\omega(B)}=
\varepsilon(\omega_*)\Delta_{j,B}$.

\end{itemize}
\end{subremarques}

Comme dans \ref{infini1}, on fixe une représentation algébrique
irréductible $V$ de $\G_\C$, et on note $\varphi$ un paramètre de Langlands
du $L$-paquet associé à $V$. On note $\lambda$ le plus haut poids de
$V$ relativement à $\B$. 
Remarquons que la restriction de $\lambda$ à $\T_G(
\R)$ est égale au quasi-caractère $\chi_{\varphi,B}$ du point (1) de la
remarque ci-dessus. 

On peut supposer que $\varphi(W_\C)\subset
\widehat{\T}_G$ et que la restriction de $\varphi$ à $W_\C=\C^\times$ est
de la forme $z\fle z^{\lambda+\rho}\overline{z}^\mu$, où $\mu\in X^*(\T_G)
\otimes_\Z\C$ est tel que $\lambda+\rho-\mu\in X^*(\T_G)$.
Comme dans la définition \ref{def:facteur_transfert} ci-dessus, on note
$a:W_\R\fl\widehat{\T}_G$ le $1$-cocycle tel que $\eta\circ\eta_{B_H}\circ
\widehat{j}$ et $\eta_B.a$ soient conjugués par $\widehat{\G}$.

La proposition ci-dessous est immédiate.

\begin{subproposition}\label{prop:Phi_H_varphi} On suppose que le
quasi-caractère $\chi_B$ de $\T_G(\R)$ associé au $1$-cocycle $a$ est
algébrique (d'après le (2) de la remarque \ref{rq:facteur_transfert} ci-dessus,
cette condition est indépendante du choix de $\B$).
Pour tout $\omega_*\in\Omega_*$, on note $V_{H,\omega_*}$ la représentation
algébrique irréductible de $\H_\C$ de plus haut poids
$(\omega_*\lambda-\chi_{\omega_*B})\circ j$ relativement à $\B_H$
et $\varphi_{H,\omega_*}$
un paramètre de Langlands du $L$-paquet de la série discrète de $\H(\R)$
associé à $V_{H,\omega_*}$. Alors $\Phi_H(\varphi)=\{\varphi_{H,\omega_*},
\omega_*\in\Omega_*\}$. Si de plus la restriction de $\chi_B$ à $\A_G$ est
triviale (par exemple, si $\eta$ envoie ${}^L\H_{ad}$ dans ${}^L\G_{ad}$),
alors $\A_G$ agit par le même caractère sur $V$ et sur les $V_{H,\omega_*}$,
$\omega_*\in\Omega_*$ (on identifie $\A_G$ à un sous-groupe de $\H$ en
utilisant $j$).

\end{subproposition}

On remarque que les conditions de la proposition
sont vérifiées pour les groupes symplectiques si on choisit les morphismes
$\eta$ comme dans \ref{endoscopie1} (sous la proposition
\ref{prop:groupes_endoscopiques_GSp}) et \ref{stabilisation2} (cf l'exemple
\ref{ex:groupes_symplectiques} ci-dessous).
En fait, comme
$(\T_G/\A_G)$ est anisotrope sur $\R$, $\chi_B$ est algébrique si et
seulement si sa restriction à $\A_G(\R)$ est algébrique (ce qui est vrai
en particulier si cette restriction est triviale).

\begin{subexemple}\label{ex:groupes_symplectiques} Soient $n\in\Nat^*$ et
$n_1,n_2\in\Nat$ tels que $n_2$ soit pair et $n=n_1+n_2$. On pose
$\G=\GSp_{2n}$, et on prend pour $(\H,s,\eta_0)$ le triplet endoscopique
elliptique de $\G$ associé à $n_1$ comme dans la proposition
\ref{prop:groupes_endoscopiques_GSp} (donc $\H=\G(\Sp_{2n_1}\times\Sp_{2n_2})
$). On choisit le prolongement évident $\eta$ de $\eta_0$.

On a défini dans
\ref{infini1} des tores
maximaux elliptiques $\T_G=\T_{G,ell}$ et $\T_H=\T_{H,ell}$ de $\G$ et $\H$,
et des isomorphismes $\beta_G:\T_G\iso\G(\U(1)^n)$ et $\beta_H:\T_H\iso\G(
\U(1)^n)$.
On prend $j=\beta_G^{-1}\circ\beta_H:\T_H\iso\T_G$. On a aussi
défini $u_G\in\G(\C)$ tel que $\Int(u_G^{-1})$ envoie $\T_{G,\C}$ sur le 
tore diagonal de $\G_\C$, et on utilise la conjugaison par $u_G$ pour
identifier $\Omega_G$ et $\{\pm 1\}^n\rtimes\Sgoth_n$.
Avec cette identification, on a
$\Omega_H=\Omega_1\times\Omega_2$, où $\Omega_1=\{\pm 1\}^{n_1}\rtimes\Sgoth_
{n_1}$ et $\Omega_2=\{(\varepsilon_1,\dots,\varepsilon_{n_2})
\in\{\pm 1\}^{n_2}|\varepsilon_1\dots\varepsilon_{n_2}=1\}\rtimes\Sgoth_{n_2}$
(le groupe $\Omega_1\times\Omega_2$
se plonge de manière évidente dans $\{\pm 1\}^n\rtimes\Sgoth_n$).

Soit
\[\B=\Int(u_G)\left(\begin{array}{ccc}* & & * \\ & \ddots & \\ 0 & & *
\end{array}\right)\in\Bo_\G(\T_G).\]
Alors $\Omega_*$ est l'ensemble des $(\varepsilon_1,\dots,\varepsilon_n)\rtimes
\sigma$ dans $\{\pm 1\}^n\rtimes\Sgoth_n$ tels que $\sigma^{-1}_{|\{1,\dots,n_1
\}}$ et $\sigma^{-1}_{|\{n_1+1,\dots,n\}}$ soient croissants et que
$\varepsilon_1=\dots=\varepsilon_{n-1}=1$. On pourrait aussi
donner une description explicite de
la bijection $\Phi_H(\varphi)\iso\Omega_*$, comme dans la remarque 3.3.3
de \cite{M3}.

Calculons le $1$-cocyle $a$ de la définition \ref{def:facteur_transfert}.
Soit $\B_H=j^*\B$. On choisit des plongements $\widehat{\T}_G\subset
\widehat{\G}$ et $\widehat{\T}_H\subset\widehat{\H}$ tels que le
diagramme suivant commute :
\[\xymatrix{\widehat{\H}\ar[r]^-{\eta_0} & \widehat{\G} \\
\widehat{\T}_H\ar[u] & \widehat{\T}_G\ar[u]\ar[l]^-{\widehat{j}}_\sim}\]
On peut supposer que les restrictions de $\eta_B$  et $\eta_{B_H}$ à
$\widehat{\T}_G$ et $\widehat{\T}_H$ sont données par ces plongements.
Les restrictions de $\eta_B$ et $\eta_{B_H}$ à $W_\C$ sont déterminées par
la définition des ces morphismes. Enfin, on a $\eta_B(\tau)=\Phi_G\times\tau$
(resp. $\eta_{B_H}(\tau)=\Phi_H\times\tau$), où $\Phi_G$ (resp. $\Phi_H$)
est un relèvement dans $\widehat{\G}$ (resp. $\widehat{\T}$) de
l'élément $(-1,\dots,-1)\rtimes 1$ (resp. $((-1,\dots,-1)\rtimes 1)\times((-1,
\dots,-1)\rtimes 1)$) du
groupe de Weyl $W(\widehat{\T}_G,\widehat{\G})\simeq\{\pm 1\}^n\rtimes
\Sgoth_n$ (resp. $W(\widehat{\T}_H,\widehat{\H})\simeq\Omega_1\times\Omega_2$).
On peut choisir $\Phi_G$ et $\Phi_H$ tels que $\eta_0(\Phi_H)=\Phi_G$. On a
alors $a(\tau)=1$, et la restriction de $a$ à $W_\C$ est $z\fle z^{\rho_{B_H}
-\rho_B}\overline{z}^{\rho_B-\rho_{B_H}}$. Donc le quasi-caractère $\chi_B$
de $\T_G(\R)$ associé à $a$ est le caractère (algébrique)
$(\rho_{B_H}\circ j^{-1})\rho_B
^{-1}$ (on remarque que $\rho_B$ et $\rho_{B_H}$ sont des caractères).
Autrement dit, on a
\[\chi_B\circ\Int(u_\G)(\lambda_1,\dots,\lambda_{2n})=(\lambda_1\dots\lambda_
{n_1})^{n_2}\lambda_{n_1+1}\dots\lambda_n.\]

\end{subexemple}

\vspace{.5cm}

On suppose à nouveau que $\G$ est quelconque.
Soit $\M$ un sous-groupe de Levi cuspidal de
$\G$, et soit $(\M',s_M,\eta_{M,0})\in
\Ell_\G(\M)$ (cf \ref{endoscopie2}) dont l'image dans $\Ell(\G)$ est
$(\H,s,\eta_0)$. On a une classe de conjugaison de
sous-groupes de Levi de $\H$ associée à $(\M',s_M,\eta_{M,0})$, et on note
$\M_H$ un élément de cette classe; on suppose que $\M_H$ est cuspidal.
On fixe des tores maximaux elliptiques
$\T_M$ et $\T_{M_H}$ de $\M$ et $\M_H$. On fixe $u_M\in\G(\C)$ et
$u_{M_H}\in\H(\C)$ tels que $u_M^{-1}\T_M(\C)u_M=\T_G(\C)$ et
$u_{M_H}^{-1}\T_{M_H}(\C)u_{M_H}=\T_H(\C)$, et on suppose qu'il existe un
isomorphisme admissible $j_M:\T_{M_H}\iso\T_M$ tel que le diagramme suivant
soit commutatif
\[\xymatrix@C=40pt{\T_{M,\C}\ar[r]^-{\Int(u_M^{-1})} & \T_{G,\C}
 \\
\T_{M_H,\C}\ar[r]_-{\Int(u_{M_H}^{-1})}\ar[u]^{j_M} & \T_{H,\C}\ar[u]_j}\]

L'isomorphisme $j_M$ induit comme
ci-dessus des applications $j_{M*}:\Phi(\T_{M_H},\H)\fl\Phi(\T_M,\G)$
et $j_M^*:\Bo_G(\T_M)\fl\Bo_H(\T_{M_H})$ (et des applications similaires
si on remplace $\G$ par $\M$ et $\H$ par $\M_H$).
On utilise la conjugaison par $u_M$ (resp. $u_{M_H}$) pour identifier 
$\Omega_G$ (resp. $\Omega_H$) et $W(\T_M(\C),\G(\C))$ (resp.
$W(\T_{M_H}(\C),\H(\C))$). Si $\B\in\Bo_G(\T_M)$, on lui associe le 
sous-ensemble $\Omega_*\subset\Omega_G$ et la bijection $\Phi_H(\varphi)\iso
\Omega_*$ définis par $\Int(u_M^{-1})(\B)\in\Bo_G(\T_G)$.

On note $R$ (resp. $R_H$) l'ensemble des racines réelles de $\Phi(\T_M,\G)$
(resp. $\Phi(\T_{M_H},\H)$). On a $j_{M*}(R_H)\subset R$.
Si $\gamma\in\T_M(\R)$ et $\gamma_H\in\T_{M_H}(\R)$, on définit des signes
$\varepsilon_R(\gamma)$ et $\varepsilon_{R_H}(\gamma_H)$ comme dans
\ref{infini1}.
Soit $\gamma\in
\T_{M,reg}(\R)$. On définit des entiers
$n_H(\gamma,\B)$, $\B\in\Bo(\T_M)$,
de manière similaire aux entiers $n(\gamma,\B)$ de \ref{infini1},
mais en utilisant $j_{M*}(R_H)$ au lieu de $R$
(si $j_M^{-1}(\gamma)$ n'est pas dans $Z(\H)(\R)Im(\H_{sc}(\R)\fl
\H(\R))$, on prend $n_H(\gamma,\B)=0$ pour tout $\B\in\Bo(\T_M)$).
On définit une fonction $S\Theta_\varphi^H$ sur $\T_{M,reg}(\R)$
en utilisant la formule du fait \ref{fait:formule_Phi_M_G} pour $S\Theta_
\varphi$ où on a remplacé $n(\gamma,\B)$ par $n_H(\gamma,\B)$.
La preuve du lemme 4.1 de \cite{GKM} s'adapte facilement à ce cas et
implique que la fonction $\Phi_M^G(.,S\Theta_\varphi^H):=|D_M^G|^{1/2}
S\Theta_\varphi^H$ se prolonge par continuité à $\T_M(\R)$.

D'après \cite{K-NP} p 23, le morphisme $\eta$ détermine un $L$-morphisme
$\eta_M:{}^L\M_H={}^L\M'\fl{}^L\M$, unique à conjugaison par $\widehat{\M}$
près, et qui prolonge $\eta_{M,0}$.
\footnote{Si $\G=\GSp_{2n}$, il existe un choix évident pour ce morphisme,
puisque $\M$ et $\M'$ sont déployés sur $\Q$.}
On utilise ce morphisme $\eta_M$ pour
définir des facteurs de transfert $\Delta_{j_M,\B_M}$, pour tout $\B_M\in\Bo_M
(\T_M)$ : pour tout $\gamma_H\in\T_{M_H}(\R)$, on pose
\[\Delta_{j_M,\B_M}(\gamma_H,\gamma)=
(-1)^{q(\G)+q(\H)}
\chi_{\B_M}(\gamma)
\prod_{\alpha\in\Phi(\T_M,\B_M)-j_{M*}(\Phi(\T_{M_H},\B_{M_H}))}
(1-\alpha(\gamma^{-1}))\]
(noter le signe),
où $\gamma=j_M(\gamma_H)$, $\B_{M_H}=j_M^*(\B_M)$ et $\chi_{B_M}$ est le 
quasi-caractère de $\T_M(\R)$ associé au $1$-cocyle $a_M:W_\R\fl\widehat{\T}_M$
tel que
$\eta_M\circ\eta_{B_{M_H}}\circ\widehat{j}_M$ et $\eta_{B_M}.a_M$ soient
conjugués par $\widehat{\M}$.

La proposition suivante est une généralisation des calculs de \cite{K-SVLR}
p 186.

\begin{subproposition}\label{prop:transfert_caracteres} On suppose que
la condition de la proposition \ref{prop:Phi_H_varphi} sur $\eta$ est vérifiée.
On fixe $\B\in\Bo_G(\T_M)$ (qui détermine $\Omega_*$
et $\Phi_H(\varphi)\iso\Omega_*$), et on note $\B_M=\B\cap\M$.
Soient $\gamma_H\in\T_{M_H}(\R)$ et $\gamma=j_M(\gamma_H)$. Alors :
\[\varepsilon_R(\gamma^{-1})\varepsilon_{R_H}(\gamma_H^{-1})
\Delta_{j_M,B_M}(\gamma_H,
\gamma)\Phi_M^G(\gamma^{-1},S\Theta_\varphi^H)=\sum_{\varphi_H\in\Phi_H
(\varphi)}\varepsilon(\omega_*(\varphi_H))\Phi_{M_H}^H(\gamma_H^{-1},S\Theta_
{\varphi_H}).\]

\end{subproposition}

\begin{subremarque}\label{rq:M_M_H_regulier} Comme $\Delta_{j_M,B_M}(\gamma_H,
\gamma)=0$ si $\gamma_H$ n'est pas $(\M,\M_H)$-régulier, le terme de droite
dans l'égalité de la proposition est non nul seulement si $\gamma_H$
est $(\M,\M_H)$-régulier.

\end{subremarque}

\begin{subremarque} Pour les groupes unitaires de \cite{M3} 2.1, on a
toujours $j_{M*}(R_H)=R$. La proposition 3.3.4 de \cite{M3} est donc un
cas particulier de la proposition ci-dessus.

\end{subremarque}

\begin{preuvep} Soient $\B_H=j_M^*\B$ et $\B_{M_H}=j_M^*\B_M$. Soit
$\Pa$ (resp. $\Pa_H$) le sous-groupe parabolique de $\G$ (resp. $\H$)
de sous-groupe de Levi $\M$ (resp. $\M_H$) et contenant $\B$ (resp.
$\B_H$).
On note $\Phi^+=\Phi(\T_M,\B)$, $\Phi_M^+=\Phi(\T_M,\B_M)$,
$\Phi_H^+=\Phi(\T_{M_H},\B_H)$ et $\Phi_{M_H}^+=\Phi(\T_{M_H},\B_{M_H})$.
Soit $\B_0=\Int(u_M)\B\in\Bo_G(\T_G)$. D'après les
hypothèses, le caractère $\chi_{B_0}$ de $\T_G(\R)$ est algébrique. On note
$\chi_B$ le caractère $\chi_{B_0}\circ j_M$ de $\T_M$, et on définit de
même un caractère $\chi_{\omega(B)}$, pour tout $\omega\in\Omega$.
Il résulte des définitions que la quasi-caractère $\chi_{B_M}$ de 
$\T_M(\R)$ est égal à $\chi_B\prod\limits_{\alpha\in\Phi^+-(\Phi_M^+\cup
\Phi_H^+)}|\alpha|^{1/2}=\chi_B\delta_{\Pa(\R)}^{1/2}(\delta_{\Pa_H(\R)}^{-1/2}
\circ j_M^{-1})$.

Comme $j_M^{-1}(\T_{M,reg}(\R))$ est dense dans $\T_{M_H}(\R)$
et que les deux côtés de l'égalité à prouver sont des fonctions continues,
on peut supposer que $\gamma$ est régulier dans $\G$. Le facteur de
transfert $\Delta_{j_M,\B_M}(\gamma_H,\gamma)$ est alors égal à
\[(-1)^{q(\G)+q(\H)}\chi_B(\gamma)
\delta_{\Pa(\R)}^{1/2}(\gamma)\delta_{\Pa_H(\R)}^{-1/2}(\gamma_H)
\Delta_{B_M}(\gamma^{-1})\Delta_{B_{M_H}}(\gamma_H^{-1})^{-1}.\]

Soient $\varphi_H\in\Phi_H(\varphi)$ et $\omega_*=\omega_*(\varphi_H)$. D'après
le fait \ref{fait:formule_Phi_M_G}, on a
\begin{flushleft}$\displaystyle{\Phi_{M_H}^H(\gamma_H^{-1},S\Theta_{\varphi_H})
=(-1)^{q(\H)}\varepsilon_{R_H}(\gamma_H^{-1})\delta_{\Pa_H(\R)}^{-1/2}(\gamma_
H)\Delta_{B_{M_H}}(\gamma_H^{-1})^{-1}}$\end{flushleft}
\begin{flushright}$\displaystyle{\sum_{\omega_H\in\Omega_H}\varepsilon(
\omega_H)n(\gamma_H^{-1},\omega_H\B_H)((\omega_H\omega_*\lambda)(\omega_H\chi_
{\omega_*B})^{-1})(\gamma^{-1})\prod_{\alpha\in\Phi_H(\omega_H)}\alpha
(\gamma_H).}$\end{flushright}

Or 
$n(\gamma_H^{-1},\omega_H\B_H)=n_H(\gamma^{-1},\omega_H\omega_*\B)$ et
\[(\omega_H\chi_{\omega_*B})(\gamma)=\chi_{\omega_H\omega_*B}(\gamma)=\chi_B
(\gamma)\prod_{\alpha\in\Phi(\omega_H\omega_*)}\alpha(\gamma)\prod_{\alpha\in
\Phi_H(\omega_H)}\alpha^{-1}(\gamma_H)\]
(cf la remarque \ref{rq:facteur_transfert}), donc
$\varepsilon_R(\gamma^{-1})\varepsilon_{R_H}(\gamma_H^{-1})
\Delta_{j_M,\B_M}(\gamma_H,\gamma)^{-1}
\varepsilon(\omega_*)\Phi_{M_H}^H(\gamma_H^{-1},S\Theta_{\varphi_H})$ est
égal à
\[(-1)^{q(\G)}\varepsilon_R(\gamma^{-1})\delta_{\Pa(\R)}^{-1/2}(\gamma)
\sum_{\omega_H\in\Omega_H}\varepsilon(\omega_H\omega_*)n_H(\gamma^{-1},
\omega_H\omega_*\B)(\omega_H\omega_*\lambda)(\gamma^{-1})\prod_{\alpha\in
\Phi(\omega_H\omega_*)}\alpha(\gamma).\]
Ceci implique le résultat cherché.

\end{preuvep}

\begin{subexemple}\label{ex:groupes_symplectiques_suite}
On utilise les notations du lemme
\ref{lemme:Levi_endoscopiques_GSp}. On a $\M\simeq\Gr_m^r\times\GL_2^t\times
\GSp_{2m}$, avec $r+2t+m=n$. Soit $\T_M$ un tore maximal elliptique
de $\M_\R$. On a $\T_M=\Gr_{m,\R}^r\times\T_1\times\dots\times\T_t\times
\T$, où $\T_1,\dots,\T_t$ sont des tores maximaux elliptiques de $\GL_{2,\R}$
et $\T$ est un tore maximal elliptique de $\GSp_{2m,\R}$.
Pour tout $i\in\{1,\dots,r\}$, on note
$e_i:\T_M\fl\Gr_{m,\R}$ la projection sur le $i$-ième facteur $\Gr_{m,\R}$.
Pour tout $j\in\{1,\dots,t\}$, on note $\alpha_j:\T_M\fl\Gr_m$ le composé de
la projection sur le facteur $\T_j$ et du morphisme $\det:\T_j\fl\Gr_{m,\R}$.
Alors l'ensemble $R$ des racines réelles de $\T_M$ dans $\GSp_{2n,\R}$ est égal
à $\{\pm(e_i-e_{i'}),1\leq i<i'\leq r\}\cup\{\pm(e_i+e_{i'}+c),1\leq i\leq i'
\leq r\}\cup\{\pm(\alpha_1+c)\}\cup\dots\cup\{\pm(\alpha_t+c)\}$.
Donc $R$ est de type $C_r\times A_1^t$.
Soit $j\in\{1,\dots,t\}$. Le quotient de $\T_j$ par le centre
$\Gr_{m,\R}I_2$ de $\GL_{2,\R}$ est anisotrope sur $\R$, donc
$(\T_j/\Gr_{m,\R}I_2)(\R)$ est connexe. On en déduit que $\alpha_j(g)>0$
pour tout $g\in\T_M(\R)$.

On suppose que $s_M$ est égal à l'élément $(1,s_{A,B,m_1,m_2})$ défini dans
le lemme \ref{lemme:Levi_endoscopiques_GSp} ($A\subset\{1,\dots,r\}$,
$B\subset\{1,\dots,t\}$, $m_1,m_2\in\Nat$ tels que $m_2\not=1$ et $m=m_1+m_2$).
Alors $\M_H\simeq\Gr_m^r\times\GL_2^t\times\G(\Sp_{2m_1}\times\SO_{2m_2})$.
Comme $\M_H$ est cuspidal, $m_2$ est forcément pair. On peut supposer que
$\T_{M_H}=\Gr_{m,\R}^r\times\T_1\times\dots\times\T_t\times\T'$, où $\T'$ est
un tore maximal elliptique de $\G(\Sp_{2m_1}\times\SO_{2m_2})$, et que
$j_M:\T_{M_H}\iso\T_M$ est l'identité sur $\Gr_{m,\R}^r\times\T_1\times\dots
\times\T_t$. On note $A_1=\{1,\dots,r\}-A$ et $A_2=A$.
Alors l'image par $j_{M*}$ de l'ensemble $R_H$ des racines
réelles de $\T_{M_H}$ dans $\H$ est
$\{\pm(e_i-e_{i'}),i<i',i,i'\in A_1\}\cup\{\pm(e_i+e_{i'}+c),i,i'\in A_1\}
\cup\{\pm(e_i-e_{i'}),i<i',i,i'\in A_2\}\cup\{\pm(e_i+e_{i'}+c),i<i',i,i'\in
A_2\}\cup\{\pm(\alpha_1+c)\}
\cup\dots\cup\{\pm(\alpha_t+c)\}$.
Donc $j_{M*}(R_H)$ est de type $C_{r_1}\times
D_{r_2}\times A_1^t$, où $r_1=|A_1|$ et $r_2=|A_2|$. On voit aussi que
$j_{M*}(R_H)=R$ si et seulement si $r_2=0$.

\end{subexemple}

On se place dans la situation de l'exemple ci-dessus. On suppose de plus
que $\M$ est standard et que
l'isomorphisme $\M\simeq\Gr_m^r\times\GL_2^t\times\GSp_{2m}$ est, à l'ordre des
facteurs dans la partie linéaire près, celui défini dans
\ref{formule_points_fixes1}.
On note $\Omega^\M$ le groupe d'automorphismes de $\M$ engendré par les
automorphismes suivants :
\begin{itemize}
\item[-] pour tout $i\in\{1,\dots,r\}$, le morphisme
\[u_i:(\lambda_1,\dots,\lambda_r,g_1,\dots,g_t,g)\fle(\lambda_1,\dots,\lambda_
{i-1},c(g)^{-1}\lambda_i^{-1},\lambda_{i+1},\dots,\lambda_r,g_1,\dots,g_t,g),\]
\item[-] pour tout $j\in\{1,\dots,t\}$, le morphisme
\[v_j:(\lambda_1,\dots,\lambda_r,g_1,\dots,g_t,g)\fle(\lambda_1,\dots,\lambda_
r,g_1,\dots,g_{j-1},c(g)^{-1}A_2{}^tg_j^{-1}A_2,g_{j+1},\dots,,g_t,g),\]

\end{itemize}
($\lambda_1,\dots,\lambda_r\in\Gr_m$, $g_1,\dots,g_t\in\GL_2$, $g\in\GSp_{2m}
$, $A_2=\left(\begin{array}{cc}0 & 1 \\ 1 & 0\end{array}\right)\in\GL_2(\Z)$).
Alors $\Omega^\M$ agit sur $\M$ par conjugaison par
des éléments de $\G(\Q)$. On fait agir $\Omega^\M$ sur $\M_H\simeq
\Gr_m^r\times\GL_2^t\times\G(\Sp_{2m_1}\times\SO_{2m_2})$ par les mêmes
formules. On note $\varepsilon_\kappa:\Omega^\M\fl\{\pm 1\}$ le morphisme
qui envoie les $v_1,\dots,v_t$ sur $1$ et tel que, pour tout $i\in\{1,\dots,
r\}$, $\varepsilon_\kappa(u_i)=-1$ si et seulement si $i\in A$. 
On voit facilement qu'un élément $\omega$ de $\Omega^\M$ agit sur
$\M_H$ par conjugaison par un élément de $\H(\Q)$ si et seulement si
$\varepsilon_\kappa(\omega)=1$.

\begin{subcorollaire}\label{cor:transfert_caracteres_GSp} On se place dans
la situation de la proposition \ref{prop:transfert_caracteres}, et on définit
une fonction $\psi:\M_H(\R)\fl\C$ par la formule suivante : pout tout
$\gamma_H\in\M_H(\R)$,
\[\psi(\gamma_H)=\sum_{\varphi_H\in\Phi_H
(\varphi)}\varepsilon(\omega_*(\varphi_H))\Phi_{M_H}^H(\gamma_H^{-1},S\Theta_
{\varphi_H}).\]

Alors, pour tous $\gamma_H\in\M_H(\R)$ et $\omega\in\Omega^\M$, on a
\[\psi(\omega(\gamma_H))=\varepsilon_\kappa(\omega)\psi(\gamma_H).\]

\end{subcorollaire}

\begin{preuve} Comme $\Ker(\varepsilon_\kappa)$ agit sur $\M_H$ par
conjugaison par des éléments de $\H(\Q)$ et que $\psi$ est clairement
invariante par conjugaison par $\H(\R)$, il suffit de montrer l'égalité
du corollaire pour un $\omega\in\Omega^\M$ tel que $\varepsilon_\kappa(\omega)=
-1$, par exemple $\omega=u_i$, avec $i\in A$.

Soit $\gamma_H\in\M_H(\R)$. Si $\gamma_H$ n'est pas elliptique ou n'a pas
d'image dans $\M(\R)$, alors il en est de meme de $\omega(\gamma_H)$, et
$\psi(\gamma_H)=\psi(\omega(\gamma_H))=0$. On peut donc supposer qu'il
existe des tores maximaux elliptiques $\T_M$ et $\T_{M_H}$ de
$\M_\R$ et $\M_{H,\R}$ et un isomorphisme $j_M:\T_{M_H}\iso\T_M$ comme
ci-dessus tels que
$\gamma_H\in\T_{M_H}(\R)$. On note $\gamma=j_M(\gamma_H)$;
il est clair que $u_i(\gamma)=j_M(u_i(\gamma_H))$.
D'après la proposition \ref{prop:transfert_caracteres},
\[\psi(\gamma_H)=\varepsilon_R(\gamma^{-1})\varepsilon_{R_H}(\gamma_H^{-1})
\Delta_{j_M,B_M}(\gamma_H,\gamma)\Phi_M^G(\gamma^{-1},
S\Theta_\varphi^H).\]
On voit facilement, en utilisant le fait que $j_{M*}(R_H^+)$
est stable par $u_i$, que $\Phi_M^G(\gamma^{-1},S\Theta_\varphi^H)=
\Phi_M^G(u_i(\gamma)^{-1},S\Theta_\varphi^H)$. 
Il est tout aussi facile de vérifier que
$\Delta_{j_M,B_M}(u_i(\gamma_H),u_i(\gamma))=\Delta_{j_M,B_M}(\gamma_H,
\gamma)$.
Il suffit donc de montrer que :
\[\varepsilon_R(u_i(\gamma)^{-1})\varepsilon_{R_H}(u_i(\gamma_H)
^{-1})=-\varepsilon_R(\gamma^{-1})\varepsilon_{R_H}(\gamma_H^{-1}).\]
Il est facile de voir que $u_i$ fait changer de signe exactement une
racine de $R_\gamma^+$ ($2e_i+c$ ou $-(2e_i+c)$, avec les notations
de l'exemple \ref{ex:groupes_symplectiques_suite}), donc
$\varepsilon_R(u_i(\gamma)^{-1})=-\varepsilon_R(\gamma^{-1})$.
Comme $i\in A$, cette racine n'est pas dans $j_{M*}(R_H)$, donc $u_i$ préserve
$R_{H,\gamma_H}^+$, et $\varepsilon_{R_H}
(u_i(\gamma_H)^{-1})=\varepsilon_{R_H}(\gamma_H^{-1})$. Ceci finit la
preuve.

\end{preuve}

\subsection{Calcul de certains $\Phi_M(\gamma,\Theta)$}
\label{infini3}

Soit $n\in\Nat^*$.
On note $\G=\GSp_{2n}$. Soient $\M$ un sous-groupe de Levi cuspidal standard de
$\G$ et $\Pa$ le sous-groupe parabolique standard
de $\G$ de sous-groupe de Levi $\M$.
On note $\B$ le sous-groupe des matrices triangulaires supérieures de $\G$
(c'est un sous-groupe de Borel contenu dans $\Pa$).
On note $\M_l$ la partie linéaire de $\M$ et $\M_h$ sa partie hermitienne. Soit
$\T_M$ un ($\R$)-tore maximal ($\R$-)elliptique de $\M_\R$. On note
$R$ l'ensemble des racines réelles de $\T_M$ dans $\G$, et $R^+$ l'ensemble
des racines de $R$ qui sont positives pour l'ordre donné par $\Pa$.

Soit $\gamma_M\in\T_M(\R)$. On note $I=\M_{\gamma_M}$. On fixe $\kappa_M\in
\Kgoth_\M(I/\R)$, et on note $(\M',s_M,\eta_{M,0},\gamma')$ le
$\M$-quadruplet endoscopique de $(\M,\gamma_M)$ associé à $\kappa_M$
(cf la fin de \ref{endoscopie2};
on remarque que $\gamma'$ est $(\M,\M')$-régulier
par construction). Pour tout $\kappa\in\Kgoth_\G(I/\R)$, on note
$(\M'_\kappa,s_{M,\kappa},\eta_{M,\kappa,0},\gamma'_\kappa)$ le
$\G$-quadruplet endoscopique de $(\M,\gamma_M)$ associé à $\kappa$, et $(\H_
\kappa,s_\kappa,\eta_{\kappa,0})$ le triplet endoscopique de $\G$ associé à
$(\M'_\kappa,s_{M,\kappa},\eta_{M,\kappa,0})$ (même référence).
Si l'image de $\kappa$ par l'application naturelle
$\Kgoth_\G(I/\R)\fl\Kgoth_\M(I/\R)$ est $\kappa_M$, alors il existe un
isomorphisme $\M'_\kappa\simeq\M'$ qui envoie $\gamma'_\kappa$ sur
$\gamma'$, et on peut supposer que $\eta_{M,0}=\eta_{M,\kappa,0}$ et
$s_MZ(\widehat{\M})=s_{M,\kappa}Z(\widehat{\M})$. 
On note $\Kgoth_\G(I/\R)_e$ l'ensemble des $\kappa\in\Kgoth_\G(I/\R)$ tels
que :
\begin{itemize}
\item[$\bullet$] les triplets endoscopiques $(\H_\kappa,s_\kappa,\eta_{\kappa,
0})$ et $(\M'_\kappa,s_{M,\kappa},\eta_{M,\kappa,0})$ sont elliptiques;
\item[$\bullet$] les groupes $\M'_\kappa$ et $\H_\kappa$ admettent des
$\R$-tores maximaux elliptiques (sur $\R$).

\end{itemize}

Soient $r,t,m\in\Nat$ tels que 
$\M_l\simeq\Gr_m^r\times\GL_2^t$ et $\M_h\simeq\GSp_{2m}$; on a
$n=r+2t+m$. Pour tout $i\in\{1,\dots,r\}$, on note $e_i:\M\fl
\Gr_m$ la projection sur le $i$-ième facteur $\Gr_m$. Pour tout $j\in\{1,\dots,
t\}$, on note $\alpha_j:\M\fl\Gr_m$ le composé de la projection sur le
$j$-ième facteur $\GL_2$ et du morphisme $\det:\GL_2\fl\Gr_m$.

Soient $m_1,m_2\in\Nat$ tels que $m=m_1+m_2$,
$m_2\not=1$ et $\M'=\M_l\times\G(\Sp_{2m_1}\times\SO_{2m_2})$. Comme
$\M'$ contient un $\R$-tore maximal elliptique sur $\R$, $m_2$ est forcément
pair.
On peut supposer que $s_M=
(1,s_{\varnothing,\varnothing,m_1,m_2})$, où les notations sont celles du
lemme \ref{lemme:Levi_endoscopiques_GSp}. Alors il résulte facilement de
ce lemme et de la remarque \ref{rq:Kgoth_G_M} que $\{\kappa\in\Kgoth_\G
(I/\R)_e|\kappa\fle\kappa_M\}$ s'identifie à l'ensemble des $(1,s_{A,B,m_1,
m_2})$, avec $A\subset\{1,\dots,r\}$ tel que $|A|$ soit pair et $B\subset\{1,
\dots,t\}$. Pour tout $s=(s_1,\dots,s_n)\in\{\pm 1\}^n$,
il existe une unique permutation
$\sigma\in\Sgoth_n$ telle que $\sigma^{-1}_{|\{i\in\{1,\dots,n\}|s_i=1\}}$
et $\sigma^{-1}_{|\{i\in\{1,\dots,n\}|s_i=-1\}}$ soient croissantes et
que $\sigma.s:=(s_{\sigma^{-1}(1)}\dots,s_{\sigma^{-1}(n)})$ soit de la
forme $(1,\dots,1,-1,\dots,-1)$; on la note $\sigma(s)$.
Si $\kappa\in\Kgoth_\G(I/\R)_e$ et
$(1,s_{A,B,m_1,m_2})$ est associé à $\kappa$ comme ci-dessus, on note
$\varepsilon(\kappa)$ la signature de $\sigma(s_{A,B,m_1,m_2})$ et, pour
tout $C\subset\{1,\dots,t\}$, on note $\varepsilon_C(\kappa)=
(-1)^{|B\cap C|}$.

\begin{subproposition}\label{prop:calcul_Phi_M_GSp} 
Soit $V$ une représentation
algébrique irréductible de $\G$. Soit $\varphi:W_\R\fl{}^L\G$ un paramètre de
Langlands du $L$-paquet de la série discrète de $\G(\R)$ associé à la
contragrédiente de $V$. Si $c(\gamma_M)<0$ (ce qui peut arriver seulement si
$\M_h$ est un tore), alors $\Phi_{\M'_\kappa}^{\H_\kappa}({\gamma'_\kappa}
^{-1},S\Theta_{\varphi_H})=0$ pour tout $\kappa\in\Kgoth_\G(I/\R)_e$ et tout
paramètre
de Langlands elliptique $\varphi_H$ de $\H_\kappa$. On suppose que
$c(\gamma_M)>0$, et on fixe $\kappa_M\in\Kgoth_M(I/\R)$ et $(\M',s_M,
\eta_{M,0},\gamma')$ comme plus haut.
\begin{itemize}
\item[(i)] Soit $C\subset\{1,\dots,t\}$ non vide. Alors
\[\sum_{{\kappa\in\Kgoth_\G(I/\R)_e}\atop{\kappa\fle\kappa_M}}
\Delta_{\M',s_{M,\kappa},\infty}^\M(\gamma',\gamma_M)^{-1}\varepsilon(\kappa)
\varepsilon_C(\kappa)\sum_{\varphi_H\in
\Phi_{H_\kappa}(\varphi)}\varepsilon(\omega_*(\varphi_H))\Phi_{\M'}^{\H_\kappa}
({\gamma'}^{-1},S\Theta_{\varphi_H})=0.\]

\item[(ii)] Si $e_i(\gamma_M)^2c(\gamma_M)\leq 1$ et $\alpha_j(\gamma_M)
c(\gamma_M)\leq 1$ pour tous
$i\in\{1,\dots,r\}$ et $j\in\{1,\dots,t\}$, alors :
\begin{flushleft}$\displaystyle{
2^{r+t}\frac{k(\M)}{k(\G)}(n_M^G)^{-1}(-1)^{q(\G)+
\dim(\A_M/\A_G)}\sum_{{\kappa\in\Kgoth_\G(I/\R)_e}\atop{\kappa\fle\kappa_M}}
\Delta_{\M',s_{M,\kappa},\infty}^\M(\gamma',\gamma_M)^{-1}\varepsilon(\kappa)
}$\end{flushleft}
\begin{flushright}$\displaystyle{
\sum_{\varphi_H\in
\Phi_{H_\kappa}(\varphi)}\varepsilon(\omega_*(\varphi_H))\Phi_{\M'}^{\H_\kappa}
({\gamma'}^{-1},S\Theta_{\varphi_H})=L_M(\gamma_M),
}$\end{flushright}
où $L_M(\gamma_M)$ est défini dans \ref{formule_points_fixes2}
(on utilise $\B$ pour former les bijections $\varphi_H\fle\omega_*(\varphi_H)$
dans le membre de gauche).
\footnote{Si $e_i(\gamma_M)^2c(\gamma_M)\geq 1$ et $\alpha_j(\gamma_M)
c(\gamma_M)\geq 1$ pour tous
$i\in\{1,\dots,r\}$ et $j\in\{1,\dots,t\}$,
on a un énoncé
similaire en remplaçant tout les $\Ho^*(Lie(\N_Q),V)_{>0}$ dans la définition
de $L_M(\gamma_M)$ par des $\Ho^*(Lie(\N_Q),V)_{<0}$.}

\end{itemize}

\end{subproposition}

On utilise la normalisation des facteurs de transfert de
\ref{infini2}. Les notations $\Phi(\varphi)$,
$\omega_*(\varphi_H)$, $S\Theta_{\varphi_H}$ et $\Phi_{\M'}^{\H_\kappa}$ sont
définies dans \ref{infini1} et \ref{infini2}, et
les notations $k(\M)$, $k(\G)$ et $n_M^G$ sont définies dans
\ref{stabilisation1}.

\begin{preuve} Si $c(\gamma_M)<0$, alors, pour tout $\kappa\in\Kgoth_\G(I/\R)
_e$, on a $c(\gamma'_\kappa)=c(\gamma_M)<0$, donc $\gamma'_k\not\in
Z(\H_\kappa)(\R)Im(\H_{\kappa,sc}(\R)\fl\H_\kappa(\R))$. Ceci prouve
l'annulation de tous les $\Phi_{\M'_\kappa}^{\H_\kappa}({\gamma'_\kappa}^{-1},
S\Theta_{\varphi_H})$. On écrit $\gamma_M=(\gamma_l,\gamma_h)$, avec
$\gamma_l\in\M_l(\R)$ et $\gamma_h\in\M_h(\R)$. Alors $c(\gamma_M)=c(\gamma_h)
$. Soit $\T$ un $\R$-tore maximal elliptique de $\M_h$ tel que
$\gamma_h\in\T(\R)$. Alors $\T/(\Gr_{m,\R}I_{2m})$ est anisotrope sur $\R$,
donc le groupe de ses $\R$-points est connexe.
On suppose que $\M_h$ n'est pas un tore, c'est-à-dire que
$m>0$.
Alors la restriction de
$c$ à $\Gr_{m,\R}I_{2m}$ est $\lambda I_{2m}\fle\lambda^2$, donc envoie
le groupe des $\R$-points de $\Gr_{m,\R}I_{2m}$ dans $\R^{+*}$, et on a
$c(\T(\R))\subset\R^{+*}$. Donc $c(\gamma_M)>0$.

On suppose à partir de maintenant que $c(\gamma_M)>0$ (avec $m$ quelconque),
et on utilise librement les notations de \ref{formule_points_fixes2}.
Soit $\T_M$ (resp. $\T_{M'}$) un $\R$-tore maximal elliptique de
$\M$ (resp. $\M'$) tel que $\gamma_M\in\T_M(\R)$ (resp. $\gamma'\in\T_{M'}
(\R)$). Comme dans l'exemple \ref{ex:groupes_symplectiques_suite}, on
écrit $\T_M=\Gr_{m,\R}^r\times\T_1\times\dots\times\T_t\times\T$ et
$\T_{M'}=\Gr_{m,\R}^r\times\T_1\times\dots\times\T_t\times\T'$, et on
choisit $j_M:\T_{M'}\iso\T_M$ qui est l'identité sur $\Gr_{m,\R}^r\times\T_1
\times\dots\times\T_t$.
On peut supposer que $\gamma_M=j_M(\gamma')$.

Soient $A\subset\{1,\dots,r\}$ et $B\subset\{1,\dots,t\}$ tels que $|A|$ soit
pair. On note $s_{A,B}=(1,s_{A,B,m_1,m_2})$ et, si $\kappa$
correspond à $s_{A,B}$, on note $\H_{A,B}=\H_\kappa$.
Soit $\sigma\in\Sgoth_r$ l'unique permutation telle que $\sigma^{-1}_{|A}$ et
$\sigma^{-1}_{|\{1,\dots,r\}-A}$ soient croissantes et $\sigma^{-1}(A)=
\{n+1-|A|,\dots,n\}$; on note $\varepsilon(A)$ la signature de $\sigma$. On
voit facilement (en utilisant le fait que $|A|$ est pair) que
$\varepsilon(\kappa)=\varepsilon(A)$.
Donc le membre de gauche de l'égalité du point (i) (resp. (ii)) de
la proposition est égal à (resp. au produit de $2^{r+t}k(\M)k(\G)^{-1}
(n_M^G)^{-1}(-1)^{q(\G)+\dim(\A_M/\A_G)}$ et de)
\[\sum_{A,B}\Delta_{\M',s_{A,B},\infty}^\M(\gamma',\gamma_M)^{-1}
\varepsilon(A)(-1)^{|B\cap C|}\sum
_{\varphi_H\in\Phi_{H_{A,B}}(\varphi)}\varepsilon(\omega_*(\varphi_H))\Phi_
{\M'}^{\H_{A,B}}({\gamma'}^{-1},S\Theta_{\varphi_H})\]
(resp. la même somme mais sans le facteur $(-1)^{|B\cap C|}$),
où $A$ parcourt l'ensemble des sous-ensembles de cardinal pair de
$\{1,\dots,r\}$ et $B$ parcourt l'ensemble des sous-ensembles de $\{1,\dots,
t\}$.
En particulier, les membres de gauche dans (i) et (ii) sont des fonctions
continues de $\gamma_M\in\T_M(\R)$.
Il est clair que les membres de droite sont aussi des fonctions continues
de $\gamma_M$.
On peut donc supposer que $\gamma_M$ est $\G$-régulier, et dans (ii),
on peut supposer que
$e_i(\gamma_M)^2c(\gamma_M)<1$ pour tout $i\in\{1,\dots,r\}$ et que $
\alpha_j(\gamma_M)c(\gamma_M)<1$ pour tout $j\in\{1,\dots,t\}$.

On note $MG^{(i)}=MG^{(i)}(\gamma_M,V,C)$ (resp. $MG=MG(\gamma_M,V)$, resp.
$MD=MD(\gamma_M,V)$) le membre de gauche de l'égalité de (i) (resp. le
membre de gauche de l'égalité de (ii), resp. le membre de droite de
l'égalité de (ii)).

Calculons $MG^{(i)}$ et $MG$.
On remarque d'abord que $MG^{(i)}(\gamma_M,V,C)=MG^{(i)}(\gamma_M^{-1},V^*,
C)$ (où $V^*$ est
la contragrédiente de $V$). Donc, quitte à remplacer $\varphi$ par un
paramètre de Langlands du $L$-paquet de la série discrète de $\G(\R)$ associé
à $V$, on peut remplacer $\gamma_M$ et $\gamma'$ par $\gamma_M^{-1}$ et
${\gamma'}^{-1}$ dans le membre de gauche. 

On a $n_M^G=2^rr!2^tt!$ et, d'après le lemme \ref{lemme:calcul_k}, $k(\G)=
2^{n-1}$ et $k(\M)=2^{m-1}$ si $m\geq 1$ (si $m=0$, on a $\M=\Gr_m^r\times\GL_
2^t\times\Gr_m$, donc $k(\M)=1$). Donc
\[(n_M^G)^{-1}\frac{k(\M)}{k(\G)}=m_P2^{-2r-3t}(r!)^{-1}(t!)^{-1}\]
(d'après la définition de $m_P$ dans \ref{formule_points_fixes2},
on a $m_P=1$ si $m\geq 1$ et $m_P=2$ si $m=0$), donc
\[MG=m_P2^{-r-2t}(r!)^{-1}(t!)^{-1}(-1)^{q(\G)+\dim(\A_M/\A_G)}MG^{(i)}
(\gamma_M,V,\varnothing).\]

On a (cf l'exemple \ref{ex:groupes_symplectiques_suite})
$R=\{\pm(e_i-e_j),1\leq i<j\leq r\}\cup\{\pm(e_i+e_j+c),1\leq i\leq j\leq r\}
\cup\{\pm(\alpha_i+c),1\leq i\leq t\}$ et
$R^+=\{e_i-e_j,1\leq i<j\leq r\}\cup
\{e_i+e_j+c,1\leq i\leq j\leq r\}\cup\{\alpha_i+c,1\leq i\leq t\}$. 
Comme $MG^{(i)}$ et $MD$ ne changent pas si on remplace $\gamma_M$ par
$g\gamma_M g^{-1}$ avec $g\in\Nor'_\G(\M)(\Q)$, on peut supposer que :
\begin{itemize}
\item[$\bullet$] il existe $r'\in\{1,\dots,r\}$ tel que $e_i(\gamma_M)>0$
pour $i\in\{1,\dots,r'\}$, et $e_i(\gamma_M)<0$ pour $i\in\{r'+1,\dots,r\}$;
\item[$\bullet$] si $1\leq i<i'\leq r'$ ou $r'+1\leq i<i'\leq r$, alors 
$(0<)e_i(\gamma_M)e_{i'}(\gamma_M)^{-1}<1$.

\end{itemize}

Soient $A\subset
\{1,\dots,r\}$ et $B\subset\{1,\dots,t\}$; on suppose que $|A|$ est pair.
On note $R_{A,B}=j_{M*}(R_{H_{A,B}})$, $R_{A,B}^+=j_{M*}(R_{H_{A,B}}^+)=
R^+\cap R_{A,B}$.
Une formule explicite pour $R_{A,B}$ est donnée dans l'exemple
\ref{ex:groupes_symplectiques_suite}.
D'après la proposition \ref{prop:transfert_caracteres}, on a
\[\Delta_{\M',s_{A,B},\infty}^\M({\gamma'}^{-1},\gamma_M^{-1})\sum_
{\varphi_H\in\Phi_{H_{A,B}}(\varphi)}\varepsilon(\omega_*(\varphi_H))
\Phi_{\M'}^{\H_{A,B}}(\gamma',S\Theta_{\varphi_H})=
\varepsilon_R(\gamma_M)\varepsilon_{R_H}(\gamma')
\Phi_\M^\G(\gamma_M,S\Theta_\varphi^{A,B}),\]
où $S\Theta_\varphi^{A,B}=S\Theta_\varphi^{H_{A,B}}$.
D'après le fait \ref{fait:formule_Phi_M_G} (dont on utilise les notations),
$\Phi_\M^\G(\gamma_M,S\Theta_\varphi^{A,B})$ est égal à
\[(-1)^{q(\G)}\varepsilon_R(\gamma
_M)\delta_{\Pa(\R)}^{1/2}(\gamma_M)\Delta_{\B_M}(\gamma_M)^{-1}\sum_{\omega\in
\Omega}\varepsilon(\omega)n_{A,B}(\gamma_M,\omega\B)(\omega\lambda)(\gamma_M)
\prod_{\alpha\in\Phi(\omega)}\alpha^{-1}(\gamma_M),\]
où $\lambda$ est le plus haut poids de $V$ relativement à $\B$ et
$n_{A,B}=n_{H_{A,B}}$.

Calculons $R_{\gamma_M}$.
On sait (cf l'exemple
\ref{ex:groupes_symplectiques_suite}) que $\alpha_i(\gamma_M)>0$
pour tout $i\in\{1,\dots,t\}$ (et on a supposé que $c(\gamma_M)>0$).
Donc $\{\pm(\alpha_1+c),\dots,\pm(\alpha_t+c)\}\subset R_{\gamma_M}$.
On note $I^+=\{1,\dots,r'\}=\{i\in\{1,\dots,r\}|e_i(\gamma_M)>0\}$
et $I^-=\{r'+1,\dots,r\}=\{i\in\{1,\dots,r\}|e_i(\gamma_M)<0\}$. Alors
$R_{\gamma_M}\cap\{\pm(e_i-e_j),1\leq i<j\leq r\}=\{\pm(e_i-e_j)|i<j
\mbox{ et }i,j\in I^+\}\cup\{\pm(e_i-e_j)|i<j\mbox{ et }i,j\in I^-\}$,
et $R_{\gamma_M}\cap\{\pm (e_i+e_j+c),1\leq i\leq j\leq r\}=\{\pm (e_i+e_j+c),
i,j\in I^+\}\cup\{\pm(e_i+e_j+c),i,j\in I^-\}$.

Si de plus $\gamma_M$ vérifie les conditions de (ii), alors
$0<\alpha(\gamma_M)<1$ pour tout
$\alpha\in R^+\cap R_{\gamma_M}$, donc $R_{\gamma_M}^+=R_{\gamma_M}\cap(-R^+)$.
Donc, dans ce cas,
$\Phi^+\cap(-R_{\gamma_M}^+)=\Phi^+\cap R_{\gamma_M}$ est de cardinal
$t+|I^+|^2+|I^-|^2$, et
\[\varepsilon_R(\gamma_M)=(-1)^{t+|I^+|^2+|I^-|^2}=(-1)^{t+|I^+|+|I^-|}
=(-1)^{r+t}.\]
On voit de même que, si $\gamma_M$ vérifie les hypothèses
de (ii), alors $\varepsilon_{R_H}(\gamma')=(-1)^{r_1+t_1+t_2}$, où $r_1,r_2,
t_1, t_2$ sont définis comme dans le lemme \ref{lemme:Levi_endoscopiques_GSp}.
Comme $r_2=|A|$ est pair, on en déduit que $\varepsilon_R(\gamma_M)=
\varepsilon_{R_H}(\gamma')$ dans ce cas.

D'autre part, on a $\dim(\A_M/\A_G)=r+t$.

On obtient donc :
\begin{flushleft}$\displaystyle{
MG^{(i)}=(-1)^{q(\G)}\varepsilon_{R_H}(\gamma')\delta_{\Pa(\R)}^{1/2}(\gamma_M)
\Delta_{\B_M}(\gamma_M)^{-1}
}$\end{flushleft}
\begin{flushright}$\displaystyle{
\sum_{\omega\in\Omega}\varepsilon(\omega)(\omega
\lambda)(\gamma_M)\prod_{\alpha\in\Phi(\omega)}\alpha^{-1}(\gamma_M)\left(
\sum_{A,B}\varepsilon(A)(-1)^{|B\cap C|}n_{A,B}(\gamma_M,\omega\B)\right)
}$\end{flushright}
et
\begin{flushleft}$\displaystyle{
MG=m_P2^{-r-2t}(r!)^{-1}(t!)^{-1}\delta_{\Pa(\R)}^{1/2}(\gamma_M)
\Delta_{\B_M}(\gamma_M)^{-1}
}$\end{flushleft}
\begin{flushright}$\displaystyle{
\sum_{\omega\in\Omega}\varepsilon(\omega)(\omega
\lambda)(\gamma_M)\prod_{\alpha\in\Phi(\omega)}\alpha^{-1}(\gamma_M)
\left(\sum_{A,B}\varepsilon(A)n_{A,B}(\gamma_M,\omega\B)\right),
}$\end{flushright}
où, dans la dernière somme de ces deux expressions,
$A$ parcourt l'ensemble des sous-ensembles de
cardinal pair de $\{1,\dots,r\}$, et $B$ parcourt l'ensemble des
sous-ensembles de $\{1,\dots,t\}$.

On remarque que
$R_{A,B}$ ne dépend pas de $B$, donc
\[\sum_{A,B}\varepsilon(A)(-1)^{|B\cap C|}n_{A,B}(\gamma_M,\omega\B)=
\left(\sum_B(-1)^{|B\cap C|}\right)\left(\sum_A\varepsilon(A)
n_{A,\varnothing}(\gamma_M,\omega\B)\right).\]
Si $C\not=\varnothing$, alors $\sum\limits_B(-1)^{|B\cap C|}=0$.
Ceci finit la preuve
de (i). On s'intéresse maintenant à l'égalité de (ii), donc au cas où
$C=\varnothing$. Alors $\sum\limits_B(-1)^{|B\cap C|}=2^t$, c'est-à-dire que
\[\sum_{A,B}\varepsilon(A)n_{A,B}(\gamma_M,\omega\B)=
2^t\sum_A\varepsilon(A)
n_{A,\varnothing}(\gamma_M,\omega\B).\]
On utilise les formules de l'article \cite{Her} de Herb pour calculer
$\sum\limits_A\varepsilon(A)n_{A,\varnothing}(\gamma_M,\omega\B)$.
On note $R_{A,\varnothing}=R_A$ et $n_{A,\varnothing}=n_A$.
On note $R_{A,\gamma_M}^+=R_A\cap R_{\gamma_M}
^+=R_A\cap R_{\gamma_M}\cap(-R^+)$.
On note $A_1=\{1,\dots,r\}-A$,
$A_2=A$, $A_1^+=A_1\cap I^+$, $A_1^-=A_1\cap I^-$, $A_2^+=A_2\cap I^+$,
$A_2^-=A_2\cap I^-$.
Alors $R_{A,\gamma_M}$ est le produit direct du système de racines de type
$C_{|A_1^+|}$ sur les racines $e_i+c/2$, $i\in A_1^+$, du système de racines de
type $C_{|A_1^-|}$ sur les racines $e_i+c/2$, $i\in A_1^-$,
du système de racines de type $D_{|A_2^+|}$ sur les racines $e_i+c/2$,
$i\in A_2^+$, du système de racines de type $D_{|A_2^-|}$ sur les racines $e_i
+c/2$, $i\in A_2^-$ et des $t$ systèmes
de racines de type $A_1$ sur les racines $\alpha_i+c$, $1\leq i\leq t$.
Donc $-1$ est dans le groupe de Weyl de $R_{A,\gamma_M}$ si et seulement si
$|A_2^+|$ et $|A_2^-|$ sont pairs, ce qui revient à demander que
$|A_2^-|$ soit pair (car $|A_2|$ est pair); on voit facilement que ceci est
équivalent au fait que $\gamma'$ soit dans $Z(\H_{A,B})(\R)Im(\H_{A,B,sc}(\R)
\fl\H_{A,B}(\R))$, pour un $B\subset\{1,\dots,t\}$ (ou pour tout $B\subset
\{1,\dots,t\}$; cette condition ne dépend pas de $B$). Si cette condition
n'est pas vérifiée, alors $n_A(\gamma_M,\omega\B)=0$ pour tout $\omega\in
\Omega$. On suppose donc à partir de maintenant que que $|A_2^-|$ est pair.

Pour tous $a,b\in\R$, on note
\[c_1(a)=\left\{\begin{array}{ll}1 & \mbox{ si }a>0 \\
0 & \mbox{ sinon}\end{array}\right.\]
\[c_{2,C}(a,b)=\left\{\begin{array}{ll}1 & \mbox{ si }
0<a<b\mbox{ ou }0<-b<a \\
0 & \mbox{ sinon}\end{array}\right.\]
\[c_{2,D}(a,b)=\left\{\begin{array}{ll}1 & \mbox{ si }
a>|b| \\
0 & \mbox{ sinon}\end{array}\right.\]
et
\[c_2(a,b)=c_{2,C}(a,b)+c_{2,D}(a,b)=
\left\{\begin{array}{ll}0 & \mbox{ si }a+b\leq 0\mbox{ ou }a\leq 0 \\
1 & \mbox{ si }a>0\mbox{ et }b>0 \\
2 & \mbox{ sinon}\end{array}\right..\]
Comme dans la section
\ref{lemmes_combinatoires}, on note, pour tout ensemble fini $I$,
$\Par^0_{\leq 2}(I)$ l'ensemble
des partitions (non ordonnées) $\{I_z,z\in Z\}$ de $I$
telles que $|I_z|\leq 2$ pour tout $z\in Z$ et qu'au
plus un des $I_z$ soit de cardinal $1$. Si $I=\{1,\dots,r\}$, on note
$\Par^0_{\leq 2}(r)=\Par^0_{\leq 2}(I)$.
Soit $p=\{I_z,z\in Z\}
\in\Par^0_{\leq 2}(r)$. On choisit une énumération $z_1,\dots,z_k$
de $Z$, et on note $\sigma$ l'élément de $\Sgoth_r$ tel que :
\begin{itemize}
\item[$\bullet$] si $1\leq i<i'\leq k$, alors, pour tous $s\in\sigma(I_{z_
i})$ et $s'\in\sigma(I_{z_{i'}})$, on a $s<s'$;
\item[$\bullet$] pour tout $z\in Z$ tel que $|I_z|=2$, si
$I_z=\{s_1,s_2\}$ avec $s_1<s_2$, alors $\sigma(s_1)<\sigma(s_2)$.

\end{itemize}
Alors la signature de $\sigma$ ne dépend pas de l'énumération de $Z$ choisie,
et on la note $\varepsilon(p)$. (Cette définition est la même que celle
de la section \ref{lemmes_combinatoires}.)
Si $I$ est un ensemble fini totalement
ordonné et $p\in\Par^0_{\leq 2}(I)$, on définit
$\varepsilon(p)$ en utilisant l'ordre pour identifier $I$ à
$\{1,\dots,|I|\}$.
D'autre part, soient $\mu=(\mu_1,\dots,\mu_r)\in\R^r$ et $I$ un
sous-ensemble de $\{1,\dots,r\}$ de cardinal $1$ ou $2$. Si $|I|=1$ et
$I=\{s\}$, on note $c_I(\mu)=c_{I,C}(\mu)=c_1(\mu_s)$;
si $|I|=2$ et $I=\{s_1,s_2\}$
avec $s_1<s_2$, on note $c_{I,C}(\mu)=c_{2,C}(\mu_{s_1},\mu_{s_2})$, 
$c_{I,D}(\mu)=c_{2,D}(\mu_{s_1},\mu_{s_2})$ et $c_I(\mu)=c_2(\mu_{s_1},
\mu_{s_2})$. Si $\mu\in\R^r$,
$I\subset\{1,\dots,r\}$ et $p=\{I_z,z\in Z\}\in
\Par^0_{\leq 2}(Z)$, on pose
\[c_C(p,\mu)=\prod_{z\in Z}c_{I_z,C}(\mu).\]
On définit de même $c_D(p,\mu)$ (si $|I|$ est pair) et $c(p,\mu)$.

On remarque que $\A_M=\Gr_m^r\times\Gr_m^t\times\Gr_m$, et que le
sous-groupe $\A_G$ de $\A_M$ est le dernier facteur $\Gr_m$. Donc
$(e_1,\dots,e_r,\alpha_1,\dots,\alpha_t)$ est une base de
$X^*(\A_M/\A_G)\otimes_\Z\R$. Soit $\chi\in(X^*(\A_M/\A_G)\otimes_\Z\R)_{reg}$.
On écrit $\chi=\sum_{i=1}^r\mu_ie_i+\sum_{j=1}^t\nu_j\alpha_j$, avec
$\mu=(\mu_1,\dots,\mu_r)\in\R^r$ et $\nu_1,\dots,\nu_t\in\R$. On écrit,
comme dans \ref{infini1}, $\gamma_M=\exp(x)\gamma_1$, avec
$x\in X_*(\A_M)\otimes_\Z\R$ et $\gamma_1\in\T_M(\R)_1$. Alors, d'après
le théorème 2.12 de \cite{Her} (et les hypothèses sur $\gamma_M$), le
coefficient $\overline{c}_{R_{A,\gamma_M}}(x,\chi)$ est égal au produit de
$2^{r+t}c_1(\nu_1)\dots c_1(\nu_t)$ et de
\[\sum_{p_1^+}\sum_{p_1^-}\sum_{p_2^+}\sum_{p_2^-}\varepsilon(p_1^+)
\varepsilon(p_2^+)\varepsilon(p_1^-)\varepsilon(p_2^-)c_C(p_1^+,\mu)c_C(p_1^-,
\mu)c_D(p_2^+,\mu)c_D(p_2^-,\mu),\]
où $p_i^+$ parcourt $\Par^0_{\leq 2}(A_i^+)$ et $p_i^-$ parcourt $\Par^0_
{\leq 2}(A_i^-)$.
On en déduit que $\sum\limits_A\varepsilon(A)\overline{c}_{R_{A,
\gamma_M}}(x,\chi)$ (où $A$ parcourt l'ensemble des sous-ensembles de cardinal
pair de $\{1,\dots,r\}$) est le produit de $2^{r+t}$
et de
\[c(x,\chi):=
c_1(\nu_1)\dots c_1(\nu_t)
\sum_{p^+\in\Par^0_{\leq 2}(I^+)}\sum_{p^-\in\Par^0_{\leq 2}(I^-)}
\varepsilon(p^+)\varepsilon(p^-)c(p^+,\mu)c(p^-,\mu).\]

On revient au calcul de $\sum\limits_{A,B}\varepsilon(A)n_{A,B}
(\gamma_M,\omega\B)$.
On fixe $\omega\in\Omega$. Par définition (cf \ref{infini1}),
on a $n_{A,B}(\gamma_M,\omega\B)=\overline{c}_{R_{A,\gamma_M}}(x,p(\omega(
\lambda+\rho-\lambda_0)))$, où $\rho=\rho_B$, $\lambda$ est toujours le plus
haut poids de $V$ relativement à $\B$ et $\lambda_0$ est comme
dans \ref{infini1} le caractère par lequel $\A_G$ agit sur $V$.
D'après les calculs ci-dessus, on a donc
\[\sum_{A,B}\varepsilon(A)n_{A,B}(\gamma_M,\omega(\B))=2^{r+2t}c(x,p(\omega(
\lambda+\rho-\lambda_0))).\]
Finalement :
\begin{flushleft}$\displaystyle{
MG=m_P(r!)^{-1}(t!)^{-1}\delta_{\Pa(\R)}^{1/2}(\gamma_M)
\Delta_{\B_M}(\gamma_M)^{-1}
}$\end{flushleft}
\begin{flushright}$\displaystyle{
\sum_{\omega\in\Omega}\varepsilon(\omega)c(x,p(\omega(\lambda+\rho-\lambda_0)))
(\omega\lambda)(\gamma_M)\prod_{\alpha\in\Phi(\omega)}\alpha^{-1}(\gamma_M).
}$\end{flushright}

On calcule maintenant le membre de droite $MD$ de l'égalité de (ii).
Pour tout $j\in\{1,\dots,t\}$, soit $u_j\in\GL_2(\C)$ tel que $\Int(u_j)$
envoie $\T_{j,\C}$ sur le tore diagonal de $\GL_{2,\C}$. Soit $v\in
\GSp_{2m}(\C)$ tel que $\Int(v)$ envoie $\T_\C$ sur le tore diagonal
de $\GSp_{2m,\C}$. On note
$u=(1,\dots,1,u_1,\dots,u_t,v)\in\M(\C)=(\C^\times)^r\times\GL_2(\C)^t\times
\GSp_{2m}(\C)$. Alors $\Int(u)$ envoie $\T_{M,\C}$ sur le tore diagonal
$\T_{0,\C}$ de $\GSp_{2n,\C}$, et $\Int(u)(\Pa)\supset\B$. On
utilise la conjugaison par $u$ pour identifier $X^*(\T_0)$ et
$X^*(\T_M)$, et on
note $\B_M=\M_\C\cap\Int(u)^{-1}(\B_\C)$ (c'est un sous-groupe de Borel
de $\M_\C$ contenant $\T_{M,\C}$). 

Comme dans la section \ref{lemmes_combinatoires},
on note $\Par_{ord}(r,t)$ l'ensemble des
partitions ordonnées $P$ de $\{1,\dots,r+2t\}$ telles que,
pour tout $i\in\{1,\dots,t\}$, $r+2i-1$ et $r+2i$ soient dans le même
ensemble de $P$. Soit $P=(I_1,\dots,I_k)\in\Par_{ord}(r,t)$. On note
$|P|=k$ et $s_i=|I_i|$, pour tout $i\in\{1,\dots,k\}$. Soit $\sigma_P$
l'unique permutation de $\Sgoth_{r+2t}$ telle que, pour tout $i\in\{0,\dots,
k-1\}$, $\sigma_P^{-1}$ soit croissante sur $\{s_1+\dots+s_i+1,\dots,
s_1+\dots+s_{i+1}\}$ et envoie $\{s_1+\dots+s_i+1,\dots,s_1+\dots,s_{i+1}\}$
sur $I_i$. On note $\varepsilon(P)$ la signature de $\sigma_P$.

Soit $P=(I_1,\dots,I_k)\in\Par_{ord}(r,t)$. On va associer à $P$ un
élément $(\QP,g)=(\QP_P,g_P)$ de  $\Par(\M)$ (les notations sont celles
de \ref{formule_points_fixes2}).
Pour tout $i\in\{1,\dots,k\}$,
on note $s_i=|I_i|$, $r_i=|I_i\cap\{1,\dots,r\}|$
et $t_i=\frac{1}{2}(s_i-r_i)$;
on a $t_1,\dots,t_k\in\Nat$ d'après la définition de $\Par_{ord}(r,t)$.
On prend $\QP=\Pa_S$, où $S=\{s_1,s_1+s_2,\dots,s_1+\dots+s_k\}$. 
Donc
$\Le_Q=\GL_{s_1}\times\dots\times\GL_{s_k}$ et $\G_Q=\GSp_{2m}$. 
On note $P_r=(I_1\cap\{1,\dots,r\},\dots,I_k\cap\{1,\dots,r\})$ et
$P_t=(J_1,\dots,J_k)$ avec, pour tout $i\in\{1,\dots,k\}$,
$J_i=\{l\in\{1,\dots,t\}|r+2l\in I_i\}$. On associe à $P_r$ et $P_t$ des
permutations $\sigma_r\in\Sgoth_r$ et $\sigma_t\in\Sgoth_t$ comme on l'a fait
pour $P$ (le fait que certains des ensembles composant $P_r$ et $P_t$
peuvent être vides n'a pas d'importance pour cette construction).
On choisit $g$ tel que :
\begin{itemize}
\item[$\bullet$] La restriction à $\M_h$ de $\Int(g)$ est l'identité
$\M_h\fl\G_Q$ (comme $\M_Q$ et $\M$ sont standard, on a $\M_h=\G_Q$);
\item[$\bullet$] $\Int(g)$ envoie $\M_l$ sur un sous-groupe de Levi
standard de $\Le_Q$.
\item[$\bullet$] Pour tout $i\in\{0,\dots,k-1\}$, le composé du morphisme
$\M_l\fl\Le_Q$ ci-dessus et de la projection de $\Le_Q$ sur le
facteur $\GL_{s_{i+1}}$ est égal au morphisme qui envoie
$(x_1,\dots,x_r,g_1,\dots,g_t)\in\Gr_m^r\times\GL_2^t$ sur
\[diag(x_{\sigma_r^{-1}(r_1+\dots+r_i+1)},\dots,x_{\sigma_r^{-1}(r_1+\dots+
r_{i+1})},g_{\sigma_t^{-1}(t_1+\dots+t_i+1)},\dots,g_{\sigma_t^{-1}(t_1+\dots+
t_{i+1})}).\]
\item[$\bullet$] Soit $u_P\in\Le_Q(\C)$ tel que, pour tout $i\in\{0,\dots,
k-1\}$,
le projeté de $u_P$ sur le facteur $\GL_{s_{i+1}}(\C)$ soit
$diag(1,\dots,1,u_{\sigma_t^{-1}(t_1+\dots+t_i+1)},\dots,u_{\sigma_t^{-1}(t_1+
\dots+t_{i+1})})$.
Si on utilise la conjugaison par $u$ pour identifier $(\T_M\cap\M_l)_\C$ à
$(\T_0\cap\M_l)=\Gr_{m,\C}^{r+2t}$ et la conjugaison par $u_P$ pour identifier
$\Int(g)(\T_M\cap\M_l)_\C$ à $(\T_0\cap\Le_Q)_\C=(\T_0\cap\M_l)_\C$, alors
$\Int(g)$ est l'application qui envoie $(z_1,\dots,z_{r+2t})$ sur
$(z_{\sigma_P^{-1}(1)},\dots,z_{\sigma_P^{-1}(r+2t)})$.

\end{itemize}
On note $\gamma=\Int(g)(\gamma_M).$
On remarque que, pour tout $i\in\{1,\dots,k\}$,
l'intersection de $g\B g^{-1}$ avec le facteur $\GL_{s_i}$ de $\Le_Q$ est
le sous-groupe des matrices triangulaires supérieures.
On note $\M_0=g\M g^{-1}=\Le_0\times\G_Q$, $\T_{M_0}=\Int(g)(\T_M)$
(le tore maximal de $\M_{0,\R}$ qui contient $\gamma$)
et $\Pa_0\subset\QP$ le
sous-groupe parabolique standard de sous-groupe de Levi $\M_0$.
On utilise la conjugaison par $(u_P,v)$ pour identifier $X^*(\T_{M_0})$
et $X^*(\T_0)$.

On utilise le théorème de Kostant (voir par exemple \cite{GHM} \S11) pour
calculer $\Ho^*(Lie(\N_Q),V)_{>0}$.
On note $\Omega$ le groupe de Weyl de $\T_0(\C)$ dans $\G(\C)$, $\Omega_{M_Q}$
le groupe de Weyl de $\T_0(\C)$ dans $\M_Q(\C)$ et $\Phi^+=\Phi(\T_0,\B)$.
Pour tout $\omega\in\Omega$, on note comme avant $\Phi(\omega)=\Phi^+\cap
(-\omega\Phi^+)$. Alors $\Omega'_Q:=\{\omega\in\Omega|\Phi(\omega)\subset
\Phi(\T_0,Lie(\N_Q))\}$ est un ensemble de représentants de $\Omega_{M_Q}\sous
\Omega$. D'après le théorème de Kostant, on a, pour tout $i\in\Nat$,
\[\H^i(Lie(\N_Q),V)\simeq\bigoplus_{\omega'\in\Omega'_Q,\ell(\omega')=i}
V_{\M_Q,\omega'(\lambda+\rho)-\rho},\]
où $\ell:\Omega\fl\Nat$ est la fonction longueur, $\lambda\in X^*(\T_0)$ est le
plus haut poids de $V$ relativement à $\B$, $\rho=\frac{1}{2}\sum\limits_
{\alpha\in\Phi^+}\alpha$ et, pour tout $\mu\in X^*(\T_0)$ dominant,
$V_{\M_Q,\mu}$ est
la représentation algébrique de $\M_Q$ de plus haut poids $\mu$
relativement à $\B\cap\M_Q$. Pour tout $s\in\{1,\dots,n\}$, on note
\[\varpi_s:\Gr_m\fl\T_0,\quad\lambda\fle\left(\begin{array}{ccc}
\lambda I_s & & 0 \\ & I_{2(n-s)} & \\ 0 & & \lambda^{-1}I_s
\end{array}\right).\]
Alors, par définition de $\Ho^i(\dots)_{>0}$ dans \ref{formule_points_fixes2}
(et de l'opération de troncature dans \cite{M2} 4.2),
on a pour tout $i\in\Nat$ un isomorphisme
\[\H^i(Lie(\N_Q),V)_{>0}\simeq\bigoplus_{\omega'} V_{\M_Q,\omega'(\lambda+\rho)
-\rho},\]
où $\omega'$ parcourt l'ensemble des éléments de $\Omega'_Q$ vérifiant
$\ell(\omega')=i$ et $\langle \omega'(\lambda+\rho)-\rho,\varpi_s\rangle >-ns+\frac{1}{2}
s(s-1)$ pour tout $s\in S$.
Comme $\omega(\lambda_0)=\lambda_0$ pour tout $\omega\in\Omega$ et
$\langle\rho,\varpi_s\rangle =ns-\frac{1}{2}s(s-1)$ et
$\langle \lambda_0,\varpi_s\rangle =0$ pour tout $s\in\{1,\dots,n\}$, on peut remplacer la
deuxième condition par la condition : pour tout $s\in S$, $\langle \omega'(\lambda+
\rho+\lambda_0),\varpi_s\rangle >0$.

D'après la formule du caractère
de Weyl (c'est-à-dire la première formule du fait \ref{fait:formule_Phi_M_G}),
on a, pour tout $\omega'\in\Omega'_Q$,
\[\Tr(\gamma,V_{\M_Q,\omega'(\lambda+\rho)-\rho})=\Delta_{\B\cap\M_Q}(\gamma)
^{-1}\sum_{\omega_M\in\Omega_{M_Q}}\varepsilon(\omega_M)(\omega_M(\omega'(
\lambda+\rho)-\rho))(\gamma)\prod_{\alpha\in\Phi_{M_Q}(\omega_M)}\alpha^{-1}
(\gamma),\]
où $\Phi_{M_Q}=\Phi(\T_0,\M_Q)$.
Soient $\omega'\in\Omega'_Q$ et $\omega_M\in\Omega_{M_Q}$. Comme
$\varpi_s$ est invariant par $\Omega_M$ pour tout $s\in S$, on a, pour tout
$s\in S$ :
\[\langle \omega_M\omega'(\lambda+\rho+\lambda_0),\varpi_s\rangle >0\Longleftrightarrow
\langle \omega'(\lambda+\rho+\lambda_0),\varpi_s\rangle >0.\]
D'autre part, on voit facilement que
\[\omega_M\omega'(\rho)-\omega_M(\rho)-\sum_{\alpha\in\Phi_{M_Q}(\omega_M)}
\alpha=-\sum_{\alpha\in\Phi(\omega_M\omega')}\alpha.\]
On note $\Phi_{M_Q}^+=\Phi^+\cap\Phi_{M_Q}$ et $\Phi_0^+=\Phi^+\cap\Phi
(\T_0,\M_0)$. On a
\[|D_{M_0}^{M_Q}(\gamma)|^{1/2}=\prod_{\alpha\in\Phi^+_{M_Q}-\Phi_0^+}
|\alpha(\gamma)|^{1/2}|1-\alpha^{-1}(\gamma)|,\]
\[\delta_{\QP(\R)}^{1/2}(\gamma)=\prod_{\alpha\in\Phi^+-\Phi_{M_Q}^+}
|\alpha(\gamma)|^{1/2},\]
\[\delta_{\Pa_0(\R)}^{1/2}(\gamma)=\prod_{\alpha\in\Phi^+-\Phi_0^+}
|\alpha(\gamma)|^{1/2},\]
\[\Delta_{\B\cap\M_Q}(\gamma)=\prod_{\alpha\in\Phi_{M_Q}^+}(1-\alpha^{-1}
(\gamma))\]
et
\[\Delta_{\B\cap\M_0}(\gamma)=\prod_{\alpha\in\Phi_0^+}(1-\alpha^{-1}(\gamma))
.\]
Donc
\[|D_{M_0}^{M_Q}(\gamma)|^{1/2}\delta_{\QP(\R)}^{1/2}(\gamma)\Delta_{\B\cap
\M_Q}(\gamma)^{-1}=\eta_Q(\gamma)\delta_{\Pa_0(\R)}^{1/2}(\gamma)\Delta_{\B\cap
\M_0}(\gamma)^{-1},\]
avec
\[\eta_Q(\gamma)=\prod_{\alpha\in\Phi^+_{M_Q}-\Phi_0^+}\frac{|1-\alpha^{-1}
(\gamma)|}{1-\alpha^{-1}(\gamma)}.\]
Comme $\gamma$ est elliptique dans $\M_0(\R)$, les racines de
$\Phi^+_{M_Q}-\Phi^+_0$ sont réelles ou complexes quand elles sont vues
comme caractères sur $\T_{M_0}(\C)$. Si $\alpha\in\Phi^+_{M_Q}-\Phi^+_0$ est
complexe, alors il existe $\alpha'\in\Phi^+_{M_Q}-\Phi^+_0$ telle que
$\alpha'\not=\alpha$ et $\alpha'(\gamma)=\overline{\alpha(\gamma)}$, et on a
\[\frac{|1-\alpha^{-1}(\gamma)|}{1-\alpha^{-1}(\gamma)}\frac{|1-{\alpha'}^{-1}
(\gamma)|}{1-{\alpha'}^{-1}(\gamma)}=1.\]
Si $\alpha\in\Phi^+_{M_Q}-\Phi^+_0$ est une racine réelle, alors
$|1-\alpha^{-1}(\gamma)|(1-\alpha^{-1}(\gamma))^{-1}$ est égal à $1$ si
$\alpha(\gamma)<0$ ou $\alpha(\gamma)>1$, et à $-1$ sinon. Donc, si on
note comme avant $R_\gamma^+$ l'ensemble des racines réelles $\alpha$ de
$\Phi$ telles que $\alpha(\gamma)>1$, alors
\[\eta_Q(\gamma)=(-1)^{|\Phi^+_{M_Q}\cap(-R_\gamma^+)|}.\]
On en déduit que
$|D_{M_0}^{M_Q}(\gamma)|^{1/2}\delta_{\QP(\R)}^{1/2}(\gamma)
Tr(\gamma,Lie(\N_Q,V)_{>0})$ est égal à
\[\eta_Q(\gamma)\delta_{\Pa_0(\R)}^{1/2}(\gamma)
\Delta_{\B\cap\M_0}(\gamma)^{-1}\sum_\omega\varepsilon(\omega)
(\omega\lambda)(\gamma)\prod_{\alpha\in\Phi(\omega)}\alpha^{-1}(\gamma),\]
où $\omega$ parcourt l'ensemble des éléments de $\Omega$ tels que
$\langle \omega(\lambda+\rho+\lambda_0),\varpi_s\rangle >0$ pour tout $s\in S$.

On plonge $\Sgoth_{r+2t}$ dans $\Sgoth_n$ en faisant agir $\Sgoth_{r+2t}$
de manière habituelle sur $\{1,\dots,r+2t\}$, et trivialement sur
$\{r+2t+1,\dots,n\}$. Comme $\Omega\simeq\{\pm 1\}^n\rtimes\Sgoth_n$, on en
déduit un plongement de $\Sgoth_{r+2t}$ dans $\Omega$, qu'on utilise
pour identifier $\Sgoth_{r+2t}$ à un sous-groupe de $\Omega$.
On remarque que, pour tout $\omega\in\Omega$, $\varepsilon(\omega)=
\varepsilon(P)\varepsilon(\sigma_P^{-1}\omega)$ et
$(\omega\lambda)(\gamma)=(\sigma_P^{-1}\omega\lambda)(\gamma_M)$.
Soit $\mu\in X^*(\T_0)\otimes_\Z\R$.
On lui associe l'élément $y_\mu=(y_1,\dots,y_{r+2t})$
de $\R^{r+2t}$ défini de la manière suivante : $y_1=
\langle \mu,\varpi_1\rangle $ et,
pour tout $i\in\{2,\dots,r+2t\}$,
$y_i=\langle \mu,\varpi_i\rangle -\langle \mu,\varpi_{i-1}\rangle $.
On écrit $\mu>_P0$ si, pour tout $i\in\{1,\dots,k\}$,
$\sum\limits_{j\in I_1\cup\dots\cup I_i}y_j>0$ (c'est-à-dire si et seulement
si, avec les notations de la section
\ref{lemmes_combinatoires}, $y_{\mu,P}>0$).
Alors il est facile de voir que, pour tout $\omega\in\Omega$, on a
$\langle \omega(\lambda+\rho+\lambda_0),\varpi_s\rangle >0$ pour tout $s\in S$
si et seulement si $\sigma_P^{-1}\omega(\lambda+\rho+\lambda_0)>_P0$.
Calculons $\eta_Q(\gamma)$. On a $\eta_Q(\gamma)=(-1)^{|\Phi^+_{M_Q}
\cap(-R_\gamma^+)|}$, et $\Phi_{M_Q}^+\cap(-R_\gamma^+)$ est l'ensemble des
racines $\alpha\in\Phi^+_{M_Q}$ qui sont réelles sur $\Int(g_P)(\T_M)$ et
telles que $0<\alpha(\gamma)<1$. L'ensemble $\Int(g_P)^{-1}(\Phi^+_{M_Q}\cap
(-R_\gamma^+))$ est inclus dans l'ensemble des racines de la forme $e_i-e_{i'}
$, avec $i<i'$ et $i,i'$ dans le même ensemble de $P$; en fait, c'est le
sous-ensemble des racines $\alpha$ de cette forme et vérifiant
$0<\alpha(\gamma_M)<1$. D'après les hypothèses sur $\gamma_M$, on a donc
\[|\Phi^+_{M_Q}\cap(-R_\gamma^+)|=\frac{1}{2}\sum_{i=1}^k(|I_i^+|(|I_i^+|-1)+
|I_i^-|(|I_i^-|-1)),\]
où, pour tout $i\in\{1,\dots,k\}$, $I_i^+=I_i\cap I^+$ et $I_i^-=I_i\cap I^-$.
Donc le signe $\eta_Q(\gamma)$ ne dépend que de $P$. On note $P^+=(I_1^+,
\dots,I_k^+)$, $P^-=(I_1^-,\dots,I_k^-)$ et $\varepsilon'(P^\pm)=(-1)^{\frac{1}
{2}\sum\limits_{i=1}^k|I_i^\pm|(|I_i^\pm|-1)}$. Alors
\[\eta_Q(\gamma)=\varepsilon'(P^+)\varepsilon'(P^-).\]
De plus, on a
\[\begin{array}{rcl}\displaystyle{\delta_{\Pa_0(\R)}^{1/2}(\gamma)\prod_
{\alpha\in\Phi(\omega)}\alpha^{-1}(\gamma)} & = & \displaystyle{\prod_{\alpha
\in\Phi^+}|\alpha(\gamma)|^{1/2}\prod_{\alpha\in\Phi_0^+}|\alpha(\gamma)|^
{-1/2}\prod_{\alpha\in\Phi^+\cap(-\omega\Phi^+)}\alpha^{-1}(\gamma)} \\
 & = & \displaystyle{\prod_{\alpha\in\sigma_P^{-1}\Phi^+}|\alpha(\gamma_M)|
^{1/2}\prod_{\alpha\in\Phi_0^+}|\alpha(\gamma)|^{-1/2}\prod_{\alpha\in
\sigma_P^{-1}\Phi^+\cap(-\sigma_P^{-1}\omega\Phi^+)}\alpha^{-1}(\gamma_M)}
\end{array}\]
et
\[\delta_{\Pa(\R)}^{1/2}(\gamma_M)\prod_{\alpha\in\Phi(\sigma_P^{-1}\omega)}
\alpha^{-1}(\gamma_M)=\prod_{\alpha\in\Phi^+}|\alpha(\gamma_M)|^{1/2}
\prod_{\alpha\in\Phi_M^+}|\alpha(\gamma_M)|^{-1/2}\prod_{\alpha\in\Phi^+\cap
(-\sigma_P^{-1}\omega\Phi^+)}\alpha^{-1}(\gamma_M),\]
où $\Phi_M^+=\Phi(\T_M,\B_M)$. Il est facile de voir que
$\prod\limits_{\alpha\in\Phi_0^+}\alpha(\gamma)=\prod\limits_{\alpha\in
\Phi_M^+}\alpha(\gamma_M)$. Donc
\[\delta_{\Pa_0(\R)}^{1/2}(\gamma)\prod_{\alpha\in\Phi(\omega)}\alpha^{-1}
(\gamma)\delta_{\Pa(\R)}^{-1/2}(\gamma_M)\prod_{\alpha\in\Phi(\sigma_P^{-1}
\omega)}\alpha(\gamma_M)=\prod_{\alpha\in\Phi^+\cap(-\sigma_P^{-1}\Phi^+)}
|\alpha(\gamma_M)|\alpha^{-1}(\gamma_M).\]
On rappelle qu'on a défini un élément $u\in\G(\C)$ tel que $\Int(u)$ envoie
$\T_{M,\C}$ sur le tore diagonal $\T_{0,\C}$. Pour tout $i\in\{1,\dots,n\}$,
on note $e_i:\T_{M,\C}\fl\Gr_{m,\C}$ le composé de $\Int(u)$ et du morphisme
$\T_{0,\C}\fl\Gr_{m,\C},diag(x_1,\dots,x_{2n})\fle x_i$ (ceci est cohérent
avec la définition de $e_1,\dots,e_r$ plus haut). Alors
l'ensemble $\Phi^+\cap(-\sigma_P^{-1}\Phi^+)=\Phi(\sigma_P^{-1})$ est inclus
dans l'ensemble des
$e_i-e_j$ avec $1\leq i<j\leq r+2t$ et $j\not=i+1\mod 2$ si $i\in\{r+1,\dots,
r+2t\}$ (ceci résulte directement du fait que $P\in\Par_{ord}(r,t)$). Donc
les racines de $\Phi(\sigma_P^{-1})$ sont réelles ou complexes
(mais pas imaginaires) sur $\T_M(\R)$. De plus, comme $\sigma_P$, vu comme
élément du groupe de Weyl de $\T_M(\C)$ dans $\G(\C)$, est dans le groupe de
Weyl de $\T_M(\R)$ dans $\G(\R)$ (et même dans $(\Nor_\G(\T_M)/\T_M)(\Q)$),
si $\alpha\in\Phi(\sigma_P^{-1})$ est complexe, alors sa
racine conjuguée est aussi dans $\Phi(\sigma_P^{-1})$. On en
déduit que
\[\prod_{\alpha\in\Phi(\sigma_P^{-1})}|\alpha(\gamma_M)|\alpha^{-1}(\gamma_M)
=(-1)^{|\{\alpha\in\Phi(\sigma_P^{-1})|\alpha(\gamma_M)\in\R\mbox{ et }
\alpha(\gamma_M)<0\}|}.\]
Si $\alpha=e_i-e_j$ avec $1\leq i<j\leq r+2t$, on a $\alpha(\gamma_M)\in\R$
si et seulement si $j\leq r$, et, si cette condition est vérifiée, on a
$\alpha(\gamma_M)<0$ si et seulement si $i\in I^+$ et $j\in I^-$. Donc le
signe ci-dessus ne dépend que de $P$; on le note $\varepsilon''(P)$. On a
\[\varepsilon''(P)=(-1)^{|\{(i,j)\in I^+\times I^-|i<j\mbox{ et }\sigma_P^{-1}
(i)>\sigma_P^{-1}(j)\}|}.\]
Enfin, on remarque que, comme $P\in\Par_{ord}(r,t)$, on a
$\varepsilon(P)=\varepsilon(P_r)$.
En utilisant les ordres sur $I^+$ et $I^-$
hérités de l'ordre sur $\{1,\dots,r\}$, on peut définir des signes
$\varepsilon(P^+)$ et $\varepsilon(P^-)$, et on voit facilement que
\[\varepsilon(P_r)\varepsilon''(P)=\varepsilon(P^+)\varepsilon(P^-).\]
Les calculs ci-dessus montrent que 
$|D_{M_0}^{M_Q}(\gamma)|^{1/2}\delta_{\QP(\R)}^{1/2}(\gamma)
Tr(\gamma,Lie(\N_Q,V)_{>0})$ est égal à
\begin{flushleft}$\displaystyle{
\delta_{\Pa(\R)}^{1/2}(\gamma_M)\Delta_{\B_M}
(\gamma_M)^{-1}\sum_{\omega\in\Omega}\varepsilon(\omega)(\omega\lambda)
(\gamma_M)\prod_{\alpha\in\Phi(\omega)}\alpha^{-1}(\gamma_M)
}$\end{flushleft}
\begin{flushright}$\displaystyle{
(-1)^{|P|}\varepsilon(P^+)\varepsilon(P^-)\varepsilon'(P^+)\varepsilon'(P^-)
\ungras_{\omega(\lambda+\rho+\lambda_0)>_P0}.
}$\end{flushright}

D'autre part, on a par définition (cf \ref{formule_points_fixes2}) $m_Q=1$ si
$m>0$, et $m_Q=2$ si $m=0$; donc $m_Q=m_P$.
De plus, 
\[\dim(\A_{M_0}/\A_{M_Q})=r+t-k=r+t-|S|\]
\[n_{M_0}^{M_Q}=r_1!t_1!\dots r_k!t_k!\]
\[|(\Nor'_\G(\M)/(\Nor'_\G(\M)\cap g^{-1}\M_Qg))(\Q)|=r!t!(r_1!t_1!
\dots r_k!t_k!)^{-1}.\]
Donc
\begin{flushleft}$\displaystyle{
|(\Nor'_\G(\M)/(\Nor'_\G(\M)\cap g^{-1}\M_Qg))(\Q)|^{-1}m_Q(-1)^{\dim(\A_{
M_0}/\A_{M_Q})}(n_{M_0}^{M_Q})^{-1}=
}$\end{flushleft}
\begin{flushright}$\displaystyle{
m_P(-1)^{r+t}(r!t!)^{-1}(-1)^{|S|},
}$\end{flushright}
et seul le facteur $(-1)^{|S|}$ dépend de $\QP$.
Donc le terme de $MD$ associé à $(\QP,g)$ est égal à
\begin{flushleft}$\displaystyle{
m_P(-1)^{r+t}(r!t!)^{-1}\delta_{\Pa(\R)}^{1/2}(\gamma_M)\Delta_{\B_M}
(\gamma_M)^{-1}\sum_{\omega\in\Omega}\varepsilon(\omega)(\omega\lambda)
(\gamma_M)\prod_{\alpha\in\Phi(\omega)}\alpha^{-1}(\gamma_M)
}$\end{flushleft}
\begin{flushright}$\displaystyle{
(-1)^{|P|}\varepsilon(P^+)\varepsilon(P^-)\varepsilon'(P^+)\varepsilon'(P^-)
\ungras_{\omega(\lambda+\rho+\lambda_0)>_P0}.
}$\end{flushright}

De plus, lorsque $P$ parcourt $\Par_{ord}(r,t)$, $(\QP_P,g_P)$ parcourt un
ensemble de représentants de $\Par(\M)/\sim$. 
Finalement, $MD$ est égal à
\begin{flushleft}$\displaystyle{
m_P(-1)^{r+t}(r!t!)^{-1}\delta_{\Pa(\R)}^{1/2}(\gamma_M)\Delta_{\B_M}
(\gamma_M)^{-1}\sum_{\omega\in\Omega}\varepsilon(\omega)(\omega\lambda)
(\gamma_M)\prod_{\alpha\in\Phi(\omega)}\alpha^{-1}(\gamma_M)
}$\end{flushleft}
\begin{flushright}$\displaystyle{
\sum_{P\in\Par_{ord}(r,t)}(-1)^{|P|}\varepsilon(P^+)\varepsilon(P^-)
\varepsilon'(P^+)\varepsilon'(P^-)\ungras_{\omega(\lambda+\rho+\lambda_0)>_P0}.
}$\end{flushright}

En comparant cette formule avec la formule obtenue plus haut pour
$MG$, on voit que,
pour conclure, il suffit de montrer que, pour tout $\mu\in X^*(\T_0)_\R$,
\[c(x,p(\mu))=(-1)^{r+t}\sum_{P\in\Par_{ord}(r,t)}(-1)^{|P|}\varepsilon(P^+)
\varepsilon(P^-)\varepsilon'(P^+)\varepsilon'(P^-)\ungras_{\mu>_P0}.\]
Mais ceci résulte directement du corollaire \ref{lemme_combinatoire_idiot_GSp},
appliqué à $y_\mu$ (on remarque que, si $y_\mu=(y_1,\dots,y_{r+2t})$ et
$p(\mu)=\sum\limits_{i=1}^r\mu_ie_i+\sum\limits_{j=1}^t\nu_j\alpha_j$,
alors $y_i=\mu_i$ pour tout $i\in\{1,\dots,r\}$ et $y_{r+2j-1}+y_{r+2j}=
\nu_j$ pour tout $j\in\{1,\dots,t\}$).

\end{preuve}

\section{Intégrales orbitales en $p$}
\label{partie_en_p}

Soit $p$ un nombre premier. Si $\G$ est un groupe de la forme
$\GL_{n_1}\times\dots\times\GL_{n_r}\times\G(\Sp_{2m_1}\times\SO_{2m_2})$,
on note, pour toute extension $L$ de $\Q_p$,
$\Hecke_{\G(L)}$ l'algèbre de Hecke des fonctions $\G(L)\fl\C$ à support
compact et bi-invariantes par $\G(\Of_L)$ (munie du produit de convolution).
Cette algèbre est commutative.

Soit $a\in\Nat^*$. On fixe une clôture algébrique
$\overline{\Q}_p$ de $\Q_p$ et on note $L$ l'extension non ramifiée de
degré $a$ de $\Q_p$ dans $\overline{\Q}_p$.
Soit $\sigma\in\Gal(L/\Q_p)$ le relèvement du Frobenius arithmétique. 
On fixe $n\in\Nat$, et on
note $\G=\GSp_{2n}$. Soient $\M$ un sous-groupe de Levi cuspidal standard
de $\G$ et $\Pa$ le sous-groupe parabolique standard de $\G$ de
sous-groupe de Levi $\M$.
On écrit $\M=\M_l\times\M_h\simeq\Gr_m^r\times\GL_2^t\times
\GSp_{2m}$ et on note $\Omega^\M$ le groupe d'automorphismes de $\M$
défini à la fin de \ref{infini2}.

Soient $(\M',s_M,\eta_{M,0})
\in\Ell_\G(\M)$, $(\H,s,\eta_0)$ son image dans $\Ell(\G)$ et $\M_H$ un
sous-groupe de Levi de $\H$ associé à $\M'$.
On peut supposer que $s=s_M=s_{A,B,m_1,m_2}$ (avec les notations
du lemme \ref{lemme:Levi_endoscopiques_GSp}), où $A\subset\{1,\dots,r\}$,
$B\subset\{1,\dots,t\}$ et $m_1+m_2=m$. On suppose que $m_2$ et
$|A|$ sont pairs (c'est-à-dire que $\H$ et $\M'$ sont cuspidaux). 
On fait agir $\Omega^\M$ sur $\M_H\simeq
\Gr_m^r\times\GL_2^t\times\G(\Sp_{2m_1}\times\SO_{2m_2})$ par les mêmes
formules que celles qui définissent son action sur $\M$.
On note $\varepsilon_{s_M}:\Omega^\M\fl\{\pm 1\}$ le morphisme
qui était noté $\varepsilon_\kappa$ dans \ref{infini2}.

On
note $s'_M=s_{\varnothing,\varnothing,m_1,m_2}\in Z(\widehat{\M}')$
(donc les composantes dans
$\widehat{\M}_h$ de $s_M$ et $s'_M$ sont les mêmes, et la composante dans
$\widehat{\M}_l$ de $s'_M$ est triviale). 
Alors $(\M',s'_M,\eta_{M,0})$ est un triplet endoscopique elliptique de $\M$,
qui est équivalent à $(\M',s_M,\eta_{M,0})$ en tant que $\M$-triplet
endoscopique.

Dans \cite{K-SVLR} p 179-180, Kottwitz a expliqué comment associer
à la donnée endoscopique $(\H,s,\eta_0)$ de $\G$ et à $L$ un morphisme de
``transfert tordu'' $b:\Hecke_{\G(L)}\fl\Hecke_{\H(\Q_p)}$ (attention, ce
morphisme dépend de $s$, et pas seulement de son image dans
$Z(\widehat{\G})/Z(\widehat{\H})$). Cette construction est rappelée
dans la section 9.1 de \cite{M3}. On notera
$b_M:\Hecke_{\M(L)}\fl\Hecke_{\M'(\Q_p)}$ le morphisme de transfert tordu
obtenu en utilisant le triplet endoscopique $(\M',s'_M,\eta_{M,0})$ de $\M$.

Soit $\mu:\Gr_m\fl\G$, $\lambda\fle\left(\begin{array}{cc}\lambda I_n & 0 \\
0 & I_n \end{array}\right)$. Le cocaractère $\mu$ se factorise par
$\M$, et on note encore $\mu$ le cocaractère de $\M$ déduit de $\mu$.
Soit $\varpi_L$ une uniformisante de $L$. On pose
\[\phi=\ungras_{\G(\Of_L)\mu(\varpi_L^{-1})\G(\Of_L)}\in\Hecke_{\G(L)}\]
\[\phi^\M=\ungras_{\M(\Of_L)\mu(\varpi_L^{-1})\M(\Of_L)}\in\Hecke_{\M(L)}.\]
On note $f^\H=b(\phi)\in\Hecke_{\H(\Q_p)}$, $f^\H_{\M_H}\in\Hecke_{\M_H
(\Q_p)}$ le terme constant de $f^\H$ en $\M_H$,
$\phi_\M\in\Hecke_{\M(L)}$ le terme constant de $\phi$ en $\M$,
$\psi^{\M'}=b_M(\phi_\M)\in\Hecke_{\M'(\Q_p)}=\Hecke_{\M_H(\Q_p)}$ et
$f^{\M'}=b_M(\phi^\M)\in\Hecke_{\M'(\Q_p)}=\Hecke_{\M_H(\Q_p)}$.

Enfin, on écrit $\M'=\M_H=\M_{H,l}\times\M_{H,h}$, avec
$\M_{H,l}=\M_l$ et $\M_{H,h}\simeq\G(\Sp_{2m_1}\times\SO_{2m_2})$.
Pour tout $\gamma_H=(\lambda_1,\dots,\lambda_r,h_1,\dots,h_t,g)\in\M_H(\Q_p)
=(\Q_p^\times)^r\times\GL_2(\Q_p)^t\times\M_{H,h}(\Q_p)$ d'image
$\gamma$ dans $\M(\Q_p)$ (un tel $\gamma$ existe toujours, car $\M$ est
déployé), on note
\[I(\gamma_H)=I(\gamma)=\{i\in\{1,\dots,r\},|\lambda_i|_p=p^{-a}\}.\]

\begin{proposition}\label{prop:identite_en_p_GSp} Soit
$\gamma_H\in\M_H(\Q_p)$. Soit $\gamma$ une image de $\gamma_H$ dans
$\M(\Q_p)$.
Si $O_{\gamma_H}(f^\H_{\M_H})\not=0$, alors
$|c(\gamma_H)|_p=|c(\gamma)|_p=p^a$. De plus :
\begin{itemize}
\item[(i)] Si $\gamma\in\M_l(\Z_p)\M_h(\Q_p)$, alors
\[O_{\gamma_H}(f^\H_{\M_H})=\delta_{\Pa(\Q_p)}^{1/2}(\gamma)O_{\gamma_H}
(f^{\M'}).\]
\item[(ii)] Il existe $C\subset\{1,\dots,t\}$, indépendant du choix de $A$ et
$B$ (mais dépendant de $\gamma_H$), tel que
\[O_{\gamma_H}(f^\H_{\M_H})=(-1)^{|I(\gamma_H)\cap A|}
\varepsilon_C(s_M)O_{\gamma_H}(\psi^{\M'}),\]
où $\varepsilon_C(s_M)=(-1)^{|C\cap B|}$.
De plus, $C$ est non vide si et seulement si
$\Omega^\M(\gamma)\cap\M_l(\Z_p)\M_h(\Q_p)=\varnothing$.
\item[(iii)] Pour tout $\omega\in\Omega^\M$,
\[O_{\omega(\gamma_H)}(f^\H_{\M_H})=\varepsilon_{s_M}(\omega)
O_{\gamma_H}(f^\H_{\M_H}).\]

\end{itemize}

\end{proposition}

Pour tout $\delta\in\M(L)$ $\sigma$-semi-simple,
on définit $\alpha_p(\gamma,\delta)$ comme l'article \cite{K-SVLR} de
Kottwitz (\S 7, p 180).
En appliquant le lemme fondamental tordu dont l'énoncé est rappelé dans la
section 5.3 de \cite{M3}
(et en utilisant le
fait que, avec la normalisation de \ref{stabilisation0}, les facteurs
de transfert sont les mêmes pour $s_M$ et $s'_M$), on obtient donc le
corollaire suivant :

\begin{corollaire}\label{cor:identite_en_p_GSp}
\begin{itemize}
\item[(i)] Si $\gamma\in\M_l(\Z_p)\M_h(\Q_p)$, alors
\[SO_{\gamma_H}(f^\H_{\M_H})=\delta_{\Pa(\Q_p)}^{1/2}(\gamma)
\sum_\delta<\alpha_p(\gamma_,\delta),s'_M>\Delta_{\M_H,s_M,p}^\M(\gamma_H,
\gamma)e(\delta)TO_\delta(\phi^\M).\]
\item[(ii)] Il existe $C\subset\{1,\dots,t\}$, indépendant du choix de $A$ et
$B$ (mais dépendant de $\gamma_H$), tel que
\[SO_{\gamma_H}(f^\H_{\M_H})=\varepsilon_C(s_M)(-1)^{|I(\gamma_H)\cap A|}
\sum_\delta<\alpha_p(\gamma,
\delta),s'_M>\Delta_{\M_H,s_M,p}^\M(\gamma_H,\gamma)
e(\delta)TO_\delta(\phi_\M).\]
De plus, $C$ est non vide si et seulement si
$\Omega^\M(\gamma)\cap\M_l(\Z_p)\M_h(\Q_p)=\varnothing$.

\end{itemize}
Dans les deux sommes ci-dessus, les facteurs de transfert sont normalisés
comme dans \ref{stabilisation0}, $\delta$ parcourt l'ensemble des classes de
$\sigma$-conjugaison de $\M(L)$ telles que $\gamma$ soit conjugué dans
$\G(\overline{\Q}_p)$ à $N\delta:=\delta\sigma(\delta)\dots\sigma^{a-1}
(\delta)$, et $e(\delta)=e(\G_\delta^\sigma)$, où $\G_\delta^\sigma$ est
le $\sigma$-centralisateur de $\delta$ dans $\G$ (cf \cite{K-SVLR} p 181, ou
\ref{formule_points_fixes2}) et $e$ est le signe de \cite{K-SC}.

\end{corollaire}

\begin{preuvep} Soient $\T$ le tore diagonal de $\G$, $\B$ le sous-groupe
de Borel standard de $\G$, $\rho=\frac{1}{2}\sum\limits_{\alpha\in\Phi(\T,
\B)}\alpha$ et $\rho_M=\frac{1}{2}\sum\limits_{\alpha\in\Phi(\T,
\B\cap\M)}\alpha$. Alors 
\[<\rho,\mu>=\frac{1}{2}n(n+1)\]
\[<\rho_M,\mu>=\frac{1}{2}m(m+1).\]
On identifie $\C[X_*(\T)]$ à $\C[X^{\pm 1},X_1^{\pm 1},\dots,X_n^{\pm 1}]$
de la manière suivante ; on envoie $X$ sur le cocaractère $\mu$ et, pour
tout $i\in\{1,\dots,n\}$, on envoie $X_i$ sur le cocaractère
$\lambda\fle diag(a_1(\lambda),\dots,a_{2n}(\lambda))$, où
\[a_j(\lambda)=\left\{\begin{array}{ll}\lambda & \mbox{ si }j=i \\
\lambda^{-1} & \mbox{ si }j=2n+1-i \\
1 & \mbox{ sinon}\end{array}\right..\]
Soit $\Omega$ le groupe de Weyl de $\T$ dans $\G$ (absolu ou relatif, cela
ne change rien dans ce cas). On a $\Omega\simeq\{\pm 1\}^n\rtimes\Sgoth_n$, et
ce groupe agit sur $\C[X^{\pm 1},X_1^{\pm 1},\dots,X_n^{\pm 1}]$ de la manière
suivante :
\begin{itemize}
\item[$\bullet$] Soit $\sigma\in\Sgoth_n$. Alors $\sigma(X)=X$ et, pour
tout $i\in\{1,\dots,n\}$, $\sigma(X_i)=X_{\sigma^{-1}(i)}$.
\item[$\bullet$] Soit $e=(e_1,\dots,e_n)\in\{\pm 1\}^n$. On note
$I=\{i\in\{1,\dots,n\}|e_i=-1\}$.
Alors $e(X)=X\prod\limits_{i\in I}X_i^{-1}$ et, pour tout $i\in\{1,\dots,n\}$,
$e(X_i)=X_i^{e_i}$.

\end{itemize}
D'après le théorème 2.1.3 de \cite{K-SVTOI} (et la reformulation de ce
théorème sous la formule (2.3.4) de cet article),
la transformée de Satake de $\phi$ est
\[p^{an(n+1)/2}X^{-1}\sum_{I\subset\{1,\dots,n\}}\prod_{i\in I}X_i.\]

On identifie $\T$ aux tores diagonaux de $\M$, $\H$ et $\M_H$. On peut choisir
ces identifications (simplement en respectant l'ordre naturel sur les
coordonnées de $\T$) de manière à ce que :
\begin{itemize}
\item[$\bullet$] La transformée de Satake de $\phi_\M$ soit égale
à celle de $\phi$.
\item[$\bullet$] La transformée de Satake de $\phi^\M$ soit égale à
\[p^{am(m+1)/2}X^{-1}\sum_{I\subset\{r+2t+1,\dots,n
\}}\prod_{i\in I}X_i.\]
\item[$\bullet$] La transformée de Satake de $f^{\M'}$ soit égale à
\[p^{am(m+1)/2}X^{-a}\sum_{I\subset\{r+2t+1,\dots,n
\}}(-1)^{|I\cap K'|}\prod_{i\in I}X_i^{a},\]
où $K'=\{r+2t+m_1+1,\dots,n\}$.
\item[$\bullet$] Les transformées de Satake de $f^\H$ et $f^\H_{\M_H}$ soient
toutes les deux égales à
\[p^{an(n+1)/2}X^{-a}\sum_{I\subset\{1,\dots,n\}}(-1)^{|I\cap K|}
\prod_{i\in I}X_i^{a},\]
où $K=A\cup\{r+2j-1,r+2j,j\in B\}\cup K'$.
\item[$\bullet$] La transformée de Satake de $\psi^{\M'}$ soit égale à
\[p^{an(n+1)/2}X^{-a}\sum_{I\subset\{1,\dots,n\}}(-1)^{|I\cap K'|}
\prod_{i\in I}X_i.\]

\end{itemize}
Pour tout $I\subset\{1,\dots,r\}$, on note $\psi_I$ la fonction de transformée
de Satake $\prod\limits_{i\in I}X_i^{a}$ sur le facteur $(\Gr_m^r)(\Q_p)$
de $\M_{H,l}(\Q_p)$. Pour tout $j\in\{1,\dots,t\}$, on note $\psi_j^{(0)}$
(resp. $\psi_j^{(1)}$, resp. $\psi_j^{(2)}$) la fonction de transformée
de Satake $1$ (resp. $X_{r+2j-1}^{a}+X_{r+2j}^{a}$, resp. $X_{r+2j-1}^{a}
X_{r+2j}^{a}$) sur le $j$-ième facteur $\GL_2(\Q_p)$ de $\M_{H,l}(\Q_p)$.
Enfin, on note $\psi_h$ la fonction de transformée de Satake
$X^{-a}\sum\limits_{I\subset\{r+2t+1,\dots,n\}}(-1)^{|I\cap K'|}\prod\limits
_{i\in I}X_i^{a}$ sur $\M_{H,h}(\Q_p)$.
Alors, d'après les calculs ci-dessus :
\[f^{\M'}=p^{am(m+1)/2}\psi_\varnothing\psi_1^{(0)}\dots\psi_t^{(0)}\psi_h\]
\[\psi^{\M'}=p^{an(n+1)/2}\left(\sum_{I\subset\{1,\dots,r\}}\psi_I\right)
\left(\prod_{j=1}^t(\psi_j^{(0)}+\psi_j^{(1)}+\psi_j^{(2)}\right)\psi_h\]
\[f^\H_{\M_H}=p^{an(n+1)/2}\left(\sum_{I\subset\{1,\dots,r\}}(-1)^{|I\cap A|}
\psi_I\right)\left(\prod_{j=1}^t(\psi_j^{(0)}-\psi_j^{(1)}+\psi_j^{(2)})\right)
\psi_h.\]

Soient $\gamma_H$ et $\gamma$ comme dans l'énoncé de la proposition. Il est
clair que $O_{\gamma_H}(f^\H_{\M_H})\not=0$ seulement si $|c(\gamma_H)|_p=
p^a$, et que $c(\gamma)=c(\gamma_H)$. Dans la suite de la preuve, on suppose
que $O_{\gamma_H}(f^\H_{\M_H})\not=0$ (sinon, toutes les égalités à prouver
sont triviale).

Montrons (i). On suppose que $\gamma\in\M_l(\Z_p)\M_h(\Q_p)$. Alors
$\gamma_H\in\M_{H,l}(\Z_p)\M_{H,h}(\Q_p)$. Soient $I\subset\{1,\dots,r\}$
et $a_1,\dots,a_t\in\{0,1,2\}$. Alors
\[O_{\gamma_H}(\psi_I(\prod_{j=1}^t\psi_j^{(a_j)})\psi_h)=0\]
sauf si $I=\varnothing$ et $a_1=\dots=a_t=0$. Donc
\[O_{\gamma_H}(f^\H_{\M_H})=p^{an(n+1)/2}O_{\gamma_H}(\psi_\varnothing\psi_1^
{(0)}\dots\psi_t^{(0)}\psi_h).\]
Pour conclure, il suffit de remarquer que, comme $\gamma\in\M_l(\Z_p)\M_h
(\Q_p)$, on a $\delta_{\Pa(\Q_p)}^{1/2}=|c(\gamma)|_p^{(n(n+1)-m(m+1))/2}
=p^{an(n+1)/2}p^{-am(m+1)/2}$.

Montrons (ii).
On sait que :
\begin{itemize}
\item[-] les supports des fonctions des fonctions $\psi_I$,
$I\subset\{1,\dots,r\}$, dont deux à deux disjoints;
\item[-] pour tout $j\in\{1,\dots,t\}$, les supports des fonctions $h\fle
O_h(\psi_j^{(0)})$, $h\fle O_h(\psi_j^{(1)})$ et $h\fle O_h(\psi_j^{(2)})$
sont deux à deux disjoints (cela résulte par exemple du fait qu'un
$h\in\GL_2(\Q_p)$ qui est dans le support de $h'\fle O_{h'}(\psi_j^{(k)})$
vérifie $|\det(h)|_p=p^{-ak}$).

\end{itemize}
On en déduit qu'il existe un unique $I\subset\{1,\dots,r\}$ et d'uniques
$a_1,\dots,a_t\in\{0,1,2\}$ tels que $O_{\gamma_H}(\psi_I\psi_1^
{(a_1)}\dots\psi_t^{(a_t)}\psi_h)\not=0$. On note $C=\{j\in\{1,\dots,t\}|
a_t=1\}$. Il est clair que $I(\gamma_H)=I$. L'égalité de (ii)
résulte des formules ci-dessus pour $f^\H_{\M_H}$ et $\psi^{\M_H}$.
Enfin, il est clair d'après la définition de l'action de $\Omega^\M$ que
$C\not=\varnothing$ si et seulement si
$\Omega^\M(\gamma)\cap\M_l(\Z_p)\M_h(\Q_p)=\varnothing$.

Montrons (iii). Il est clair d'après la définition de $\psi^{\M'}$ que
$O_{\omega(\gamma_H)}(\psi^{\M'})=O_{\gamma_H}(\psi^{\M'})$ pour tout
$\omega\in\Omega^\M$. D'après (ii), il suffit donc de voir comment
$C$ et $I(\gamma_H)$ varient. Il est facile de voir que l'ensemble $C$
est le même pour tous les $\omega(\gamma_H)$, $\omega\in\Omega^\M$.
On utilise les notations $u_1,\dots,u_r,v_1,\dots,v_t$ de la fin
de \ref{infini2} pour les générateurs de $\Omega^\M$. Si $j\in\{1,\dots,t\}$,
alors $I(v_j(\gamma_H))=I(\gamma_H)$. Si $i\in I(\gamma_H)$, alors
$I(u_i(\gamma_H))=I(\gamma_H)-\{i\}$; si $i\in\{1,\dots,r\}-I(\gamma_H)$,
alors $I(u_i(\gamma_H))=I(\gamma_H)\cup\{i\}$. Le point (iii) résulte
de ces observations et de la définition de $\varepsilon_{s_M}:\Omega^\M
\fl\{\pm 1\}$.

\end{preuvep}

\section{Stabilisation}
\label{stabilisation}

On commence par rappeler les normalisations des mesures de Haar et des facteurs
de transfert que l'on utilise (ce sont celles du chapitre 5 de \cite{M3}).
On renvoie à \cite{M3} 5.3
pour le rappel de ce qui est su sur les lemmes fondamentaux et conjectures
de transfert qui apparaissent dans cet article (en résumé, tous les
résultats nécessaires sont maintenant démontrés grâce aux travaux de
Kottwitz, Clozel, Labesse, Hales, Waldspurger, Laumon-Ngo et Ngo).

\subsection{Normalisations}
\label{stabilisation0}

\subsubsection*{Mesures de Haar}

On utilise les règles suivantes pour normaliser les mesures de Haar :
\begin{itemize}
\item[(1)] Lorsque l'on est dans la situation du théorème
\ref{th:points_fixes_moi}
on utilise les normalisations définies juste avant ce théorème.
\item[(2)] Soit $\G$ un groupe réductif connexe sur $\Q$. On choisit
toujours sur $\G(\Af)$ une mesure de Haar telle que les volumes des 
sous-groupes compacts ouverts soient des nombres rationnels. Soit $p$ un nombre
premier en lequel $\G$ est non ramifié, et soit $L$ une extension finie non 
ramifiée de $\Q_p$; alors on choisit sur $\G(L)$ la mesure de Haar telle
que le volume des sous-groupes compacts hyperspéciaux soit $1$. Si on a
fixé une mesure de Haar $dg_f$ sur $\G(\Af)$, alors on choisit la mesure de
Haar $dg_\infty$ sur $\G(\R)$ telle que $dg_fdg_\infty$ soit la mesure de
Tamagawa sur $\G(\Ade)$ (cf \cite{O}).
\item[(3)] (cf \cite{K-STF:EST} 5.2)
Soient $F$ un corps local de caractéristique $0$, $\G$ un
groupe réductif connexe sur $F$ et $\gamma\in\G(F)$ semi-simple. 
On note $I=\G_\gamma$, et on choisit des mesures de Haar sur $\G(F)$ 
et $I(F)$. Si $\gamma'\in\G(F)$ est stablement conjugué à $\gamma$, alors
$I'=\G_{\gamma'}$ est une forme intérieure de $I$, donc la mesure sur
$I(F)$ donne une mesure sur $I'(F)$. Lorsque l'on forme l'intégrale orbitale
stable en $\gamma$ d'une fonction de $C_c^\infty(\G(F))$, on utilise ces
mesures sur les centralisateurs des éléments stablement conjugués à $\gamma$.
\item[(4)] Soient $F$ un corps local ou global de caractéristique $0$,
$\G$ un groupe réductif connexe sur $F$ et $(\H,s,\eta_0)$ un triplet
endoscopique elliptique de $\G$. Soit $\gamma_H\in\H(F)$ un élément
semi-simple $(\G,\H)$-régulier. On suppose qu'il existe une image $\gamma\in
\G(F)$ de $\gamma_H$. Alors $I=\G_\gamma$ est une forme intérieure
de $I_H=\H_{\gamma_H}$ (\cite{K-STF:EST} 3.1). On choisit toujours des
mesures de Haar qui se correspondent sur $I(F)$ et $I_H(F)$.
\item[(5)] Soit $\G$ un groupe réductif connexe sur $\Q$ comme
dans \ref{formule_points_fixes2}, et soit $(\gamma_0;\gamma,\delta)$ un triplet
vérifiant les conditions (C) de \ref{formule_points_fixes2} et tel que
l'invariant $\alpha(\gamma_0;\gamma,\delta)$ soit trivial. Alors on peut
associer à $(\gamma_0;\gamma,\delta)$ un groupe (réductif connexe sur $\Q$) $I$
comme dans \ref{formule_points_fixes2} tel que $I_\R$ soit une forme intérieure
de $I(\infty):=\G_{\R,\gamma_0}$. Si on a choisi une mesure de Haar sur
$I(\R)$ (par exemple en suivant la règle (2), si on a déjà une mesure de
Haar sur $I(\Af)$), alors on prend sur
$I(\infty)(\R)$ la mesure de Haar correspondante.

\end{itemize}

\subsubsection*{Facteurs de transfert}

On renvoie à \cite{K-STF:EST} pour l'énoncé des propriétés des facteurs de
transfert que l'on utilisera ici.
Noter que les facteurs de transfert ont été définis en toute
généralité (pour l'endoscopie ordinaire) par Langlands et Shelstad, cf
\cite{LS1} et \cite{LS2}. La relation de \cite{K-STF:EST} 5.6 est prouvée
dans \cite{LS2} 4.2, et la conjecture 5.3 de \cite{K-STF:EST} est prouvée
dans la proposition 1 (section 3) de \cite{K-TN}.

Soit $\G$ un groupe réductif connexe sur $\Q$, et soit
$(\H,s,\eta_0)$ un triplet endoscopique elliptique de $\G$. On choisit un
$L$-morphisme $\eta:{}^L\H\fl{}^L\G$ qui prolonge $\eta_0$. Les facteurs de
transfert locaux associés à $\eta$ ne sont définis qu'à un scalaire près. 
On suppose désormais que $\G=\GSp_{2n}$, $n\in\Nat$, et on fixe une
normalisation des facteurs de transfert.

En la place infinie, on normalise le facteur de transfert comme dans
\cite{K-SVLR} \S7 p 184-185 (cette normalisation est rappelée dans 
\ref{infini2}), en utilisant le morphisme $j$ de \ref{infini2}
et le sous-groupe de Borel standard.

Soit $p$ un nombre premier tel que $\G$ et $\H$ soient
non ramifiés en $p$. On 
normalise le facteur de transfert en $p$ comme dans \cite{K-SVLR} \S7
p 180-181. Si
$\eta$ est non ramifié en $p$, alors cette normalisation est celle donnée
par les $\Z_p$-structures sur $\G$ et $\H$, définie par Hales (cf \cite{H} II
7 et \cite{Wa3} 4.6).

On choisit les facteurs de transfert aux autres places de manière à ce que
la propriété \cite{K-STF:EST} 6.10 (b) soit vérifiée. On note $\Delta_{\H,v}^
\G$ les facteurs de transfert ainsi normalisés.

Soit $\M$ un sous-groupe de Levi de $\G$, et soit $(\M',s_M,
\eta_{M,0})\in\Ell_\G(\M)$ d'image $(\H,s,\eta_0)$ dans $\Ell(\G)$. Comme dans
\ref{endoscopie2}, \ref{infini2} et \ref{partie_en_p}), on associe à $(\M',s_M,
\eta_{M,0})$ un sous-groupe de Levi $\M_H\simeq\M'$ de $\H$
et un $L$-morphisme $\eta_M:{}^L\M_H={}^L\M'\fl{}^L\M$ qui prolonge
$\eta_{M,0}$. On va définir une normalisation des facteurs de
transfert pour $\eta_M$ associée à ces données.

En la place infinie, on normalise le facteur de transfert comme dans
\ref{infini2}, pour le sous-groupe de Borel de $\M$ intersection avec
$\M$ du sous-groupe de Borel de $\G$ choisi ci-dessus. 

Si $v$ est une place finie de $\Q$, on normalise le facteur de transfert en
$v$ pour que
\[\Delta_{\M_H,v}^{\M}(\gamma_H,\gamma)=|D_{\M_H}^\H(\gamma_H)|^{1/2}
|D_\M^\G(\gamma)|^{-1/2}\Delta_{\H,v}^\G(\gamma_H,\gamma),\]
si $\gamma_H\in\M_H(\Q_v)$ est semi-simple $\G$-régulier et $\gamma\in\M(\Q_v)$
est une image de $\gamma_H$ (cf \cite{K-NP} lemme 7.5).

On note $\Delta_{\M_H,s_M,v}^\M$ les facteurs de transfert ainsi normalisés.
On remarque que, si $p$ est une place finie où $\G$ et $\H$ sont
non ramifiés,
alors $\Delta_{\M_H,s_M,p}^\M$ ne dépend que de l'image de $(\M',s_M,
\eta_{M,0})$ dans $\Ell(\M)$ (car c'est simplement le facteur de transfert en
$p$ pour $\eta_M$ normalisé comme dans \cite{K-SVLR} p 180-181, c'est-à-dire,
si $\eta_M$ est non ramifié en $p$, avec la normalisation donnée par les
$\Z_p$-structures sur $\M$ et $\M_H$). En revanche, les facteurs de
transfert $\Delta_{\M_H,s_M,v}^\M$ en la place infinie, en les
places finies où $\G$ ou $\H$ est ramifié dépendent en général de
$(\M',s_M,\eta_{M,0})\in\Ell_\G(\M)$, et pas uniquement de son image dans
$\Ell(\M)$ (d'où le ``$s_M$'' dans la notation).

\subsection{Le côté géométrique de la formule des traces stables, d'après
Kottwitz}
\label{stabilisation1}

On rappelle ici la formule de Kottwitz sur le côté géométrique de la
formule des traces stable pour une fonction cuspidale stable en la place
infinie.
La référence est \cite{K-NP}.

On utilise les notations de \ref{infini1}.
Soit $\G$ un groupe algébrique sur $\Q$. On suppose que $\G$ est cuspidal
(c'est-à-dire que $(\G/\A_\G)_\R$ a un tore maximal anisotrope).
Soit $\K_\infty$ un 
sous-groupe compact maximal de $\G(\R)$.
Soient $\G^*$ une forme
intérieure de $\G$ sur $\Q$ qui est quasi-déployée, $\overline{\G}$ une
forme intérieure de $\G$ sur $\R$ telle que $(\overline{\G}/\A_{\G,\R})(\R)$
soit compact et $\T_e$ un tore maximal elliptique de $\G_\R$. On
note
\[\overline{v}(\G)=e(\overline{\G})\vol(\overline{\G}(\R)/\A_\G(\R)^0)\]
($e(\overline{\G})$ est le signe associé à $\overline{\G}$ dans
\cite{K-SC}), et
\[k(\G)=|Im(\H^1(\R,\T_e\cap\G_{der})\fl\H^1(\R,\T_e))|.\]
Pour tout sous-groupe de Levi $\M$ de $\G$, on note comme dans
\ref{endoscopie2}
\[n_M^G=|(\Nor_\G(\M)/\M)(\Q)|;\]
de plus, si $\gamma\in\M(\Q)$, on note
\[\overline{\iota}^M(\gamma)=|(Cent_M(\gamma)/Cent_M(\gamma)^0)(\Q)|.\]
Soit $\nu$ un quasi-caractère de $\A_G(\R)^0$. On note $\Pi_{temp}(\G(\R),\nu)$
(resp. $\Pi_{disc}(\G(\R),\nu)$) l'ensemble des éléments $\pi$ de $\Pi_{temp}
(\G(\R))$ (resp. $\Pi_{disc}(\G(\R))$) tels que la restriction à
$\A_G(\R)^0$ du caractère central de $\pi$ est $\nu$. On note
$C_c^\infty(\G(\R),\nu^{-1})$ l'ensemble des fonctions $f_\infty:\G(\R)\fl\C$
lisses à support compact modulo $\A_G(\R)^0$ et telles que, pour tout
$(z,g)\in\A_\G(\R)^0
\times \G(\R)$, $f_\infty(zg)=\nu^{-1}(z)f_\infty(g)$. On dit que $f_\infty\in 
C_c^\infty(\G(\R),\nu^{-1})$ est \emph{cuspidale stable} si $f_\infty$ est 
$\K_\infty$-finie à droite et à gauche et si la fonction
\[\Pi_{temp}(\G(\R),\nu)\fl\C,\pi\fle Tr(\pi(f_\infty))\]
est nulle en dehors de $\Pi_{disc}(\G(\R),\nu)$
et constante sur les $L$-paquets de $\Pi_{disc}(\G(\R),\nu)$.

Soit $f_\infty\in C_c^\infty(\G(\R),\nu^{-1})$. Pour tout
$L$-paquet $\Pi$ de $\Pi_{disc}(\G(\R),\nu)$, on note
$Tr(\Pi(f_\infty))=\sum\limits_{\pi\in\Pi}Tr(\pi(f_\infty))$ et
$\Theta_\Pi=\sum\limits_{\pi\in\Pi}\Theta_\pi$.
Pour tout sous-groupe de Levi cuspidal $\M$ de $\G$, on définit une fonction
$S\Phi_\M(.,f_\infty)=S\Phi_M^G(.,f_\infty)$ sur $\M(\R)$ par la formule
\[S\Phi_\M(\gamma,f_\infty)=(-1)^{\dim(\A_\M/\A_\G)}k(\M)k(\G)^{-1}\overline{v}
(\M_\gamma)^{-1}\sum_\Pi\Phi_\M(\gamma^{-1},\Theta_\Pi)Tr(\Pi(f_\infty)),\]
où $\Pi$ parcourt l'ensemble des $L$-paquets de $\Pi_{disc}(\G(\R),\nu)$ et
$\M_\gamma=\Cent_\M(\gamma)$. On a bien sûr $S\Phi_\M(\gamma,f_\infty)=0$
si $\gamma$ n'est pas semi-simple elliptique dans $\M(\R)$.
Si $\M$ n'est pas cuspidal, on pose $S\Phi_M^G=0$.

Soit $f:\G(\Ade)\fl\C$. On suppose que $f=f^\infty f_\infty$, avec $f^\infty\in
C_c^\infty(\G(\Af))$ et $f_\infty\in C_c^\infty(\G(\R),\nu^{-1})$.
Pour tous sous-groupe de Levi $\M$ de $\G$, on note
\[ST_M^G(f)=\tau(\M)\sum_\gamma\overline{\iota}^M(\gamma)^{-1}
SO_\gamma(f_\M^\infty)S\Phi_\M(\gamma,f_\infty),\]
où $\gamma$ parcourt l'ensemble des classes
de conjugaison stables dans $\M(\Q)$ qui sont semi-simples et elliptiques dans
$\M(\R)$, et $f^\infty_\M$ est le terme constant de $f^\infty$ en un
sous-groupe parabolique de $\G$ de sous-groupe de Levi $\M$.
On pose
\[ST^G(f)=\sum_\M (n_M^G)^{-1}ST_M^G(f),\]
où $\M$ parcourt l'ensemble des classes de conjugaison (sous $\G(\Q)$) de 
sous-groupes de Levi de $\G$.

Dans le cas où le groupe dérivé de $\G$ est simplement connexe et $\G$ est
quasi-déployé, Kottwitz a
montré (\cite{K-NP}, théorème 5.1) que $ST^G(f)$ est le côté géométrique
de la formule des traces stable pour $\G$ si $f_\infty$ est stable cuspidale.
Il y a deux problèmes si l'on veut utiliser ce
résultat. D'abord, \cite{K-NP} n'est pas publié. Ensuite, il
n'est connu que dans le cas où $\G^{der}$ est simplement connexe, alors
que l'on voudrait l'utiliser pour les groupes endoscopiques des
groupes symplectiques, qui ne vérifient pas cette propriété en général
(ceci n'est sans doute pas sérieux, car, dans la stabilisation de la
formule des points fixes, les termes correspondant aux éléments de
centralisateur non connexe sont automatiquement nuls). Cependant, rien
ne nous empêche d'utiliser la distribution $ST^G$ (qui est bien définie)
au lieu du
côté géométrique de la formule des traces stable, et c'est ce que nous ferons
dans la suite. 

On peut tout de même, en utilisant le même raisonement que dans l'article
\cite{Lau} de Laumon, montrer le résultat ci-dessous, qui est essentiellement
une conséquence triviale du théorème 5.1 de \cite{K-NP} et qui permettra
d'obtenir des applications de la stabilisation de la formule des points
fixes pour $\GSp_4$ et $\GSp_6$.

\begin{subproposition}\label{prop:ST=T} Notons, pour $\G$ réductif connexe
sur $\Q$, $T^G$ la distribution de la formule des traces invariante d'Arthur.
Alors, si $\G=\GSp_{2n}$ ou $\G(\Sp_{2n}\times\SO_4)$ et si
$f=f^\infty f_\infty:\G(\Ade)\fl\C$ est comme ci-dessus, avec $f_\infty$
cuspidale stable, alors $T^G(f)=ST^G(f)$.

\end{subproposition}

\begin{preuve} Supposons d'abord que $\G=\GSp_{2n}$. Alors, d'après la
preuve du lemme 3.1 de \cite{Lau}, les $\kappa$-intégrales orbitales
de $f_\infty$ sont nulles si $\kappa\not=1$ (cela résulte du fait que
$f_\infty$ est cuspidale stable et que $\Ho^1(\R,\GSp_{2n})=\{1\}$).
Le même raisonnement montre que, pour tout sous-groupe de Levi cuspidal
$\M$ de $\G$, les analogues pour $\Phi_M^G(.,f_\infty)$ des $\kappa$-intégrales
orbitales (définies comme dans le lemme 8.3.2 de \cite{M3}) s'annulent
si $\kappa\not=1$. Grâce à ces observations, il suffit pour obtenir le
résultat d'appliquer à tout Levi cuspidal $\M$ le
processus de préstabilisation de la partie elliptique de la formule des
traces, qui est décrit par exemple dans le début de la preuve du théorème
9.6 de \cite{K-STF:EST} (cf en particulier la formule (9.6.5) de loc. cit.).

Supposons maintenant que $\G=\G(\Sp_{2n}\times\SO_4)$. Remarquons d'abord
que $\GSO_4$ est isomorphe au groupe $\H=\GL_1\sous(\GL_2\times\GL_2)$ de
\cite{Lau} (défini dans (2.6) de cet article), et qu'on peut choisir un
isomorphisme qui envoie $c:\GSO_4\fl\Gr_m$ sur le morphisme
$\GL_1\sous(\GL_2\times\GL_2)$, $(g_1,g_2)\fle\det(g_1)\det(g_2)$.
\footnote{Indiquons comment construire un tel isomorphisme. (Cette
construction est celle de \cite{Di}, Chapitre IV, section 8, point 7.)
Le groupe
$\GSO_4$ est isomorphe au groupe $\GSO(M_2(\Q),\langle.,.\rangle)$, où
$M_2(\Q)$ est l'algèbre des matrices carrées de taille $2$ sur $\Q$, vue ici
simplement comme un $\Q$-espace vectoriel de dimension $4$, et
$\langle.,.\rangle$ est la forme bilinéaire symétrique sur $M_2(\Q)$ qui
envoie $(X,Y)$ sur $\Tr({}^tXJYJ^{-1})$, avec $J=\left(\begin{array}{cc}
0 & 1 \\ -1 & 0\end{array}\right)$. On définit un morphisme
$u:\GL_2\times\GL_2\fl\GL(M_2(\Q))$ en envoyant $(g_1,g_2)\in\GL_2\times\GL_2$
sur l'automorphisme linéaire $X\fle g_1Xg_2$ de $M_2(\Q)$. Alors le
noyau de $u$ est $\GL_1$ (plongé dans $\GL_2\times\GL_2$ par $\lambda\fle
(\lambda I_2,\lambda^{-1}I_2)$ comme dans \cite{Lau} (2.6)), et l'image
de $u$ est $\GSO(M_2(\Q),\langle.,.\rangle)$.
}
Soient $F$ un corps local ou global, $\M$ un sous-groupe de Levi de
$\G$ et $\gamma\in\M(F)$. On note $I=\M_\gamma$, $\G_1=\GSp_{2n}$,
$\M_1=\M\cap\G_1$ (un sous-groupe de Levi de $\G_1$) et $I_1=I\cap\M_1$.
Si $\gamma=(\gamma_1,\gamma_2)$ avec $\gamma_1\in\GSp_{2n}$ et
$\gamma_2\in\GSO_4$, alors il est clair que $\gamma_1\in\M_1(F)$ et
$I_1=\M_{1,\gamma_1}$.
En utilisant le raisonnement de le preuve du lemme 3.2 de \cite{Lau}
(avec la suite exacte $1\fl\G_1\fl\G\fl\PGL_2\times\PGL_2\fl 1$ au lieu
de la suite exacte $1\fl\Gr_m\fl\H\fl\overline{\H}\fl 1$ de loc. cit.), on
montre facilement que le morphisme évident
$\Kgoth_{\G_1}(I_1/F)\fl\Kgoth_{\G}(I/F)$ est un isomorphisme. Une fois
cette égalité connue, on peut finir la preuve exactement comme dans le
cas $\G=\GSp_{2n}$.

\end{preuve}

Enfin, on calcule $k(\G)$ pour $\G$ un groupe orthogonal-symplectique.

\begin{sublemme}\label{lemme:calcul_k}
Soient $n_1,n_2\in\Nat$ avec $n_2$ pair.
On note $n=n_1+n_2$ et $\G=\G(\Sp_{2n_1}\times
\SO_{2n_2})$. Alors
\[k(\G)=\left\{\begin{array}{ll}1 & \mbox{ si }n=0 \\
2^{n-1} & \mbox{ si }n\geq 1\mbox{ et }n_2=0 \\
2^{n-2} & \mbox{ si }n\geq 1\mbox{ et }n_2\geq 1\end{array}\right..\]

\end{sublemme}

En particulier, $k(\Gr_m)=k(\GSp_0)=1=k(\GL_2)=k(\GSp_2)=1$.

\begin{preuve} On note $\Gamma(\infty)=\Gal(\C/\R)$. Soit $\T$ un tore
maximal elliptique de $\G_\R$.
D'après \cite{K-NP} (7.9.3) (cf la preuve du lemme 5.4.2 de \cite{M3}), on a
\[k(\G)=|\pi_0(\widehat{\T}^{\Gamma(\infty)})||\pi_0(Z(\widehat{\G})^
{\Gamma(\infty)})|^{-1}.\]

On utilise le tore maximal $\T=\T_{ell}$ de $\G$ défini dans
\ref{infini1}. Comme $\G$ est déployé, $\Gamma(\infty)$ agit trivialement
sur $Z(\widehat{\G})$, donc
$|\pi_0(Z(\widehat{\G})^{\Gamma(\infty)})|=|\pi_0(Z(\widehat{\G}))|$. Or on a
déjà vu (dans la preuve du lemme \ref{lemme:Tamagawa}) que
$|\pi_0(Z(\widehat{\G}))|$ est égal à $1$ si $n_2=0$, et à $2$ si $n_2\geq 1$.
De plus, on a vu dans la remarque \ref{rq:T_Sp_SO} que $\T\simeq\G(\U(1)^n)$,
où $\U(1)$ est le groupe des éléments de norme $1$ de $E=\Q[\sqrt{-1}]$.
Donc $\widehat{\T}^{\Gamma(\infty)}=\widehat{\T}^
{\Gal(E/\Q)}$ est calculé dans le point (i) du lemme 2.3.3 de \cite{M3},
et la conclusion résulte de ce calcul.

\end{preuve}

\begin{sublemme}\label{lemme:k_tau_k_tau}(\cite{K-NP} (7.8.1))
Soient $\M$ un sous-groupe de Levi de $\G:=\GSp_{2n}$, $(\M',s_M,\eta_{M,0})
\in\Ell_\G(\M)$ et $(\H,s,\eta_0)$ le triplet endoscopique elliptique
de $\G$ associé. On suppose que $\M$, $\M'$ et $\H$ sont cuspidaux. Alors
\[\frac{\tau(\G)}{\tau(\H)}\frac{\tau(\M')}{\tau(\M)}=\frac{k(\H)}{k(\G)}
\frac{k(\M)}{k(\M')}.\]

\end{sublemme}

Il s'agit d'un résultat général, qui est démontré dans \cite{K-NP} 8.1. Ici,
il est facile de le vérifier directement en utilisant les lemmes
\ref{lemme:Tamagawa}, \ref{lemme:Levi_endoscopiques_GSp} et
\ref{lemme:calcul_k}.

\subsection{Stabilisation de la formule des points fixes}
\label{stabilisation2}

Soit $n\in\Nat$. On note $\G=\GSp_{2n}$, et on considère la donnée
de Shimura de \ref{formule_points_fixes1}.
On fixe un nombre premier $p$, un entier
naturel non nul $j$, une représentation algébrique $V$ de $\G$ et
une fonction $f^{p,\infty}=\prod\limits_{v\not=p,\infty}f_v\in
C_c^\infty(\G(\Af^p))$.
On note $\varphi:W_\R\fl\widehat{\G}\times W_\R$ un paramètre de Langlands
du $L$-paquet de la série discrète de $\G(\R)$ associé à la contragrédiente
de $V$ (cf \ref{infini1}, après la remarque \ref{rq:Phi_M}).
On identifie l'ensemble $\Ell(\G)$ de \ref{endoscopie2}
à l'ensemble des données
déterminées comme dans la proposition \ref{prop:groupes_endoscopiques_GSp}
par les entiers $n_1$ tels que $0\leq n_1\leq n$ et $n_1\not=n-1$,
et on note $\Ell^0(\G)$
le sous-ensemble des données déterminées par un $n_1$ tel que $n-n_1$ soit
pair.
Pour tout
sous-groupe de Levi cuspidal $\M$ de $\G$, on note $\Ell^0_\G(\M)$
l'ensemble des éléments de $\Ell_\G(\M)$ dont l'image dans $\Ell(\G)$ 
est équivalente à un élément de $\Ell^0(\G)$. On identifie $\Ell_\G(\M)$ et
$\Ell^0_\G(\M)$ avec les ensembles de représentants donnés par le lemme
\ref{lemme:Levi_endoscopiques_GSp}.

Soit $(\H,s,\eta_0)\in\Ell^0(\G)$. On choisit comme prolongement
$\eta:{}^L\H\fl {}^L\G$ de $\eta_0$ le morphisme $\eta_0\times id_{W_\Q}$.
On définit une fonction $f_\H^{(j)}=f_{\H,p}^{(j)}\prod\limits_{v\not=p}
f_{\H,v}\in C^\infty(\H(\Ade))$ de la manière suivante (cf la section 7 de
\cite{K-SVLR}): 
\begin{itemize}
\item[$\bullet$] Pour toute place $v\not=p,\infty$ de $\Q$, $f_{\H,v}$
est un transfert de $f_v$ (au sens habituel, qui est rappelé dans
\cite{M3} 5.3).
\item[$\bullet$] $f_{\H,p}^{(j)}$ est la fonction de $\Hecke
(\H(\Q_p),\H(\Z_p))$ obtenue par transfert tordu à partir de 
$\phi_j^{\G}$ en utilisant $\eta_0\times id_{W_{\Q_p}}$ comme prolongement de
$\eta_0$ (cf la section \ref{partie_en_p}, et cf la preuve de la proposition
\ref{prop:identite_en_p_GSp} pour le calcul explicite de la transformée de
Satake de $f_{\H,p}^{(j)}$).
\item[$\bullet$] On utilise les notations de la section \ref{infini}.
Pour tout paramètre de Langlands elliptique $\varphi_H:W_\R\fl\widehat{\H}
\rtimes W_\R$, on pose
\[f_{\varphi_H}=d(\H)^{-1}\sum_{\pi\in\Pi(\varphi_H)}f_\pi.\]
On note $\B$ le sous-groupe de Borel standard de $\G_\C$ (c'est-à-dire le
groupe des matrices triangulaires supérieures). Il détermine comme dans
\ref{infini2} un sous-ensemble $\Omega_*\subset\Omega_G$ et
une bijection $\Phi_H(\varphi)\iso\Omega_*$, 
$\varphi_H\fle\omega_*(\varphi_H)$. On prend
\[f_{\H,\infty}=(-1)^{q(\G)}\sum_{\varphi_H\in\Phi_H
(\varphi)}\det(\omega_*(\varphi_H))f_{\varphi_H}.\]

\end{itemize}

\begin{subremarque}\label{rq:trivialite_mu_s} Dans \cite{K-SVLR} \S7, on
a un facteur $<\mu,s>$ dans $f_{\H,\infty}$, où $\mu$ est le cocaractère
de $\G_\C$ déterminé par la donnée de Shimura comme dans
\ref{formule_points_fixes1}.
Ici, comme $(\H,s,\eta_0)\in\Ell^0(\G)$, on a $<\mu,s>=1$.

\end{subremarque}

On note $f^\infty=f^{\infty,p}\ungras_{\G(\Z_p)}$.

\begin{subtheoreme}\label{th:adieu_partie_lineaire_GSp} 
Soit $\M$ un sous-groupe de Levi cuspidal standard de $\G$. Alors,
pour $j$ assez grand,
\[\Tr_M(f^\infty,j)=(n_M^G)^{-1}\sum_{(\M',s_M,\eta_{M,0})\in\Ell^0_\G
(\M)}\tau(\G)\tau(\H)^{-1}|\Lambda_\G(\M',s_M,\eta_{M,0})|^{-1}
ST_{M'}^H(f_\H^{(j)}),\]
où $\Tr_M(f^\infty,j)$ est défini au-dessus du théorème
\ref{th:points_fixes_moi} et, pour tout $(\M',s_M,\eta_{M,0})\in\Ell_\G^0(\M)$,
$(\H,s,\eta_0)$
est l'élément de $\Ell^0(\G)$ correspondant.

\end{subtheoreme}

Pour $\M=\G$, ce théorème est dû à Kottwitz (\cite{K-SVLR} théorème
7.2).

\begin{subcorollaire}\label{cor:stab_FT_IC} Si $j$ est assez grand, alors
\[\Tr(Frob_p^j\times f^\infty,\Ho^*(\overline{S}^\K_{\overline{\Q}},
IC^{\K}V))=
\sum_{(\H,s,\eta_0)\in\Ell^0(\G)}\iota(\G,\H)ST^{H}(f_{\H}^{(j)}).\]

\end{subcorollaire}

\begin{preuve}
Comme $\G$ est déployé, le corollaire résulte immédiatement du
théorème ci-dessus et du lemme \ref{lemme:Levi_endoscopiques}.

\end{preuve}

\begin{subremarque} Comme dans \cite{M3} 7.1, on pourrait donner un sens
au corollaire (ou même au théorème) pour tout $j\in\Z$ et prouver que
la validité du corollaire pour $j>>0$ implique sa validité pour tout
$j\in\Z$.

\end{subremarque}

\begin{preuvet} Comme $j$ est fixé, on omet l'exposant $(j)$ dans cette
preuve. Soit $\Pa$ le sous-groupe parabolique standard de $\G$ de sous-groupe
de Levi $\M$.
On note $\M_l$ la partie linéaire de $\M$ et $\M_h$ sa partie
hermitienne. Comme dans \ref{formule_points_fixes1}, on a
$m,r_1,\dots,r_k\in\Nat$
tels que $n=m+r_1+\dots+r_k$, $\M_h\simeq\GSp_{2m}$ et $\M_l\simeq\GL_{r_1}
\times\dots\times\GL_{r_k}$.
Comme $\M$ est cuspidal, $r_i\leq 2$ pour tout $i\in\{1,
\dots,k\}$; on a donc, quitte à changer l'ordre des facteurs, $\M_l\simeq
\Gr_m^r\times\GL_2^t$, avec $r,t\in\Nat$ tels que $r+2t=n-m$.

Comme le théorème est déjà connu pour $\M=\G$, on peut supposer que
$\M\not=\G$.
D'après le lemme \ref{lemme:Levi_endoscopiques_GSp}, on a
$|\Lambda_\G(\M',s_M,\eta_{M,0})|=1$ pour tout
$(\M',s_M,\eta_{M,0})\in\Ell_\G(\M)$.
Notons, avec les conventions de l'énoncé,
\[\Tr'_M=(n_M^G)^{-1}
\sum_{(\M',s_M,\eta_{M,0})\in\Ell^0_\G(\M)}\tau(\G)\tau(\H)^
{-1}ST_{M'}^H(f_\H).\]
D'après la définition de $ST_{M'}^H$ dans \ref{stabilisation1}, on
a
\[\Tr'_M=(n_M^G)^{-1}
\tau(\G)\sum_{(\M',s_M,\eta_{M,0})\in\Ell_\G^0(\M)}\tau(\H)^{-1}
\tau(\M')\sum_{\gamma'}SO_{\gamma'}((f^\infty_\H)_{\M'})S\Phi_{M'}^H(\gamma',
f_{\H,\infty}),\]
où $\gamma'$ parcourt l'ensemble des classes de conjugaison stable semi-simples
de $\M'(\Q)$ qui sont elliptiques dans $\M'(\R)$. 
D'après la proposition \ref{prop:transfert_caracteres}
(et la remarque \ref{rq:M_M_H_regulier}), seuls les termes indexés par
une classe de conjugaison $\gamma'$ qui est $(\M,\M')$-régulière peuvent
être non nuls.
D'après le lemme \ref{lemme:param_groupes_endoscopiques}, on a
\[\Tr'_M=(n_M^G)^{-1}
\tau(\G)\sum_{\gamma_{0,M}}\sum_{\kappa\in\Kgoth_\G(I_M/\Q)_e}\tau(\M')
\tau(\H)^{-1}\psi(\gamma_{0,M},\kappa),\]
où :
\begin{itemize}
\item[$\bullet$] $\gamma_{0,M}$ parcourt l'ensemble des classes de conjugaison
stable semi-simples de $\M(\Q)$ qui sont elliptiques dans $\M(\R)$. 
\item[$\bullet$] $I_M=\M_{\gamma_{0,M}}$ et $\Kgoth_\G(I_M/\Q)$ est défini 
au-dessus du lemme \ref{lemme:param_groupes_endoscopiques}.
\item[$\bullet$] Soient $\gamma_{0,M}$ comme ci-dessus et $\kappa\in\Kgoth_\G(I_M/
\Q)$. Soit $(\M',s_M,\eta_{M,0},\gamma')$ un $\G$-quadruplet
endoscopique associé à $\kappa$ par le lemme
\ref{lemme:param_groupes_endoscopiques}; on choisit toujours
$(\M',s_M,\eta_{M,0})$ dans l'ensemble de représentants du lemme
\ref{lemme:Levi_endoscopiques_GSp} (donc $s_M$ est uniquement déterminé par
$\kappa$).
Le sous-ensemble $\Kgoth_\G(I_M/\Q)_e$
de $\Kgoth_\G(I_M/\Q)$ est l'ensemble des $\kappa$ tels que
$(\M',s_M,\eta_{M,0})\in\Ell^0_\G(\M)$. Si $\kappa\in\Kgoth_\G
(I_M/\Q)_e$, on note
\[\psi(\gamma_{0,M},\kappa)=SO_{\gamma'}
((f_\H^\infty)_{\M'})S\Phi_{M'}^H(\gamma',f_{\H,\infty}),\]
où $\H$ a la même signification qu'avant.

\end{itemize}
Donc
\[\Tr'_M=(n_M^G)^{-1}
\sum_{\gamma_{0,M}}\sum_{\kappa_M\in\Kgoth_M(I_M/\Q)}\Psi(\gamma_{0,M},\kappa_M
),\]
où, avec les mêmes notations qu'avant,
\[\Psi(\gamma_{0,M},\kappa_M)=\tau(\G)\sum_{{\kappa\in\Kgoth_\G(I_M/\Q)_e}\atop
{\kappa\fle\kappa_M}}\tau(\M')\tau(\H)^{-1}\psi(\gamma_{0,M},\kappa).\]
(Il y a un morphisme canonique $\Kgoth_\G(I_M/\Q)\fl\Kgoth_\M(I_M/\Q)$, cf la
remarque \ref{rq:Kgoth_G_M}.)

On définit le groupe d'automorphismes $\Omega^M$ de $\M$ comme
dans \ref{infini2} et \ref{partie_en_p}. Rappelons que ce groupe
agit sur $\M$ par conjugaison par des éléments de $\G(\Q)$.
On a clairement $\Omega^\M\simeq\{\pm 1\}^r\times\{\pm 1\}^t$, donc
$|\Omega^\M|=2^{r+t}$.

On fixe $\gamma_{0,M}$ et $\kappa\in\Kgoth_\G(I_M/\Q)_e$ comme ci-dessus, et on
note $\kappa_M$ l'image de $\kappa$ dans $\Kgoth_\M(I_M/\Q)$. On écrit
$\gamma_{0,M}=\gamma_{0,L}\gamma_0$, avec $\gamma_{0,L}\in\M_l(\Q)$ et
$\gamma_0\in\M_h(\Q)$.
Calculons
$\psi(\gamma_{0,M},\kappa)$. Soit $(\M',s_M,\eta_{M,0},\gamma')$
un $\G$-quadruplet endoscopique associé à $\kappa$ par le lemme
\ref{lemme:param_groupes_endoscopiques}.
On note $(\H,s,\eta_0)$ l'élément de $\Ell^0(\G)$
associé à $(\M',s_M,\eta_{M,0})$, et on définit $s'_M$ comme dans
\ref{partie_en_p} (c'est-à-dire égal à $s_M$ dans la composante
$\widehat{\M}_h$, et trivial dans la composante $\widehat{\M}_l$); on
remarque que $s'_M$ ne dépend pas de $\kappa$, mais uniquement de $\kappa_M$.
D'après le (ii) du lemme
\ref{lemme:transfert_terme_constant}, on a
\[SO_{\gamma'}((f^{p,\infty}_\H)_{\M'})=\sum_{\gamma_M}\Delta_{\M',s_M}^
{\M,\infty,p}(\gamma',\gamma_M)e(\gamma_M)O_{\gamma_P}(f^{\infty,p}_\M),\]
où ${\gamma_M}$ parcourt l'ensemble des classes de conjugaison semi-simples
de $\M(\Af^p)$ et $e({\gamma_M})=\prod\limits_{v\not=p,\infty}e(\M_{\Q_v,
\gamma_{M,v}})$ ($e$ est le signe de \cite{K-SC}).
D'après \cite{K-STF:EST} 5.6, cette dernière somme est égale à
\[\Delta_{\M',s_M}^{\M,\infty,p}(\gamma',\gamma_{0,M})\sum_{\gamma_M}<\alpha(\gamma_{0,M},
{\gamma_M}),\kappa>e({\gamma_M})O_{\gamma_M}(f^{\infty,p}_\M),\]
où $\gamma_M$ parcourt l'ensemble des classes de conjugaison de $\M(\Af^p)$
qui sont stablement conjuguées à $\gamma_{0,M}$ en toute place $v\not=p,\infty$
et $\alpha(\gamma_{0,M},{\gamma_M})$
est noté $inv(\gamma_{0,M},{\gamma_M})$ dans \cite{K-STF:EST} (l'article
\cite{K-STF:EST} ne fait qu'énoncer une conjecture, mais cette conjecture a
été démontrée depuis, voir \cite{M3} 5.3 pour des explications).

Pour tout $\gamma_M\in\M(\Af^p)$, on écrit $\gamma_M=\gamma_L\gamma$, avec
$\gamma_L\in\M_l(\Af^p)$ et $\gamma\in\M_h(\Af^p)$. Rappelons que
la partie linéaire de $\M$ est isomorphe à $\Gr_m^r\times
\GL_2^t$, et elle est aussi isomorphe à la partie linéaire de $\M'$.
Donc, si $\gamma_{0,M}$ et $\gamma_M$ sont stablement conjugués,
alors $\gamma_{0,L}$ et $\gamma_L$ sont conjugués, et la formule ci-dessus
devient
\[SO_{\gamma'}((f^{p,\infty}_\H)_{\M'})=\Delta_{\M',s_M}^{\M,\infty,p}(\gamma',
\gamma_{0,M})\sum_{\gamma}<\alpha(\gamma_0,\gamma),s'_M>e(\gamma)O_{\gamma_{0,L}
\gamma}
(f^{\infty,p}_\M),\]
où $\gamma$ parcourt l'ensemble des classes de conjugaison de $\M_h(\Af^p)$
qui sont stablement conjuguées à $\gamma_0$ en toute place $v\not=p,\infty$
et $\alpha(\gamma_0,\gamma)$, $e(\gamma)$ sont définis comme ci-dessus.
En particulier, l'expression ci-dessus ne dépend que de
$\kappa_M$ (et pas de $\kappa$).

D'autre part, d'après le corollaire \ref{cor:identite_en_p_GSp} (et avec les
notations de ce corollaire), on a
\[SO_{\gamma'}((f_{\H,p})_{\M'})
=\varepsilon_C(s_M)(-1)^{|I(\gamma_{0,M})\cap A|}
\sum_{\delta_M}<\alpha_p({\gamma_M},
{\delta_M}),s'_M>\Delta_{\M_H,s_M,p}^\M(\gamma_H,{\gamma_M})
e({\delta_M})TO_{\delta_M}(\phi_\M),\]
où ${\delta_M}$ parcourt l'ensemble des classes
de $\sigma$-conjugaison de $\M(L)$ telles que $\gamma_{0,M}\in\Norme{\delta_M}$.
Si de plus $\gamma_{0,M}\in\M_l(\Z_p)\M_h(\Q_p)$ (ie $\gamma_{0,L}\in
\M_l(\Z_p)$), alors, d'après le même corollaire,
\[SO_{\gamma'}((f_{\H,p})_{\M'})=\delta_{\Pa(\Q_p)}^{1/2}(\gamma_{0,M})\sum_
{\delta_M}<\alpha_p(\gamma_{0,M},{\delta_M}),s'_M>\Delta_{\M',s_M,p}^{\M}(\gamma',
\gamma_{0,M})e({\delta_M})TO_{\delta_M}(\phi_j^{\M}),\]
où ${\delta_M}$ parcourt le même ensemble d'indices; gr{\^a}ce au résultat
principal de \cite{K-BC} et au fait que $\phi_j^{\M}=\ungras_{\M_h(\Of_L)}
\phi_j^{\M_h}$ (cette dernière fonction est définie dans
\ref{formule_points_fixes2}), cette expression est égale à
\[\delta_{\Pa(\Q_p)}^{1/2}(\gamma_{0,M})\sum_\delta<\alpha_p(\gamma_0,\delta),
s'_M>\Delta_{\M',s_M,p}^{\M}(\gamma',\gamma_{0,M})e(\delta)O_{\gamma_{0,L}}
(\ungras_{\M_h(\Z_p)})TO_{\delta}(\phi_j^{\M_h}),\]
où $\delta$ parcourt l'ensemble des classes
de $\sigma$-conjugaison de $\M_h(L)$ telles que $\gamma_0\in\Norme{\delta}$
et $\alpha_p(\gamma_0,\delta)$, $e(\delta)$ sont définis de manière analogue
à $\alpha_p(\gamma_{0,M},\delta_M)$, $e(\delta_M)$.
On remarque que la seule partie de cette dernière expression qui pourrait
dépendre de $\kappa$ est $\Delta_{\M',s_M,p}^{\M}(\gamma',\gamma_{0,M})$.

Enfin, on a $<\alpha_\infty(\gamma_0),s'_M>=1$, où $\alpha_\infty(\gamma_0)$
est comme dans \cite{K-SVLR} p 167 (cf
la remarque \ref{rq:trivialite_mu_s}).
On en déduit que $\psi(\gamma_{0,M},\kappa)$ est égal au produit de
\[\psi_\infty(\gamma_{0,M},\kappa_M):=
\Delta_{\M',s_M,\infty}^{\M}(\gamma',\gamma_{0,M})^{-1}
S\Phi_{\M'}^{\H}(\gamma',f_{\H,\infty})\]
et de
\begin{flushleft}$\displaystyle{
\psi^\infty(\gamma_{0,M},\kappa):=\varepsilon_C(s_M)(-1)^{|I(\gamma_{0,M})\cap A|}
\sum_{\gamma}\sum_{\delta_M}<\alpha(\gamma_0,\gamma),s'_M><\alpha(\gamma_{0,M},
\delta_M),s'_M>
}$\end{flushleft}
\begin{flushright}$\displaystyle{
e(\gamma)e({\delta_M})O_{\gamma_{0,L}\gamma}(f^{\infty,p}_\M)
TO_{\delta_M}(\phi_\M),
}$\end{flushright}
où $\gamma$ et $\delta_M$ parcourent les mêmes ensembles d'indices
qu'avant. De plus, si $\gamma_{0,L}\in\M_l(\Z_p)$, alors
$\psi^\infty(\gamma_{0,M},\kappa)$ est égal à
\[\delta_{P(\Q_p)}^{1/2}(\gamma_{0,M})
\sum_{(\gamma,\delta)}<\alpha(\gamma_0;\gamma,\delta),s'_M>
e(\gamma)e(\delta)O_{\gamma_{0,L}\gamma}(f^{\infty,p}_\M)
O_{\gamma_{0,L}}(\ungras_{\M_l(\Z_p)})TO_{\delta}(\phi_j^{\M_h}),\]
où $\gamma$ et $\delta$ parcourent les mêmes ensembles d'indices
qu'avant et
$\alpha(\gamma_0;\gamma,\delta)$ est comme dans
\ref{formule_points_fixes2} (ie défini dans \cite{K-SVLR} \S2).

Par définition de $S\Phi_{\M'}^\H$ (cf \ref{stabilisation1})
et avec les notations de \ref{infini3},
$S\Phi_{\M'}^\H(\gamma',f_{\H,\infty})$ est égal à
\[(-1)^{\dim(\A_M/\A_G)+q(\G)}
k(\M')k(\H)^{-1}\overline{v}(I_M)^{-1}\varepsilon(\kappa)\sum_{\varphi_H\in
\Phi_H(\varphi)}\det(\omega_*(\varphi_H))\Phi_{\M'}^\H({\gamma'}^{-1},
S\Theta_{\varphi_H}).\]
Le signe $\varepsilon(\kappa)$ est là car les morphismes $\eta:{}^L\H\fl
{}^L\G$ (donc les ensembles $\Omega_*$ de
\ref{infini2}) ne sont pas les mêmes selon que l'on voit
$\H$ comme le premier élément d'un triplet endoscopique de
$\Ell^0(\G)$ (ce que l'on fait pour définir $f_{\H,\infty}$),
ou comme le premier élément du triplet endoscopique déterminé par
$s_M$ (ce que l'on fait dans la formule ci-dessus).

Supposons que $\Omega^\M(\gamma_{0,M})\cap\M_l(\Z_p)\M_h(\Q_p)=\varnothing$.
Alors l'ensemble $C$ qui apparaît dans le point (ii) du corollaire
\ref{cor:identite_en_p_GSp} (et dans la formule pour $\psi^\infty(\gamma_{0,M},
\kappa)$) est non vide. D'après le point (i) de la proposition
\ref{prop:calcul_Phi_M_GSp} et le lemme \ref{lemme:k_tau_k_tau},
$\Psi(\gamma_{0,M},\kappa_M)=0$ (noter que le signe qui est noté
$\varepsilon_C(s_M)$ ici et dans \ref{partie_en_p} est le signe
$\varepsilon_C(\kappa)$ de \ref{infini3}).

D'autre part, il résulte du point (ii) du lemme
\ref{lemme:transfert_terme_constant} et de la normalisation des facteurs
de transfert dans \ref{stabilisation0}, si $v\not=p,\infty$, alors
$SO_{\gamma'}((f_{\H,v})_{\M'})$ ne change pas si on remplace
$\gamma'$ par $\omega(\gamma')$, $\omega\in\Omega^\M$
(remarquer que, bien que le groupe $\Omega^\M$ n'agisse pas
sur $\M'$ par conjugaison par des éléments de $\H(\Q)$, la
fonction $\gamma'\fle |D_{M'}^H(\gamma')|_v$ est invariante par
$\Omega^\M$).
En appliquant le point (iii) de la proposition
\ref{prop:identite_en_p_GSp} et le corollaire
\ref{cor:transfert_caracteres_GSp}, on en déduit que
$\psi(\omega(\gamma_{0,M}),\kappa)=\psi(\gamma_{0,M},\kappa)$, pour tout
$\omega\in\Omega^\M$ (quel que soit $\gamma_{0,M}\in\M(\Q)$ semi-simple).

On déduit des calculs ci-dessus que
\[\Tr'_M=(n_M^G)^{-1}2^{r+t}\sum_{\gamma_{0,M}}\sum_{\kappa_M\in\Kgoth_M(I_M/\Q)}
\Psi(\gamma_{0,M},\kappa_M),\]
où, dans la première somme, on se restreint aux $\gamma_{0,M}$ qui sont dans
$\M_l(\Z_p)\M_h(\Q_p)$.

On fixe $\gamma_{0,M}$ et $\kappa_M\in\Kgoth_\M(I_M/\Q)$ comme ci-dessus, avec
$\gamma_{0,M}\in\M_l(\Z_p)\M_h(\Q_p)$. Alors, d'après le calcul des
$\psi(\gamma_{0,M},\kappa)$ ci-dessus, $\Psi(\gamma_{0,M},\kappa_M)$ est égal au
produit de
\[\Psi^\infty(\gamma_{0,M},\kappa_M):=\delta_{\Pa(\Q_p)}^{1/2}(\gamma_{0,M})
\sum_{(\gamma,\delta)}<\alpha(\gamma_0;\gamma,\delta),s'_M>
e(\gamma)e(\delta)O_{\gamma_{0,L}\gamma}(f^{\infty,p}_\M)O_{\gamma_{0,L}}
(\ungras_{\M_l(\Z_p)})TO_{\delta}(\phi_j^{\M_h})\]
(où $\gamma$ et $\delta$ parcourent les mêmes ensembles d'indices
qu'avant) et de
\[\Psi_\infty(\gamma_{0,M},\kappa_M):=
\tau(\G)\sum_{{\kappa\in\Kgoth_\G(I_M/\Q)_e}\atop{\kappa\fle\kappa_M}}\tau(\M')
\tau(\H)^{-1}\Delta_{\M',s_M,\infty}^{\M}(\gamma',\gamma_{0,M})^{-1}
S\Phi_{\M'}^{\H}(\gamma',f_{\H,\infty}).\]

On veut utiliser la proposition \ref{prop:calcul_Phi_M_GSp} pour calculer
$\Psi_\infty(\gamma_{0,M},\kappa_M)$, donc il faut vérifier les conditions
de cette proposition. Soient $e_1,\dots,e_r,\alpha_1,\dots,
\alpha_t:\M\fl\Gr_m$ comme dans l'exemple
\ref{ex:groupes_symplectiques_suite}. Comme $\gamma_{0,M}\in\M_l(\Z_p)\M_h(\Q_p)$,
on a $|e_1(\gamma_{0,M})|_p=\dots=|e_r(\gamma_{0,M})|_p=|\alpha_1(\gamma_{0,M})|_p=\dots=
|\alpha_t(\gamma_{0,M})|_p=1$.
D'après la remarque 1.7.5 de \cite{M3}, la fonction ${\gamma_M}\fle O_{{\gamma_M}}(
(f^{\infty,p})_{\M})$ sur $\M(\Af^p)$ est à support compact modulo
conjugaison. Donc il existe $D\in\R^{+*}$ (dépendant uniquement de
$\M$) tel que, pour tout ${\gamma_M}\in\M(\Af^p)$ vérifiant $O_{{\gamma_M}}
((f^{\infty,p})_{\M})\not=0$, on ait 
\[\inf(\inf_{1\leq i\leq r}|e_i({\gamma_M})^2c({\gamma_M})|_{\Af^p},\inf_{1\leq j
\leq t}|\alpha_j({\gamma_M})c({\gamma_M})|_{\Af^p})\geq D.\]

On suppose que $j$ est assez grand pour que $p^jD\geq 1$. 
Alors, si $\gamma_{0,M}$ est tel que $\Psi^\infty(\gamma_{0,M},\kappa_M)\not=0$
(et $\gamma_{0,M}\in\M_l(\Z_p)\M_h(\Q_p)$), on a, pour tout $i\in\{1,\dots,r\}$,
\[|e_i(\gamma_{0,M})^2c(\gamma_{0,M})|_\infty=|e_i(\gamma_{0,M})^2c(\gamma_{0,M})|_p^{-1}
|e_i(\gamma_{0,M})^2c(\gamma_{0,M})|_{\Af^p}^{-1}\leq p^{-j}D^{-1}\leq 1,\]
et, pour tout $j\in\{1,\dots,t\}$,
\[|\alpha_j(\gamma_{0,M})c(\gamma_{0,M})|_\infty=|\alpha_j(\gamma_{0,M})c(\gamma_{0,M})|_p^{-1}
|\alpha_j(\gamma_{0,M})c(\gamma_{0,M})|_{\Af^p}^{-1}\leq p^{-j}D^{-1}\leq 1,\]
car $|c(\gamma_{0,M})|_p=p^j$ d'après la proposition
\ref{prop:identite_en_p_GSp}.
Donc, d'après la proposition \ref{prop:calcul_Phi_M_GSp}, le lemme
\ref{lemme:kgoth_Q_R} et le lemme
\ref{lemme:k_tau_k_tau}, on a
\[\Psi_\infty(\gamma_{0,M},\kappa_M)=2^{-(r+t)}(n_M^G)\tau(\M)e(I_M(\infty))
L_M(\gamma_{0,M})\vol(I_M(\infty)
(\R)/\A_M(\R)^0)^{-1}\ungras_{c(\gamma_{0,M})>0},\]
où $I_M(\infty)$ est une forme intérieure sur $\R$ de $I_M$ telle que
$I_M(\infty)/\A_M$ soit anisotrope.

On pose $I_L=\M_{l,\gamma_{0,L}}$ et $I=\M_{h,\gamma_0}$. Soit $I_L(\infty)$
(resp. $I(\infty)$) une forme intérieure sur $\R$ de $I_L$ (resp. $I$)
qui est anisotrope modulo $\A_{I_L}$ (resp. $A_I$). 
D'après \cite{GKM} 7.10 (et le fait que, avec les notations de cet
article, $\tau(I_L)={\mathcal D}(I_L)=1$), on a
\[\chi(I_L)=e(I_L)\vol(I_L(\infty)(\R)/\A_{I_L}(\R)^0)^{-1}.\]
Donc $\Psi(\gamma_{0,M},\kappa_M)$ est égal à
\begin{flushleft}$\displaystyle{
2^{-(r+t)}(n_M^G)\tau(\M)\delta_{\Pa(\Q_p)}^{1/2}(\gamma_{0,M})L_M(\gamma_{0,M})
\vol^{-1}\ungras_{c(\gamma_{0,M})>0}
}$\end{flushleft}
\begin{flushright}$\displaystyle{
\sum_{(\gamma,\delta)}<\alpha(\gamma_0;\gamma,\delta),s'_M>e(\gamma)
e(\delta)e(I(\infty))O_{\gamma_{0,L}\gamma}(f^{\infty,p}_\M)
O_{\gamma_{0,L}}(\ungras_{\M_l(\Z_p)})TO_{\delta}(\phi_j^{\M_j}),
}$\end{flushright}
où on a écrit $\vol$ pour $\vol(I(\infty)(\R)/\A_I(\R)^0)$.

Finalement, on trouve que $\Tr'_M$ est égal à
\begin{flushleft}$\displaystyle{
\tau(\M_h)\sum_{\gamma_{0,M}}\delta_{\Pa(\Q_p)}^{1/2}(\gamma_{0,M})\chi(I_L)
L_M(\gamma_{0,M})\vol^{-1}\ungras_{c(\gamma_0)>0}
}$\end{flushleft}
\begin{flushright}$\displaystyle{
\sum_{\kappa_M}
\sum_{(\gamma,\delta)}<\alpha(\gamma_0;\gamma,\delta),s'_M>e(\gamma)
e(\delta)e(I(\infty))O_{\gamma_{0,L}\gamma}(f^{\infty,p}_\M)
O_{\gamma_{0,L}}(\ungras_{\M_l(\Z_p)})TO_{\delta}(\phi_j^\M),
}$\end{flushright}
où les ensembles d'indices sont les mêmes que plus haut (on peut
omettre la condition $\gamma_{0,M}\in\M_l(\Z_p)\M_h(\Q_p)$, car elle est
nécessaire pour que le terme associé à $\gamma_{0,M}$ dans la somme ci-dessus
soit non nul). Cette formule résulte des calculs ci-dessus et du fait que
$\tau(\M)=\tau(\M_h)$ (car $\tau(\M_l)=1$) et $c(\gamma_{0,M})=c(\gamma_0)$.
De plus, en utilisant à nouveau le fait que $\M_l$ est isomorphe à un
produit direct de groupes $\Gr_m$ et $\GL_2$, il est facile de voir que,
pour tout $\gamma_{0,M}$, $\Kgoth_\M(I_M/\Q)=\Kgoth_{\M_h}(I/\Q)$.
Donc, en appliquant le
raisonnement de \cite{K-SVLR} \S4 à $\M_h$, on trouve que $\Tr'_M$ est égal à
\[\sum_{\gamma_{0,L}}
\sum_{(\gamma_0;\gamma,\delta)\in C'_{\M_h,j}}\chi(I_L)c(\gamma_0;\gamma,\delta)
\delta_{\Pa(\Q_p)}^{1/2}(\gamma_{0,L}\gamma_0)O_{\gamma_{0,L}\gamma}(f^{\infty,p}_\M)
O_{\gamma_{0,L}}(\ungras_{\M_l(\Z_p)})TO_{\delta}(\phi_j^\M)L_M(\gamma_{0,L}
\gamma_0),\]
où $\gamma_{0,L}$ parcourt l'ensemble des classes de conjugaison de
$\M_l(\Q)$ qui sont semi-simples et elliptiques dans $\M_l(\R)$.
Ceci est l'expression pour $\Tr_M(f^\infty,j)$ donnée dans
\ref{formule_points_fixes2}.

\end{preuvet}

Comme dans le théorème ci-dessus, soit $\M$ un sous-groupe de Levi cuspidal
standard de $\G$. Soit $\gamma_M\in\M(\Q)$ semi-simple et elliptique dans
$\M(\R)$. On note $I_M=\M_{\gamma_M}$, et on utilise la notation $\Kgoth_
\G(I_M/\Q)_e$ de la preuve du théorème.
Rappelons qu'on a défini dans \ref{infini3} un sous-ensemble
$\Kgoth_\G(I_M/\R)_e$ de $\Kgoth_\G(I_M/\R)$.

\begin{sublemme}\label{lemme:kgoth_Q_R} Le morphisme canonique
$\Kgoth_\G(I_M/\Q)\fl\Kgoth_\G(I_M/\R)$ (obtenu en restreignant l'inclusion
$(Z(\widehat{I_M})/Z(\widehat{\G}))^{\Gal(\overline{\Q}/\Q)}\subset
(Z(\widehat{I_M})/Z(\widehat{\G}))^{\Gal(\C/\R)}$) induit une bijection
$\Kgoth_\G(I_M/\Q)_e\iso\Kgoth_\G(I_M/\R)_e$.

\end{sublemme}

Ce lemme résulte facilement de la description explicite des $\G$-triplets
endoscopiques elliptiques de $\M$ (sur un corps local ou global $F$)
dans le lemme
\ref{lemme:Levi_endoscopiques_GSp} et de la remarque \ref{rq:Kgoth_G_M}
(cf le début de \ref{infini3} pour plus de détails sur la manière d'en
déduire la description de $\Kgoth_\G(I_M/F)_e$ dans le cas $F=\R$; le point
est que cette description ne dépend pas de $F$).

\begin{sublemme}\label{lemme:transfert_terme_constant}(cf \cite{K-NP} 7.10
et lemme 7.6) On fixe une place $v$ de $\Q$.
Soient $\M$ un sous-groupe de Levi de $\G$, $(\M',s_M,\eta_{M,0})\in\Ell_\G
(\M)$ et $(\H,s,\eta_0)$ l'image de $(\M',s_M,\eta_{M,0})$ dans $\Ell(\G)$.
Comme le lemme \ref{lemme:Levi_endoscopiques}, on identifie $\M'$ à un
sous-groupe de Levi de $\H$.
On choisit des prolongements compatibles $\eta:{}^L\H\fl{}^L\G$ et
$\eta_M:{}^L\M'\fl{}^L\M$ de $\eta_0$ et $\eta_{M,0}$, et on normalise les
facteurs de transfert comme dans \ref{stabilisation0}.

\begin{itemize}
\item[(i)] Soit $f\in C_c(\G(\Q_v))$. Alors, pour tout $\gamma\in\M(\Q_v)$
semi-simple et $\G$-régulier, on a
\[SO_\gamma(f_\M)=|D_M^G(\gamma)|_v^{1/2}SO_\gamma(f).\]
(On rappelle que $D_M^G(\gamma)=\det(1-\Ad(\gamma),Lie(\G)/Lie(\M)$.)
\item[(ii)] Soit $f\in C_c^\infty(\G(\Q_v))$, et soit $f^\H\in C_c^\infty
(\H(\Q_v))$ un transfert de $f$. Alors, pour tout $\gamma_H\in\M'(\Q_v)$
semi-simple $(\M,\M')$-régulier, on a
\[SO_{\gamma_H}((f^\H)_{\M'})=\sum_\gamma\Delta_{\M',s_M}^\M(\gamma_H,\gamma)
e(\M_\gamma)O_\gamma(f_\M),\]
où $\gamma$ parcourt l'ensemble des classes de conjugaison semi-simples
de $\M(\Q_v)$ qui sont des images de $\gamma_H$.

\end{itemize}

\end{sublemme}

Le lemme ci-dessus est le lemme 6.3.3 de \cite{M3}, et la preuve de Kottwitz
est rappelée dans \cite{M3} loc. cit.

\section{Applications}
\label{applications}

Dans cette section, on donne une application de la formule des points fixes
stabilisée au calcul des composantes isotypiques de la cohomologie
d'intersection, dans le cas où $\G=\GSp_4$ ou $\G=\GSp_6$. On a choisi de
se restreindre à ces cas car la stabilisation du côté géométrique de la
formule des traces est alors inconditionnelle (et triviale), cf la proposition
\ref{prop:ST=T}. Si l'on accepte d'utiliser les résultats de \cite{K-NP}
(et leur extension au cas où $\G^{der}$ n'est pas simplement connexe),
alors les résultats de cette section ont des analogues (plus compliqués)
pour $\GSp_{2n}$ (voir la section 7.1 de \cite{M3}).
\footnote{Bien sûr, si on admet les conjectures d'Arthur sur
le spectre discret des groupes $\G(\Sp_{2n_1}\times\SO_{2n_2})$,
alors on peut faire encore mieux,
comme cela est expliqué dans les sections 8 à 10 de \cite{K-SVLR}.}
Notons
cependant qu'on ne peut pas utiliser la récurrence sur les groupes
endoscopiques de la section 7.2 de \cite{M3} dans le cas des groupes
symplectiques.

On suppose donc que $\G=\GSp_4$ ou $\GSp_6$. Si $\G=\GSp_4$, on note
$\H=\GSO_4$; si $\G=\GSp_6$, on note $\H=\G(\Sp_2\times\SO_4)$. Alors,
d'après la proposition \ref{prop:groupes_endoscopiques_GSp}, $\G$ et
$\H$ sont les seuls groupes endoscopiques elliptiques de $\G$ qui
apparaissent dans la stabilisation de la formule des points fixes.
On a $\iota(\G,\H)=\frac{1}{4}$ d'après les calculs de \ref{endoscopie1}.
On note $\eta:{}^L\H\fl{}^L\G$ le prolongement évident du morphisme
$\eta_0$ de \ref{endoscopie1}.

Soit $\K$ un sous-groupe compact ouvert net de $\G(\Af)$.
On fixe, comme dans \ref{formule_points_fixes2}, une représentation
algébrique irréductible $V$ de $\G$, un nombre premier $\ell$ et un
isomorphisme $\overline{\Q}_\ell\simeq\C$.
On note $\Hecke_K=\Hecke(\G(\Af),\K):=C_c^\infty(\K\sous\G(\Af)/\K)$,
et on définit un objet $W_\ell$
du groupe de Grothendieck des représentations de $\Hecke_K\times\Gal(
\overline{\Q}/\Q)$ dans un $\Q_\ell$-espace vectoriel de dimension finie par
\[W_\ell=\sum_{i\geq 0}(-1)^i[\Ho^i(\overline{S}^\K_{\overline{\Q}},IC^{\K}
V_{\overline{\Q}})].\]
On a la décomposition isotypique
de $W_\ell\otimes_{\Q_\ell}\overline{\Q}_\ell$ en tant que $\Hecke_\K$-module,
\[W_\ell\otimes_{\Q_\ell}\overline{\Q}_\ell=
\sum_{\pi_f}W_\ell(\pi_f)\otimes\pi_f^{\K},\]
où $\pi_f$ parcourt l'ensemble des classes d'équivalence de représentations
admissibles irréductibles de $\G(\Af)$ telles que $\pi_f^{\K}\not=0$ et où
les $W_\ell(\pi_f)$ sont des représentations virtuelles
de $\Gal(\overline{\Q}/\Q)$ dans des $\overline{\Q}_\ell$-espaces vectoriels
de dimension finie. Comme il n'existe qu'un nombre fini de $\pi_f$ telles que
$W_\ell(\pi_f)\not=0$, les représentations virtuelles $W_\ell(\pi_f)$ sont
définies sur une extension finie de $\Q_\ell$.

\begin{notation}
Soient $\G'$ un groupe algébrique réductif connexe sur
$\Q$ et $\xi$ un quasi-caractère sur $\A_{G'}(\R)^0$.
On note $\Pi(\G'(\Ade),\xi)$ l'ensemble des classes
d'isomorphisme de représentations admissibles irréductibles de $\G'(\Ade)$ sur
lesquelles $\A_{G'}(\R)^0$ agit par $\xi$. Pour toute $\pi\in\Pi(\G'(\Ade),
\xi)$, on note $m_{disc}(\pi)$ la multiplicité avec laquelle $\pi$
appara{\^i}t comme facteur direct dans $L^2(\G'(\Q)\sous\G'(\Ade),\xi)$
(cf \cite{A-L2}, \S2).
On note $\Pi_{disc}(\G'(\Ade),\xi)$ l'ensemble des $\pi\in\Pi(\G'(\Ade),\xi)$
telles que $m_{disc}(\pi)\not=0$.

\end{notation}

Soit $\xi_G$ le quasi-caractère par lequel
le groupe $\A_G(\R)^0$ agit sur la contragrédiente 
de $V$. On considère le morphisme
\[\varphi:W_\R\stackrel{j}{\fl}{}^L\H_\R\stackrel{\eta_\infty}{\fl}{}^L
\G_\R\stackrel{p}{\fl}{}^L(\A_G)_\R,\]
où $j$ est l'inclusion évidente, $\eta_\infty$ est induit par $\eta$
et $p$ est le dual de l'inclusion $\A_G\fl\G$. Le morphisme
$\varphi$ est le paramètre de Langlands d'un quasi-caractère sur $\A_G(\R)$,
et on note $\chi$ la restriction à $\A_G(\R)^0$
de ce quasi-caractère. Comme $\A_H=\A_G$, on peut définir un quasi-caractère
$\xi_H$ sur $\A_H(\R)^0$ par
\[\xi_H=\xi_G\chi^{-1}.\]
Ce quasi-caractère vérifie la propriété suivante : si $\varphi_H:
W_\R\fl{}^L\H_\R$ est un paramètre de Langlands correspondant à un
$L$-paquet de représentations de $\H(\R)$ de caractère central $\xi_H$
sur $\A_H(\R)^0$, alors $\eta_\infty\circ\varphi_H:W_\R\fl{}^L\G_\R$
correspond à un
$L$-paquet de représentations
de $\G(\R)$ de caractère central $\xi_G$ sur $\A_G(\R)^0$. (Cette
construction est celle de \cite{K-NP} 5.5).

On note $\Pi_G=\Pi_{disc}(\G(\Ade),\xi_G)$ et $\Pi_H=\Pi_{disc}(\H(\Ade),
\xi_H)$.

Soit $M_H$ l'ensemble des classes de $\H(\Q)$-conjugaison de cocaractères
$\mu_H:\Gr_m\fl\H$ tels que $\mu_H$, vu comme un cocaractère de $\G$
grâce à l'isomorphisme évident entre le tore diagonal de $\H$ et
celui de $\G$, soit conjugué (par $\G(\Q)$) au cocaractère $\mu:\Gr_m\fl
\G$ de \ref{formule_points_fixes1}. Pour tout $\mu_H\in M_H$, on note
$r_{-\mu_H}$ la représentation de ${}^L\H$ associée à $-\mu_H$ par
Kottwitz (cf \cite{K-SVTOI} 2.1.2, \cite{K-SVLR} p 193); la représentation
$r_{-\mu_H}$ est triviale sur $W_\Q$, et sa restriction à $\widehat{\H}$
est algébrique de plus haut poids $-\mu_H$. On a de même une représentation
$r_{-\mu}$ de ${}^L\G$ associée à $-\mu$.
On va définir une fonction $\varepsilon:M_H\fl\{\pm 1\}$. 
Si $\G=\GSp_4$, alors un système de représentants de $M_H$ est
$z\fle diag(z,z,1,1)$ et $z\fle diag(z,1,z,1)$; on envoie le premier
cocaractère sur $1$, et le deuxième sur $-1$. Si $\G=\GSp_6$, alors
un système de représentants de $M_H$ est $z\fle diag(z,z,z,1,1,1)$,
$z\fle diag(z,1,1,z,z,1)$, $z\fle diag(z,1,z,1,z,1)$ et $z\fle diag(z,z,1,z,
1,1)$; on envoie les deux premiers cocaractères sur $1$ et les deux derniers
sur $-1$.

Comme dans \ref{stabilisation2}, on associe à $V$ des fonctions cuspidales
stables $f_{\G,\infty}=f_\infty\in C^\infty(\G(\R))$ et $f_{\H,\infty}
\in C^\infty(\H(\R))$.

Pour toute représentation irréductible admissible $\pi_f$ de $\G(\Af)$,
on note
\[c_\G(\pi_f)=\sum_{{\pi_\infty\in\Pi(\G(\R)),}\atop{\pi_f\otimes\pi_\infty
\in\Pi_\G}}m_{disc}(\pi_f\otimes\pi_\infty)\Tr(\pi_\infty(f_\infty))\]
(cette somme n'a qu'un nombre fini de termes non nuls, car il n'existe
qu'un nombre fini de $\pi_\infty\in\Pi(\G(\R))$ telles que $\Tr(\pi_\infty(
f_\infty))\not=0$).
Pour toute représentation irréductible admissible $\pi_{H,f}$ de
$\H(\Af)$, on note
\[c_\H(\pi_{H,f})=\sum_{{\pi_{H,\infty}\in\Pi(\H(\R)),}\atop{\pi_{H,f}\otimes
\pi_{H,\infty}\in\Pi_\H}}m_{disc}(\pi_{H,f}\otimes\pi_{H,\infty})\Tr(\pi_{H,
\infty}(f_{\H,\infty}))\]
(cette somme est finie pour la même raison que dans le cas de $c_\G(\pi_f)$).

Soit $\pi_f=\bigotimes\limits_p\pi_p$ une représentation irréductible 
admissible de $\G(\Af)$. Alors on note $R_H(\pi_f)$ l'ensemble des classes
d'équivalence de représentations irréductibles admissibles $\pi_H=
\bigotimes\limits_p\pi_{H,p}$ de $\H(\Af)$ telles que, pour presque
tout nombre premier $p$ où $\pi_f$ et $\pi_{H,f}$ sont non ramifiées,
le morphisme $\eta:{}^L\H\fl{}^L\G$ envoie un 
paramètre de Langlands de $\pi_{H,p}$ sur un paramètre de Langlands de 
$\pi_p$.

Soit $p$ un nombre premier. On a fixé des plongements $\Q\subset\overline{\Q}
\subset\overline{\Q}_p$, qui déterminent
un morphisme $\Gal(\overline{\Q}_p/\Q_p)\fl\Gal(\overline{\Q}/\Q)$.
On note $Frob_p\in\Gal(\overline{\Q}_p/\Q_p)$ un relèvement du Frobenius 
géométrique, et on utilise la m{\^e}me notation pour son image dans
$\Gal(\overline{\Q}/\Q)$. Si $\G'$ est un groupe réductif non ramifié sur
$\Q_p$ et $\pi_p$ est une représentation non ramifiée de $\G'(\Q_p)$, on
note $\varphi_{\pi_p}:W_{\Q_p}\fl{}^L\G'$ un paramètre de Langlands de $\pi_p$.

\begin{subtheoreme}\label{th:composantes_isotypiques}
Soit $\pi_f$ une représentation
irréductible admissible de $\G(\Af)$ telle que $\pi_f^\K\not=0$. Alors
il existe une fonction $f^\infty\in C_c^\infty(\G(\Af))$ telle que,
pour presque tout nombre premier $p$ et pour tout $m\in\Z$, on ait
\begin{flushleft}$\displaystyle{\Tr(Frob_p^m,W_\ell(\pi_f))=p^{md/2}
c_\G(\pi_f)\dim(\pi_f^\K)\Tr(r_{-\mu}\circ\varphi_
{\pi_p}(Frob_p^m))}$\end{flushleft}
\begin{flushright}$\displaystyle{+\iota(\G,\H)p^{md/2}
\sum_{\pi_{H,f}\in R_H(\pi_f)}
c_\H(\pi_{H,f})\Tr(\pi_{H,f}((f^\infty)^\H))\sum_{\mu_H\in M_H}
\varepsilon(\mu_H)
\Tr(r_{-\mu_H}\circ\varphi_{\pi_{H,p}}(Frob_p^m)),}$
\end{flushright}
où $(f^\infty)^\H$ est un transfert de $f^\infty$ et
$d=\dim S^\K$ (donc $d=3$ si $\G=\GSp_4$ et $d=6$ si $\G=\GSp_6$).

\end{subtheoreme}

On pourrait facilement déduire de ce théorème des résultats sur la valeur
absolue des valeurs propres de Hecke de $\pi_f$ en presque toute place,
comme dans le théorème 7.5 de \cite{Lau} ou la section 6.2 de \cite{M3}.
Nous ne le ferons pas ici.

\begin{preuve} Il suffit de prouver l'égalité du théorème pour $m$ assez grand 
(où la signification de ``assez grand'' peut dépendre de $p$).

Pour tout ensemble fini $S$ de nombres premiers, on note
$\Ade_S=\prod\limits_{p\in S}\Q_p$, $\Ade_f^S=\prod\limits_{p\not\in S}
{}\!\!'\,\Q_p$ et $\K^S=\prod\limits_{p\not\in S}\G(\Z_p)$.
Pour tout $\K_S\subset\G(\Ade_S)$, on note $\Hecke(\G(\Ade_S),\K_S)$ et
$\Hecke(\G(\Ade_f^S),\K^S)$ les algèbres de Hecke des fonctions
à support compact sur $\G(\Ade_S)$ (resp. $\G(\Ade_f^S)$) qui sont
bi-invariantes par $\K_S$ (par $\K^S$).
Si $\pi'_f=\bigotimes{}'\pi'_p$ est une représentation irréductible
admissible de $\G(\Af)$, on note $\pi'_S=\bigotimes\limits_{p\in S}\pi'_p$
et ${\pi'}^S=\bigotimes\limits_{p\not\in S}{}\!\!'\,\pi'_p$.

Notons $R'$ l'ensemble des classes d'isomorphisme de représentations 
irréductibles admissibles $\pi'_f$ de $\G(\Af)$ telles que $\pi'_f\not\simeq
\pi_f$, que
$(\pi'_f)^\K\not=0$, et que $W_\ell(\pi'_f)\not=0$ ou
$c_\G(\pi'_f)\not=0$. Alors
$R'$ est fini, donc il existe une fonction $h$ dans $\Hecke_\K$
telle que $\Tr(\pi_f(h))=1$ et $\Tr(\pi'_f(h))=0$ pour toute $\pi'_f\in 
R'$.

Soit $T$ un ensemble fini de nombres premiers tel que $\ell\in T$,
que toutes les
représentations de $R'$ soient non ramifiées en dehors de $T$,
que $\K=\K_T\K^T$ avec $\K_T\subset\G(\Ade_T)$ et que
$h=h_T\ungras_{\K^T}$ avec $h_T\in\Hecke(\G(\Ade_T),\K_T)$. Alors, pour
toute fonction $g^T$ de $\Hecke(\G(\Ade_f^T),\K^T)$, on a
$\Tr(\pi_f(h_Tg^T))=\Tr(\pi^T(g^T))$ et $\Tr(\pi'_f(h_Tg^T))=0$ si
$\pi'_f\in R'$.

Soit $R'_H$ l'ensemble des classes d'équivalence
de représentations irréductibles admissibles $\rho_f$ de $\H(\Af)$ telles 
que $\rho_f\not\in R_H(\pi_f)$ et que
$c_\H(\rho_f)\not=0$. Alors $R'_H$ est fini, donc
il existe $g^T\in\Hecke(\G(\Ade_f^T),\K^T)$ telle que $\Tr(\pi^T(g^T))=1$ et,
pour tout $\rho_f\in R'_H$, si $k^T$ est la fonction sur $\H(\Ade_f^T)$
obtenue par transfert non ramifié à partir de $g^T$ en utilisant le morphisme 
$\eta$, on ait $\Tr(\rho^T(k^T))=0$.

Soit $S\supset T$ un ensemble fini de nombres premiers tel que $g^T=g_{S-T}
\ungras_{\K^S}$, avec $g_{S-T}$ une fonction sur $\G(\Ade_{S-T})$. 
On pose
\[f^\infty=h_Tg^T.\]
Soit $p\not\in S\cup\{\ell\}$ un nombre premier.
Alors $f^\infty=f^{\infty,p}\ungras_{\G(\Z_p)}$, donc on peut, pour
tout $m\in\Nat^*$, associer à
$f^\infty$ des fonctions $f_\G^{(m)}=f_{\G,p}^{(m)}\prod\limits
_{v\not=p}f_{\G,v}
\in C^\infty(\G(\Ade))$ et $f_\H^{(m)}=f_{\H,p}^{(m)}\prod\limits
_{v\not=p}f_{\H,v}
\in C^\infty(\H(\Ade))$ comme dans \ref{stabilisation2}. On remarque que
$\prod\limits_{v\not=\infty,p}f_{\G,v}=f^{\infty,p}$.

D'après le corollaire \ref{cor:stab_FT_IC} et la proposition
\ref{prop:ST=T}, pour $m$ assez grand,
\[\Tr(Frob_p^m\times f^\infty,W_\ell)=
T^G(f_\G^{(m)})+\iota(\G,H)T^H(f_\H^{(m)}).\]

D'après les calculs de \cite{A-L2} p 267-268, on a
\[T^G(f_\G^{(m)})=\sum_{\rho\in\Pi_\G}m_{disc}(\rho)\Tr(\rho(f^{(m)}))\]
\[T^H(f_\H^{(m)})=\sum_{\rho_H\in\Pi_\H}m_{disc}(\rho_H)\Tr(\rho_H(f_\H^{
(m)})).\]
Soit $\rho_H=\rho_H^{\infty,p}\otimes\rho_{H,p}\otimes\rho_{H,\infty}
=\rho_{H,f}\otimes\rho_{H,\infty}\in\Pi_\H$. On a
\[\Tr(\rho_H(f_\H^{(m)}))=\Tr(\rho_H^{\infty,p}(f_\H^{\infty,p}))
\Tr(\rho_{H,p}(f_{\H,p}^{(m)}))\Tr(\rho_{H,\infty}(f_{\H,\infty})).\]
Comme $f_{H,p}^{(m)}\in\Hecke(\H(\Q_p),\H(\Z_p))$, on voit que la trace
ci-dessus s'annule si $\rho_H$ est ramifiée en $p$.
Supposons que $\rho_H$ est non ramifiée en $p$. Alors on a
\[\Tr(\rho_H^{\infty,p}(f_\H^{\infty,p}))=\Tr(\rho_{H,f}(f_\H^{\infty,p}
\ungras_{\H(\Z_p)}))=\Tr(\rho_{H,f}((f^\infty)^\H))\]
(car $\rho_{H,p}^{\H(\Z_p)}$ est de dimension $1$), et, en utilisant le
calcul de $f_{\H,p}^{(m)}$ dans la preuve de la proposition
\ref{prop:identite_en_p_GSp} (où cette fonction est notée $f^\H$),
on voit que
\[\Tr(\rho_{H,p}(f_{\H,p}^{(m)}))=p^{md/2}\sum_{\mu_H\in M_H}\varepsilon(\mu_H)
\Tr(r_{-\mu_H}\circ\varphi_{\rho_{H,p}}(Frob_p^m)).\]
Enfin, d'après le choix de $f^\infty$, on a 
\[c_\H(\rho_{H,f})\Tr(\rho_{H,f}((f^\infty)^\H))=0\]
si $\rho_{H,f}\not\in R_H(\pi_f)$.

De même, pour toute $\rho=\rho_f\otimes\rho_\infty\in\Pi_\G$, on a
$\Tr(\rho(f_\G^{(m)}))=0$ si $\rho$ est ramifiée en $p$ et, si $\rho$ est non
ramifiée en $p$, alors
\[c_\G(\rho_f)\Tr(\rho_f((f_\G^{(m)})^\infty))=c_\G(\rho_f)\Tr(\rho_f(
f^\infty))=0\]
si $\rho_f\not\simeq\pi_f$ et
\[\Tr(\rho(f_\G^{(m)}))=\dim(\pi_f^{\K})\Tr(\rho_\infty(f_\infty))
p^{md/2}\Tr(r_{-\mu}\circ\varphi_{\pi_p}(Frob_p^m))\]
si $\rho_f\simeq\pi_f$.

Le théorème résulte immédiatement de ces calculs.

\end{preuve}

En appliquant les méthodes de la preuve ci-dessus à toute la cohomologie
d'intersection, on obtient la proposition suivante :

\begin{proposition}\label{prop:fonction_L} Pour tout nombre premier $p$
tel que $\K=\G(\Z_p)\K^p$ avec $\K^p\subset\G(\Af^p)$, on a
\begin{flushleft}$\displaystyle{
\log L_p(s,\Ho^*(\overline{S}^\K_{\overline{\Q}},IC^{\K}V_{\overline{\Q}}))=
\sum_{\pi_f}c_\G(\pi_f)\dim(\pi_f^\K)\log
L(s-\frac{d}{2},\pi_p,r_{-\mu})
}$\end{flushleft}
\begin{flushright}$\displaystyle{
+\sum_{\pi_{H,f}}c_\H(\pi_{H,f})\Tr(\pi_{H,f}
((\ungras_\K)^\H))
\sum_{\mu_H\in M_H}\varepsilon(\mu_H)\log L(s-\frac{d}{2},\pi_{H,p},r_{-\mu_H})
,
}$\end{flushright}
où la première (resp. deuxième) somme est sur les classes d'isomorphisme de
représentations irréductibles admissibles $\pi_f$ (resp. $\pi_{H,f}$) de
$\G(\Af)$ (resp. $\H(\Af)$), et $(\ungras_\K)^\H)$ est un transfert de
$\ungras_\K$.

\end{proposition}

\section{Appendice : lemmes combinatoires}
\label{lemmes_combinatoires}

Cet appendice contient les lemmes combinatoires qui sont utilisés dans
la preuve de la proposition \ref{prop:calcul_Phi_M_GSp}.

Soit $n\in\Nat$. On fait agir $\Sgoth_n$ sur $\R^n$ par $(\sigma,
(\lambda_1,\dots,\lambda_n))\fle (\lambda_{\sigma^{-1}(1)},\dots,\lambda_{
\sigma^{-1}(n)})$.
On note $\tau$ l'élément de $\Sgoth_n$ qui envoie $i\in\{1,\dots,n-1\}$ sur
$i+1$ et $n$ sur $1$.

Soit $\lambda=(\lambda_1,\dots,\lambda_n)\in\R^n$. On dit que $\lambda>0$
si $\lambda_1>0,\lambda_1+\lambda_2>0,\dots,\lambda_1+\dots+\lambda_n>0$,
et on note $\Sgoth(\lambda)=\{\sigma\in\Sgoth_n|\sigma(\lambda)>0\}$.

Soient $I$ un ensemble fini et $\lambda=(\lambda_i)_{i\in I}\in\R^I$.
Pour tout $J\subset I$, on note $s_J(\lambda)=\sum\limits_{i\in J}
\lambda_i$. On note $\delta(\lambda)$ le minimum des $s_J(\lambda)/|J|$, pour
$J$ parcourant l'ensemble des sous-ensembles de $I$ tels que
$s_J(\lambda)>0$ (s'il n'existe pas de tel $J$, alors $\delta(\lambda)=
-\infty$); si $\delta(\lambda)>0$, on note $N(\lambda)$ le minimum des
$|J|$, pour $J\subset I$ tel que $s_J(\lambda)/|J|=\delta(
\lambda)$.

On note $\Par(I)$ l'ensemble des partitions de $I$, et
$\Par_{ord}(I)$ l'ensemble des partitions ordonnées de
$I$. On a une application évidente d'oubli de l'ordre
$\oubli:\Par_{ord}(I)\fl\Par(I)$, et une partition $p\in\Par(I)$ a
$|p|!$ antécédents par cette application, où $|p|$ est le nombre d'ensembles
composant $p$. Pour $P\in\Par_{ord}(I)$, on note $|P|=|\oubli(P)|$.
Pour tout $k\in\Nat$, on note $\Par^k(I)$ l'ensemble des
partitions de $\Par(I)$ qui ont exactement $2k$ ou $2k+1$ ensembles
de cardinal impair, et $\Par_{ord}^k(I)=\oubli^{-1}(\Par^k(I))$.
On note $\Par^0_{\leq 2}(I)$ l'ensemble des partitions de $\Par^0(I)$ dont
tous les ensembles sont de cardinal $\leq 2$ (donc, si $p\in\Par^0_{\leq 2}
(I)$, tous les ensembles de $p$ sont de cardinal $2$, sauf peut-être un qui est
de cardinal $1$). Enfin, on note $\Dcal(I)$ l'ensemble des couples
$(I_1,I_2)$, où $I_1$ et $I_2$ sont des sous-ensembles disjoints
(éventuellement vides) de $I$ tels que $I=I_1\cup I_2$.
Si $I=\{1,\dots,n\}$, on note $\Par(I)=\Par(n)$, etc.

Soit $J$ un sous-ensemble de $I$. Si $P$ (resp. $p$) est une partition ordonnée
(resp. une partition) de $I$, on note $P\cap J$ (resp. $p\cap J$) la partition
ordonnée (resp. la partition) de $J$ que l'on obtient en intersectant les
ensembles de $P$ (resp. $p$) avec $J$ et en omettant les intersections vides.

Soit $n\in\Nat$. Soit $P=(I_1,\dots,I_k)\in\Par_{ord}(n)$. Pour tout
$i\in\{1,\dots,k\}$, on note $n_i=|I_i|$. Il existe une unique
permutation $\sigma\in\Sgoth_n$ telle que, pour tout $i\in\{0,\dots,k-1\}$,
la restriction de $\sigma^{-1}$ à $\{n_1+\dots+n_i+1,\dots,n_1+\dots+n_{i+1}\}$
soit croissante, et $\sigma^{-1}(\{n_1+\dots+n_i+1,\dots,n_1+\dots+n_{i+1}\})=
I_{i+1}$. On note cette permutation $\sigma_P$, et on pose $\varepsilon(P)=
sgn(\sigma_P)$ (où $sgn:\Sgoth_n\fl\{\pm 1\}$ est le morphisme
signature). Si $I$ est un ensemble fini totalement
ordonné et $P\in\Par_{ord}(I)$, on
utilise l'ordre sur $I$ pour définir $\varepsilon(P)$.
On remarque de plus que, si $p\in
\Par^0(I)$, alors la restriction à $\oubli^{-1}(p)$ de l'application
$\varepsilon:\Par_{ord}(I)\fl\{\pm 1\}$ est constante; on note $\varepsilon(p)$
sa valeur.
Si $(J,K)\in\Dcal(I)$, on peut lui associer une permutation $\sigma_{(J,K)}$
de la même façon (le fait que $J$ et $K$ puissent être vides n'a aucune
importance pour cette construction); on note $\varepsilon(J,K)=
sgn(\sigma_{(J,K)})$.

Soit à nouveau $I$ un ensemble fini quelconque. Pour toute $p=\{I_\alpha,\alpha
\in A\}\in\Par(I)$, on pose
\[\varepsilon'(p)=(-1)^{\frac{1}{2}\sum_{\alpha\in A}|I_\alpha|
(|I_\alpha|-1)}.\]
Si $P\in\Par_{ord}(I)$, on note $\varepsilon'(P)=\varepsilon'(\oubli(P))$.
On remarque que, pour tout $k\in\Nat$, l'application $\varepsilon'$ est
constante sur $\Par^k(I)$, de valeur $(-1)^{m-k}$, où $m$ est la valeur
entière de $n/2$.
\footnote{Je remercie le referee qui m'a signalé ce fait utile.}

Soient $I^+,I^-$ deux sous-ensembles disjoints (éventuellement vides) de
$I$.
Si $P\in\Par_{ord}(I)$, on note $P^+=P\cap I^+$ et $P^-=P\cap I^-$.
On suppose que $I$ est totalement
ordonné, et on munit $I^+$ et $I^-$ des ordres hérités de l'ordre sur
$I$. On peut alors définir, pour $P\in\Par_{ord}(I)$, des signes
$\varepsilon(P^+)$ et $\varepsilon(P^-)$.

Soit $\lambda=(\lambda_i)_{i\in I}\in\R^I$.
Pour toute $P=(I_1,\dots,I_k)\in\Par_{ord}(I)$, on note
$\lambda_P=(s_{I_1}(\lambda),\dots,s_{I_k}(\lambda))\in\R^k$. On note
\[\Par_{ord}(\lambda)=\{P=(I_1,\dots,I_k)\in\Par_{ord}(I)|\lambda_P>0\}\]
\[\Par(\lambda)=\{p=\{I_\alpha,\alpha\in A\}\in\Par(I)|\forall\alpha\in A,
s_{I_\alpha}(\lambda)>0\}.\]
On note
$\Par^0_{ord}(\lambda)=\Par_{ord}(\lambda)\cap\Par^0_{ord}(I)$ et
$\Par^0_{\leq 2}(\lambda)=\Par^0_{\leq 2}(I)\cap\Par(\lambda).$

On définit des fonctions $c_1:\R\fl\Nat$ et $c_2:\R^2\fl\Nat$ par les formules
suivantes :
\[c_1(a)=\left\{\begin{array}{ll}0 & \mbox{ si }a\leq 0 \\
1 & \mbox{ si }a>0\end{array}\right.\]
\[c_2(a,b)=\left\{\begin{array}{ll}0 & \mbox{ si }a+b\leq 0\mbox{ ou }a\leq 0\\
1 & \mbox{ si }a>0\mbox{ et }b>0 \\
2 & \mbox{ sinon}\end{array}\right..\]

On suppose que $I$ est muni d'un ordre total
(par exemple, $I\subset\{1,\dots,n\}$,
avec l'ordre habituel sur $\{1,\dots,n\}$).
Soit $p=\{I_\alpha,\alpha\in A\}\in\Par^0_{\leq 2}(I)$. Soit $\alpha\in A$. Si
$I_\alpha$ a un seul élément, on écrit $I_\alpha=\{i\}$ et on pose $c_{I_
\alpha}=c_1(\lambda_i)$; si $I_\alpha$ a deux éléments, on écrit $I_\alpha=
\{i_1,i_2\}$
avec $i_1<i_2$ et on pose $c_{I_\alpha}(\lambda)=c_2(\lambda_{i_1},\lambda_
{i_2})$. On note
\[c(p,\lambda)=\prod_{\alpha\in A}c_{I_\alpha}(\lambda)\]
(bien sûr, on a $c(p,\lambda)=0$ si $p\not\in\Par^0_{\leq 2}(\lambda)$).

Soient $n,m\in\Nat$. On note $\Par(n,m)$ le sous-ensemble de
$\Par(n+2m)$ formé des partitions $p$ de
$\{1,\dots,n+2m\}$ telles que, pour tout $i\in\{1,\dots,m\}$, $n+2i-1$ et
$n+2i$ soient dans le même ensemble de $p$. On note $\Par_{ord}(n,m)=
\oubli^{-1}(\Par(n,m))$.
Si $\lambda\in\R^{n+2m}$, on note $\Par_{ord}(n,m,\lambda)=\Par_{ord}(n,m)\cap
\Par_{ord}(\lambda)$.

\begin{proposition}\label{prop_combinatoire_idiote_GSp1}
Soient $R$ un anneau commutatif et $I$ un ensemble fini totalement
ordonné. On se donne, pour tout sous-ensemble $I'$ de $I$, des fonctions
$a_{I'},b_{I'}:\Dcal(I')\fl R$ et $c_{I'},d_{I'}:\Par_{ord}(I')\fl R$ vérifiant
la condition suivante : pour tous $P=(I_1,\dots,I_r)\in\Par_{ord}(I')$ et
$(J,K)\in\Dcal(I')$ tels qu'il existe $k\in\{1,\dots,r\}$ avec $J=I_1\cup
\dots\cup I_k$, on a
\[c_{I'}(P)=a_{I'}(J,K)c_J(P\cap J)c_K(P\cap K)\]
\[d_{I'}(P)=b_{I'}(J,K)d_J(P\cap J)d_K(P\cap K).\]

Soient $\lambda=(\lambda_i)_{i\in I}\in\R^I$ et $(I^+,I^-)\in\Dcal(I)$.
On note
$\lambda^+=(\lambda_i)_{i\in I^+}\in\R^{I^+}$ et $\lambda^-=(\lambda_i)_
{i\in I^-}\in\R^{I^-}$. Alors
\[\sum_{P\in\Par_{ord}(\lambda)}(-1)^{|P|}c_{I^+}(P\cap I^+)d_{I^-}(P\cap I^-)=
\sum_{P^+\in\Par_{ord}(\lambda^+)}(-1)^{|P^+|}c_{I^+}(P^+)\sum_{P^-\in
\Par_{ord}(\lambda^-)}(-1)^{|P^-|}d_{I^-}(P^-).\]

\end{proposition}

\begin{exemple}\label{ex_LCI}
Donnons des exemples de fonctions $a_{I'}$ et $c_{I'}$
vérifiant la condition de la proposition ci-dessus.
\begin{itemize}
\item[(1)] $a_{I'}=1$ et $c_{I'}=1$.
\item[(2)] $a_{I'}=1$
et $c_{I'}$ égale à la fonction $\Par_{ord}(I')\fl\Z$, $P\fle\varepsilon'(P)$.
\item[(3)] $a_{I'}$ égale à la fonction $\Dcal(I')\fl\Z$, $(J,K)\fle
\varepsilon(J,K)$, et $c_{I'}$ égale à la fonction $\Par_{ord}(I')\fl\Z$,
$P\fle\varepsilon(P)$.
\item[(4)] $a_{I'}$ égale à la fonction $\Dcal(I')\fl\Z$, $(J,K)\fle
\varepsilon(J,K)$, et $c_{I'}$ égale à la fonction $\Par_{ord}(I')\fl\Z$,
$P\fle\varepsilon(P)\varepsilon'(P)$.
\item[(5)] En général, il est clair que, si l'on a des fonctions $a_{I'},
b_{I'},c_{I'},d_{I'}$ comme dans la proposition ci-dessus, alors les fonctions
$a_{I'}b_{I'}$ et $c_{I'}d_{I'}$ vérifient la condition de cette proposition.
\item[(6)] Soient $a_{I'},b_{I'},c_{I'},d_{I'}$ comme dans la proposition
ci-dessus. Soient $I^+,I^-$ deux sous-ensembles disjoints de $I$. Pour
tout sous-ensemble $I'$ de $I$, on définit des fonctions $ab_{I'}:\Dcal(I')
\fl R$ et $cd_{I'}:\Par_{ord}(I')\fl R$ par les formules suivantes : pour
tous $(J,K)\in\Dcal(I')$ et $P\in\Par_{ord}(I')$,
\[ab_{I'}(J,K)=a_{I'\cap I^+}(J\cap I^+,K\cap I^+)b_{I'\cap I^-}(J\cap I^-,
K\cap I^-)\]
\[cd_{I'}(P)=c_{I'\cap I^+}(P\cap I^+)d_{I'\cap I^-}(P\cap I^-).\]
Alors il est facile de voir que ces fonctions satisfont la condition de la
proposition ci-dessus.

\end{itemize}

\end{exemple}

\begin{preuvepn}{\ref{prop_combinatoire_idiote_GSp1}}
On raisonne par récurrence sur le couple $(|I|,|\Par_{ord}(\lambda)|)$ (on
utilise l'ordre lexicographique). Si $I=\varnothing$, le résultat est évident.
Si $\Par_{ord}(\lambda)=\varnothing$, alors $\sum\limits_{i\in I}\lambda_i\leq
0$, donc $\sum\limits_{i\in I^+}\lambda_i\leq 0$ ou $\sum\limits_{i\in I^-}
\lambda_i\leq 0$, donc $\Par_{ord}(\lambda^+)=\varnothing$ ou
$\Par_{ord}(\lambda^-)=\varnothing$, et le résultat est vrai aussi. On suppose
donc que $I\not=\varnothing$ et $\Par_{ord}(\lambda)\not=\varnothing$. On
suppose aussi que $I^+\not=\varnothing$ et $I^-\not=\varnothing$, car sinon le
résultat est trivial.
On distingue deux cas : $N(\lambda)=|I|$ et $N(\lambda)<|I|$.

Traitons d'abord le cas où $N(\lambda)<|I|$. On utilise les notations $\delta$,
$J$, $\lambda'$,
$\mu$, $\nu$, etc du lemme \ref{LCI_de_base3}. Si cela est possible, on
choisit $J$ tel que $J\subset I^+$ ou $J\subset I^-$.
On note de plus
$K=I-J$, $J^\pm=J\cap I^\pm$, $K^\pm=K\cap I^\pm$, et on définit
${\lambda'}^\pm$, $\mu^\pm$, $\nu^\pm$ de manière similaire à $\lambda^\pm$.
D'après le lemme \ref{LCI_de_base3}, on a $\Par_{ord}(\lambda)=\Par_{ord}
(\lambda')\sqcup\Par'(\lambda)$, où $\Par'(\lambda)=\Par_{ord}'(I)\cap
\Par_{ord}(\lambda)$, et on a de plus une bijection naturelle
$\Par'(\lambda)\iso\Par_{ord}(\mu)\times\Par_{ord}(\nu)$. Comme
$s_J(\lambda)>0$, on a $s_{J^+}(\lambda)>0$ ou $s_{J^-}(\lambda)>0$. Quitte
à intervertir $I^+$ et $I^-$, on peut supposer que $s_{J^+}(\lambda)>0$
(donc en particulier $J^+\not=\varnothing$).
Alors, d'après le lemme \ref{LCI_de_base}, on a $s_{J^-}(\lambda)\leq 0$.
Supposons d'abord $J\cap I^-\not=\varnothing$. D'après le choix
de $J$, pour tout $L\subset I^\pm$ tel que $s_L(\lambda)>0$, on a
$s_L(\lambda)/|L|\geq\delta$, et $s_L(\lambda)/|L|>\delta$ si $L\subset J$
(on ne peut pas avoir $L=J$, car $J^+,J^-\not=\varnothing$),
donc $s_L(\lambda')=s_L(\lambda)-\delta|L\cap J|>0$.
On en
déduit que $\Par_{ord}(\lambda^\pm)=\Par_{ord}({\lambda'}^\pm)$. De plus, comme
$s_{J^-}(\lambda)\leq 0$, on a $\Par_{ord}(\mu^-)=\varnothing$.
En appliquant l'hypothèse de récurrence à $\mu\in\R^J$, on en déduit que
\[\sum_{P_1\in\Par_{ord}(\mu)}(-1)^{|P_1|}c_{J^+}(P_1\cap J^+)d_{J^-}(P_1\cap
J^-)=0,\]
donc
\begin{flushleft}$\displaystyle{
\sum_{P\in\Par'(\lambda)}(-1)^{|P|}c_{I^+}(P\cap I^+)d_{I^-}(P\cap I^-)=
a_{I^+}(J^+,K^+)b_{I^-}(J^-,K^-)\sum_{P_1\in\Par_{ord}(\mu)}(-1)^{|P_1|}
}$\end{flushleft}
\begin{flushright}$\displaystyle{
c_{J^+}(P_1\cap J^+)d_{J^-}(P_1\cap J^-)\sum_{P_2\in\Par_{ord}(\nu)}(-1)
^{|P_2|}c_{K^+}(P_2\cap K^+)d_{K^-}(P_2\cap K^-)=0
}$\end{flushright}
(on applique les hypothèse sur $c_{I^+}$ et $d_{I^-}$ aux partitions
$(J^+,K^+)$ et $(J^-,K^-)$), donc
\[\sum_{P\in\Par_{ord}(\lambda)}(-1)^{|P|}c_{I^+}(P\cap I^+)d_{I^-}(P\cap
I^-)=\sum_{P\in\Par_{ord}(\lambda')}(-1)^{|P|}c_{I^+}(P\cap I^+)d_{I^-}(P\cap
I^-).\]
L'égalité cherchée résulte de l'hypothèse de récurrence appliquée à $\lambda'$
(et du fait que $\Par_{ord}({\lambda'}^\pm)=\Par_{ord}(\lambda^\pm)$).
Supposons maintenant que $J\cap I^-=\varnothing$, c'est-à-dire que
$J\subset I^+$. En appliquant le lemme \ref{LCI_de_base3} à $I^+$, on voit
que $\Par_{ord}(\lambda^+)=\Par_{ord}({\lambda'}^+)\sqcup\Par'(\lambda^+)$,
où $\Par'(\lambda^+)=\Par'_{ord}(I^+)\cap\Par_{ord}(\lambda^+)$ et
$\Par'_{ord}(I^+)$ est défini comme $\Par'_{ord}(I)$,
mais en utilisant $J^+\subset
I^+$ (au lieu de $J\subset I$). D'autre part, on a $\mu=\mu^+$ et
$\lambda^-={\lambda'}^-=\nu^-$. Comme $\Par_{ord}(\lambda)=\Par_{ord}(\lambda')
\sqcup\Par'(\lambda)$, on a
\begin{flushleft}$\displaystyle{
\sum_{P\in\Par_{ord}(\lambda)}(-1)^{|P|}c_{I^+}(P\cap I^+)d_{I^-}(P\cap I^-)=
}$\end{flushleft}
\begin{flushright}$\displaystyle{
\sum_{P\in\Par_{ord}(\lambda')}(-1)^{|P|}c_{I^+}(P\cap I^+)d_{I^-}(P\cap I^-)
+\sum_{P\in\Par'(\lambda)}(-1)^{|P|}c_{I^+}(P\cap I^+)d_{I^-}(P\cap I^-).
}$\end{flushright}
D'après l'hypothèse sur $c_{I^+}$ (et le fait que $I^-=K^-$), on a
\begin{flushleft}$\displaystyle{
\sum_{P\in\Par'(\lambda)}(-1)^{|P|}c_{I^+}(P\cap I^+)d_{I^-}(P\cap I^-)=
}$\end{flushleft}
\begin{flushright}$\displaystyle{
a_{I^+}(J,K^+)\sum_{P_1\in\Par_{ord}(\mu)}(-1)^{|P_1|}c_J(P_1)\sum_{
P_2\in\Par_{ord}(\nu)}(-1)^{|P_2|}c_{K^+}(P_2\cap K^+)d_{I^-}(P_2\cap I^-).
}$\end{flushright}
En appliquant l'hypothèse de récurrence à $\nu$, on trouve que ceci est égal à
\[a_{I^+}(J,K^+)\sum_{P_1\in\Par_{ord}(\mu)}(-1)^{|P_1|}c_J(P_1)\sum_{
P_2^+\in\Par_{ord}(\nu^+)}(-1)^{|P_2^+|}c_{K^+}(P_2^+)\sum_{P^-\in\Par_{ord}
(\lambda^-)}(-1)^{|P^-|}d_{I^-}(P^-).\]
D'autre part, en appliquant l'hypothèse de récurrence à $\lambda'$, on trouve
\begin{flushleft}$\displaystyle{
\sum_{P\in\Par_{ord}(\lambda')}(-1)^{|P|}c_{I^+}(P\cap I^+)d_{I^-}(P\cap I^-)=
}$\end{flushleft}
\begin{flushright}$\displaystyle{
\sum_{P^+\in\Par_{ord}({\lambda'}^+)}(-1)^{|P^+|}c_{I^+}(P^+)
\sum_{P^-\in\Par_{ord}(\lambda^-)}(-1)^{|P^-|}d_{I^-}(P^-).
}$\end{flushright}
L'égalité cherchée résulte de ces calculs et du fait que $\Par_{ord}(\lambda^+)
=\Par_{ord}({\lambda'}^+)\sqcup\Par'(\lambda^+)$, avec $\Par'(\lambda^+)\iso
\Par_{ord}(\mu)\times\Par_{ord}(\nu^+)$.

Il reste à traiter le cas où $N(\lambda)=|I|$. 
D'après le lemme
\ref{LCI_de_base}, on a $s_{I^+}(\lambda)\leq 0$ ou $s_{I^-}(\lambda)\leq 0$.
Quitte à échanger $I^+$ et $I^-$, on peut supposer que $s_{I^+}(\lambda)\leq
0$. Donc le membre de droite de l'égalité de la proposition est nul, et
il s'agit de montrer que $\sum\limits_{P\in\Par_{ord}(\lambda)}(-1)^{|P|}
c_{I^+}(P\cap I^+)d_{I^-}(P\cap I^-)=0$.
On note $\Par'$ le sous-ensemble de $\Par_{ord}(\lambda)$ formé des
$P=(I_1,\dots,I_k)$ telles qu'il existe $r\in\{1,\dots,k\}$ vérifiant l'une des
deux conditions suivantes :
\begin{itemize}
\item[(a)] $r\leq k-1$, $I_r\subset I^-$ et $I_{r+1}\subset I^+$;
\item[(b)] $I_r\cap I^+\not=\varnothing$, $I_r\cap I^-\not=\varnothing$, et
$(I_1,\dots,I_{r-1},I_r\cap I^-,I_r\cap I^+,I_{r+1},\dots,I_k)\in
\Par_{ord}(\lambda)$.

\end{itemize}
On définit une application $\iota:\Par'\fl\Par'$ de la manière suivante :
Soit $P=(I_1,\dots,I_k)\in\Par'$. Soit $r$ le plus petit élément de
$\{1,\dots,k\}$ qui vérifie (a) ou (b). Si $r$ vérifie (a), on pose
$\iota(P)=(I_1,\dots,I_{r-1},I_r\cup I_{r+1},I_{r+2},\dots,I_k)$. Si $r$
vérifie (b), on pose $\iota(P)=(I_1,\dots,I_{r-1},I_r\cap I^-,I_r\cap I^+,
I_{r+1},\dots,I_k)$. Il est clair que, pour tout $P\in\Par'$, $(-1)^{|P|}=
-(-1)^{|\iota(P)|}$, $\iota(P)\cap I^\pm=P\cap I^\pm$ et $\iota(
\iota(P))=P$. Donc
\[\sum_{P\in\Par'}(-1)^{|P|}c_{I^+}(P\cap I^+)d_{I^-}(P\cap I^-)=0.\]
Il suffit donc de montrer que $\Par'=\Par_{ord}(\lambda)$. Soit
$P=(I_1,\dots,I_k)\in\Par_{ord}(\lambda)-\Par'$. On va montrer que
$P\cap I^+\in\Par_{ord}(\lambda^+)$; ceci contredit le fait que
$s_{I^+}(\lambda)\leq 0$, donc finit la preuve. On note $P\cap I^+=
(J_1,\dots,J_l)$. Montrons par récurrence sur $r$ que, pour tout
$r\in\{1,\dots,l\}$, $\sum\limits_{i=1}^rs_{J_i}(\lambda)>0$. Soit
$r\in\{1,\dots,l\}$; on suppose que l'hypothèse de récurrence est vérifiée
pour tout $r'<r$. Soit $m\in\{1,\dots,k\}$ tel que $J_r\subset I_m$.

On se place d'abord dans le cas où, pour tout $i\in\{1,\dots,m\}$,
$I_i\subset I^+$ ou $I_i\subset I^-$; en particulier, comme $I_m\cap I^+=
J_r\not=\varnothing$, on a $I_m\subset I^+$, donc $I_m=J_r$.
Comme $P\not\in\Par'$, on a
$I_i\subset I^+$ pour tout $i\in\{1,\dots,m-1\}$; donc $r=m$ et $I_i=J_i$
pour tout $i\in\{1,\dots,r\}$. Comme $P\in\Par_{ord}(\lambda)$,
$\sum\limits_{i=1}^r s_{J_i}(\lambda)=\sum\limits_{i=1}^r s_{I_i}(\lambda)>0$.

On traite maintenant le cas où il existe $i\in\{1,\dots,m\}$
tel que $I_i\cap I^+\not=\varnothing$ et $I_i\cap I^-\not=\varnothing$.
Soit $n$ le plus grand élément de $\{1,\dots,m\}$ tel que 
$I_n\cap I^+\not=\varnothing$ et $I_n\cap I^-\not=\varnothing$. Soit
$s\in\{1,\dots,r\}$ tel que $J_s=I_n\cap I^+$. On note $K=I_n-J_s=I_n\cap I^-$.
D'après la définition de $n$, on a $I_i\subset I^+$ ou $I_i\subset I^-$
pour tout $i\in\{n+1,\dots,m-1\}$. Si $m=n$, alors $I_i\subset I^+$ pour
tout $i\in\{n+1,\dots,m-1\}$ (trivialement). Si $m>n$, alors $I_m\subset
I^+$ ou $I_m\subset I^-$ (par définition de $n$), donc $I_m\subset I^+$
(car $I_m\cap I^+=J_r\not=\varnothing$); comme $P\not\in\Par'$,
on a donc forcément $I_i\subset I^+$ pour tout $i\in\{n+1,\dots,m-1\}$.
Dans les deux cas, on en déduit que $r-s=m-n$
et $I_{n+i}=J_{s+i}$ pour tout $i\in\{1,\dots,r-s\}$.
D'après l'hypothèse de récurrence (si $s\geq 2$) ou trivialement (si $s=1$),
on a $\sum\limits_{i=1}^{s-1}s_{J_i}(\lambda)\geq 0$. De plus, on a
\[\left(\sum_{i=1}^{n-1}s_{I_i}(\lambda)\right)+s_{J_s}(\lambda)+s_K(\lambda)+
\sum_{i=s+1}^rs_{J_i}(\lambda)=\sum_{i=1}^ms_{I_i}(\lambda)>0\]
(car $P\in\Par_{ord}(\lambda)$), et
\[\left(\sum_{i=1}^{n-1}s_{I_i}(\lambda)\right)+s_K(\lambda)\leq 0\]
(car $P\not\in\Par'$), donc $\sum\limits_{i=s}^{r}s_{J_i}(\lambda)>0$, et
on en déduit que $\sum\limits_{i=1}^rs_{J_i}(\lambda)>0$.

\end{preuvepn}

\begin{corollaire}
\label{prop_combinatoire_idiote_GU}
Soient $n\in\Nat^*$ et $\lambda=(\lambda_1,\dots,\lambda_n)\in\R^n$.
Alors 
\[\sum_{P\in\Par_{ord}(\lambda)}(-1)^{|P|}=\left\{
\begin{array}{ll}(-1)^n & \mbox{ si }
\lambda_r>0\mbox{ pour tout }r\in\{1,\dots,n\} \\
0 & \mbox{ sinon}\end{array}\right..\]

\end{corollaire}

\begin{preuve} On raisonne par récurrence sur $n$. Le résultat est évident
si $n=1$. On suppose que $n\geq 2$, et que le résultat du corollaire est
connu pour $n-1$. On applique la proposition
\ref{prop_combinatoire_idiote_GSp1}, avec les choix suivants pour
$a_I,b_I,c_I,d_I$ : pour tout $I$, on prend $a_I=b_I=c_I=d_I=1$.
La proposition \ref{prop_combinatoire_idiote_GSp1}, pour $I=\{1,\dots,n\}$,
$I^+=\{1,\dots,n-1\}$ et $I^-=\{n\}$, dit que :
\[\sum_{P\in\Par_{ord}(\lambda)}(-1)^{|P|}=\sum_{P\in\Par_{ord}
(\lambda_1,\dots,\lambda_{n-1})}(-1)^{|P|}\sum_{Q\in\Par_{ord}
(\lambda_n)}(-1)^{|Q|}.\]
L'égalité cherchée résulte alors de l'hypothèse de récurrence appliquée à
$(\lambda_1,\dots,\lambda_{n-1})$ (et du cas $n=1$).

\end{preuve}

\begin{proposition}\label{prop_combinatoire_idiote_GSp2}
Soient $n\in\Nat$ et $\lambda\in\R^n$. Alors
\[\sum_{P\in\Par_{ord}(\lambda)}(-1)^{|P|}\varepsilon(P)\varepsilon'(P)=
(-1)^n\sum_{p\in\Par^0_{\leq 2}(n)}\varepsilon(p)c(p,\lambda).\]

\end{proposition}

\begin{corollaire}\label{lemme_combinatoire_idiot_GSp}
Soient $n,m\in\Nat$, $\lambda=(\lambda_1,\dots,\lambda_{n+2m})\in\R^{n+2m}$
et $I^+,I^-$ des sous-ensembles disjoints de $\{1,\dots,n\}$ tels que
$\{1,\dots,n\}=I^+\cup I^-$. Pour tout $i\in\{1,\dots,m\}$, on note
$\lambda'_i=\lambda_{n+2i-1}+\lambda_{n+2i}$. Alors
\begin{flushleft}$\displaystyle{
\sum_{P\in\Par_{ord}(n,m,\lambda)}(-1)^{|P|}\varepsilon(P^+)\varepsilon(P^-)
\varepsilon'(P^+)\varepsilon'(P^-)=
}$\end{flushleft}
\begin{flushright}$\displaystyle{
(-1)^{n+m}c_1(\lambda'_1)\dots c_1(\lambda'_m)
\sum_{p^+\in\Par^0_{\leq 2}(I^+)}\sum_{p^-\in\Par^0_{\leq 2}(I^-)}\varepsilon
(p^+)\varepsilon(p^-)c(p^+,\lambda)c(p^-,\lambda),}$\end{flushright}
où, pour toute $P\in\Par_{ord}(n,m)$, $P^+=P\cap I^+$ et $P^-=P\cap I^-$.

\end{corollaire}

\begin{preuve} Soit $u:\{1,\dots,n+2m\}\fl\{1,\dots,n+m\}$ définie par les
formules suivantes : $u(i)=i$ si $1\leq i\leq n$, et, pour tout
$i\in\{1,\dots,m\}$, $u(n+2i-1)=u(n+2i)=n+i$.
Soit $\varphi:\Par_{ord}(n,m)\fl\Par_{ord}(n+m)$ l'application
qui envoie $(I_1,\dots,I_k)$ sur $(u(I_1),\dots,u(I_k))$. Il est clair que
$\varphi$ est bijective et que, pour toute $P\in\Par_{ord}(n,m)$, on a
$|P|=|\varphi(P)|$, $P\cap I^\pm=\varphi(P)\cap I^\pm$ et
$P\in\Par_{ord}(n,m,\lambda)$ si et seulement si $\varphi(P)\in\Par_{ord}
(\mu)$, où $\mu=(\lambda_1,\dots,\lambda_n,\lambda'_1,\dots,\lambda'_m)$.
Donc le membre de gauche de l'égalité de la proposition est égal à
\[\sum_{P\in\Par_{ord}(\mu)}(-1)^{|P|}\varepsilon(P\cap I^+)\varepsilon(P\cap
I^-)\varepsilon'(P\cap I^+)\varepsilon'(P\cap I^-).\]
On applique la proposition \ref{prop_combinatoire_idiote_GSp1} à
$(\{1,\dots,n\},\{n+1,\dots,n+m\})\in\Dcal(\{1,\dots,n+m\})$,
avec, pour tout $I\subset\{1,\dots,n+m\}$, $a_I:(J,K)\in\Dcal(I)\fle\varepsilon
(J\cap I^+,K\cap I^+)\varepsilon(J\cap I^-,K\cap I^-)$,
$c_I:P\in\Par_{ord}(I)\fle\varepsilon(P\cap I^+)\varepsilon(P\cap I^-)
\varepsilon'(P\cap I^+)\varepsilon'(P\cap I^-)$ (cf les points (4) et (6) de
l'exemple \ref{ex_LCI}), $b_I=1$ et $d_I=1$.
On trouve que la somme ci-dessus est égale à
\[\sum_{P\in\Par_{ord}(\lambda_1,\dots,\lambda_n)}(-1)^{|P|}\varepsilon
(P\cap I^+)\varepsilon(P\cap I^-)\varepsilon'(P\cap I^+)\varepsilon'(P\cap I^-)
\sum_{Q\in\Par_{ord}(\lambda'_1,\dots,\lambda'_m)}(-1)^{|Q|}.\]
D'après le corollaire \ref{prop_combinatoire_idiote_GU},
\[\sum_{Q\in\Par_{ord}(\lambda'_1,\dots,\lambda'_m)}(-1)^{|Q|}=(-1)^mc_1
(\lambda'_1)\dots c_1(\lambda'_m).\]
On applique une deuxième fois la proposition
\ref{prop_combinatoire_idiote_GSp1}, cette fois à $(I^+,I^-)\in\Dcal(\{1,\dots,
n\})$, avec, pour tout $I\subset\{1,\dots,n\}$, $a_I=b_I:(J,K)\in\Dcal(I)\fle
\varepsilon(J,K)$ et $c_I=d_I:P\in\Par_{ord}(I)\fle\varepsilon(P)\varepsilon'
(P)$. On trouve que $\sum\limits_{P\in\Par_{ord}(\lambda_1,\dots,\lambda_n)}
(-1)^{|P|}\varepsilon(P\cap I^+)\varepsilon(P\cap I^-)\varepsilon'(P\cap I^+)
\varepsilon'(P\cap I^-)$ est égal à
\[\sum_{P^+\in\Par_{ord}(\lambda^+)}(-1)^{|P^+|}\varepsilon(P^+)\varepsilon'
(P^+)\sum_{P^-\in\Par_{ord}(\lambda^-)}(-1)^{|P^-|}\varepsilon(P^-)\varepsilon'
(P^-),\]
où $\lambda^\pm=(\lambda_i)_{i\in I^\pm}$. L'égalité du corollaire résulte
donc de la proposition \ref{prop_combinatoire_idiote_GSp2}, appliquée
à $\lambda^+$ et $\lambda^-$.

\end{preuve}

\begin{preuvepn}{\ref{prop_combinatoire_idiote_GSp2}}
Comme dans la preuve de la proposition \ref{prop_combinatoire_idiote_GSp1},
on raisonne par récurrence sur
le couple $(n,|\Par_{ord}(\lambda)|)$. Si $n\leq 1$ ou si $\Par_{ord}(\lambda)=
\varnothing$, le résultat est évident. 

Supposons que $n=2$ et $\Par_{ord}(\lambda)\not=\varnothing$ (donc
$\lambda_1+\lambda_2>0$). Alors $\Par^0_{\leq 2}(\lambda)=\{\{1,2\}\}$, donc
$\sum\limits_{p\in\Par^0_{\leq 2}(\lambda)}\varepsilon(p)c(p,\lambda)=
c(\{\{1,2\}\},\lambda)$.
Si $\lambda_1>0$ et $\lambda_2>0$, alors $\Par_{ord}(\lambda)=
\{\{1,2\},(\{1\},\{2\}),(\{2\},\{1\})\}$, donc $\sum\limits_{P\in\Par_{ord}
(\lambda)}(-1)^{|P|}\varepsilon(P)\varepsilon'(P)=1$.
Si $\lambda_1>0$ et $\lambda_2\leq 0$,
alors $\Par_{ord}(\lambda)=\{\{1,2\},(\{1\},\{2\})\}$, donc $\sum\limits_{P\in
\Par_{ord}(\lambda)}(-1)^{|P|}\varepsilon(P)\varepsilon'(P)=2$.
Si $\lambda_1\leq 0$ et
$\lambda_2>0$, alors $\Par_{ord}(\lambda)=\{\{1,2\},(\{2\},\{1\})\}$, donc
$\sum\limits_{P\in\Par_{ord}(\lambda)}(-1)^{|P|}\varepsilon(P)\varepsilon'(P)
=0$. Dans tous les cas, l'égalité cherchée est vraie.

Supposons que $n\geq 3$ et $\Par_{ord}(\lambda)\not=\varnothing$.
Comme $\Par_{ord}(\lambda)\not=\varnothing$, on a $\lambda_1+\dots+\lambda_n
>0$ et $\delta:=\delta(\lambda)>0$.
Soient $J$ et $\lambda'$ comme dans le lemme \ref{LCI_de_base3}, dont on
utilise les notations. On note $\Par_{ord}'=\Par_{ord}(\lambda)\cap\Par'_{ord}
(n)$ et $\Par'=\Par^0_{\leq 2}(\lambda)\cap\Par'(n)$.
D'après le (ii) du lemme \ref{LCI_de_base3}, on a $\Par_{ord}(\lambda)=\Par_
{ord}(\lambda')\sqcup\Par'_{ord}$, $\Par^0_{\leq 2}(\lambda)=\Par^0_{\leq 2}
(\lambda')\sqcup\Par'$ et $|\Par_{ord}(\lambda')|<|\Par_{ord}(\lambda)|$.
Soit $p=\{I_\alpha,\alpha\in A\}\in\Par^0_{\leq 2}(\lambda')$. Soit $\alpha
\in A$. Si $J\cap I_\alpha=\varnothing$, on a évidemment $c_{I_\alpha}
(\lambda')=c_{I_\alpha}(\lambda)$. Si $I_\alpha\subset J$, on a $I_\alpha\not=
J$ car $p\in\Par(\lambda')$, donc, d'après le (i) du lemme \ref{LCI_de_base3},
$c_{I_\alpha}(\lambda')=c_{I_\alpha}(\lambda)$. Si $|J\cap I_\alpha|=1$ et
$J\cap I_\alpha\not= J$, alors, toujours d'après le (i) du lemme
\ref{LCI_de_base3}, $c_{I_\alpha}(\lambda')=c_{I_\alpha}(\lambda)$.
Finalement, on trouve que $c_p(\lambda')=c_p(\lambda)$ sauf si
$|J|=1$, et, si $|J|=1$, le seul élément $\alpha$ de $A$ tel que $c_{I_\alpha}
(\lambda')\not=c_{I_\alpha}(\lambda)$ est celui qui vérifie $J\subset
I_\alpha$.

Supposons que $|J|\geq 3$. Alors $\Par'=\varnothing$ (d'après le (i) du
lemme \ref{LCI_de_base}), donc
\[\sum_{p\in\Par^0_{\leq 2}(\lambda)}\varepsilon(p)c(p,\lambda)=
\sum_{p\in\Par^0_{\leq 2}(\lambda')}\varepsilon(p)c(p,\lambda)=\sum
_{p\in\Par^0_{\leq 2}(\lambda')}\varepsilon(p)c(p,\lambda').\]
D'autre part, si $\mu\in\R^{|J|}$ est défini comme dans le lemme
\ref{LCI_de_base3}, on a $N(\mu)=|J|$, donc $\sum\limits_{P\in\Par_{ord}
(\mu)}(-1)^{|P|}\varepsilon(P)\varepsilon'(P)=0$ d'après les lemmes
\ref{sous_LCI_GSp_pair} et \ref{sous_LCI_GSp_impair}. En utilisant les
points (iii) et (iv) du lemme \ref{LCI_de_base3}, on en déduit que
$\sum\limits_{P\in\Par'_{ord}}(-1)^{|P|}\varepsilon(P)\varepsilon'(P)=0$,
donc que $\sum\limits_{P\in\Par_{ord}(\lambda)}(-1)^{|P|}\varepsilon(P)
\varepsilon'(P)=\sum\limits_{P\in\Par_{ord}(\lambda')}(-1)^{|P|}\varepsilon(P)
\varepsilon'(P)$. L'égalité cherchée résulte
alors de l'hypothèse de récurrence, appliquée à $\lambda'$.

Supposons que $|J|=2$. On utilise toujours les notations du lemme
\ref{LCI_de_base3}. D'après les points (iii) et (iv) de ce lemme, 
$\sum\limits_{P\in\Par_{ord}(\lambda)}(-1)^{|P|}\varepsilon(P)\varepsilon'(P)$
est égale à
\[\sum_{P\in\Par_{ord}(\lambda')}(-1)^{|P|}\varepsilon(P)\varepsilon'(P)+
\varepsilon_0\sum_{P_1\in\Par_{ord}(\mu)}(-1)^{|P_1|}\varepsilon(P_1)
\varepsilon'(P_1)\sum_{P_2\in\Par_{ord}(\nu)}(-1)^{|P_2|}\varepsilon(P_2)
\varepsilon'(P_2).\]
D'autre part, comme $|J|$ est pair, l'application $\varphi$ du lemme
\ref{LCI_de_base3} induit une bijection $\Par'\iso\Par^0_{\leq 2}(\mu)\times
\Par^0_{\leq 2}(\nu)$. Donc, d'après les points (iii) et (iv) de ce lemme
et le fait que $c(p,\lambda)=c(p,\lambda')$ pour toute $p\in\Par^0_{\leq 2}
(\lambda')$, $\sum\limits_{p\in\Par^0_{\leq 2}(\lambda)}
\varepsilon(p)c(p,\lambda)$ est égale à
\[\sum_{p\in\Par^0_{\leq 2}(\lambda')}\varepsilon(p)c(p,\lambda')+
\varepsilon_0\sum_{p_1\in\Par^0_{\leq 2}(\mu)}\varepsilon(p_1)
c(p_1,\mu)\sum_{p_2\in\Par^0_{\leq 2}(\nu)}\varepsilon(p_2)c(p_2,
\nu).\]
L'égalité cherchée résulte alors de l'hypothèse de récurrence, appliquée à
$\lambda'$, $\mu$ et $\nu$.

Supposons que $|J|=1$. On a comme avant, d'après les points (iii) et
(iv) du lemme \ref{LCI_de_base} et l'hypothèse de récurrence appliquée à
$\lambda'$, $\mu$ et $\nu$,
\begin{flushleft}$\displaystyle{
\sum_{P\in\Par_{ord}(\lambda)}(-1)^{|P|}
\varepsilon(P)\varepsilon'(P)
}$\end{flushleft}
\begin{flushright}$\displaystyle{
=\sum_{P\in\Par_{ord}(\lambda')}(-1)^{|P|}
\varepsilon(P)\varepsilon'(P)+\varepsilon_0\sum_{P_1\in\Par_{ord}(\mu)}(-1)^
{|P_1|}\varepsilon(P_1)\varepsilon'(P_1)\sum_{P_2\in\Par_{ord}(\nu)}(-1)^
{|P_2|}\varepsilon(P_2)\varepsilon'(P_2)
}$\end{flushright}
\begin{flushright}$\displaystyle{
=(-1)^n\sum_{p\in\Par^0_{\leq 2}(\lambda')}\varepsilon(p)c(p,\lambda')+
(-1)^n\varepsilon_0\sum_{p_1\in\Par^0_{\leq 2}(\mu)}
\varepsilon(p_1)c(p_1,\mu)\sum_{p_2\in\Par^0_{\leq 2}(\nu)}
\varepsilon(p_2)c(p_2,\nu).
}$\end{flushright}
On traite d'abord le cas où $n$ est impair. Si $n$ est impair, alors
$\varphi$ induit une bijection $\Par'\iso\Par^0_{\leq 2}(\mu)\times
\Par^0_{\leq 2}(\nu)$, donc la somme ci-dessus est égale à
\[(-1)^n\sum_{p\in\Par^0_{\leq 2}(\lambda')}\varepsilon(p)c(p,\lambda')+
(-1)^n\sum_{p\in\Par'}\varepsilon(p)c(p,\lambda),\]
et il reste à montrer que
\[\sum_{p\in\Par^0_{\leq 2}(\lambda')}\varepsilon(p)c(p,\lambda')=
\sum_{p\in\Par^0_{\leq 2}(\lambda')}\varepsilon(p)c(p,\lambda).\]
Soit $p\in\Par^0_{\leq 2}(\lambda')$. Soient $\alpha_0,\alpha_1\in A$ tels que
$|I_{\alpha_0}|=1$ et $J\subset I_{\alpha_1}$. On écrit $I_{\alpha_0}=\{i_0\}$
et $I_{\alpha_1}=\{i_1,i_2\}$, avec $J=\{i_1\}$. On définit $p'\in
\Par(n)$ par $p'=\{I_\alpha,\alpha\in A-\{\alpha_0,\alpha_1\}\}\cup\{\{i_0,i_1
\},\{i_2\}\}$; comme $\lambda'_{i_1}=0$ et $s_{I_{\alpha_1}}(\lambda')>0$, on
a $\lambda'_{i_2}>0$, donc $p'\in\Par^0_{\leq 2}(\lambda')$.
L'application $p\fle p'$ est donc une involution (sans points fixes)
de $\Par^0_{\leq 2}
(\lambda')$, et on vérifie facilement (en examinant toutes les possibilités
pour les positions relatives de $i_0$, $i_1$ et $i_2$) que
\[\varepsilon(p)c(p,\lambda)+\varepsilon(p')c(p',\lambda)=
\varepsilon(p)c(p,\lambda')+\varepsilon(p')c(p',\lambda').\] 
L'égalité cherchée en résulte.

On traite enfin le cas où $n$ est pair (et $|J|=1$). Dans ce cas, on a
$\Par'=\varnothing$, donc il faut montrer que $\sum\limits_{p\in\Par^0_{\leq 2}
(\lambda)}\varepsilon(p)c(p,\lambda)$ est égale à
\[\sum_{p\in\Par^0_{\leq 2}(\lambda')}
\varepsilon(p)c(p,\lambda')+\varepsilon_0\sum_{p_1\in\Par^0_{\leq 2}(\mu)}
\varepsilon(p_1)c(p_1,\mu)\sum_{p_2\in\Par^0_{\leq 2}(\nu)}
\varepsilon(p_2)c(p_2,\nu).\]
Soit $p=\{I_\alpha\in A\}\in\Par^0_{\leq 2}(\lambda')$. Soit $\alpha_0\in A$
tel que $J\subset I_{\alpha_0}$. On écrit $I_{\alpha_0}=\{i_1,i_2\}$, avec
$J=\{i_1\}$. Comme $s_{I_{\alpha_0}}(\lambda')>0$ et $\lambda'_{i_1}=0$, on
a $\lambda'_{i_2}>0$. On pose $p_1=\{\{i_1\}\}$ et $p_2=\{I_\alpha,\alpha\in A-
\{\alpha_0\}\}\cup\{\{i_2\}\}$. Alors l'application $p\fle (p_1,p_2)$ induit
une bijection $\Par^0_{\leq 2}(\lambda')\iso\Par^0_{\leq 2}(\mu)\times\Par^0_
{\leq 2}(\nu)$. De plus, on voit facilement (en distinguant les cas
$i_1<i_2$ et $i_1>i_2$) que
\[\varepsilon(p)c(p,\lambda)=\varepsilon(p)c(p,\lambda')+\varepsilon_0
\varepsilon(p_1)c(p_1,\mu)\varepsilon(p_2)c(p_2,\nu).\]
L'égalité cherchée en résulte.

\end{preuvepn}

\begin{lemme}\label{LCI_GSp_pair}
Soient $n\in\Nat$ et $\lambda=(\lambda_1,\dots,\lambda_n)\in\R^n$.

Alors
\[\sum_{P\in\Par^0_{ord}(\lambda)}(-1)^{|P|}\varepsilon(P)\varepsilon'(P)=
(-1)^n\sum_{p\in\Par^0_{\leq 2}(\lambda)}\varepsilon(p).\]

\end{lemme}

\begin{preuve}
Comme dans la preuve de la proposition \ref{prop_combinatoire_idiote_GSp1}, on
raisonne par récurrence sur $(n,|\Par_{ord}(\lambda)|)$. Le résultat est
immédiat si $n\leq 2$ ou si $\Par_{ord}(\lambda)=\varnothing$.
On suppose
donc que $n\geq 3$ et que $\Par_{ord}(\lambda)\not=\varnothing$ (c'est-à-dire
que $\lambda_1+\dots+\lambda_n>0$). On note $\delta=\delta(\lambda)$. Soient
$J$ et $\lambda'$ comme dans le lemme \ref{LCI_de_base3}, dont on utilise les
notations. On note $\Par'_{ord}=\Par'_{ord}(n)\cap\Par^0_{ord}(\lambda)$ et
$\Par'=\Par'(n)\cap\Par^0_{\leq 2}(\lambda)$.
D'après le (ii) du lemme \ref{LCI_de_base3}, on a $|\Par_{ord}(\lambda')|<
|\Par_{ord}
(\lambda)|$, $\Par^0_{ord}(\lambda)=\Par^0_{ord}(\lambda')\sqcup\Par'_{ord}$
et $\Par^0_{\leq 2}(\lambda)=\Par^0_{\leq 2}(\lambda')\sqcup\Par'$. En
appliquant l'hypothèse de récurrence à $\lambda'$, on trouve
$\sum\limits_{P\in\Par^0_{ord}(\lambda')}(-1)^{|P|}\varepsilon(P)
\varepsilon'(P)=(-1)^n\sum\limits_{p\in\Par^0_{\leq 2}(\lambda')}
\varepsilon(p)$. Il suffit donc
de montrer que $\sum\limits_{P\in\Par'_{ord}}(-1)^{|P|}\varepsilon(P)
\varepsilon'(P)=(-1)^n\sum\limits_{p\in\Par'}\varepsilon(p)$.

On rappelle qu'on utilise les notations du lemme \ref{LCI_de_base3}; en
particulier, on note $m=N(\lambda)(=|J|)$. Si $m=n$ (c'est-à-dire
$J=\{1,\dots,n\}$), on a $\Par'_{ord}=\Par^0_{ord}(\lambda)$ (car $\Par'
_{ord}(n)=\Par_{ord}(n)$ par définition) et $\Par'=\Par^0_{\leq 2}(\lambda)=
\varnothing$ (car $\Par(\lambda)=\{\{1,\dots,n\}\}$ d'après le (i) du lemme
\ref{LCI_de_base}, et $n\geq 3$), donc l'égalité cherchée résulte du lemme
\ref{sous_LCI_GSp_pair} ci-dessous. Si $m$ et $n-m$ sont tous les deux impairs,
alors il est clair que $\Par'_{ord}=\Par'=\varnothing$. On peut donc supposer
que $m<n$, et que $m$ et $n-m$ ne sont pas tous les deux impairs. Alors le
(iv) du lemme \ref{LCI_de_base3} implique que $\varphi_{ord}(\Par'_{ord})=
\Par^0_{ord}(\mu)\times\Par^0_{ord}(\nu)$ et $\varphi(\Par')=\Par^0_{\leq 2}
(\mu)\times\Par^0_{\leq 2}(\nu)$. L'égalité cherchée résulte donc de
l'hypothèse de récurrence appliquée à $\mu$ et à $\nu$ et du point (iii) du
lemme \ref{LCI_de_base3}.

\end{preuve}

\begin{lemme}\label{sous_LCI_GSp_pair} On suppose que $n\geq 3$.
Soit $\lambda=(\lambda_1,\dots,\lambda_n)\in\R^n$ tel que $\lambda_1+\dots+
\lambda_n>0$ et $N(\lambda)=n$. Alors
\[\sum_{P\in\Par^0_{ord}(\lambda)}(-1)^{|P|}\varepsilon(P)\varepsilon'(P)=0.\]

\end{lemme}

\begin{preuve} On note $m$ la partie entière de $n/2$. D'après la remarque
sous la définition de $\varepsilon'$, on a $\varepsilon'(P)=(-1)^m$ pour
toute $P\in\Par^0_{ord}(n)$. Il suffit donc de montrer que
\[\sum_{P\in\Par^0_{ord}(\lambda)}(-1)^{|P|}\varepsilon(P)=0.\]
Soit
$\Par_{ord}'$ le sous-ensemble de $\Par^0_{ord}(\lambda)$ formé
des partitions ordonnées $P$ telles que, pour tout $i\in\{1,\dots,m\}$, $2i$ et
$2i-1$ soient dans le même ensemble de $P$. On définit par récurrence
descendante sur $i\in\{1,\dots,m\}$ des sous-ensembles $\Par''_{ord,i}$
de $\Par^0_{ord}(n)$ de la manière suivante :
$\Par''_{ord,i}$ est l'ensemble des $P\in\Par^0_{ord}(n)-\bigcup\limits
_{j=i+1}^m\Par''_{ord,j}$ telles que $2i$ et $2i-1$ soient dans des ensembles
différents de $P$. Pour tout $i\in\{1,\dots,m\}$, on note $\Par''_{ord,i}
(\lambda)=\Par''_{ord,i}\cap\Par_{ord}(\lambda)$. Alors $\Par^0_{ord}(\lambda)-
\Par'_{ord}$ est union disjointe des $\Par''_{ord,i}(\lambda)$,
$1\leq i\leq m$.

D'après les lemmes \ref{LCI_de_base} et \ref{LCI_de_base2},
pour toute $p\in\Par(n)$, le
cardinal de $\oubli^{-1}(p)\cap\Par_{ord}(\lambda)$ est $(|p|-1)!$,
c'est-à-dire $|\oubli^{-1}(p)|/|p|$.

Soit $i\in\{1,\dots,m\}$. On note $\Par''_i=\oubli(\Par''_{ord,i})$. D'après
la remarque ci-dessus, pour tout $k\in\Nat^*$ :
\[\sum_{P\in\Par''_{ord,i}(\lambda),|P|=k}(-1)^{|P|}\varepsilon(P)
=(-1)^k(k-1)!\sum_{p\in\Par''_i,|p|=k}\varepsilon(p).\]
On définit une involution de $\Par''_i$ en envoyant une partition
$p$ sur la partition $p'$ obtenue en échangeant $2i$ et $2i-1$. On a
alors $\varepsilon(p')=-\varepsilon(p)$ et
$|p'|=|p|$; donc, pour tout $k\in\Nat^*$,
\[\sum_{p\in\Par''_i,|p|=k}\varepsilon(p)=0.\]

On en déduit que, pour tout $i\in\{1,\dots,m\}$, $\sum\limits_{P\in\Par''_{ord,
i}(\lambda)}(-1)^{|P|}\varepsilon(P)=0$, donc que
\[\sum_{P\in\Par^0_{ord}(\lambda)}(-1)^{|P|}\varepsilon(P)
=\sum_{P\in\Par'_{ord}}(-1)^{|P|}\varepsilon(P).\]
D'autre part, on remarque que $\varepsilon(P)=1$ pour toute $P\in\Par'_{ord}$.
Donc
\[\sum_{P\in\Par^0_{ord}(\lambda)}(-1)^{|P|}\varepsilon(P)
=\sum_{P\in\Par'_{ord}}(-1)^{|P|}.\]
Notons $\mu=(\lambda_1+\lambda_2,\dots,\lambda_{n-1}+\lambda_n)$ si $n$ est
pair, et $\mu=(\lambda_1+\lambda_2,\dots,\lambda_{n-2}+\lambda_{n-1},
\lambda_n)$
si $n$ est impair. Alors $\Par'_{ord}$ est de manière évidente en bijection
avec $\Par_{ord}(\mu)$ (et cette bijection est compatible avec les applications
$P\fle |P|$), donc le lemme résulte du corollaire
\ref{prop_combinatoire_idiote_GU}
(on utilise ici l'hypothèse $n\geq 3$ pour voir que
les coordonnées de $\mu$ ne peuvent pas être toutes strictement positives).

\end{preuve}

\begin{remarque}\label{rq:LCI_GSp_pair} Le lemme \ref{LCI_GSp_pair} reste
valable si l'on remplace tous les ``$>0$'' par des ``$\geq 0$'' dans les
définitions de $\Par(\lambda)$ et $\Par_{ord}(\lambda)$ (c'est-à-dire si l'on
remplace $\Par(\lambda)$ par l'ensemble des $p=\{I_\alpha,\alpha\in A\}\in
\Par(n)$ telles que $s_{I_\alpha}(\lambda)\geq 0$ pour tout $\alpha\in A$ et
$\Par_{ord}(\lambda)$ par l'ensemble des $P=(I_1,\dots,I_r)\in\Par_{ord}(n)$
telles que $s_{I_1}(\lambda)+\dots+s_{I_i}(\lambda)\geq 0$ pour tout
$i\in\{1,\dots,r\}$).

En effet, si l'on définit $\lambda'\in\R^n$ par $\lambda'=(\lambda_1+\eta,
\dots,\lambda_n+\eta)$ avec $\eta>0$ assez petit, alors, pour tout
$I\subset\{1,\dots,n\}$, $s_I(\lambda)\geq 0$ si et seulement si $s_I(\lambda')
>0$. Il suffit alors d'appliquer le lemme \ref{LCI_GSp_pair} à $\lambda'$.

\end{remarque}

On rappelle que, pour tout $k\in\Nat$, on a noté
$\Par^k(n)$ l'ensemble des $p=\{I_\alpha,\alpha\in A\}$
dans $\Par(n)$ telles que le cardinal de l'ensemble
$\{\alpha\in A,|I_\alpha|\ est\ impair\}$ soit $2k$ ou $2k+1$,
et $\Par^k_{ord}(n)=\oubli^{-1}(\Par^k(n))$.
Pour tous $\lambda\in\R^n$ et $k\in\Nat$, on note $\Par^k_{ord}(\lambda)=
\Par^k_{ord}(n)\cap \Par_{ord}(\lambda)$.

\begin{lemme}\label{sous_LCI_GSp_impair}
On suppose que $n\geq 3$. Soit
$\lambda=(\lambda_1,\dots,\lambda_n)\in\R^n$ tel que $\lambda_1+\dots+\lambda_n
>0$ et $N(\lambda)=n$.
Alors, pour tout $k\geq 1$,
\[\sum_{P\in\Par^k_{ord}(\lambda)}(-1)^{|P|}\varepsilon(P)\varepsilon'(P)=0.\]

\end{lemme}

L'égalité du lemme est bien sûr vraie aussi si $k=0$, puisque c'est le
résultat du lemme \ref{sous_LCI_GSp_pair} dans ce cas. (On a séparé les deux
lemmes, car le lemme \ref{LCI_GSp_pair}, qui utilise le lemme
\ref{sous_LCI_GSp_pair}, sert dans la démonstration du lemme
\ref{sous_LCI_GSp_impair}.)

\begin{preuve}
On note $m$ la partie entière de $n/2$. D'après la remarque sous la
définition de $\varepsilon'$, on a $\varepsilon'(P)=(-1)^{m-k}$ pour toute
$P\in\Par^k_{ord}(n)$. Il suffit donc de prouver que
\[\sum_{P\in\Par^k_{ord}(\lambda)}(-1)^{|P|}\varepsilon(P)=0.\]

On suppose d'abord que $n$ est impair. On va montrer un résultat plus précis.
Soit $p=\{I_\alpha,\alpha\in A\}\in\Par(n)$ telle que deux au moins des
$I_\alpha$ aient un cardinal impair. Montrons que 
\[\sum_{P\in\Par_{ord}(\lambda)\cap\oubli^{-1}(p)}\varepsilon(P)=0.\]

Soit $P=(I_1,\dots,I_k)\in\oubli^{-1}(p)$. D'après
le lemme \ref{LCI_de_base} (appliqué à $\lambda_P$), il existe un unique
$l\in\{1,\dots,k\}$ tel que $P':=(I_l,\dots,I_k,I_1,\dots,I_{l-1})$ soit
dans $\Par_{ord}(\lambda)$. De plus, on a $\sigma_{P}=\tau^{|I_1|+\dots+
|I_{l-1}|}\sigma_{P'}$, donc $\varepsilon(P')=\varepsilon(P)$ (comme $n$ est
impair, $sgn(\tau)=1$).
On en déduit que
\[\sum_{P\in\Par_{ord}(\lambda)\cap\oubli^{-1}(p)}\varepsilon(P)
=\frac{1}{|p|}\sum_{P\in\oubli^{-1}(p)}\varepsilon(P).\]
On écrit $p=\{I_\alpha,\alpha\in A\}$. Soient $\alpha_1,\alpha_2\in A$
distincts tels que $|I_{\alpha_1}|$ et $|I_{\alpha_2}|$ soient impairs. On
considère l'involution $\iota$ de $\oubli^{-1}(p)$ qui échange les places des
ensembles $I_{\alpha_1}$ et $I_{\alpha_2}$. Alors $\varepsilon(\iota(P))=
-\varepsilon(P)$
pour tout $P\in\oubli^{-1}(p)$, donc $\sum\limits_{P\in\oubli^{-1}(p)}
\varepsilon(P)=0$.

On suppose maintenant que $n$ est pair. Soit $k\in\Nat^*$. Soit $P=(I_1,\dots,
I_r)\in\Par^k_{ord}(\lambda)$. On note $\mu=\lambda_P\in\R^r$ (donc
$\mu_i=s_{I_i}(\lambda)$).
Comme $n$ est pair, il y a exactement $2k$
indices $i\in\{1,\dots,r\}$ tels que $|I_i|$ soit impair. On les note
$i_1,\dots,i_{2k}$, avec $i_1<\dots<i_{2k}$. 
On utilise les notions introduites au-dessus du lemme \ref{LCI_blocs}, et
on note $P=Q_1\dots Q_{2k}$ la décomposition de $P$ en blocs de centres
$I_{i_1},\dots,I_{i_{2k}}$ (cf le (iii) du lemme \ref{LCI_blocs}).
Dans la suite de la preuve, on utilisera toujours cette décomposition en blocs,
et on l'appellera la décomposition en blocs de $P$. On remarque que les
blocs sont tous de taille impaire.

On note $\Par'$ l'ensemble des $P\in\Par^k_{ord}(\lambda)$ telles que, pour
tout $m\in\{1,\dots,k\}$, $Q_{2m-1}$ soit positif et $Q_{2m}$ négatif.
On définit par récurrence sur $l\in\{1,\dots,2k-1\}$ des sous-ensembles
$\Par''_l$ de $\Par^k_{ord}(\lambda)$ de la manière
suivante : 
\begin{itemize}
\item[$\bullet$] pour tout $m\in\{1,\dots,k\}$, $\Par''_{2m-1}$ est l'ensemble
des $P\in\Par^k_{ord}(\lambda)-\bigcup\limits_{1\leq l\leq 2m-3}\Par''_l$ tels
que $Q_{2m-1}$ et $Q_{2m}$ soient tous les deux positifs;
\item[$\bullet$] pour tout $m\in\{1,\dots,k-1\}$, $\Par''_{2m}$ est l'ensemble
des $P\in\Par^k_{ord}(\lambda)-\bigcup\limits_{1\leq l\leq 2m-2}\Par''_l$ tels
que $Q_{2m+1}$ et $Q_{2m}$ soient tous les deux négatifs.

\end{itemize}
Alors on a $\Par^k_{ord}(\lambda)=\Par'\sqcup\coprod\limits_{1\leq l\leq 2k-1}
\Par''_l$. 

Soit $l\in\{1,\dots,2k-1\}$. Soit $P\in\Par''_l$. Soient
$Q_1,\dots,Q_{2k}$ les blocs de $P$. On note $\iota(P)$ la partition
ordonnée qu'on obtient en échangeant les blocs $Q_l$ et $Q_{l+1}$.
Alors, d'après le (ii) du lemme \ref{LCI_blocs}, $\iota(P)$ est encore dans
$\Par_{ord}(\lambda)$, $|\iota(P)|=|P|$, $\varepsilon(\iota(P))=-
\varepsilon(P)$.
De plus, il est clair que $\iota(P)\in\Par''_l$.
Donc l'application $P\fle\iota(P)$
est une involution de $\Par''_l$. On en déduit que $\sum\limits_{P\in
\Par''_l}(-1)^{|P|}\varepsilon(P)=0$. Finalement,
\[\sum_{P\in\Par^k_{ord}(n)}(-1)^{|P|}\varepsilon(P)=
\sum_{P\in\Par'}(-1)^{|P|}\varepsilon(P).\]

Il reste à montrer que $\sum\limits_{P\in\Par'}(-1)^{|P|}\varepsilon(P)=0$.
Supposons d'abord que $k\geq 2$. Soit $P\in\Par'$, et soit $P=Q_1\dots Q_{2k}$
la décomposition en blocs de $P$.
Soit $\Sgoth_{2k}'$ l'ensemble des permutations
$\sigma\in\Sgoth_{2k}$ qui envoient les entiers impairs entre $1$ et
$2k$ sur des entiers impairs (donc $\sigma$ envoie les entiers pairs sur des
entiers pairs). On note $\Qpar$ l'ensemble des partitions ordonnées de
$\{1,\dots,n\}$ de la forme $Q_{\sigma(1)}\dots Q_{\sigma(2k)}$, avec $\sigma
\in\Sgoth'_{2k}$.
Il suffit de montrer que $\sum\limits_{P'\in\Par'\cap\Qpar}
(-1)^{|P'|}\varepsilon(P')=0$
(quelle que soit la partition ordonnée $P$ de
départ). Soit $P'=Q_{\sigma(1)}\dots Q_{\sigma(2k)}$, avec $\sigma\in\Sgoth'_
{2k}$. D'après le (iv) du lemme \ref{LCI_blocs}, il existe un unique
$a\in\{1,\dots,2k\}$ tel que $P'':=Q_{\sigma(a)}\dots Q_{\sigma(2k)}Q_{\sigma
(1)}\dots Q_{\sigma(a-1)}\in\Par_{ord}(\lambda)$. Comme une partition
de $\Par_{ord}(\lambda)$ ne peut pas commencer par un bloc négatif, $a$ est
pair; donc $P''\in\Qpar$ et $\varepsilon(P'')=\varepsilon(P)$.
On en déduit que
\[\sum_{P'\in\Qpar\cap\Par'}(-1)^{|P'|}\varepsilon(P')
=\frac{1}{k}\sum_{P'\in\Qpar}(-1)^{|P'|}\varepsilon(P').\]
Notons, pour toute $P'=Q_{\sigma(1)}\dots Q_{\sigma(2k)}\in\Qpar$,
$\iota(P')$ la partition ordonnée obtenue à partir de $P'$ en échangeant
$Q_{\sigma(1)}$ et $Q_{\sigma(3)}$. Alors $\iota$ est une involution de
$\Qpar$, et $|\iota(P')|=|P'|$ et $\varepsilon(\iota(P'))=-\varepsilon(P')$
pour toute $P'\in\Qpar$. Donc $\sum\limits_{P'\in\Qpar}(-1)^{|P'|}\varepsilon
(P')=0$.

Traitons enfin le cas où $k=1$. 
Soit $J\subset\{1,
\dots,n\}$ de cardinal impair. On note $\Par_J$ l'ensemble des partitions
ordonnées $P\in\Par^1_{ord}(\lambda)$ telles que $J$ soit le
deuxième ensemble de cardinal impair de $P$.
On a $\Par_J=\coprod\limits_{(J_1,J_2)}\Par_
{J,J_1,J_2}$, où :
\begin{itemize}
\item[$\bullet$] $(J_1,J_2)$ parcourt l'ensemble des partitions en deux
ensembles de
$\{1,\dots,n\}-J$ telles que $|J_1|$ est impair (donc $|J_2|$ est pair),
$s_{J_1}(\lambda)>0$ et $s_{J_2}(\lambda)\leq 0$ (noter qu'alors on a
forc{\'e}ment $s_{J_1\cup J}(\lambda)>0$);
\item[$\bullet$] $\Par_{J,J_1,J_2}$ est l'ensemble des $P=(I_1,\dots,I_k)
\in\Par_J$ telles que, si le deuxième ensemble de cardinal impair de $P$
est $I_r$ (i.e., $I_r=J$), alors $J_1=I_1\cup\dots\cup I_{r-1}$ et
$J_2=I_{r+1}\cup\dots\cup I_k$.

\end{itemize}

Soient $\mu=(\lambda_i)_{i\in J_1}$ et $\nu\in\R^{|J_2|}$ l'{\'e}l{\'e}ment
obtenu {\`a} partir de $(-\lambda_i)_{i\in J_2}$ en inversant l'ordre sur les
indices. Alors se donner un {\'e}l{\'e}ment de $P$ de $\Par_{J,J_1,J_2}$
revient {\`a} se donner un {\'e}l{\'e}ment $P_1$ de $\Par^0_{ord}(\mu)$ et
un {\'e}l{\'e}ment $P_2$ de l'analogue de $\Par^0_{ord}(\nu)$ qu'on
obtient en rempla\c{c}ant toutes les in{\'e}galit{\'e}s strictes par des
in{\'e}galit{\'e}s larges dans la d{\'e}finition.

En appliquant le lemme \ref{LCI_GSp_pair} (et la remarque
\ref{rq:LCI_GSp_pair}) à $\mu$ et $\nu$,
on trouve, pour tout $(J_1,J_2)$ comme ci-dessus :
\[\sum_{P\in\Par_{J,J_1,J_2}}(-1)^{|P|}\varepsilon(P)\varepsilon'(P)=
(-1)^{n+1-|J|+\frac{1}{2}|J|(|J|-1)}
\sum_{p\in\Par'_{J,J_1,J_2}}\varepsilon(p),\]
où :
\begin{itemize}
\item[$\bullet$] $\Par'_{J,J_1,J_2}$ est l'ensemble des partitions $p=\{I_
\alpha,\alpha\in A\}\in\Par(n)$ telles qu'il existe une partition
$A=A_1\sqcup A_2\sqcup\{\alpha_0\}$ de $A$ avec :
\begin{itemize}
\item[-] $J=I_{\alpha_0}$, $J_1=\bigcup\limits_{\alpha\in A_1}I_{\alpha}$,
$J_2=\bigcup\limits_{\alpha\in A_2}I_\alpha$;
\item[-] pour tout $\alpha\in A_1$ (resp. $A_2$), $|I_\alpha|\leq 2$ (resp.
$|I_\alpha|=2$) et $s_{I_\alpha}(\lambda)>0$ (resp. $\leq 0$);
\item[-] il existe un unique $\alpha\in A_1$ tel que $|I_\alpha|=1$;

\end{itemize}
\item[$\bullet$] pour toute $p\in\Par'_{J,J_1,J_2}$, on note $\varepsilon(p)=
\varepsilon(P)$, où $P\in\oubli^{-1}(p)$ est obtenue en choisissant un
ordre sur les ensembles de $p$ qui place $J$ après l'autre ensemble de cardinal
impair ($\varepsilon(P)$ ne dépend que de l'ordre des ensembles de
cardinal impair de $P$).

\end{itemize}
Soit $\Par'_J$ l'ensemble des $p=\{I_\alpha,\alpha\in A\}\in\Par(n)$ telles
que :
\begin{itemize}
\item[$\bullet$] il existe $\alpha_0\in A$ tel que $I_{\alpha_0}=J$;
\item[$\bullet$] il existe $\alpha_1\in A-\{\alpha_0\}$ tel que
$|I_{\alpha_1}|=1$ et $s_{I_{\alpha_1}}(\lambda)>0$;
\item[$\bullet$] pour tout $\alpha\in A-\{\alpha_0,\alpha_1\}$, $|I_\alpha|=2$.

\end{itemize}
Alors $\Par'_J=\coprod\limits_{(J_1,J_2)}\Par'_{J,J_1,J_2}$, où $(J_1,J_2)$
parcourt le même ensemble d'indices que ci-dessus. De plus, $n+1-|J|$ est pair
et $|J|$ est impair, donc $n+1-|J|+\frac{1}{2}|J|(|J|-1)=\frac{1}{2}(|J|-1)
\mod 2$.
Donc $\sum\limits_{P\in
\Par_J}(-1)^{|P|}\varepsilon(P)\varepsilon'(P)=(-1)^{\frac{1}{2}(|J|-1)}
\sum\limits_{p\in\Par'_J}\varepsilon(p)$.

Soit $k\in\{1,\dots,n\}$. On note $\Par''_k$ l'ensemble
des partitions $p=\{I_\alpha,\alpha\in A\}\in\Par(n)$ telles que :
\begin{itemize}
\item[$\bullet$] il existe $\alpha_0\in A$ tel que $I_{\alpha_0}=\{k\}$;
\item[$\bullet$] il existe $\alpha_1\in A-\{\alpha_0\}$ tel que
$|I_{\alpha_1}|$ soit impair;
\item[$\bullet$] pour tout $\alpha\in A-\{\alpha_0,\alpha_1\}$, $|I_\alpha|=2$.

\end{itemize}
Pour toute $p\in\Par''_k$, on pose $\varepsilon(p)=\varepsilon(P)$, où
$P\in\oubli^{-1}(p)$ est obtenue en choisissant un ordre sur les ensembles de
$p$ qui place $\{k\}$ avant l'autre ensemble de cardinal impair, et on note
$\varepsilon''(p)=(-1)^{\frac{1}{2}(m-1)}$, où $m$ est le cardinal de
l'ensemble de cardinal impair de $p$ qui n'est pas $\{k\}$.

Alors 
\[\sum_J(-1)^{\frac{1}{2}(|J|-1)}\sum_{p\in\Par'_J}\varepsilon(p)
=\sum_k\sum_{p\in\Par''_k}\varepsilon(p)\varepsilon''(p),\] 
où, dans la première somme, $J$ parcourt l'ensemble des sous-ensembles de
cardinal impair de $\{1,\dots,n\}$ et, dans la deuxième somme,
$k$ parcourt l'ensemble des éléments de $\{1,\dots,n\}$ tels que
$\lambda_k>0$. Donc :
\[\sum_{P\in\Par^1_{ord}(\lambda)}(-1)^{|P|}\varepsilon'(P)\varepsilon(P)
=\sum_J\sum_{P\in\Par_J}
(-1)^{|P|}\varepsilon'(P)
\varepsilon(P)=\sum_J(-1)^{\frac{1}{2}(|J|-1)}\sum_{p\in\Par'_J}
\varepsilon(p)=\sum_k\sum_{p\in\Par''_k}\varepsilon(p)\varepsilon''(p),\]
où $J$ et $k$ parcourent les mêmes ensembles que plus haut.

Soit $k\in\{1,\dots,n\}$. Montrons que $\sum\limits_{p\in\Par''_k}
\varepsilon(p)\varepsilon''(p)=0$ (ce qui finit la preuve du lemme).
On voit facilement que
ceci résulte du lemme \ref{sous_sous_LCI_GSp_impair} ci-dessous (appliqué
à $n-1$, qui est bien impair et $\geq 3$).

\end{preuve}

\begin{lemme}\label{sous_sous_LCI_GSp_impair} Soit $n\geq 3$ impair.
On note $\Par$ l'ensemble des partitions
$p\in\Par(n)$ telles que l'un des ensembles de $p$ soit de cardinal impair et
tous les autres ensembles de $p$ soient de cardinal $2$. Pour toute $p\in\Par$,
on note $\varepsilon''(p)=(-1)^{\frac{1}{2}(|I|-1)}$, où $I$ est
l'ensemble de cardinal impair de $p$. Alors
\[\sum\limits_{p\in\Par}\varepsilon(p)\varepsilon''(p)=0\]
($\varepsilon(p)$ est défini pour
$p\in\Par$, car $p$ a un seul ensemble de cardinal impair).

\end{lemme}

\begin{preuve} Le raisonnement ressemble beaucoup à celui de la preuve
du lemme \ref{sous_LCI_GSp_pair}.
On note $m=(n-1)/2$.
Soit $\Par'$ le sous-ensemble de $\Par$ formé des partitions $p$ telles que,
pour tout $i\in\{1,\dots,m\}$, $2i$ et $2i-1$ soient dans le même ensemble de
$p$. On définit par récurrence descendante sur $i\in\{1,\dots,m\}$ des
sous-ensembles $\Par''_i$ de $\Par$ de la manière suivante :
$\Par''_i$ est l'ensemble des $p\in\Par-\bigcup\limits_{j=i+1}^m\Par''_j$
telles que $2i$ et $2i-1$ soient dans des ensembles
différents de $p$.

Soit $i\in\{1,\dots,m\}$.
On définit une involution de $\Par''_i$ en envoyant une partition
$p$ sur la partition $p'$ obtenue en échangeant $2i$ et $2i-1$. On a
alors $\varepsilon(p')=-\varepsilon(p)$ et $\varepsilon''(p')=\varepsilon''
(p)$; donc
\[\sum_{p\in\Par''_i}\varepsilon(p)\varepsilon''(p)=0.\]
On en déduit que $\sum\limits_{p\in\Par}\varepsilon(p)\varepsilon''(p)
=\sum\limits_{p\in\Par'}\varepsilon(p)\varepsilon''(p)$.
D'autre part, on remarque que $\varepsilon(p)=1$ pour toute $p\in\Par'$.
Donc
\[\sum_{p\in\Par}\varepsilon(p)\varepsilon''(p)=\sum_{p\in\Par'}\varepsilon''
(p).\]
On remarque qu'une partition de $\Par'$ est entièrement déterminée par la
donnée de son ensemble de cardinal impair. Soit $p\in\Par'$, et soit $I$ son
ensemble de cardinal impair; alors :
\begin{itemize}
\item[-] $n\in I$;
\item[-] $I$ est entièrement déterminé par son intersection avec l'ensemble
$J:=\{2i-1,1\leq i\leq m\}$, et tous les sous-ensembles de $J$ apparaissent
de cette manière;
\item[-] $|I\cap J|=\frac{1}{2}(|I|-1)$.

\end{itemize}
Finalement,
\[\sum_{p\in\Par'}\varepsilon''(p)=\sum_{k=0}^m(-1)^kC_m^k=0.\]

\end{preuve}

\begin{lemme}\label{LCI_de_base} Soit
$\lambda=(\lambda_1,\dots,\lambda_n)\in\R^n$ tel que $\lambda_1+\dots+\lambda_n
>0$. Alors :
\begin{itemize}
\item[(i)] Les conditions suivantes sont équivalentes :
\begin{itemize}
\item[(a)] $N(\lambda)=n$;
\item[(b)] il n'existe pas de partition $(I_1,I_2)$ de $\{1,\dots,n\}$ telle
que $s_{I_1}(\lambda)>0$ et $s_{I_2}(\lambda)>0$.
\end{itemize}

\item[(ii)] Il existe $k\in\Z$ tel que $\tau^k(\lambda)>0$. De plus, s'il
n'existe pas de partition $(I_1,I_2)$ de $\{1,\dots,n\}$ telle que $s_{I_1}
(\lambda)>0$ et $s_{I_2}(\lambda)>0$, alors un tel entier $k$ est uniquement
déterminé modulo $n$.

\end{itemize}
\end{lemme}

\begin{preuve} Montrons (i). Supposons que $N(\lambda)=n$ et qu'il existe une
partition $(I_1,I_2)$ de $\{1,\dots,n\}$ telle que $s_{I_1}(\lambda)>0$ et
$s_{I_2}(\lambda)>0$.
D'après la définition de $\delta(\lambda)$, on a
$s_{I_1}(\lambda)>\delta(\lambda)|I_1|$ et $s_{I_2}(\lambda)>
\delta(\lambda)|I_2|$, donc $\lambda_1+\dots+\lambda_n=s_{I_1}(\lambda)+
s_{I_2}(\lambda)> n\delta(\lambda)$, contradiction. Donc (a) implique
(b).

Réciproquement, supposons que $N(\lambda)<n$, et soit $I_1\subset\{1,\dots,n\}$
tel que $s_{I_1}(\lambda)/|I_1|=\delta(\lambda)$ et $|I_1|=N(\lambda)$.
On note $I_2=\{1,\dots,n\}
-I_1$. Alors $s_{I_2}(\lambda)=\lambda_1+\dots+\lambda_n-s_{I_1}(\lambda)
\geq (n-|I_1|)\delta(\lambda)>0$. Donc (b) implique (a).

Montrons (ii).
On note $s=\min\{\lambda_1+\dots+\lambda_l,1\leq
l\leq n\}$. Soit $k$ le plus grand élément de $\{1,\dots,n\}$ tel que 
$\lambda_1+\dots+\lambda_k=s$. Si $l\in\{k+1,\dots,n\}$, on a $\lambda_1+\dots
+\lambda_l >\lambda_1+\dots+\lambda_k$, donc $\lambda_{k+1}+\dots+\lambda_l>0$.
Si $l\in\{1,\dots,k\}$, on a 
\[\begin{array}{rcl}\lambda_{k+1}+\dots+\lambda_n+\lambda_1+\dots+\lambda_l &
= & (\lambda_1+\dots+\lambda_n)-(\lambda_1+\dots+\lambda_k)+(\lambda_1+\dots+
\lambda_l) \\
 & > & -(\lambda_1+\dots+\lambda_k)+(\lambda_1+\dots+\lambda_l) \\
& \geq & 0.\end{array}\]
Ceci prouve que $\tau^k(\lambda)=(\lambda_{k+1},\dots,
\lambda_n,\lambda_1,\dots,\lambda_k)>0$.

Montrons la dernière assertion de (ii). On suppose qu'il existe $k,l\in\{1,
\dots,n\}$ tels que $k<l$, $\tau^k(\lambda)>0$ et $\tau^l(\lambda)>0$. On note
$I_1=\{k+1,\dots,l\}$ et $I_2=\{1,\dots,n\}-I_1$. Alors $s_{I_1}(\lambda)=
\lambda_{k+1}+\dots+\lambda_l>0$ car $\tau^k(\lambda)>0$, et
$s_{I_2}(\lambda)=\lambda_{l+1}+\dots+\lambda_n+\lambda_1+\dots+\lambda_k>0$
car $\tau^l(\lambda)>0$. 

\end{preuve}

\begin{remarque}\label{rq:LCI_de_base} L'entier $k$ défini dans la
preuve de la première partie de (ii) du lemme ci-dessus est l'unique élément de
$\{1,\dots,n\}$ vérifiant les deux propriétés suivantes : 
\begin{itemize}
\item[(a)] pour tout $l\in\{k+1,\dots,n\}$, $\lambda_{k+1}+\dots+\lambda_l>0$;
\item[(b)] pour tout $l\in\{2,\dots,k\}$, $\lambda_l+\dots+\lambda_k\leq 0$.

\end{itemize}
(Un tel $k$ existe même si $\lambda_1+\dots+\lambda_n\leq 0$.)

\end{remarque}

\begin{lemme}\label{LCI_de_base2} Soit $\lambda=(\lambda_1,
\dots,\lambda_n)\in\R^n$. On suppose que $\lambda_1+\dots+\lambda_n>0$ et
qu'il n'existe pas de partition $\{I_1,I_2\}$ de $\{1,\dots,n\}$ telle que
$s_{I_1}(\lambda)>0$ et $s_{I_2}(\lambda)>0$. Alors $|\Sgoth(\lambda)|=(n-1)!$.

\end{lemme}

\begin{preuve} D'après le lemme \ref{LCI_de_base},
$\Sgoth_n$ est union disjointe des
sous-ensembles $\tau^k\Sgoth(\lambda)$, $0\leq k\leq n-1$. La conclusion du
lemme en résulte.

\end{preuve}

\begin{lemme}\label{LCI_de_base3} Soit $\lambda=(\lambda_1,\dots,\lambda_n)$
tel que $\lambda_1+\dots+\lambda_n>0$. On note $\delta=\delta(\lambda)$.
Soit $J\subset\{1,\dots,n\}$ tel que $s_J(\lambda)/|J|=\delta$ et $|J|=
N(\lambda)$.
On définit $\lambda'=(\lambda'_1,\dots,
\lambda'_n)\in\R^n$ par :
\[\lambda'_i=\left\{\begin{array}{ll}\lambda_i & \mbox{ si }i\not\in J \\
\lambda_i-\delta & \mbox{ si }i\in J\end{array}\right..\]
On note $\Par'_{ord}(n)$ l'ensemble des $P=(I_1,\dots,I_k)\in\Par_{ord}(n)$
telles qu'il existe $r\in\{1,\dots,k\}$ avec $J=I_1\cup\dots\cup I_r$, et
$\Par'(n)=\oubli(\Par'_{ord}(n))$. Alors :

\begin{itemize}
\item[(i)] Pour tout $K\subset\{1,\dots,n\}$ tel que $K\not=J$, on a
$s_K(\lambda')>0$ si et seulement si $s_K(\lambda)>0$.
\item[(ii)] On a
\[\Par_{ord}(\lambda)=\Par_{ord}(\lambda')\sqcup(\Par'_{ord}(n)\cap\Par_{ord}
(\lambda))\]
\[\Par(\lambda)=\Par(\lambda')\sqcup(\Par'(n)\cap\Par(\lambda)).\]
En particulier, $|\Par_{ord}(\lambda')|<|\Par_{ord}(\lambda)|$ et
$|\Par(\lambda')|<|\Par(\lambda)|$.
\end{itemize}

Notons $m=N(\lambda)(=|J|)$. 
On écrit $J=\{i_1,\dots,i_m\}$ avec $i_1<\dots<i_m$ et $K:=\{1,\dots,n\}-J=
\{j_1,\dots,j_{n-m}\}$ avec $j_1<\dots<j_{n-m}$. On a des bijections
$u_J:J\iso\{1,\dots,m\},i_r\fle r$ et $u_K:K\iso\{1,\dots,n-m\},j_r\fle r$.
Soit $\sigma\in\Sgoth_n$ défini par $\sigma(i_r)=r$ pour tout $r\in
\{1,\dots,m\}$ et $\sigma(j_r)=m+r$ pour tout $r\in\{1,\dots,n-m\}$
(autrement dit, $\sigma=\sigma_{(J,K)}$). On note $\varepsilon_0=sgn(\sigma)
(=\varepsilon(J,K))$.
On définit une application $\varphi_{ord}:\Par'_{ord}(n)\fl\Par_{ord}(m)\times
\Par_{ord}(n-m)$ de la manière suivante : si $P=(I_1,\dots,I_k)\in\Par'_{ord}
(n)$ et si $r\in\{1,\dots,k\}$ est tel que $I_1\cup\dots\cup I_r=J$, on pose
$\varphi_{ord}(P)=((u_J(I_1),\dots,u_J(I_r)),(u_K(I_{r+1}),\dots,u_K(I_k)))$.
On définit de manière similaire une application $\varphi:\Par'(n)\fl
\Par(m)\times\Par(n-m)$. Enfin, on définit $\mu=(\mu_1,\dots,\mu_m)\in\R^m$
et $\nu=(\nu_1,\dots,\nu_{n-m})\in\R^{n-m}$ par $\mu_r=\lambda_{u_J^{-1}(r)}
=\lambda_{i_r}$
et $\nu_r=\lambda_{u_K^{-1}(r)}=\lambda_{j_r}$.
Alors :

\begin{itemize}
\item[(iii)] Pour toute $P\in\Par'_{ord}(n)$, si $\varphi_{ord}(P)=(P_1,P_2)$
alors $\varepsilon(P)=\varepsilon_0\varepsilon(P_1)\varepsilon(P_2)$
et $\varepsilon'(P)=\varepsilon'(P_1)\varepsilon'(P_2)$.

\item[(iv)] L'application $\varphi_{ord}$ induit une bijection $\Par'_{ord}(n)
\cap\Par_{ord}(\lambda)\iso\Par_{ord}(\mu)\times\Par_{ord}(\nu)$, et
l'application $\varphi$ induit une bijection $\Par'(n)\cap\Par(\lambda)\iso
\Par(\mu)\times\Par(\nu)$.

\end{itemize}
\end{lemme}

\begin{preuve} Montrons (i). Soit $K\subset\{1,\dots,n\}$. Alors, par
définition de $\lambda'$, $s_K(\lambda')=s_K(\lambda)-\delta|K\cap J|$.
Il est clair que $s_K(\lambda)>0$ si $s_K(\lambda')>0$. On suppose que
$s_K(\lambda)>0$. 
Si $K\not\subset J$, alors $s_K(\lambda)\geq\delta|K|$ et
$|K|>|K\cap J|$, donc $s_K(\lambda')>0$. Si
$K\subset J$ et $K\not=J$, alors $s_K(\lambda)>\delta|K|$
(car $|K|<|J|=N(\lambda)$), donc $s_K(\lambda')>0$.

Montrons (ii).
Comme $\lambda'_i\leq\lambda_i$ pour tout $i\in\{1,\dots,n\}$
et $s_J(\lambda')=0$, il est clair que $\Par(\lambda')\subset\Par(\lambda)$,
$\Par_{ord}(\lambda')\subset\Par_{ord}(\lambda)$, $\Par(\lambda')\cap
\Par'(n)=\varnothing$ et $\Par_{ord}(\lambda')\cap\Par'_{ord}(n)=\varnothing$.

Soit $P=(I_1,\dots,I_k)\in\Par_{ord}(\lambda)-\Par'_{ord}(n)$; montrons que
$P\in\Par_{ord}(\lambda')$.
Il suffit de montrer que $s_{I_1}(\lambda')>0$, car
$(I_1\cup\dots\cup I_r,I_{r+1},\dots,I_k)
\in\Par_{ord}(\lambda)-\Par'_{ord}(n)$ pour tout $r\in\{1,\dots,k\}$.
Or $I_1\not=J$ et $s_{I_1}(\lambda)>0$, donc ceci résulte de (i).

Soit $p=\{I_\alpha,\alpha\in A\}\in\Par(\lambda)-\Par'(n)$; montrons que
$p\in\Par(\lambda')$. Soit $\alpha\in A$. On a $I_\alpha\not=J$ et
$s_{I_\alpha}(\lambda)>0$, donc, d'après (i), $s_{I_\alpha}(\lambda')>0$.

La dernière assertion de (ii) résulte du fait que $(J,\{1,\dots,n\}-J)\in
\Par_{ord}(\lambda)\cap\Par'_{ord}(n)$ et $\{J,\{1,\dots,n\}-J\}\in\Par(
\lambda)\cap\Par'(n)$.

Le point (iii) résulte facilement des définitions.
Montrons (iv). Il est clair
que $\varphi_{ord}$ et $\varphi$ sont injectives, et que $\varphi_{ord}^{-1}
(\Par_{ord}(\mu)\times\Par_{ord}(\nu))\subset\Par_{ord}(\lambda)$ et
$\varphi^{-1}(\Par(\mu)\times\Par(\nu))=\Par(\lambda)$. 
Soit $P=(I_1,\dots,I_k)\in\Par'_{ord}(n)\cap\Par(\lambda)$; montrons que
$\varphi_{ord}(P)\in\Par_{ord}(\mu)\times\Par_{ord}(\nu)$. Soit
$r\in\{1,\dots,k\}$ tel que $J=I_1\cup\dots\cup I_r$. Il s'agit de montrer que,
pour tout $s\geq r+1$, $s_{I_{r+1}\cup\dots\cup I_s}(\lambda)>0$.
Comme on peut, sans
changer le problème, remplacer $P$ par $(I_1\cup\dots\cup I_r,I_{r+1}\cup\dots
I_s,I_{s+1},\dots,I_k)$ (qui est aussi dans $\Par'_{ord}(n)\cap\Par_{ord}
(\lambda)$), on peut supposer que $r=1$ (donc $J=I_1$) et $s=2$. On a alors
\[s_{I_2}(\lambda)=s_{J\cup I_2}(\lambda)-s_J(\lambda)=s_{J\cup I_2}(\lambda)-
\delta|J|\geq \delta(|J\cup I_2|-|J|)>0.\]

\end{preuve}

Soient $\lambda=(\lambda_1,\dots,\lambda_n)\in\R^n$ et $Q=(I_1,\dots,I_r)$
une partition ordonnée d'un sous-ensemble $I$ de $\{1,\dots,n\}$. Si
$k\in\{1,\dots,r\}$, on dit que $Q$ est un \emph{bloc de centre
$I_k$} (pour $\lambda$) si, pour tous $i\in\{1,\dots,k-1\}$ et $j\in\{k+1,
\dots,r\}$, on a $s_{I_1}(\lambda)+\dots+s_{I_i}(\lambda)>0$ et
$s_{I_j}(\lambda)+\dots+s_{I_r}(\lambda)\leq 0$. On dit que
$Q$ est un \emph{bloc} s'il existe $k\in\{1,\dots,r\}$ tel que
$Q$ soit un bloc de centre $I_k$ (noter qu'un tel $k$ n'est pas forcément
unique); la \emph{taille} du bloc $Q$
est par définition $|I|$, et le \emph{support} de $Q$ est $I$.
On dit qu'un bloc $Q$ de support $I$ est \emph{positif} si $s_{I}(\lambda)>0$,
et \emph{négatif} sinon.

Soit $P\in\Par_{ord}(n)$. On dit que $(Q_1,\dots,Q_k)$ est une
\emph{décomposition en blocs} de $P$ si :
\begin{itemize}
\item[$\bullet$] pour tout $l\in\{1,\dots,k\}$, $Q_l=(I_1^l,\dots,I_{r_l}^l)$
est une partition ordonnée d'un sous-ensemble de $\{1,\dots,n\}$, et
$Q_l$ est un bloc;
\item[$\bullet$] $P=(I_1^1,\dots,I_{r_1}^1,\dots,I_1^k,\dots,I_{r_k}^k)$.

\end{itemize}
On écrira souvent $P=Q_1\dots Q_k$.

\begin{lemme}\label{LCI_blocs}\begin{itemize}
\item[(i)] Soit $Q=(I_1,\dots,I_r)$ un bloc pour $\lambda$.
Si $Q$ est positif (resp. négatif), alors, pour
tout $i\in\{1,\dots,r\}$, $s_{I_1}(\lambda)+\dots+s_{I_i}(\lambda)>0$
(resp. $s_{I_i}(\lambda)+\dots+s_{I_r}(\lambda)\leq 0$).
\item[(ii)] Soient $P\in\Par_{ord}(n)$, $P=Q_1\dots Q_k$ une
décomposition en blocs de $P$ et $l\in\{1,\dots,k-1\}$. On suppose que
$P\in\Par_{ord}(\lambda)$ et que $Q_l$ est négatif ou $Q_{l+1}$ est positif.
Alors $P':=Q_1\dots Q_{l-1}Q_{l+1}Q_l
Q_{l+2}\dots Q_k\in\Par_{ord}(\lambda)$. De plus, on a $|P'|=|P|$ et
$\varepsilon(P')=(-1)^{n_ln_{l+1}}\varepsilon(P)$, où $n_l$ (resp. $n_{l+1}$)
est la taille de $Q_l$ (resp. $Q_{l+1}$).

\end{itemize}

Dans la suite du lemme, on suppose que $\lambda_1+\dots+\lambda_n>0$ et que
$N(\lambda)=n$. Alors :
\begin{itemize}
\item[(iii)] Soient $P=(I_1,\dots,I_r)\in\Par_{ord}(\lambda)$, soit
$k\in\{1,\dots,r\}$ et soient
$i_1,\dots,i_k\in\{1,\dots,r\}$ tels que $i_1<\dots<i_k$. Alors il existe
une unique décomposition en blocs $P=Q_1\dots Q_k$ telle que $I_{i_l}$ soit
un centre de $Q_l$ pour tout $l\in\{1,\dots,k\}$.
De plus, $Q_1$ est positif et, si $k\geq 2$, $Q_k$ est négatif.
\item[(iv)] 
Soient $P=(I_1,\dots,I_r)\in\Par_{ord}(\lambda)$ et $P=Q_1\dots Q_k$ une
décomposition en blocs de $P$. Alors il existe un unique $s\in\{1,\dots,r\}$
tel que $P':=(I_s,\dots,I_r,I_1,\dots,I_{s-1})\in\Par_{ord}(\lambda)$, et
$P'$ est de la forme $P'=Q_l\dots Q_kQ_1\dots Q_{l-1}$, pour un $l\in\{1,\dots,
k\}$ uniquement déterminé.

\end{itemize}
\end{lemme}

\begin{preuve}\begin{itemize}
\item[(i)] Soit $k\in\{1,\dots,r\}$ tel que $I_k$ soit un centre de $Q$.
On suppose que $Q$ est positif. D'après la définition d'un bloc de centre
$I_k$, on a $s_{I_1}(\lambda)+\dots+s_{I_i}
(\lambda)>0$ pour tout $i\in\{1,\dots,k-1\}$. D'autre part, pour tout
$i\in\{k,\dots,r\}$,
\[\begin{array}{rcl}s_{I_1}(\lambda)+\dots+s_{I_i}(\lambda) & = &
(s_{I_1}(\lambda)+\dots+s_{I_r}(\lambda))-(s_{I_{i+1}}(\lambda)+\dots+s_{I_r}
(\lambda)) \\
& \geq & s_{I_1}(\lambda)+\dots+s_{I_r}(\lambda)>0.\end{array}\]
Le cas où $Q$ est négatif se traite de manière similaire.

\item[(ii)] On écrit $Q_l=(I_1,\dots,I_r)$ et $Q_{l+1}=(J_1,\dots,J_s)$, et
on note $S$ la somme des $s_I(\lambda)$, pour $I$ parcourant les ensembles de
$Q_1,\dots,Q_{l-1}$. Pour tous $i\in\{1,\dots,r\}$ et $j\in\{1,\dots,s\}$, on
note $B_i=s_{I_i}(\lambda)$ et $C_j=s_{J_j}(\lambda)$.
Comme $P\in\Par_{ord}(\lambda)$, on a $S>0$.
Il s'agit de montrer que, pour tous $i\in\{1,\dots,r\}$ et
$j\in\{1,\dots,s\}$,
\[S+C_1+\dots+C_j>0\]
et
\[S+C_1+\dots+C_s+B_1+\dots+B_i>0.\]
On suppose que $Q_{l+1}$ est positif. Alors la première inégalité résulte
immédiatement du point (i), et la deuxième inégalité résulte du fait que
$C_1+\dots+C_s>0$ et que $S+B_1+\dots+B_i>0$ pour tout $i\in\{1,\dots,r\}$
(car $P\in\Par_{ord}(\lambda)$).
On suppose que $Q_l$ est négatif.
Soit $j\in\{1,\dots,s\}$. Alors
\[\begin{array}{rcl}S+C_1+\dots+C_j & = & (S+B_1+\dots+B_r+C_1+
\dots+C_j)-(B_1+\dots+B_r) \\
& \geq & S+B_1+\dots+B_r+C_1+\dots+C_j>0\end{array}\]
(la première inégalité vient du fait que $Q_l$ est négatif, et la deuxième du
fait que $P\in\Par_{ord}(\lambda)$). Soit $i\in\{1,\dots,r\}$. Alors
\begin{flushleft}$S+C_1+\dots+C_s+B_1+\dots+B_i$\end{flushleft}
\begin{flushright}$\begin{array}{cl} = & (S+B_1+\dots+B_r+C_1+\dots+C_s)-
(B_{i+1}+\dots+B_r) \\
\geq & S+B_1+\dots+B_r+C_1+\dots+C_s>0.\end{array}$\end{flushright}

Enfin, la dernière phrase de (ii) est évidente.

\item[(iii)] On note $\mu=\lambda_P\in\R^r$ (donc
$\mu_i=s_{I_i}(\lambda)$).
On remarque que :
\begin{itemize}
\item[(a)] Pour tout $i\in\{1,\dots,i_1-1\}$ (et même $\{1,\dots,r\}$),
$\mu_1+\dots+\mu_i>0$.
\item[(b)] Pour tout $i\in\{i_k+1,\dots,n\}$ (et même $\{2,\dots,r\}$),
$\mu_i+\dots+\mu_r\leq 0$.
\item[(c)] Soit $l\in\{1,\dots,k-1\}$. Il existe un unique $i\in\{i_l+1,\dots,
i_{l+1}-i\}$ tel que, pour tout $j\in\{i+1,\dots,i_{l+1}-1\}$, on ait
$\mu_{i+1}+\dots+\mu_j>0$ et, pour tout $j\in\{i_l+1,\dots,i\}$, on ait
$\mu_j+\dots+\mu_i\leq 0$.

\end{itemize}
En effet, le point (a) résulte du fait que $P\in\Par_{ord}(\lambda)$. Le point
(b) résulte du fait que $P\in\Par_{ord}(\lambda)$ et du fait, d'après
l'hypothèse $N(\lambda)=n$ et le lemme \ref{LCI_de_base}, il n'existe pas
de partition $(I,J)$ de $\{1,\dots,n\}$ telle que $s_I(\lambda)>0$ et
$s_J(\lambda)>0$. Enfin, le point (c) résulte de la remarque
\ref{rq:LCI_de_base}, appliquée à $(\mu_{i_l},\dots,\mu_{i_{l+1}-1})\in
\R^{i_{l+1}-i_l}$.

Autrement dit, on peut écrire de manière unique
\[P=(A_1^1,\dots,A_{r_1}^1,B_1,C_1^1,\dots,C_{s_1}^1,\dots,A_1^k,\dots,
A_{r_k}^k,B_k,C_1^k,\dots,C_{s_k}^k)\]
avec :
\begin{itemize}
\item[$\bullet$] pour tout $l\in\{1,\dots,k\}$, $B_l=I_{i_l}$;
\item[$\bullet$] pour tout $l\in\{1,\dots,k\}$, pour tous $i\in\{1,\dots,r_l
\}$ et $j\in\{1,\dots,s_l\}$, $s_{A_1^l}(\lambda)+\dots+s_{A_i^l}(\lambda)>0$
et $s_{C_j^l}(\lambda)+\dots+s_{C_{s_l}^l}(\lambda)\leq 0$.

\end{itemize}
On pose, pour tout $l\in\{1,\dots,k\}$, $Q_l=(A_1^l,\dots,A_{r_l}^l,B_l,
C_1^l,\dots,C_{s_l}^l)$. Il est clair (d'après les points (a) et (b) ci-dessus)
que $Q_1$ est positif et $Q_k$ négatif.

\item[(iv)] L'existence et l'unicité de $s$ résultent du lemme
\ref{LCI_de_base}. Pour tout $a\in\{1,\dots,k\}$, on note $\mu_a$ la somme
des $s_I(\lambda)$, pour $I$ parcourant les ensembles de $Q_a$. Alors, en
appliquant le lemme \ref{LCI_de_base} à $(\mu_1,\dots,\mu_k)$, on voit qu'il
existe un unique $l\in\{1,\dots,k\}$ tel que $(\mu_l,\dots,\mu_k,\mu_1,\dots,
\mu_{l-1})>0$. On note $P''=Q_l\dots Q_kQ_1\dots Q_{l-1}$. Pour prouver (iv),
il suffit de montrer que $P''\in\Par_{ord}(\lambda)$. On se ramène donc à
montrer l'énoncé suivant (où les notations ne sont plus les mêmes) :

Soient $P\in\Par_{ord}(n)$ et $P=Q_1\dots Q_k$ une décomposition en blocs de
$P$. Pour tout $a\in\{1,\dots,k\}$, on note $\mu_a$ la somme
des $s_I(\lambda)$, pour $I$ parcourant les ensembles de $Q_a$. On suppose que
$(\mu_1,\dots,\mu_k)>0$. Alors $P\in\Par_{ord}(\lambda)$.

Pour tout $a\in\{1,\dots,k\}$, on écrit $Q_a=(I_1^a,\dots,I_{r_a}^a)$. 
Soient $a\in\{1,\dots,k\}$ et $i\in\{1,\dots,r_a\}$. On veut montrer que
\[\mu_1+\dots+\mu_{a-1}+s_{I_1^a}(\lambda)+\dots+s_{I_i^a}(\lambda)>0.\]
Si le bloc $Q_a$ est positif, cela résulte de (i) (et de la positivité
de $(\mu_1,\dots,\mu_k)$). Supposons que $Q_a$ est négatif. Alors, toujours
grâce à (i) :
\[\begin{array}{rcl}\mu_1+\dots+\mu_{a-1}+s_{I_1^a}(\lambda)+\dots+s_{I_i^a}
(\lambda) & = & \mu_1+\dots+\mu_a-(s_{I_{i+1}^a}(\lambda)+\dots+s_{I_{r_a}^a}
(\lambda)) \\
& \geq & \mu_1+\dots+\mu_a>0.\end{array}\]

\end{itemize}
\end{preuve}


\begin{thebibliography}{9999}

\bibitem[A]{A-L2} J. Arthur, \emph{The $L^2$-Lefschetz numbers of Hecke
operators}, Inv. Math. 97 (1989), p 257-290

\bibitem[B]{Bo} A. Borel, \emph{Automorphic $L$-functions}, dans \emph{
Automorphic forms, representations, and $L$-functions} (Proc. Symposia
in Pure Math., volume 33, 1977), tome 2, p 26-61

\bibitem[C]{C} V.I. Chernousov, \emph{The Hasse principle for groups of type
$E_8$}, Soviet Math. Dokl. 39 (1989), p 592-596

\bibitem[CD]{CD} L. Clozel et P. Delorme, \emph{Pseudo-coefficients et 
cohomologie des groupes de Lie réductifs réels}, C.R. Acad. Sc. Paris 300, 
série I (1985), p 385-287

\bibitem[D]{Di} J. Dieudonné, \emph{La géométrie des groupes classiques.
Troisième édition}, Ergebnisse der Mathematik und ihrer Grenzgebiete, Band 5,
Springer (1971)

\bibitem[GHM]{GHM} M. Goresky, G. Harder et R. MacPherson, \emph{Weighted cohomology}, Invent. math. 166 (1994), p 139-213

\bibitem[GKM]{GKM} M. Goresky, R. Kottwitz et R. MacPherson, \emph{Discrete series characters and the Lefschetz formula for Hecke operators}, Duke Math. J. 89 (1997), p 477-554 et Duke Math. J. 92 (1998), no. 3, p 665-666

\bibitem[Ha]{H} T. Hales, \emph{A simple definition of transfer factors for
unramified groups}, in \emph{Representation theory of groups and algebras},
édité par J. Adams, R. Herb, S. Kudla, J.-S. Li, R. Lipsman et J. Rosenberg,
Contemporary Mathematics 145 (1993), p 109-134


\bibitem[H]{Her} R. Herb, \emph{Characters of averaged discrete series on
semisimple real Lie groups}, Pac. J. Math. 80 (1979), p 169-177

\bibitem[K1]{K-SC} R. Kottwitz, \emph{Sign changes in harmonic analysis on
reductive groups}, Trans. A.M.S. 278 (1983), p 289-297

\bibitem[K2]{K-SVTOI} R. Kottwitz, \emph{Shimura varieties and twisted orbital
integrals}, Math. Ann. 269 (1984), p 287-300

\bibitem[K3]{K-STF:CTT} R. Kottwitz, \emph{Stable trace formula : cuspidal
tempered terms}, Duke Math. J. 51 (1984), p 611-650

\bibitem[K4]{K-BC} R. Kottwitz, \emph{Base change for units of Hecke 
algebras}, Compositio Math. 60 (1986), p 237-250

\bibitem[K5]{K-STF:EST} R. Kottwitz, \emph{Stable trace formula : elliptic
singular terms}, Math. Ann. 275 (1986), p 365-399

\bibitem[K6]{K-TN} R. Kottwitz, \emph{Tamagawa numbers}, Ann. of Math. 127 
(1988), p 629-646

\bibitem[K7]{K-SVLR} R. Kottwitz, \emph{Shimura varieties and $\lambda$-adic
representations}, dans \emph{Automorphic forms, Shimura varieties and
$L$-functions}, partie I, Perspectives in Mathematics  vol. 10, Academic Press,
San Diego, CA (1990), p 161-209

\bibitem[K8]{K-NP} R. Kottwitz, non publié

\bibitem[Lan1]{L2} R. Langlands, \emph{Stable conjugacy : definitions and lemmas}, Can. J. Math. 31 (1979), p 700-725

\bibitem[Lan2]{L3} R. Langlands, \emph{Les débuts d'une formule des traces
stable}, Publ. Math. Univ. Paris VII vol. 13, Paris (1983)

\bibitem[LR]{LR} R. Langlands et D. Ramakrishnan (éditeurs), \emph{The zeta
function of Picard modular surfaces}, publications du CRM (1992), Montréal

\bibitem[LS1]{LS1} R. Langlands and D. Shelstad, \emph{On the definition of
transfer factors}, Math. Ann. 278 (1987), p 219-271

\bibitem[LS2]{LS2} R. Langlands and D. Shelstad, \emph{Descent for transfer
factors}, The Grothendieck Festschrift, Vol. II, Progr. Math. 87, 
Birkh{\"a}user (1990), p 485-563


\bibitem[Lau1]{Lau} G. Laumon, \emph{Sur la cohomologie à supports compacts des
variétés de Shimura pour $\GSp(4)_\Q$}, Compositio Math. 105 (1997), no. 3,
p 267-359

\bibitem[Lau2]{Lau2} G. Laumon, \emph{Fonctions zêta des variétés de Siegel
de dimension trois}, dans \emph{Formes automorphes II. Le cas du groupe
$\GSp(4)$}, Astérisque 302 (2005), p 1-66

\bibitem[LN]{LN} G. Laumon and B.-C. Ngo, \emph{Le lemme fondamental pour
les groupes unitaires}, Annals of Math. 168 (2008), no. 2, p 477-573


\bibitem[M1]{M2} S. Morel, \emph{Complexes pondérés sur les compactifications
de Baily-Borel. Le cas des variétés de Siegel}, Journal of the AMS 21 (2008),
no. 1, p 23-61

\bibitem[M2]{M3} S. Morel, \emph{On the cohomology of certain non-compact
Shimura varieties} (à paraître dans la série Annals of Mathematics Studies de
Princeton University Press),
http://www.math.ias.edu/~morel/stabilisation.pdf

\bibitem[N]{Ng} B. C. Ngo, \emph{Le lemme fondamental pour les algèbres de
Lie} (soumis), arXiv:0801.0446

\bibitem[O]{O} T. Ono, \emph{On Tamagawa numbers}, in \emph{Algebraic groups
and discontinuous subgroups}, édité par A. Borel et G. Mostow, Proceedings
of Symposia in Pure Math. 9 (1966)

\bibitem[P1]{P1} R. Pink, \emph{Arithmetical compactification of mixed Shimura varieties}, thèse, Bonner Mathematische Schriften 209 (1989)

\bibitem[P2]{P2} R. Pink, \emph{On $\ell$-adic sheaves on Shimura varieties and their higher direct images in the Baily-Borel compactification}, Math. Ann. 292 (1992), 197-240

\bibitem[R]{Ro1} J. Rogawski, \emph{Automorphic representations of
unitary groups in three variables}, Annals of Mathematics Study 123,
Princeton University Press (1990)

\bibitem[W1]{Wa1} J.-L.Waldspurger, \emph{Le lemme fondamental implique le
transfert}, Comp. Math. 105 (1997), n${}^\circ$2, p 153-236

\bibitem[W2]{Wa2} J.-L.Waldspurger, \emph{Endoscopie et changement de
caractéristique}, J. Inst. Math. Jussieu 5 (2006), n${}^\circ$3, p 423-525

\bibitem[W3]{Wa3} J.-L.Waldspurger, \emph{L'endoscopie tordue n'est pas si
tordue}, Mem. Amer. Math. Soc. 194 (2008), no. 908


\end{thebibliography}
\end{document}